\documentclass[twoside, 10pt, reqno]{article}

\usepackage{eucal}
\usepackage{amsmath}
\usepackage{amssymb}
\usepackage{theorem}
\usepackage[all,v2, rotate]{xy}
\xyoption{2cell} \UseAllTwocells

\newtheorem{prop}{Proposition}[section]
\newtheorem{lem}[prop]{Lemma}

\newtheorem{cor}[prop]{Corollary}
\newtheorem{them}[prop]{Theorem}

\newtheorem{thm}[prop]{Theorem}

\theorembodyfont{\upshape}

\newtheorem{defn}[prop]{Definition}

\newtheorem{numrmk}[prop]{Remark}

\newtheorem{numex}[prop]{Example}

\newtheorem{warning}[prop]{Warning}

\newtheorem{claim}{Claim}
\newtheorem{example}[prop]{Example}
\newtheorem{ex}[prop]{Example}

\newtheorem{rmk}[prop]{Remark}
\newtheorem{rem}[prop]{Remark}

\numberwithin{equation}{section}

\newenvironment{pf}{\begin{trivlist}\item[]{\sc Proof.}}%
            {\nolinebreak $\Box$ \end{trivlist}}

\newenvironment{proof}{\begin{trivlist}\item[]{\sc Proof.}}%
            { $\Box$ \end{trivlist}}

\newcommand{\eq}{\begin{eqnarray}}
\newcommand{\eneq}{\end{eqnarray}}
\newcommand{\eqn}{\begin{eqnarray*}}
\newcommand{\eneqn}{\end{eqnarray*}}

\newcommand{\bnum}{\begin{enumerate}[{\rm(i)}]}
\newcommand{\enum}{\end{enumerate}}
\newcommand{\banum}{\begin{enumerate}[{\rm(a)}]}
\newcommand{\eanum}{\end{enumerate}}

\newcommand{\noprint}[1]{}

\newcommand{\comment}[1]{{\marginpar{\footnotesize #1}}}
\newcommand{\unsure}[1]{{\footnotesize #1}}
\newcommand{\notsure}[1]{{ #1}}
\def\comment{\noprint}
\def\unsure{\noprint}
\renewcommand{\tilde}{\widetilde}

\newcommand{\toto}{\rightrightarrows}

\newcommand{\scom}{{\scriptscriptstyle\bullet}}
\newcommand{\com}{^{\scriptstyle\bullet}}
\newcommand{\lcom}{_{\scriptstyle\bullet}}

\newcommand{\X}{\mathop{\sf X}}
\newcommand{\XX}{{\mathfrak X}}
\renewcommand{\AA}{{\mathfrak A}}
\newcommand{\BB}{{\mathfrak B}}
\newcommand{\DD}{{\mathfrak D}}
\renewcommand{\SS}{{\mathfrak S}}
\newcommand{\TT}{{\mathfrak T}}

\newcommand{\YY}{{\mathfrak Y}}
\newcommand{\UU}{{\mathfrak U}}

\newcommand{\ZZ}{{\mathfrak Z}}

\newcommand{\EE}{{\mathfrak E}}
\newcommand{\FF}{{\mathfrak F}}
\newcommand{\GG}{{\mathfrak G}}
\newcommand{\NN}{{\mathfrak N}}
\newcommand{\MM}{{\mathfrak M}}
\newcommand{\KK}{{\mathfrak K}}
\newcommand{\RR}{{\mathfrak R}}

\newcommand{\VV}{{\mathfrak V}}

\newcommand{\Gg}{{\mathfrak g}}

\newcommand{\rr}{{\mathbb R}}

\renewcommand{\O}{{\cal O}}

\newcommand{\rk}{\mathop{\rm rk}}

\newcommand{\ev}{\mathop{\rm ev}\nolimits}

\newcommand{\pr}{\mathop{\rm pr}\nolimits}

\newcommand{\Ad}{\mathop{\rm Ad}\nolimits}
\newcommand{\rank}{\mathop{\rm rank}\nolimits}

\newcommand{\codim}{\mathop{\rm codim}\nolimits}

\newcommand{\id}{\mathop{\rm id}\nolimits}
\newcommand{\Hom}{\mathop{\rm Hom}\nolimits}

\newcommand{\coker}{\mathop{\rm Coker}\nolimits}

\newcommand{\injectlim}{\mathop{\lim\limits_{\textstyle\longrightarrow}}\limits}

\newcommand{\ldiag}[1]%
       {\makebox[0cm]{${\scriptstyle#1}\downarrow\phantom{\scriptstyle#1}$}}
\newcommand{\ldiagup}[1]%
       {\makebox[0cm]{${\scriptstyle#1}\uparrow\phantom{\scriptstyle#1}$}}
\newcommand{\rdiag}[1]%
       {\makebox[0cm]{$\phantom{\scriptstyle#1}\downarrow{\scriptstyle#1}$}}
\newcommand{\sediagr}[1]%
       {\makebox[0cm]{$\phantom{\scriptstyle#1}\searrow{\scriptstyle#1}$}}
\newcommand{\nediagr}[1]%
       {\makebox[0cm]{$\phantom{\scriptstyle#1}\nearrow{\scriptstyle#1}$}}
\newcommand{\rdiagup}[1]%
       {\makebox[0cm]{$\phantom{\scriptstyle#1}\uparrow{\scriptstyle#1}$}}
\newcommand{\swdiag}[1]%
       {\makebox[0cm]{$\phantom{\scriptstyle#1}\swarrow{\scriptstyle#1}$}}
\newcommand{\sediag}[1]%
       {\makebox[0cm]{${\scriptstyle#1}\searrow\phantom{\scriptstyle#1}$}}
\newcommand{\nediag}[1]%
       {\makebox[0cm]{${\scriptstyle#1}\nearrow\phantom{\scriptstyle#1}$}}

\newcommand{\iso}{\stackrel{\sim}{\rightarrow}}

\newcommand{\doublearrowstack}[2]%
                      {{{{\scriptstyle#1}\atop{\textstyle\longrightarrow}}\atop{{\textstyle\longrightarrow}\atop{\scriptstyle#2}}}}
\newcommand{\rightleftarrowstack}[2]%
                      {{{{\scriptstyle#1}\atop{\textstyle\longrightarrow}}\atop{{\textstyle\longleftarrow}\atop{\scriptstyle#2}}}}
\newcommand{\leftrightarrowstack}[2]%
                      {{{{\scriptstyle#1}\atop{\textstyle\longleftarrow}}\atop{{\textstyle\longrightarrow}\atop{\scriptstyle#2}}}}

\newcommand{\ra}{\rightarrow}

\newcommand{\la}{\leftarrow}

\newcommand{\lra}{\longrightarrow}
\newcommand{\llra}[1]{\stackrel{#1}{\lra}}

\newcommand{\overtoparrow}%
{\makebox[0cm]{\beginpicture
\setcoordinatesystem units <.8cm,.4cm> point at 0 0
\setplotarea x from -3 to 3, y from 0 to 1
\setquadratic
\plot -3 0 0 1 3 0 /
\put{\vector(3,-1){0}}[Bl] at 3 0
\endpicture}}

\newcommand{\underbottomarrow}%
{\makebox[0cm]{\beginpicture
\setcoordinatesystem units <.8cm,.4cm> point at 0 0
\setplotarea x from -3 to 3, y from 0 to 1
\setquadratic
\plot -3 1 0 0 3 1 /
\put{\vector(3,1){0}}[Bl] at 3 1
\endpicture}}

\newcommand{\ses}[5]%
{0\longrightarrow#1\stackrel{#2}{ \longrightarrow}#3\stackrel{#4}{
\longrightarrow}#5\longrightarrow0}

\newcommand{\dt}[6]%
{#1\stackrel{#2}{longrightarrow}#3 \stackrel{#4}{\longrightarrow}#5
\stackrel{#6}{\longrightarrow} #1[1]}

\newcommand{\cat}[1]%
{(\mbox{\rm #1})}

\newcommand{\sfT}{\mathsf{T}}
\newcommand{\Top}{\mathsf{Top}}

\def\smashedlongrightarrow{\setbox0=\hbox{$\longrightarrow$}\ht0=1pt\box0}
\def\risom{\buildrel\sim\over{\smashedlongrightarrow}}

\def\smashedst{\setbox0=\hbox{$\rightrightarrows$}\ht0=4pt\box0}
\newcommand{\sst}[1]{\stackrel{#1}{\smashedst}}

\newcommand{\bfhom}{\mathbf{Hom}}

\newcommand{\map}{\operatorname{Map}}
\newcommand{\Map}{\operatorname{Map}}

\def \BH {\mathbb{H}}
\def \Gammaa {\XX}
\def \LG {{\rm L}\Gammaa}
\def  \IG {\Lambda \Gammaa}

\def \Lo {{\rm L}}

\newcommand{\lllra}[4]{\xymatrix@C=14pt{#1 \ar[r]^{#2}_*+[o][F-]{#3} & {#4}}}
\newcommand{\al}{\alpha}
\newcommand{\be}{\beta}
\newcommand{\ep}{\epsilon}

\newcommand{\ga}{\gamma}

\newcommand{\BV}{{\bf BV}}
\newcommand{\copd}{\mathop{\coprod}}

\newcommand{\Gk}{{\mathfrak k}}

\newcommand{\gm}{\Gamma}

\def \LG {{\rm L }\Gamma}
\def  \IG {\Lambda \Gamma}
\def  \IMG {\Lambda [M/G]}
\def \ING {\Lambda [N\stackrel{i}\to \gm]}

\def \PH {\mathcal{P}}
\def \PHo {\mathcal{P}^{\rm orb}}
\def \sH {\mathbb{H}}
\def \CR {\cup_{\rm orb}}

\def \LM {{\rm L}M}

\def \LXX {{\rm L}\XX}

\def \IXX {\Lambda \XX}

\def \IWG {\widetilde{\Lambda \Gamma}}

\def \Ex {\mathfrak{E}}
\def \Ob {{\mathfrak{O}}}

\def \kor {{ k}}
\def \eu {\mathfrak{e}}

\def \cac {\mathcal{C}ac}

\def \IPG {\Lambda [*/G]}
\def \LPG {{\rm L}[*/G]}
\newcommand{\LL}{\mathcal{L}}
 
\newcommand{\TopSt}{\mathsf{TopSt}}
\newcommand{\bbZ}{\mathbb{Z}}
\newcommand{\bbX}{\mathbb{X}}
\newcommand{\bbR}{\mathbb{R}}
\newcommand{\bbC}{\mathbb{C}}

\newcommand{\bbP}{\mathbb{P}}

\newcommand{\DR}[1]{H_{\rm DR}^{#1}}
\newcommand{\DRC}[1]{H_{\rm DR,\, c}^{#1}}
\newcommand{\gy}[1]%
{ {{#1}^!} }

\newcommand{\cgy}[1]%
{ {{#1}_!} }

\def\lisom{\buildrel\sim\over{\smashedlongleftarrow}}
\def\smashedlongrightarrow{\setbox0=\hbox{$\longrightarrow$}\ht0=1pt\box0}
\def\risom{\buildrel\sim\over{\smashedlongrightarrow}}

\def\smashedlongleftarrow{\setbox0=\hbox{$\longleftarrow$}\ht0=1pt\box0}
\def\lisom{\buildrel\sim\over{\smashedlongleftarrow}}
\newcommand{\torel}[1]{\stackrel{#1}{\to}}
 \newcommand{\sfC}{\mathsf{C}}
\def \A {\AA}
\let\:=\colon
\let\oldcite\cite
\def \eb {\breve{e}}
\def \B {\BB}
\def \E {\EE}
\def \Y {\YY}
\def \V {\VV}
\def \X {\XX}
\def \F {\FF}
\def \K {\KK}
\def \Z {\ZZ}
\def \N {\NN}
\def \U {\UU}
\def \D {\DD}
\def \M {\MM}
\def \age {{\rm age}}
\newcommand{\Horb}[1]{H_{#1}^{\rm orb}}

\newcommand{\Lind}{{Lindel\"of }}
\def \proper {{bounded proper }}
\def \superproper {{strongly proper }}
\def \properpt {{bounded proper}}
\def \superproperpt {{strongly proper}}
\def \regular {{regular }}
\def \Regular {{Regular }}
\def \half {{\frac{1}{2}}}
\def \sz {(\mathbb{Z}/2\mathbb{Z})^{n+1}}
\def \hidden {hidden loop }

\author{K.~Behrend, G.~Ginot, B.~Noohi and P.~Xu}

\title{String topology for stacks}

\begin{document}
\sloppy

\maketitle

\begin{abstract}
We establish the general machinery of string
topology for differentiable stacks. This machinery allows us to
treat on equal footing free loops in stacks and hidden loops. 
%In
%particular, we give a good notion of  free loop stack, and  of 
%mapping stack $\map(Y,\XX)$, where $Y$ is a compact space and $\XX$
%a topological stack, which is  functorial both in $\XX$ and $Y$ and
%behaves well enough with respect to pushouts. 
We construct a
bivariant (in the sense of Fulton and MacPherson) theory for
topological stacks: it gives us a flexible theory of Gysin maps
which are automatically compatible with pullback, pushforward and
products. Further we prove an excess formula in this context.
We introduce oriented stacks,  generalizing oriented
manifolds, which are stacks on which we can do string topology. We
prove that the homology of the free loop stack   of an oriented
stack and the homology of hidden loops (sometimes called ghost loops) are a 
Frobenius algebra which are related by a natural morphism of Frobenius algebras. 
We also prove that the homology of free loop stack has a natural structure 
of \BV-algebra, which together with the Frobenius structure fits into an 
homological conformal field theories with closed positive boundaries.  
Using our general machinery, we construct an intersection pairing for 
(non necessarily compact) almost complex orbifolds which is in 
the same relation to the
intersection pairing for manifolds as Chen-Ruan orbifold cup-product
is to ordinary cup-product of manifolds. We show that the \hidden
product of almost complex is isomorphic to the orbifold intersection
pairing twisted by a canonical class. Finally we gave some examples 
including the case of the classifying stacks $[*/G]$ of a compact Lie group.
\end{abstract}

\tableofcontents

\section*{Introduction}
\addcontentsline{toc}{section}{Introduction}

%\subsubsection{Headline}

String topology is a term coined by Chas-Sullivan~\cite{CS} to
describe the rich algebraic structure on the homology of the free
loop manifold $\LM$ of an oriented manifold $M$. The algebraic
structure in question is induced by geometric operations on loops
such as gluing or pinching of loops. In particular, $H\lcom(\LM)$
inherits a canonical product and coproduct yielding a structure of
Frobenius algebra~\cite{CS, CoGo}. Furthermore, the canonical action
of $S^1$ on $\LM$ together with the multiplicative structure make
$H\lcom(\LM)$ a \BV-algebra~\cite{CS}. These algebraic structures,
especially the loop product, are known to be related to many
subjects in mathematics and in particular mathematical
physics~\cite{Su, CFP, Coh, AbZe, CoVo}.

Many interesting geometric objects in (algebraic or differential)
geometry or mathematical physics are \emph{not} manifolds. There are,
for instance, orbifolds, classifying spaces of compact Lie groups, or,
more generally, global quotients of a manifold by a Lie group. All
these examples belong to the realm of (geometric) stacks. A natural
generalization of smooth manifolds, including the previous examples,
is given by differentiable stacks~\cite{BX} (on which one can still do
differentiable geometry).  Roughly speaking, differential stacks are
Lie groupoids {\em up to Morita equivalence}.

One important feature of differentiable stacks is that they are {\em
non-singular}, when viewed as stacks (even though their associated
coarse spaces are typically singular).  For this reason,
differentiable stacks have an intersection product on their homology,
and a loop product on the homology of their free loop stacks.

The aim of this paper is to establish the general machinery of string
topology for differentiable stacks. This machinery allows us to treat on an
equal footing free loops in stacks and \emph{hidden} loops. The latter
are loops inside the stack, which vanish on the associated coarse
space.  The stack of hidden loops in the stack $\XX$ is the {\em
  inertia stack }of $\XX$, notation $\Lambda\XX$.  The inertia stack
$\Lambda\XX\to\XX$ is an example of a family of commutative {\em
  (sic!) }groups
over the stack $\XX$, and the theory of hidden loops generalizes to
arbitrary commutative families of groups over stacks.

In the realm of stacks several new difficulties arise whose solutions
should be of independent interest.

First, we need a good notion of \emph{free loop stack} $\LXX$ of a
stack $\XX$, and more generally of mapping stack $\map(Y,\XX)$ (the
stack of stack morphisms $Y\to \XX$). For the general theory of
mapping stacks, we do not need a differentiable structure on $\XX$;
we work with topological stacks.  This is developed in \cite{Mapping}
and is discussed in
Section~\ref{mappingstack}. The issue here is to obtain a mapping
stack with a topological structure which is functorial both in $\XX$
and $Y$ and behaves well enough with respect to pushouts in order to
get geometric operations on loops. For instance, a key point in
string topology is the identification $\map(S^1\vee S^1,\XX) \cong
\LXX\times_\XX \LXX$. Since pushouts are a delicate matter in the
realm of stacks, extra care has to be taken in finding the correct
class of  topological stacks to work with (Section~\ref{topological}
and \cite{Noohi}). For this reason, we restrict our attention to the
class of {\em Hurewicz} topological stacks. These are topological
stacks which admit an atlas with a certain fibrancy property.
Without restricting to this special class of topological stacks,
$S^1\vee S^1$ would not be the pushout of two copies of $S^1$, in
the category of stacks.

A crucial step in usual string topology is the existence of a
canonical Gysin homomorphism $H\lcom(\LM\times \LM)\to
H_{\scriptscriptstyle \bullet -d}(\LM\times_M \LM)$ when $M$ is a
$d$-dimensional manifold. In fact, the loop product is the composition
\begin{multline}
H_{p}(\LM)\otimes H_{q}(\LM)\to\\
\to H_{p+q}(\LM\times \LM)\to H_{p+q-d}(\LM\times_M \LM) \to
H_{p+q-d}(\LM)\,, \label{origdenf}
\end{multline}
where the last map is obtained by gluing two loops at their base
point.

Roughly speaking the Gysin map can be obtained as follows. The free
loop manifold is equipped with a structure of Banach manifold such
that the evaluation map $\ev: \LM \to M$ which maps a loop $f$ to
$f(0)$ is a surjective submersion. The pullback along $\ev\times\ev$ of a
tubular neighborhood of the diagonal $M\to M\times M$ in $M\times M$
yields a normal bundle of codimension $d$ for the embedding
$\LM\times_M \LM \to \LM$. The Gysin map can then be constructed using a
standard argument on Thom isomorphism and Thom collapse~\cite{CoJo}.

This approach does {not} have a straightforward generalization to
stacks. For instance, the free loop stack of a differentiable stack is
not a Banach stack in general, and neither is the inertia stack. In
order to obtain a flexible theory of Gysin maps, we construct a
\emph{bivariant theory} in the sense of Fulton-MacPherson~\cite{FuMac}
for topological stacks, whose underlying homology theory is singular
homology. A bivariant theory is an efficient tool encompassing into a
unified framework both homology and cohomology as well as many
(co)homological operations, in particular Gysin homomorphisms. The
Gysin maps of a bivariant theory are automatically compatible with
pullback, pushforward, cup and cap-products (see~\cite{FuMac}). (Our
bivariant theory is somewhat weaker than that of Fulton-MacPherson, in
that products are not always defined.) Our bivariant theory applies in particular
to all orbifolds. Further we gave an excess formula allowing to compute Gysin maps for relative regular embeddings.

In Section~\ref{S:nns} we introduce {\em oriented
stacks}. 
These are the stacks over which we are able to do string
topology. Examples of oriented stacks include: oriented manifolds,
oriented orbifolds, and quotients of oriented manifolds by compact Lie
groups (if the action is orientation preserving and of finite orbit
type). A topological stack $\XX$ is \emph{orientable} if the
diagonal map $\XX\to \XX \times \XX$ factors as
\begin{equation}\label{factornorm}
\XX\stackrel{0}{\longrightarrow}\NN\longrightarrow\EE\longrightarrow\XX\times\XX\,,
\end{equation}
where $\NN$ and $\EE$ are orientable vector bundles over $\XX$ and
$\XX\times \XX$ respectively, and $\NN\to\EE$ is an isomorphism onto
an open substack (there is also the technical assumption that $\EE$ is
metrizable, and $\XX\to\EE$ factors through the unit disk bundle). The
embedding $\NN \to \EE$ plays the role of a tubular neighborhood. The
dimension of $\XX$ is $\rk\NN-\rk\EE$.

The factorization (\ref{factornorm}) gives rise to a bivariant class
$\theta\in H(\XX\to\XX\times\XX)$, the {\em orientation }of $\XX$.

\medskip

Sections 10-15 are devoted to the string topology operations, focusing
on the Frobenius, \BV-algebra and homological conformal field theory structures. The bivariant formalism
has the following consequence: if $\XX$ is  an oriented stack of
dimension $d$, then any cartesian square
\begin{equation*}
\xymatrix@C=12pt@R=10pt@M=8pt{\YY\dto \rto & \ZZ
\dto\\ \XX \rto^{\hspace{-0.3cm}\Delta} & \XX\times
 \XX}
\end{equation*}
defines a canonical Gysin map
$\Delta^!\:H\lcom(\ZZ)\to H_{\scriptscriptstyle \bullet-d}(\YY)$.
For example, the cartesian square
\begin{equation*}
\xymatrix@C=12pt@R=10pt@M=8pt{\LXX\times_\XX\LXX\dto \rto & \LXX\times\LXX
\dto\\ \XX \rto^{\hspace{-0.3cm}\Delta} & \XX\times
 \XX}
\end{equation*}
Gives rise to a Gysin map
$\Delta^!\:H\lcom(\LXX\times\LXX)\to H_{\scriptscriptstyle
  \bullet-d}(\LXX\times_\XX\LXX)$, and we can construct a loop product
$$ \star: H_{\scriptstyle
\bullet}(L\XX)\otimes H_{\scriptstyle \bullet}(L\XX) \to
H_{\scriptstyle \bullet-d}(L\XX), $$
as in~\ref{origdenf}, or~\cite{CS, CoJo, CoGo}.

We also obtain a coproduct
$$\delta\: H\lcom (\LXX)\longrightarrow
\bigoplus_{i+j=\bullet-d}H_{i}(\LXX) \otimes H_{j}(\LXX)\,.$$
Furthermore, $\LXX$ admits a natural $S^1$-action yielding the
operator $D: H\lcom (\LXX)\to H_{\scriptstyle \bullet +1} (\LXX)$
which is the composition:
$$H\lcom(\LXX)\stackrel{\times \omega}\longrightarrow
 H_{\scriptstyle \bullet+1}(\LXX\times S^1)
\longrightarrow H_{\scriptstyle \bullet+1}(\LXX), $$
where $\omega\in H_1(S^1)$ is the fundamental class.
Thus we prove that $\left( H\lcom(\LXX),\star,\delta\right)$ is
a  Frobenius algebra and that the shifted homology
$\left(H_{\scriptstyle \bullet +d}(\LXX),\star, D\right)$ is a
 \BV-algebra. Using Sullivan's chord diagram~\cite{CoGo} and our formalism of Gysin maps given by the bivariant theory, we extend the previous $\BV$ and Frobenius structure into a homological conformal field theory (in the sense of~\cite{Cos, Seg04}) with closed positive boundaries (said otherwise non-unital and non-counital) in Section~\ref{S:HCFTbig}. Roughly, this means that to any compact Riemann surface $\Sigma$ with only closed boundaries 
(with say $n$ incoming ones and $m$ outgoing ones), and such that any connected component of $\Sigma$ has a positive number of both incoming and outgoing boundary components, and to any class $\alpha$ in the homology of the mapping class group of $\Sigma$, we associate an operation $\mu_\alpha: H(\LXX)^{\otimes n}\to H(\LXX)^{\otimes m}$ compatible with the glueing and disjoint union of surfaces.  

\smallskip

Since the inertia stack can be considered as the stack of hidden
loops, the general machinery of Gysin maps yields, for any oriented
stack $\XX$, a product and a coproduct on the homology $H\lcom(\IXX)$
of the inertia stack $\IXX$, making it a Frobenius algebra,
too. Moreover in Section~\ref{S:Frobmap}, we construct a natural map
$\Phi\: \IXX \to \LXX$ inducing a morphism of Frobenius algebras in
homology.

\smallskip

In Section~\ref{S:Brane}, we explain how to adapt the loop product to the case of spheres spaces $\bfhom(S^n,\XX)$. We obtain an analogue of the loop product, called the brane product, and also study power maps $\lambda^k:H\lcom((\bfhom(S^n,\XX)) \to H\lcom((\bfhom(S^n,\XX))$ induced by the  degree $k$ maps $S^n\to S^n$. We show that for $n\ge 2$, the maps $\lambda^k$ are maps of algebras with respect to the brane product.

\smallskip

In Section~\ref{Orbifolds}, we consider almost complex orbifolds
(not necessarily compact).  Using Gysin maps and the obstruction
bundle of Chen-Ruan~\cite{CR}, we construct the \emph{orbifold
intersection pairing} on the homology of the inertia stack. It is in
the same relation to the intersection pairing on the homology of a
manifold as the Chen-Ruan orbifold cup-product~\cite{CR} is to the
ordinary cup product on the cohomology of a manifold.

The orbifold intersection pairing defines a structure of associative,
graded commutative algebra on $\Horb{\scom}(\XX)$ for any almost
complex orbifold $\XX$. As a vector space the orbifold homology
$\Horb{\scom}(\XX)$ coincides with the homology of the inertia stack
$\IXX$, but the grading is shifted according to the age as in~\cite{CR,
FG}.

In the compact case, the orbifold intersection pairing is identified
with the Chen-Ruan product, via orbifold Poincar\'e duality.

\smallskip

We also prove that the loop product, \hidden product and intersection
pairing (for almost complex orbifolds) can be twisted by a
cohomology class in $H\lcom(\LXX\times_{\XX}\LXX)$ or $H\lcom(\IXX
\times_{\XX} \IXX)$, satisfying the 2-cocycle condition (see
Propositions~\ref{pr:twistedassociativity},~\ref{twistedstring},
and~\ref{CRcocycle}).  The notion of twisting provides a connection
between the orbifold intersection pairing and the \hidden product. In
fact, we associate to an almost complex orbifold $\XX$ a canonical
vector bundle $\O_\XX\oplus \N_\XX$ over $\IXX\times_\XX \IXX$ and
prove that the orbifold intersection pairing, twisted by the Euler
class of $\O_\XX\oplus \N_\XX$, is the \hidden product of $\XX$.

\smallskip

Parallel to our work, the \hidden product for global quotient orbifolds
was studied in~\cite{LUX, GLSU}. Furthermore, a nice interpretation of
the \hidden product in terms of the Chen-Ruan product of the cotangent
bundle was given by Gonz{\'a}lez et al.~\cite{GLSU}. A loop product
for global quotients of a manifold by a finite group was studied
in~\cite{LuUrXi, LUX}. Also purely homotopical techniques to study string topology of classifying spaces of Lie groups have been recently studied in~\cite{CM}.

We close this introduction by remarking that our construction of
string operations for stacks can in fact be extended to
generalized (co)homology theories other than singular.
For instance, in view of the 
Freed-Hopkins-Teleman's work~\cite{FHT},  using $K$-theory 
may lead to interesting consequences.
In the case  of manifolds, Cohen and Godin have already
considered such generalization in ~\cite{CoGo}. 

The key point in extending our theory to other 
(co)homology theories is to cast such a
(co)homology theory as part of a bivariant theory.
Once this is done, the formalism developed in 
Section~\ref{S:Bivariant} applies to produce
the desired Gysin maps, and these in turn give rise to 
string operations. The main input needed to make the
construction of the Gysin maps possible is to produce an orientation 
class $\theta$ in the bivariant cohomology of the diagonal $\XX \to \XX\times \XX$,
and this is done by making use of appropriate Thom classes 
(Definition~\ref{D:orientablevb}) for the given (co)homology theory.

\subsection*{Conventions}
\addcontentsline{toc}{subsection}{Conventions}

\subsubsection{Topological spaces}

All topological spaces are compactly generated.\comment{otherwise,
  there is ambiguity of what fibered products are!}
The category of topological spaces endowed with the Grothendieck
  topology of open coverings is denoted $\Top$.  This is the {\em site
  }of topological spaces.

\subsubsection{Manifolds}
All manifolds are second countable and Hausdorff. In particular they
are regular \Lind and paracompact.

\subsubsection{Groupoids}

We will commit the usual abuse of notation and abbreviate a groupoid
to $\Gamma_1\toto\Gamma_0$.  A {\em topological groupoid}, is a
groupoid $\Gamma_1\toto\Gamma_0$, where $\Gamma_1$ and $\Gamma_0$ are
topological spaces, but no further assumptions is made on the source
and target maps, except continuity. A topological groupoid is a {\em
  Lie groupoid} if $\gm_1, \gm_0$ are manifolds, all the structures
maps are smooth and, in addition, the source and target maps are
subjective submersions.

\subsubsection{Stacks}

For stacks, we use the words {\em equivalent } and {\em isomorphic
}interchangeably.  We will often omit 2-isomorphisms from the
notation.  For example, we may call morphisms equal if they are
2-isomorphic.  The stack associated to a groupoid $\Gamma_1\toto
\Gamma_0$ we denote by $[\Gamma_0/\Gamma_1]$, because we think of it
as the quotient. Also if $G$ is a Lie group acting on a space $Y$, we
simply denote $[Y/G]$ the stack associated to the transformation
groupoid $Y\times G \toto Y$.

Our terminology is different from that in \cite{Noohi}. The quotient
stack $\X$ of a topological groupoid $[\Gamma_0/\Gamma_1]$ is called
a {\em topological stack} in this paper, where as in [ibid.] these
are called {\em pretopological stacks}. If the source and target map
of $[\Gamma_0/\Gamma_1]$ are local Hurewicz fibrations, then we say
that $\X$ is a {\em Hurewicz topological stack}; see Section~\ref{S:topological}.

\begin{warning}\label{W:Hurewicz}
In Sections~\ref{Loopproduct},  \ref{Frobeniusforloops}, \ref{BVstructure}, \ref{S:HCFTbig} and~\ref{S:Brane}, unless otherwise stated, the (base) stack $\XX$ will always be assumed to be a Hurewicz stack (see Section~\ref{topological}). Note that \emph{any} differentiable stack (which are our main interest) is Hurewicz (Example~\ref{E:Hurewicz}). 
\end{warning}

\subsubsection{(Co)homology}
The coefficients of our (co)homology theories will be taken in a
commutative unital ring $\kor$. All tensors products
are over $\kor$ unless otherwise specified.

\smallskip

We will write both $H(\XX)$, $H\lcom(\XX)$ for the total homology
groups $\bigoplus H_{n}(\XX)$. We use the first notation when we deal
with ungraded elements and ungraded maps, while we use the second when
where dealings with homogeneous homology classes and graded
maps. Similarly, in Section~\ref{S:Bivariant}, we use respectively the notations
$H(\XX\stackrel{f}\to \YY)$ and $H\com(\XX\stackrel{f}\to \YY)$
for the total bivariant cohomology groups when we want to deal with
ungraded maps or with graded ones.

\subsection*{Acknowledgements}
\addcontentsline{toc}{subsection}{Acknowledgements}
The authors warmly thank Gustavo Granja and Andrew Kresch for helpful and inspiring discussion on the topological issues of this paper. The authors also thank Eckhard Meinrenken for his suggestions on the Cartan model.

\vfill
\eject

%%%%%%%%%%%%%%%%%%%%%%%%%%%%%%%%%%%%%%%%%%%%%%%%%%%%%%%%%%%%%%%%%
%%%%%%%%%%%%%%%%%%%%%%%%%%%%%%%%%%%%%%%%%%%%%%%%%%%%%%%%%%%%%%%%%

%

% -------------------------------------------------------
\section{Topological stacks}{\label{S:Topst}}

We review some basic facts about topological stacks. More details
can be found in \cite{Noohi}.

\subsection{Stacks over $\Top$}\label{TopStacks}

Throughout these notes, by a stack we mean a stack over the site
$\Top$ of compactly generated topological spaces with the standard
Grothendieck topology. This means, a stack is a category $\XX$
fibered in groupoids over $\Top$ satisfying the descent condition, see Appendix~\ref{S:FiberedCat} for more details.
Alternatively, we can think of $\XX$ as a presheaf of groupoids over
$\Top$ which satisfies descent.

We list some basic facts about stacks.

\begin{itemize}

\item[$\mathbf{1.}$] Stacks over $\Top$ form a 2-category in which
2-morphisms are invertible. Therefore, given two stacks $\XX$ and
$\YY$, we have the {\em groupoid} $\Hom(\YY,\XX)$ of morphisms
between them. In the case where the source  stack $\YY=T$ is a
topological space, we usually use the alternative notation $\XX(T)$
for the above hom-groupoid. This is sometimes referred to as the
groupoid of {\em $T$-valued points} of $\XX$.

Although in practice one may really be interested only in the
category of stacks which obtained by identifying 2-isomorphic
1-morphisms, the 2-category structure can not be ignored. For
example, when we talk about {\em fiber products of stacks}, we
exclusively mean the 2-fiber product in the 2-category of stacks.

\item[$\mathbf{2.}$] The 2-category of stacks has fiber products and
inner homs, so it is cartesian closed. The 2-fiber product
$\XX\times_{\ZZ}\YY$ is characterized by the property that, for
every topological space $T$, its groupoid of $T$-valued points is
given equivalent to
    $$\XX(T)\times_{\ZZ(T)}\YY(T).$$

Given stacks $\XX$ and $\YY$ be stacks over $\Top$, the inner hom
between them, called the {\em mapping stack} $\bfhom(\YY,\XX)$, is
defined by the rule
  $$ T\in \Top \ \ \   \mapsto  \ \ \ \Hom(T\times \YY,\XX).$$
Note that we have a natural equivalence of groupoids
    $$\bfhom(\YY,\XX)(*)\cong\Hom(\YY,\XX),$$
where $*$ is a point. The mapping stack has the exponential
property. That is, given stacks $\XX$, $\YY$, and $\ZZ$, we have a
natural equivalence of stacks
 $$\bfhom(\ZZ\times\YY,\XX)\cong\bfhom(\ZZ,\bfhom(\YY,\XX)).$$

\item[$\mathbf{3.}$] The category of topological spaces embeds
fully faithfully in the 2-category of stacks. This means, given two
topological spaces $X$ and $Y$, viewed as stacks via the functor
they represent, the hom-groupoid $\Hom(X,Y)$ is equivalent to a set,
and this set is in a natural bijection with the set of continuous
functions from $X$ to $Y$.

This way, we can think of a topological space as a stack.

This embedding preserves the closed cartesian structure on $\Top$.
This means that fiber products of spaces get sent to 2-fiber
products of the corresponding stacks, and the mapping spaces (with
the compact-open topology) get sent to mapping stacks.

\item[$\mathbf{4.}$] The embedding of the category of topological
spaces in the 2-category of stacks admits a left adjoint. That is,
to every stack $\XX$ one can associate a topological space, together
with a natural map $\pi \: \XX \to \XX_{mod}$ which is universal
among maps from $\XX$ to topological spaces. (That is, every map
from $\XX$ to a topological space $T$ factors uniquely through
$\pi$.) ; See \cite{Noohi}, $\S$4.3.

The space $\XX_{mod}$ is called the {\em coarse moduli space} of
$\XX$ and it should be thought of as the ``underlying space'' of
$\XX$.

In particular, the underlying set of $\X_{mod}$ is the set of
isomorphism classes of the groupoid $\X(*)$, where $*$ stands for a
point. In other words, the points in $\X_{mod}$ are the
2-isomorphism classes of points of $\X$, where by a {\em point} of
$\X$ we mean a morphism $x \: * \to \X$.

The underlying set of the coarse moduli space of the mapping stack
$\bfhom(\YY,\XX)$ is the set of 2-isomorphism classes of morphisms
from $\YY$ to $\XX$.

\item[$\mathbf{5.}$] To a point $x \: * \to \X$ of a stack $\X$
there is associated a group $I_x$, called the {\em inertia group} of
$\X$ at $x$. By definition, $I_x$ is the group of 2-isomorphisms
from the point $x$ to itself. An element in $I_x$ is sometimes
referred to as a {\em ghost} or {\em hidden loop}; see
(\cite{Noohi}, $\S$10) since its image under the map $\XX\stackrel{\pi}\to \XX_{mod}$ is constant.  

The groups $I_x$ assemble into a stack $\IXX \to \XX$ over $\XX$
called the {\em inertia stack}. The inertia stack is defined by the
following 2-fiber square
  $$\xymatrix{ \IXX \ar[d] \ar[r] & \XX \ar[d]^{\Delta} \\
      \XX \ar[r]_{\Delta} & \XX\times\XX}$$
and will be studied in more details in Section~\ref{S:hidden}.

%\item[$\mathbf{6.}$] More generally, given a topological groupoid
%$X_1\toto X_0$, to it we can associate a {\em quotient stack}
%$[X_0/X_1]$
\end{itemize}

\subsection{Morphisms of stacks}\label{morphisms}

A morphism $f \: \XX \to \YY$ of stacks is called {\bf
representable} if for every map $T \to \YY$ from a topological space
$T$, the fiber product $T\times_{\YY}\XX$ is a topological space.
This is, roughly speaking, saying that the fibers of $f$ are
topological spaces.

Any property $\mathbf{P}$ of morphisms of topological spaces which
is invariant under base change can be defined for an arbitrary
representable morphism of stacks. More precisely, we say that a
representable morphism $f \: \XX \to \YY$ is $\mathbf{P}$, if for
every map $T \to \YY$ from a topological space $T$, the base
extension $f_T \: T\times_{\YY}\XX \to T$ is $\mathbf{P}$ as a map
of topological spaces; see (\cite{Noohi}, $\S$4.1).

This way we can talk about {\em embeddings (closed, open, locally
closed, or arbitrary) of stacks, proper morphisms, finite
morphisms}, and so on.

We say that $f \: \XX \to \YY$ is an {\em epimorphism}, if it is an
epimorphism in the sheaf theoretic sense. In the case where $f$ is
representable, this is equivalent to saying that every base
extension $f_T$ of $f$ over a topological space $T$ admits  local
sections.

\subsection{Topological stacks}\label{S:topological}

A {\bf topological} stack (\cite{Noohi}, Definition 7.1) is a stack
$\XX$ over $\Top$ which admits a representable epimorphism $p \: X
\to \XX$ from a topological space $X$. Equivalently, $\XX$ is
equivalent to the quotient stack $[X_0/X_1]$  of a topological
groupoid $X_1 \toto X_0$. This quotient stack, by definition, is the
stack associated to the presheaf of groupoids
  $$ T \ \mapsto \ X_1(T) \toto X_0(T).$$
This stack is the equivalent to the stack of torsors for the
groupoid $X_1\toto X_0$;  see (\cite{Noohi}, $\S$12). The
groupoid $X_1 \toto X_0$ is recovered from the atlas
$p \: X \to \XX$ by setting $X_0:=X$ and $X_1:=X\times_{\XX}X$.
Under this correspondence between topological stacks and topological
groupoids, morphisms of stacks correspond to Hilsum-Skandalis bibundles.

An important example to keep in mind is the case of a topological
group $G$ acting on a topological space $X$.  We define the {\bf transformation 
groupoid} $X\rtimes G \toto X$ of this action  as follows. 
As the notation suggests, the space of objects is $X$ 
and the space of arrows is $X\times G$. The source map 
$s \: X\times G \to X$ is the first projection
and the target map is the action $X\times G \to X$. 
The composition of arrows is induced from the multiplication in $G$.

The quotient stack of
the transformation groupoid $X\rtimes G \toto X$ 
is denoted by $[X/G]$. For a topological space $T$, the
groupoid $[X/G](T)$ of $T$-points of $[X/G]$ is the groupoid of
pairs $(P,\varphi)$, where $P$ is a principal $G$-bundle over $T$,
and $\varphi \: P \to X$ is a $G$-equivariant map. In the case where
$X$ is a point, the corresponding quotient stack $[*/G]$ is called
the {\em classifying stack} of $G$. Its group of $T$-points is
precisely the groupoid of principal $G$-bundles over $T$.

We list some basic facts about topological stacks.

\begin{itemize}

\item[$\mathbf{1.}$] Topological stacks form a full sub 2-category
 of the 2-category of stacks over $\Top$.

\item[$\mathbf{2.}$] The 2-category of topological stacks is closed
  under fiber products.
  It, however, does not seem to have inner hom objects. That is, it
  does not seem to be the case in general that the
  mapping stack $\bfhom(\YY,\XX)$ of two topological stacks $\XX$
  and $\YY$ is a topological stack. This
  {\em is} the case, however, if $\YY$ is the quotient stack of a groupoid
  $Y_1 \toto Y_0$
  such that $Y_0$ and $Y_1$ are compact topological spaces; see
  Proposition \ref{P:MapSt}.

\item[$\mathbf{3.}$] The stack associated to a topological space $X$
is topological. It is, in fact, equivalent to the stack associated
to the trivial groupoid $X \toto X$. Thus, the category of
topological spaces is a full subcategory of the 2-category of
topological stacks.

\item[$\mathbf{4.}$]  Let $\XX=[X_0/X_1]$ be the quotient stack of a
topological groupoid $X_1 \toto X_0$. Then, the coarse moduli space
$\XX_{mod}$ of $\XX$ is naturally homeomorphic to the course
quotient space of the groupoid $X_1 \toto X_0$. In particular, the
coarse moduli space of the quotient stack $[X/G]$ is the orbit space
$X/G$ of the action of $G$ on $X$. The coarse moduli space of the
classifying stack $[*/G]$ of $G$ is just a single point.

\item[$\mathbf{5.}$] For a point $x \: * \to \X$ of a topological stack
$\X$, the inertia  group $I_x$ is naturally a topological group. The
inertia stack $\IXX$ is a topological stack, and the natural map
$\IXX \to \XX$ is representable.

\item[$\mathbf{6.}$] Every morphism $T \to \XX$ from a topological
space $T$ to a topological stack $\XX$ is representable.
\end{itemize}

\subsection{Substacks of a topological stack}
   Let $\XX$ be a topological stack. A representable morphism $i \: \YY \to \XX$ is  
   called an {\bf embedding} if for every map $T \to \XX$, with $T$ a topological space,
   the base extension $i_T \: T\times_{\XX} \YY \to T$ is an embedding of topological spaces
   (that is, $i_T$ maps $T\times_{\XX} \YY$ homeomorphically onto a subspace of $T$). 
   In this case, we say that $\YY$ is a {\bf substack} of   $\XX$. We can similarly define
   {\bf open}, {\bf closed}, and {\bf locally closed} substacks. With a slight abuse of
   notation, we often use the notation $\YY \subseteq \XX$ for a substack.
   
   Let $p \: X \to \XX$ be an atlas for the topological stack $\XX$, and let $X_1\toto X$ be the
   corresponding groupoid presentation. Then, taking inverse image via $p$
   induces a bijection between substacks 
   $\YY \subseteq \XX$ and invariant subspaces $Y \subset X$. Under this bijection
   open (respectively, closed, locally closed) substacks of $\XX$ correspond to 
   open (respectively, closed, locally closed) subspaces of $X$.
   
   Given a family of substacks $\YY_{\alpha}$ of $\XX$, we define their {\bf intersection} 
   $\underset{\alpha}\bigcap \YY_{\alpha}$ to be the largest substack of $\XX$ which 
   is a substack of all of $\YY_{\alpha}$. The {\bf union}  $\underset{\alpha}\bigcup \YY_{\alpha}$
   is defined similarly. The intersection and union of substacks always exist. In fact, if
   $p \: X \to \XX$ is an atlas for $\XX$, then the intersection 
   $\underset{\alpha}\bigcap \YY_{\alpha}$ corresponds to the invariant subspace
   $\underset{\alpha}\bigcap Y_{\alpha}$ of $X$. The same goes with the union.
   
   Given a substack $\YY$ of $\XX$, we define its {\bf closure} $\bar{\YY}$  to be the smallest 
   closed substack containing $\XX$. The {\bf interior}  $\YY^{\circ}$ is defined to be the largest
   open substack of $\XX$ contained in $\YY$. The {\bf complement} $\YY^c$ of a substack
   $\YY \subseteq \XX$ is the largest substack of $\XX$ whose intersection with
   $\YY$ is empty. Given two substacks $\YY$ and $\ZZ$ of $\XX$, we define the
   {\bf difference} $\YY-\ZZ$ to be the substack $\YY\cap\ZZ^c$.
   All these exist, are well-defined, and can be constructed by taking the corresponding invariant subspaces
   of an atlas $p \: X \to \XX$.

\subsection{Hurewicz topological stacks}\label{topological}

As we will see in $\S$\ref{SS:push}, in order to have nice gluing
properties for maps into a stack $\X$, we need to assume $\XX$ is a
Hurewicz stack. This will be needed later on when we work with loop
stacks. We recall the definition of a Hurewicz stack.

A {\em Hurewicz fibration} is a map having the homotopy lifting
property for all topological spaces. A map $f \: X \to Y$ of
topological spaces is a {\em local Hurewicz fibration} if for every
$x \in X$ there are opens $x \in U$ and $f(x) \in V$ such that
$f(U)\subseteq V$ and $f|_U \to V$ is a Hurewicz fibration.  The
most important example for us is the case of a topological
submersion: a map $f\:X\to Y$, such that locally $U$ is homeomorphic
to $V\times \rr^n$, for some $n$.

Dually, we have the notion of {\em local cofibration}. It is known
(\cite{StromI}), that if $A\to Z$ is a closed embedding of
topological spaces, it is a local cofibration if and only if there
exists and open neighborhood $A\subset U\subset Z$ such that $A$ is
a strong deformation retract of $U$. If $A\to Z$ is a local
cofibration, so is $A\times T\to  Z\times T$ for every topological
space $T$. Moreover, the following result is essential for our
purposes (\cite{StromII}):

Given a commutative diagram, with $A\to Z$ a local cofibration and
$X\to Y$ a local fibration
$$\xymatrix{
A\dto\rto &X\dto\\ Z\rto & Y}$$ then for every point $a\in A$ there
exists an open neighborhood $Z'$ of $a$ in $Z$, such that there
exists a lifting (the dotted arrow) giving two commutative triangles
$$\xymatrix{
A'\dto\rto & X\dto\\ Z'\rto\ar@{.>}[ur] & Y}$$ where $A'=A\cap Z'$.

\begin{defn}\label{def:topological}
   A topological stack $\X$ is called {\bf Hurewicz} if it
   is equivalent to the quotient stack $[X_0/X_1]$ of a topological
   groupoid $X_1\toto X_0$ whose source and target maps are local
   Hurewicz fibrations.
\end{defn}

\begin{numex}\label{E:Hurewicz}
A topological space is a Hurewicz topological stack.  Every substack of
a Hurewicz topological stack is a Hurewicz 
topological stack. The
topological stack underlying any differentiable stack is a Hurewicz
topological stack. In particular, any global quotient $[M/G]$ of a
manifold by a Lie group is a Hurewicz topological stack.
\end{numex}

\subsection{Pushouts in the category of stacks}{\label{SS:push}}

The following generalizes (\oldcite{Noohi}, Theorem
16.2).

\begin{prop}{\label{P:glue}}
 Let $A\to Y$ be a closed embedding of Hausdorff spaces, which is a
 local cofibration.
 Let $A\to Z$ be a finite proper map of Hausdorff spaces. Suppose we
 are given a pushout
 diagram in the category of topological spaces
     $$\xymatrix{
               A \ar@{^(->} [r] \ar@{->}[d] &  \ar[d] Y  \\
                             Z \ar[r]        &    Z\vee_A Y       }$$
Then this   diagram remains a pushout diagram in the 2-category of
Hurewicz topological stacks. In other words,
 for every Hurewicz topological stack $\X$, the morphism
$$\X(Z\vee_A Y)\longrightarrow \X(Z)\times_{\X(A)} \X(Y)$$
is an equivalence of groupoids.
\end{prop}
\begin{pf}
Let us abbreviate the pushout by $U=Z\vee_A Y$.\comment{pushout or
puhs-out?}

The fully faithful property only uses that $\X$ is a topological
stack and that  $U$ is a pushout.  Let us concentrate on essential
surjectivity. Because $\X$ is a stack and we already proved full
faithfulness, the question is local in $U$. Assume given $Z\to \X$
and $Y\to \X$, and an isomorphism over $A$. Let $X_1\toto X_0$ be a
groupoid presenting $\X$, whose source and target maps are local
fibrations.

Let us remark that both $Z\to U$ and $Y\to U$ are finite proper maps
of Hausdorff spaces.  Thus we can cover $U$ by open subsets $U_i$,
such that for every $i$, both $Z_i=U_i\cap Z$ and $Y_i=Y\cap U_i$
admit liftings to $X_0$ of their morphisms to $\X$.  We thus reduce
to the case that we have $Z\to X_0$, $Y\to X_0$, and $A\to X_1$.
Next, we need to construct the dotted arrow in
$$\xymatrix{
A\rto\dto   & Y\dto\ar@{.>}[dl]\\ X_1\rto  & X_0}$$ We can cover $Y$
by opens over which this arrow exists, because $A\to Y$ is a local
cofibration and $X_1\to X_0$ a local fibration. Then for a point
$u\in U$ we choose an open neighborhood in $U$ small enough such
that the preimage in $Y$ is a disjoint union of sets over with the
dotted arrow exists. Passing to such a neighborhood of $u$ reduces
to the case that the dotted arrow exists.  Then there is nothing
left to prove.
\end{pf}

%----------------------
\subsection{Orbifolds as topological stacks}{\label{SS:orbifolds}}
  The most familiar examples of topological stacks are the orbifolds.
 An orbifold, by definition, is a topological stack which can be covered by
 open substacks of the form $[X/G]$, with $G$  a finite group.
 Any orbifold is the quotient stack of an \'etale groupoid $[X_1\sst{} X_0]$.
 Recall that being \'{e}tale means that 
 the source (hence also the target) $X_1 \to X_0$ is a local homeomorphism --
 in particular, an orbifold is a Hurewicz topological stack.
 Moreover, it can be shown that the diagonal map $X_1 \toto X_0\times X_0$
 is a closed map (with finite fibers). In fact, the converse is also true in the locally connected case. 
 Namely, the quotient stack of a locally connected \'{e}tale groupoid whose diagonal map is closed with finite
 fibers is an orbifold (see \cite{Noohi}, Propostion 14.9).
 
 We should point out that there is some inconsistency in the literature about
 terminology: in the definition of orbifold, some authors assume that the action of $G$ on $X$ 
 is generically free. For this reason, and by analogy to their algebraic geometric counterparts, 
 in {\em loc. cit.} the term {\em Deligne-Mumford}
 has been used instead of orbifold. The orbifolds for which the above generic freeness
 condition is satisfied are sometimes called reduced orbifolds.

Although every orbifold $\XX$ is 
locally the quotient stack $[X/G]$ of a finite group action, this may not
be the case globally, i.e., $\XX$ may not be {\em good}. (For
a charcterization of good orbifolds
in terms of their fundamental group see \cite{Noohi}, Theorem 18.24.)
It is known, however, that a reduced differentiable orbifold $\XX$ can 
always be globally written as a quotient 
stack $[X/G]$, where $G$ is a Lie group acting with finite stabilizers on a 
manifold $X$. This is not known to be the case for general orbifolds though.
 
Orbifolds clearly form a small subclass of all topological stacks. For instance,
every point on an orbifold has finite stabilizer group, and this is not true for an arbitrary topological
stack. The simplest example of a topological stack which is not an orbifold is the quotient stack $[*/G]$,
where $G$ is any topological group which is not finite.

% -----------------------------------------
\subsection{Geometric stacks}

I this paper we will encounter other types of stacks as well.  A {\bf differentiable
stack} is a stack on the category of $C^\infty$-manifolds, which is
isomorphic  to the quotient stack of a Lie groupoid.  Every
differentiable stack has an underlying topological stack that
is Hurewicz.  If the
Lie groupoid $X_1\toto X_0$ represents the differentiable stack
$\XX$, the underlying topological groupoid represents the underlying
topological stack.  Often we will tacitly pass from a differentiable
stack to its underlying topological stack. For more on
differentiable stacks, see~\cite{BX}. 

%Every differentiable orbifold
%is a differentiable stack. The class of differentiable  
%stacks is, however, considerably larger. A typical example is the 
%quotient stack $[*/G]$, where $G$ is Lie group acting trivially on  
%a point. This is a differentiable stack which is an orbifold
%if and only if $G$ is finite. 

An {\bf almost complex stack }is a stack on the category of almost
complex manifolds, which is isomorphic to the quotient stack of an
almost complex Lie groupoid, i.e., a Lie groupoid $X_1\toto X_0$,
where $X_0$ and $X_1$ are almost complex manifolds, and all
structure maps respect the almost complex structure.  Every almost
complex stack has an underlying differentiable stack and hence also
an underlying topological stack.

% -------------------------------------------------------
\section{Homotopy type of a topological stack}{\label{S:Htpytypes}}

% --------------------------------
\subsection{Classifying space of a topological groupoid} 

We recall the construction of the (Haefliger-Milnor)
classifying space $B\bbX$ and the {\bf universal bundle} $E\bbX$
of a topological groupoid $\bbX=[X_1 \sst{} X_0]$ from \cite{No/cohomology}. 

An element in $E\bbX$ is a sequence
$(t_0\al_0,t_1\al_1,\cdots,t_n\al_n,\cdots)$, where $\al_i \in R$
are such that $s(\al_i)$ are equal to each other,
 and $t_i \in [0,1]$ are such that all but finitely
many of them are zero and $\sum t_i=1$. As the notation suggests, we
set $(t_0\al_0,t_1\al_1,\cdots,t_n\al_n,\cdots)=(t'_0\al'_0,t'_1\al'_1,\cdots,t'_n\al'_n,\cdots)$
if $t_i=t'_i$ for all $i$ and $\al_i=\al'_i$ if $t_i\neq 0$.

Let $t_i \: E\bbX \to [0,1]$ denote the map
$(t_0\al_0,t_1\al_1,\cdots,t_n\al_n,\cdots) \mapsto t_i$, and
let $\al_i \: t_i^{-1}(0,1] \to R$ denote the
map $(t_0\al_0,t_1\al_1,\cdots,t_n\al_n,\cdots) \mapsto \al_i$. The
topology on $E\bbX$ is the weakest topology in which $t_i^{-1}(0,1]$ are all
open and $t_i$ and $\al_i$ are all continuous.

The classifying space $B\bbX$ is defined to be the quotient of
$E\bbX$ under the following equivalence relation:  two
elements $(t_0\al_0,t_1\al_1,\cdots,t_n\al_n,\cdots)$ and
$(t'_0\al'_0,t'_1\al'_1,\cdots,t'_n\al'_n,\cdots)$ of $E\bbX$ are
equivalent, if $t_i=t'_i$  for all $i$, and if there is an element
$\gamma \in X_1$ such that $\al_i=\gamma\al'_i$. (So, in particular,
$t(\al_i)=t(\al'_i)$  for all $i$.)

% --------------------------------
\subsection{Classifying space of a topological stack}{\label{SS:classifying}}

Many facts about  topological stacks can be reduced to the case of
topological spaces by virtue of the following.

\begin{thm}{\label{T:classifying}}
  For every topological stack $\XX$, there exists a topological space
  $X$ together with a morphism $\varphi \: X \to \XX$ which has the
  property that, for every morphism $T\to \XX$ from a
  topological space $T$, the pullback $T\times_\XX X\to T$ is a
  weak homotopy equivalence.
\end{thm}

A topological space $X$ with the above property is called a {\bf
classifying space }for $\XX$. A classifying space for $\XX$ can be
constructed by taking the classifying space $B\bbX$ of a groupoid
$\bbX=X_1\toto X_0$ whose quotient stack is $\XX$ 
(see \cite{No/cohomology}, Theorem 6.3). 
The above theorem implies that the
classifying space of a topological stack
is unique up to a unique isomorphism in the weak
homotopy category of topological spaces (i.e., the category of
topological spaces with weak homotopy equivalences inverted).

In the case where $\XX=[X/G]$ is the quotient stack of a group
action, the Borel construction $X\times_G{EG}$ is a classifying
space for $\XX$. Here $EG$ is the total space of the universal
principal $G$-bundle in the sense of Milnor.

% ---------------------------------
\subsection{Paracompactness of the classifying space}

In many applications, it is important to find a classifying space
for $\XX$ which is {\em paracompact}. There are various conditions
on a groupoid $X_1 \toto X_0$ which guarantee that the fat
realization of the nerve of $X_1 \toto X_0$ is paracompact. The
following is one.

\begin{defn}\label{def:Lindelof}
   A topological stack $\X$ is called {\bf regular \Lind} if it
   is equivalent to the quotient stack $[X_0/X_1]$ of a topological
   groupoid $X_1\toto X_0$ such that $X_1, X_0$ are regular \Lind spaces.
\end{defn}

The proof of the following proposition will appear elsewhere.

\begin{prop}\label{P:Lindelof}
If $\X$ is a regular \Lind stack, there exists a classifying space
for $\X$ which is a regular \Lind space, in particular paracompact.
\end{prop}
\begin{numrmk}
Every differentiable stack is regular \Lind and hence has a
paracompact classifying space.
\end{numrmk}

% ---------------------------------
\subsection{(Co)homology theories for topological stacks}

Theorem \ref{T:classifying} allows one to extend every (generalized)
(co)homology theory $h$ to the 2-category of topological stacks. For
instance, let us explain how to define $h(\XX,\AA)$ for a pair
$(\XX,\AA)$ of topological stacks.

Choose a classifying space $\varphi \: X \to \XX$, and let
$A:=\varphi^{-1}\A$. It follows that the pair $(X,A)$ is
well-defined in the weak homotopy category of pairs (i.e., is
independent of the choice of a particular classifying space $X$).
So, we can define $h(\X,\A)$ to be $h(X,A)$. It can be easily
verified that this construction is functorial in morphisms of pairs.

The  cohomology theory thus defined on topological stacks will
maintain all natural properties that it had on spaces. For example,
it will be homotopy invariant (in particular, it will not
distinguish 2-isomorphic morphisms), it will satisfy excision, it
will maintain all the products (cap, cup, etc.) that it had on
spaces, and so on.

In the case where $\XX=[X/G]$, we will recover the usual
$G$-equivariant (co)homology of $X$ defined using the Borel
construction. That is, $h([X/G])\cong h(X\times_G{EG})$.

Every (co)homology theory for topological stacks which is invariant under weak equivalences 
is induced from one on topological spaces. This is due to existence of a classifying space
$\varphi \: X \to \XX$ (Theorem \ref{T:classifying}) which forces the (co)homology of $\XX$ to be equal to
that of its classifying space $X$.

% ---------------------------------
\subsection{Eilenberg-Steenrod axioms for topological stacks}

We recall the Eilenberg-Steenrod axioms for a homology theory, formulated in the context  of topological 
stacks. Let $H_n$ be a sequence of functors
from the category of pairs $(\XX,\AA)$ of topological stacks  to the 
category of abelian groups. (By a pair $(\XX,\AA)$ we mean a topological stack $\XX$
and a substack $\AA$.) This sequence is
equipped with  natural transformations $\partial \: H_i(\XX,\AA) \to H_{i-1}(\AA)$, called 
the boundary maps. The Eilenberg-Steenrod axioms are the following:

\begin{itemize}
  \item[$\mathbf{1.}$] {\em Homotopy.}  If $f,g \: (\XX,\AA) \to (\YY,\BB)$ are homotopic as morphisms
    of pairs of stacks (in the sense of \cite{Noohi}, Definition 17.2), then they induce the
    same map on $H_n$ for all $n$.
    
  \item[$\mathbf{2.}$] {\em Excision.} Let $(\XX,\AA)$ be a pair of topological stacks.
    Let $\UU$ be a substack of $\XX$ such that the closure of $\UU$ is contained in the interior of $\AA$. 
    Then, the inclusion map $(\XX-\UU,\AA-\UU) \hookrightarrow (\XX,\AA)$  induces an isomorphism in homology.
  
  \item[$\mathbf{3.}$] {\em Dimension.} $H_n(pt)=0$ for all $n\neq 0$.
  \item[$\mathbf{4.}$] {\em Additivity.}  For any collection $\{\XX_{\alpha}\}$ of
          topological stacks, $H_n(\underset{\alpha}\coprod \XX_{\alpha})\cong\underset{\alpha}\bigoplus H_n(\XX_{\alpha})$
  \item[$\mathbf{5.}$] {\em Exactness.} For every pair  $(\XX,\AA)$ of topological stacks, the maps $i \: \AA \to \XX$ and
      $j \: (\XX,\emptyset) \to (\XX,\AA)$ induce a long exact sequence 
        $$\cdots \longrightarrow H_n(\AA) \stackrel{i_*}\longrightarrow H_n(\XX) \stackrel{j_*}\longrightarrow 
            H_n(\XX,\AA) \stackrel{\partial}\longrightarrow H_{n-1}(\AA) \longrightarrow \cdots.$$
\end{itemize}
 In the case of singular homology with coefficients in an abelian group $A$, we have $H_n=0$ for all $n<0$ and 
 $H_0(pt)=A$.

\subsection{Singular homology and cohomology}

We will fix once and for all a coefficient ring and drop it from the
notation consistently.
We assume all topological stacks are regular Lindel\"of.

Singular homology and cohomology for spaces  lifts to topological
stacks.  The singular (co)homology of the topological stack $\XX$
can be defined to be  the singular (co)homology of its classifying
space. Alternatively,
 let $\gm :
\gm_1\toto \gm_0$ be a topological groupoid presentation of $\XX$.
We denote
$\Gamma_p=\Gamma_1\times_{\Gamma_0}\ldots\times_{\Gamma_0}\Gamma_1$
($p$-fold)  the space of composable sequences of $p$ arrows in the
groupoid $\gm$.  It yields  a simplicial space $\gm \lcom$
\begin{equation}\label{sim.ma}
\xymatrix{ \ldots \gm_2 \ar[r]\ar@<1ex>[r]\ar@<-1ex>[r] &
\gm_1\ar@<-.5ex>[r]\ar@<.5ex>[r] &\gm_0\,.}
\end{equation}

The {\em singular chain complex } of $\Gamma\lcom$ is the total
complex associated to the double complex
$C\lcom(\Gamma\lcom)$~\cite{Behrend}. Here $C_q(\Gamma_p)$ is the
linear space generated by the continuous maps $\Delta_q\to\Gamma_p$.
Its homology groups
$H_q(\Gamma\lcom)=H_q\big(C\lcom(\Gamma\lcom)\big)$ are called the
{\em homology groups } of $\Gamma$.  The {\em singular cochain
complex } of $\Gamma\lcom$ is the dual of $C\lcom(\Gamma\lcom)$. In
other words the total complex associated to the bicomplex
$C^p(\gm_q)$. These groups are Morita invariant (i.e. only depend on
the quotient stack  $[\Gamma_0/\Gamma_1]$). By definition they are
the (co)homology groups of the stack $[\Gamma_0/\Gamma_1]$. They
coincide with the cohomology of the classifying space of $\X$.

This Theory generalizes  to singular cohomology theory $H(\X,\A)$  
for pairs $(\X,\A)$ of topological stacks satisfying the
Eilenberg-Steenrod axioms. This theory comes with cup products and
coincided with the usual singular cohomology when $(\X,\A)$ is a
pair of topological spaces. In the case when $\X=[X/G]$ is the
quotient stack of a topological group action, $H$ is $G$-equivariant
cohomology. 

The K\"unneth formula also holds for singular (co)homology of topological stacks. 
We only formulate the cohomology version but the homology one holds as well (with the same proof). 

\begin{prop}[K\"unneth formula]{\label{P:Kunneth}}
In the case of field coefficients, singular cohomology of topological stacks satisfies K\"unneth formula. 
That is, we have an isomorphism of graded groups
  $$H^{\bullet}(\X,\A)\otimes H^{\bullet}(\Y,\B)\cong H^{\bullet}(\X\times\Y,\X\times\B\cup\A\times\Y).$$
\end{prop}

\begin{proof}
  Like other properties of singular cohomology, 
  this is proved by choosing classifying spaces $X \to \X$ and $Y \to \Y$ and pulling back everything along
  $X\times Y \to \X\times\Y$.
\end{proof}

\begin{numrmk}\label{R:cross}
 When the coefficient is only a ring, there are still natural cross-product homomorphisms
 $H^{\bullet}(\X,\A)\otimes H^{\bullet}(\Y,\B)\to H^{\bullet}(\X\times\Y,\X\times\B\cup\A\times\Y)$ in cohomology and in homology as well
 $H_{\bullet}(\X,\A)\otimes H_{\bullet}(\Y,\B)\to H_{\bullet}(\X\times\Y,\X\times\B\cup\A\times\Y)$ (the later being further a monomorphism). As for the proof of Proposition~\ref{P:Kunneth}, this can be seen by choosing classifying space or, alternatively, by working directly with the singular (co)chain complexes of associated groupoid presentations.
\end{numrmk}

\begin{prop}{\label{P:triple}}
  Let $\X$ be a topological stack and $\A,\B \subseteq \X$
  substacks. Then, we have a cohomology long exact sequence
  {\footnotesize $$ \cdots \to H^{n-1}(\A,\A\cap\B) \to
     H^n(\X,\A\cup\B) \to H^n(\X,\B) \to
    H^n(\A,\A\cap\B) \to H^{n+1}(\X,\A\cup\B) \cdots .$$}
\end{prop}

\begin{proof}
   By Excision %see Rotman
   $H^n(\A,\A\cap\B)\cong H^n(\A\cup\B,\B)$. The result follows from
   the long exact cohomology sequence for the triple
   $(\X,\A\cup\B,\B)$.
\end{proof}

The following less standard fact is also true about
$H$.

\begin{prop}{\label{P:dot}}
Let $\A\hookrightarrow \B \hookrightarrow \X$ be closed embeddings
of topological stacks. Then, there is a natural product
   $$H^m(\X,\X-\B)\otimes H^n(\B,\B-\A)\to H^{m+n}(\X,\X-\A)$$ which
   coincides with the cup product if $\B=\X$.
\end{prop}

\begin{proof}
 One uses the fact the classifying space is paracompact
(Proposition~\ref{P:Lindelof}). It is a general fact (for instance
see~\cite{KaSc}) that if $F$ is a sheaf over a paracompact space $X$
and $Z\subset X$ is closed, then $\injectlim_{U\supset
Z}\Gamma(U,F) \stackrel{\sim}\longrightarrow \Gamma(Z,F)$, where $U$
is open. Then the result follows from the same argument as for
topological spaces in \cite{FuMac}, $\S$\;3.
\end{proof}

% -------------------------------------------------------
\section{Vector bundles on stacks}{\label{S:VB}}

We begin with the definition of a (representable) vector bundle on a
stack.

\begin{defn}{\label{D:vb}}
 Let $\XX$ be a  (topological) stack. A real {\bf vector bundle} on $\XX$ is a
 representable morphism of stacks $\EE \to \XX$ which makes $\EE$ a
 vector space object relative to $\XX$. That is, we have an addition
 morphism $\EE\times_{\XX} \EE \to \EE$ and an $\bbR$-action
 $\bbR\times\EE \to \EE$, both relative to $\X$, which satisfy the
 usual axioms. A complex vector bundle is defined analogously.
\end{defn}

A linear map between two vector bundles is defined in the obvious
way. Vector bundles on $\X$ and linear maps between them form a
category. (Notice that since a vector bundle $\EE$ is representable
over its base $\X$, we only get a category and not a 2-category.)

There are  alternative ways of defining vector bundles over stack a
$\X$ as we will see in the next proposition. All definitions are
equivalent to the one given above.

\begin{prop}{\label{P:vb}}
The following three definitions for a vector bundle on a stack $\XX$
are equivalent to the one given in Definition \ref{D:vb}, in the
sense that the corresponding categories of vector bundles over a
given stack $\X$ are naturally equivalent as linear categories:
\begin{itemize}
 \item[$\mathbf{1.}$]
    A vector bundle on $\XX$ is a representable morphism  of stacks
    $\EE
    \to \XX$ such that, for every $f \: U \to \XX$ with $U$ a topological
    space, the pullback $E_U \to U$ is endowed with the
    structure of a vector bundle. Here, $E_U:=f^*\E=U\times_{\XX}\EE$.
   $$\xymatrix@=16pt@M=8pt{   E_U \ar[r] \ar[d] & \EE \ar[d] \\
                              U \ar[r]_{f} & \XX} $$

    We also require, for every $a \: V \to U$, that the natural
    isomorphism
        $$\varphi_a \: (f\circ a)^*\EE \to a^*(f^*\EE)$$
    be a bundle map.

 \item[$\mathbf{2.}$] A vector bundle on $\XX$ is  assignment of a
   vector bundle $E_U \to U$ for every morphism $f \: U \to \X$, with
   $U$
   a topological space, together with isomorphisms $\varphi_a \:
   a^*E_U \to E_V$ of vector bundles, for every  2-commutative triangle
     $$\xymatrix@C=0pt@R=6pt@M=6pt{   V \ar[rrrr]^a \ar[ddrr]_g & &
                                          & & U \ar[lldd]^f \\
                              & \ar@{=>}_{\alpha} [ru] && & \\
                                                     & & \XX & } $$
   We require the isomorphisms $\varphi$ to satisfy the cocycle
   condition
       $$\varphi_{a\circ b}=\varphi_b\circ b^*(\varphi_a)$$
   for
   every pair of composable triangles $a$ and $b$. (Note the abuse
   of notation: the vector bundle $E_U$ also   depends on $f$, and
   the isomorphism $\varphi_a$ also depends on the 2-morphism
   $\alpha$.)

 \item[$\mathbf{3.}$] Let $\mathbb{X}=[s,t \: X_1\sst{}X_0]$
    be a groupoid presentation for $\XX$. Then, a vector bundle on
    $\XX$ is  an $\mathbb{X}$-equivariant vector
    bundle. Recall that an $\mathbb{X}$-equivariant vector bundle
    consists of a vector bundle $E$ over $X_0$, and an isomorphism
    $\psi \: s^*E\to t^*E$ of vector bundles over $X_1$ such that
    the three restrictions of $\psi$ to $X_1\times_{X_0}X_1$
    satisfy the cocycle condition.
\end{itemize}
\end{prop}

\begin{proof}
We briefly explain how to go from one definition to the other.

Let $\EE \to \XX$ be a vector bundle in the sense of Definition
\ref{D:vb}. It is clear that the pullback vector bundles $E_U$
satisfy the conditions of ($\mathbf{1}$).

To go from ($\mathbf{1}$) to ($\mathbf{2}$) is obvious.

Given a vector bundle in the sense of ($\mathbf{2}$), we obtain a
vector bundle $E_{X_0}$ on $X_0$ corresponding to the quotient map
$p \: X_0 \to \X$. It follows from the cocycle condition on
($\mathbf{2}$) that this is an $\mathbb{X}$-equivariant bundle.

Finally, given an $\mathbb{X}$-equivariant vector bundle $E$, we
define $\E$ to be the quotient stack of the groupoid
$[E_1\sst{}E_0]$, where $E_0:=E$ and $E_1:=s^*E=X_1\times_{X_0}E$.
The source map $E_1 \to E_0$ is the projection map $\pr_2
\:X_1\times_{X_0}E \to E_0$. The target map is $\pr_2\circ \psi$. It
is easy to verify that $\EE$ is a vector bundle over $\XX$ in the
sense of Definition \ref{D:vb}.
\end{proof}

% -------------------------------------------------------
\subsection{Operations on vector bundles}{\label{SS:operations}}

The standard  operations on vector bundles on spaces (e.g., direct
sum, tensor product, exterior powers, and so on) can be carried out
on vector bundles on stacks {\em mutatis mutandis}. This is more
easily seen if we think of a vector bundle as in Proposition
\ref{P:vb} ($\mathbf{2}$). In this case, we simply perform the
desired operation simultaneously on the $E_U$, for varying $U$, and
the resulting family of vector bundles, say $F_U$, will give rise to
a vector bundle $\FF$ on $\XX$.

In view of Proposition \ref{P:vb} ($\mathbf{3}$), operations on
vector bundles on $\X$ correspond to operations on
$\mathbb{X}$-equivariant vector bundles.

Similarly, we can define a {\bf metric} on a vector bundle. More
precisely, a metric on $\EE$ is the same thing as a compatible
family of metrics on $E_U$, for varying $U$. Given a presentation
$\mathbb{X}=[X_1\sst{}X_0]$ for $\X$, a metric on $\EE$ is the same
thing as an invariant metric on the $\mathbb{X}$-equivariant vector
bundle $E$. (The latter simply means a metric on the vector bundle
$E$ over $X_0$ such that the isomorphism $\psi \: s^*E \to t^*E$ is
an isometry.)

\begin{ex}{\label{E:metric}}
  Let $X$ be a paracompact topological space (say, a manifold) and
  $G$ a compact Lie group acting on it. Set $\XX:=[X/G]$. Then every
  vector bundle $\EE$ on $\XX$ admits a metric. In fact, metrics on
  $\E$ are in bijection with $G$-invariant metrics on the vector
  bundle $E:=p^*\E$ over $X$. (Here $p \: X \to \XX$ is the
  quotient map.)
\end{ex}

% -----------------------------------------
\subsection{Tangent, normal, and excess bundles}{\label{SS:tangent}}

The main examples of vector bundles we encounter in this paper are tangent 
and normal bundles. In this section we explain how they are defined. 

%---------------
\subsubsection{Tangent bundle}

Let $\XX$ be a differentiable stack and choose
a differentiable groupoid $[X_1\toto X_0]$ presenting it. Taking  
tangent bundles gives rise to a new differentiable groupoid $[TX_1\toto TX_0]$.
The quotient stack $[TX_0/TX_1]$ is denoted by $\TT\XX$ and 
is  the {\em tangent bundle} of $\XX$. 
The base maps induce a groupoid morphism
 $[TX_1\toto TX_0] \to [X_1\toto X_0]$. After passing to quotient stacks,
this induces a morphism of stacks $\TT\XX \to \XX$ which we regard as
the base map of $\TT\XX$.
It is not hard to see that, up to isomorphism of stacks over $\XX$,
$\TT\XX$ is independent of the choice of the groupoid
presentation. The tangent bundle $\TT\XX$ is functorial in the obvious sense.
 
\begin{ex}{\label{E:tgtbundle}}
  Let $\X$ be a differentiable orbifold. Choose a presentation for it by a smooth
  \'etale groupoid  $\mathbb{X}=[s,t\: X_1\sst{}X_0]$. The tangent bundle
  $TX_0$ of $X_0$ is naturally $\mathbb{X}$-equivariant because the two
  pullbacks $s^*(TX_0)$ and $t^*(TX_0)$ are both naturally
  isomorphic to $TX_1$. The corresponding vector bundle on $\X$ is
  the tangent bundle of $\X$. 
\end{ex} 
 
\begin{warning}The above example shows that, when $\XX$ is a differentiable orbifold,
the tangent bundle $\TT\XX$ is indeed a vector bundle in the sense of
Definition \ref{D:vb}. 
This, however, is not the case for arbitrary differentiable stacks.
This is seen by observing that the fiber $\TT_x\XX$ of the map $\TT\XX \to \XX$
over a point $x$ in $\XX$ is not  a vector space in general. In fact, 
the map $\TT\XX \to \XX$ may not even be representable.

\medskip

\noindent{\bf Notation:} when $\XX$ is
an orbifold we use the notation $T\X$ for the tangent bundle. 
\end{warning}

\begin{ex}{\label{E:LieGp}}
  Let $G$ be a Lie group, and let $\XX=BG=[*/G]$ be its classifying stack.
  Let $\mathfrak{g}$ be the Lie algebra of $G$. Then,
  $\TT BG=B(\mathfrak{g}\rtimes G)$, where $G$ is acting on $\mathfrak{g}$
  by the adjoint action. The fiber of the base map
  $\TT BG \to BG$ over the point $* \to BG$ is $B\mathfrak{g}=[*/\mathfrak{g}]$.
  (Here, $\mathfrak{g}$ is regarded as a group via its vector space addition.
  The stack  $B\mathfrak{g}=[*/\mathfrak{g}]$ is a
  simple example of a ``2-vector space.'')
\end{ex}

As indicated in Example \ref{E:LieGp},  the
definition of 
the tangent bundle for a general differentiable stack $\XX$
requires a more general 
class of  vector bundles which are quite a bit subtler.
These are what some authors call `stacky
vector bundles'  or `2-vector bundles'
and have the property that their fibers are, in general, 2-vector spaces. In particular,
the structure map $\EE \to \XX$ of a stacky vector
bundle is no longer representable.  
  
Locally, the tangent  2-vector bundle of $\XX$ can be presented by a length 
2 complex of vector bundles
(``the tangent complex''). A suitable model for this
is the complex $E \stackrel{\rho}\to TX$, where $\rho$ is the anchor map
of the Lie algebroid  associated to a Lie groupoid presentation $X_1\toto X$
for $\XX$. Here $E$ is the normal bundle of the unit map $\eta\: X\to X_1$. 
It is naturally identified with the relative tangent bundle $T_t$ of the target 
map $t:X_1\to X$ and the anchor map $\rho$ is the composition 
$\rho: E\cong \eta^*(T_t) \hookrightarrow \eta^*(T_t) \oplus TX 
\cong \eta^*(TX_1) \cong \eta^*(T_s) \oplus TX\to TX$ where $T_s$ is the relative 
tangent bundle of the source map $s:X_1\to X$. In Example \ref{E:LieGp},
the tangent bundle $TBG$ can be represented by the complex  of vector spaces 
$\mathfrak{g} \to 0$ (viewed as a complex
a vector bundles on a point). 
 
In the rest of the paper, the only instances where we encounter tangent bundles 
are when $\XX$ is an orbifold.

%---------------
\subsubsection{Normal bundle}

Let $\YY$ be a differentiable stack and 
$\XX \hookrightarrow \YY$ a differentiable substack.
We would like to define the normal bundle of $\XX$ in $\YY$. 
When $X$ and $Y$ are smooth manifold,
one defines the normal bundle either as the quotient
$(TY|_{X})/TX$, or as the orthogonal complement to $TX$ in 
$TY|_{X}$ (upon fixing a metric on the latter). 

None of these definitions
are available to us in the context of stacks
(except when $\YY$ is an orbifold). Nevertheless, it {\em is} possible to define
the normal bundle as a vector bundle on $\XX$. To do so, pick an 
atlas $Y \to \YY$ and let $X \subset Y$ be the invariant submanifold
corresponding to $\XX$. Let $N_{X/Y}=(TY|_X)/TX$ be the normal bundle of $X$ in $Y$.
This is an equivariant vector bundle with respect to the induced groupoid structure
on $X$, hence, after passing to quotient, gives rise to a vector bundle
on $\XX$, which we denote by $\N_{\XX/\YY}$ and call the {\bf normal
bundle} to $\XX$ in $\YY$. 

\begin{ex}{\label{E:normalbundle1}}
  Let $\YY$ be a differentiable orbifold. Then, we
  have $\N_{\XX/\YY}=(T\YY|_{\XX})/T\XX$. In fact, since we can always choose
  a metric on $T\YY$ (because $\XX$ is paracompact), we have a direct sum
  decomposition $T\YY|_{\XX}=\N_{\XX/\YY}\oplus T\XX$.
\end{ex}

%---------------
\subsubsection{Excess bundle and transversality}
Consider a 2-cartesian diagram of differentiable stacks
         $$\xymatrix@=16pt@M=8pt{ \X' \ar@{^(->}[r]^{j} \ar[d]_{p} & \Y' \ar[d]^{q} \\
                          \X   \ar@{^(->} [r]_{i} & \Y    }$$
in which the horizontal morphisms are embeddings. (Note that
if $q$ is a submersion, $i$ being an embedding implies that $j$ is an embedding. 
When $q$ is representable, this is seen by pulling back everything
along an atlas $Y \to \YY$. The general case reduces to the representable case by pulling back $j$ along an atlas
$Y' \to \Y'$.) 
The bundle $\N_{\XX'/\YY'}$  is naturally a subbundle of 
$p^*\N_{\XX/\YY}$. (This is seen using the same pullback argument we
just gave to prove that $j$ is an embedding.) We call the quotient bundle  
$\E:=(p^*\N_{\XX/\YY})/(\N_{\XX'/\YY'})$, which is
a vector bundle on $\XX'$, the {\bf excess normal bundle} of the diagram.  We say that $q$ is
{\bf transversal} to $f$ if $\E$ is trivial. 

\begin{ex}{\label{E:normalbundle2}}
  Let $i \: \XX \hookrightarrow \YY$ and $j \: \YY \hookrightarrow \ZZ$ be embeddings of differentiable stacks. Consider
  the 2-cartesian diagram 
              $$\xymatrix@=16pt@M=8pt{ \XX \ar@{^(->}[r]^{\id} \ar@{^(->} [d]_{i} & \XX \ar@{^(->} [d]^{j\circ i} \\
                          \YY  \ar@{^(->} [r]_{j} & \ZZ    }$$
  The excess normal bundle for this diagram is   $i^*\N_{\YY/\ZZ}$. The excess normal bundle for the transpose diagram
             $$\xymatrix@=16pt@M=8pt{ \XX \ar@{^(->}[r]^{i} \ar@{^(->} [d]_{\id} & \YY\ar@{^(->} [d]^{j} \\
                          \XX   \ar@{^(->} [r]_{j\circ i} & \ZZ    }$$ 
   is also $i^*\N_{\YY/\ZZ}$, because we have a short exact sequence
                          $$ 0 \to  N_{\XX/\YY}  \to N_{\XX/\ZZ} \to i^*N_{\YY/\ZZ} \to 0.$$
   This can be checked by choosing an atlas for $\ZZ$.                       
\end{ex}

% ----------------------------------------------------------------
\section{Thom isomorphism} {\label{S:Thom}}

\begin{defn}{\label{D:orientablevb}}
 We say a vector bundle $p\: \E \to \X$ of rank $n$ on a topological stack  
 $\X$ is {\bf orientable}, if there is a class $\mu \in
 H^n(\E,\E-\X)$ such that the map
  $$\begin{array}{rcl}
     H^i(\X) & \llra{\tau} & H^{i+n}(\E,\E-\X) \\
     c & \mapsto & p^*(c)\cup \mu
                      \end{array}$$
 is an isomorphism for all $i \in \mathbb{Z}$. The class $\mu$ is
 called a {\bf Thom class}, or an {\bf orientation}, for $p\: \E \to
 \X$.
 \end{defn}

\begin{lem}{\label{L:pullback}}
Let $\E \to \X$ be an oriented vector bundle and $\mu \in
 H^n(\E,\E-\X)$ a Thom class for it. Let $f \: \Y \to \X$
 be a morphisms of stacks. Then $f^*\E \to \Y$ is an oriented vector
 bundle and $f^*(\mu)$ is a Thom class for it.
\end{lem}

%\begin{proof}
%\end{proof}

\begin{lem}{\label{L:pullback2}}
 Let $\E \to \X$ be a vector bundle. Let $f \: \Y \to \X$ be a
 trivial fibration of topological stacks, and let $\nu$ be a Thom class
 for the vector bundle $f^*\E \to \Y$. Then, there is a unique Thom
 class $\mu$ for $\E$ such that $f^*(\mu)=\nu$.
\end{lem}

%\begin{proof}
%\end{proof}

\begin{prop}{\label{P:orientablevb}}
  Let $p \: \E \to \X$ be an orientable vector bundle of rank
  $n$, and let
  $\mu \in H^n(\E,\E-\X)$ be
  a Thom class for it.
  Let $\K \subset \X$ be a closed substack.
  Then, the homomorphism
   $$\begin{array}{rcl}
     H^*(\X,\X-\K) & \llra{\tau} & H^{\scriptstyle \bullet+n}(\E,\E-\K) \\
     c & \mapsto & p^*(c)\cup \mu
                      \end{array}$$
  is an isomorphism. Here, we have identified $\K$ with a closed
  substack of $\E$ via the zero section of $\E \to \X$.
\end{prop}

\begin{proof} Let $\U=\X-\K$.
The map $c  \mapsto  p^*(c)\cup \mu$ induces a map between long
exact sequences
 {\small
  $$\xymatrix@C=8pt@R=10pt@M=8pt{ \cdots \ar[r] & H^{\scriptstyle \bullet+n}(\E|_{\U},\E|_{\U}-\U) \ar[r] &
  H^{\scriptstyle \bullet+n}(\E,\E-\X) \ar[r] & H^{\scriptstyle \bullet+n}(\E,\E-\K) \ar[r] & \cdots \\
  \cdots  \ar[r] & H^{\scriptstyle \bullet}(\U) \ar[r] \ar[u]^{\cong} &
  H^{\scriptstyle \bullet}(\X) \ar[r] \ar[u]^{\cong} & H^{\scriptstyle \bullet+n}(\X,\X-\K) \ar[r] \ar[u] &
  \cdots.
    }$$}
  (The top sequence is long exact by Proposition \ref{P:triple}.)
  The claim follows from 5-lemma.
\end{proof}

\begin{prop}{\label{P:Thomdot}}
   In Proposition \ref{P:orientablevb}, identify $\X$ with a closed
   substack of $\E$ via the zero section. Then, for every $c \in
   H^*(\X,\X-\K)$, we have $\tau(c)=c\cdot\mu$, where $\cdot$ is the
   product of Proposition \ref{P:dot}.
\end{prop}

%\begin{proof}
%\end{proof}

\begin{prop}{\label{P:metrizedorientablevb}}
 In Proposition \ref{P:orientablevb}, assume that $\E$ is metrized, and let
 $\D$ denote its disc bundle of radius $r$. Set $\LL=p^{-1}(\K)\cap\D$.
 and let $\rho \: H^{\scriptstyle \bullet}(\E,\E-\K) \to
  H^{\scriptstyle \bullet}(\E,\E-\LL)$ be the restriction homomorphism.
  Then the homomorphism
    $$\begin{array}{rcl}
     H^{\scriptstyle \bullet}(\X,\X-\K) & \llra{\tau} & H^{\scriptstyle \bullet+n}(\E,\E-\LL) \\
     c & \mapsto & \rho(p^*(c)\cup \mu)
                      \end{array}$$
    is an isomorphism.
In particular, the map $\rho$ is an isomorphism.
\end{prop}

\begin{proof}  Let $\U=\X-\K$.
 In the case where $\K=\X$, a standard deformation retraction
 argument shows that $\rho$ is an isomorphism, so the result follows
 from Proposition \ref{P:orientablevb}. The general case reduces to
 this case by considering the map of long exact sequences induced by
     $c  \mapsto \rho(p^*(c)\cup \mu)$,
{\small
 $$\xymatrix@C=8pt@R=10pt@M=8pt{ \cdots \ar[r] & H^{\scriptstyle \bullet+n}(\E|_{\U},\E|_{\U}-\D|_{\U}) \ar[r] &
  H^{\scriptstyle \bullet+n}(\E,\E-\D) \ar[r] & H^{\scriptstyle \bullet+n}(\E,\E-\LL) \ar[r] & \cdots \\
  \cdots  \ar[r] & H^{\scriptstyle \bullet}(\U) \ar[r] \ar[u]^{\cong} &
  H^{\scriptstyle \bullet}(\X) \ar[r] \ar[u]^{\cong} & H^{\scriptstyle \bullet+n}(\X,\X-\K) \ar[r] \ar[u] &
  \cdots,
    }$$}
 and applying 5-lemma.
\end{proof}

The following lemma strengthens Proposition \ref{P:orientablevb}.

\begin{lem}{\label{L:section}}
    Let $p \: \E \to \X$ be an orientable vector bundle of rank
  $n$, and let
  $\mu \in H^n(\E,\E-\X)$ be
  a Thom class for it.
  Let $\K \subset \X$ be a closed substack, and $\K' \subset \E$ a closed substack
  of $\E$ mapping isomorphically to $\K$ under $p$.
  Then, we have a natural isomorphism $H^{\scriptstyle \bullet}(\X,\X-\K) \cong
  H^{\scriptstyle \bullet+n}(\E,\E-\K')$.
\end{lem}

%\begin{proof}
%\end{proof}

\begin{lem}{\label{L:summand}}
  Let $p\: \E \to \X$ and $q\: \F \to \X$ be  vector bundles over
 $\X$, and assume that $\E$ is
  oriented. Then, an orientation for $\F$ determines an orientation
  for $\E\oplus\F$, and vice versa. Indeed, if $\mu$ is an
  orientation for $\E$, and $\nu$ an orientation for $\F$, then
  $\mu\cdot p^*(\nu)=\nu\cdot q^{*}(\mu)$ is an orientation for $\E\oplus\F$. Here,
  $\cdot$ is the product of Proposition \ref{P:dot}.
\end{lem}

\begin{proof} We only prove one of the statements, namely, the case
 where $\E$ and  $\E\oplus\F$ are  oriented. We show that $\F$ is
 also oriented.
  Assume $\E$ and $\F$ have rank
  $m$ and $n$,   respectively, and let
  $\mu \in H^m(\E,\E-\X)$ and $\nu \in  H^{m+n}(\E\oplus\F,\E\oplus\F-\X)$
  be orientations for $\E$ and $\E\oplus\F$.  The class
  $q^*(\mu) \in H^{m}(\E\oplus\F,\E\oplus\F-\F)$ is an
  orientation for the pullback bundle $q^*(\E)\cong\E\oplus\F$ over $\F$; note
  that the bundle map $q^*(\E) \to \F$ can be naturally identified
  with the second projection map $\pi \: \E\oplus\F \to \F$. By Proposition
  \ref{P:orientablevb}, applied to the vector bundle  $\pi \: \E\oplus\F \to \F$,
  we have an isomorphism
         $$\begin{array}{rcl}
     H^n(\F,\F-\X) & \ra & H^{n+m}(\E\oplus\F,\E\oplus\F-\X) \\
     c & \mapsto & \pi^*(c)\cup q^*(\mu).
                      \end{array}$$
   The inverse image of $\nu$ under this isomorphism is the desired
   orientation class in $H^n(\F,\F-\X)$.
\end{proof}

In Lemma \ref{L:summand}, we call the orientation  on $\E\oplus\F$
the {\bf sum} of the orientations of $\E$ and $\F$, and the
orientation on $\F$ the {\bf difference} of the orientations on
$\E\oplus\F$ and $\E$.

\begin{lem}{\label{L:summand2}}
   Let $0 \to \E \to \M \to \F \to 0$ be a short exact sequence of
   vector bundles over a topological stack $\X$. Then, the choice of
   orientations on two of the three vector bundles uniquely
   determines an orientation on the third one. Moreover, we have the
   relation following relation between the orientation classes:
     $$\mu_{\E}\cdot p^*(\mu_{\F})=\mu_{\M}.$$
   Here, $p$ stands for the morphism of pairs $(\M,\M-\E) \to (\F,\F-\X)$, and 
   $\cdot$ is the product of Proposition \ref{P:dot}.
\end{lem}

\begin{proof}
  Apply Lemma \ref{L:pullback2} to the trivial fibration $f \: \M \to
  \X$ to reduce the problem to the split case and then apply Lemma
  \ref{L:summand}.
\end{proof}

\begin{lem}{\label{L:section2}}
  In Lemma \ref{L:section}, assume we are given another oriented
  vector bundle $\F \to \X$ of rank $m$, and endow $\E\oplus\F$ with the sum
  orientation. Let $\K'' \subset \E\oplus\F$ be a closed substack mapping
  isomorphically to $\K'$ under the projection
  $\E\oplus\F \to \E$. Then, the diagram
    $$\xymatrix@C=-40pt@R=10pt@M=8pt{ H^{\scriptstyle \bullet}(\X,\X-\K)
    \ar[rr]^{\cong} \ar[rd]_{\cong} & & H^{\scriptstyle \bullet+n}(\E,\E-\K') \ar[ld]^{\cong} \\
    & H^{\scriptstyle \bullet+n+m}(\E\oplus\F,\E\oplus\F-\K'') & }$$
 commutes. (All the isomorphisms in this diagram are the ones of Lemma
 \ref{L:section}. So,  in the case where $\K=\K'=\K''$, the isomorphisms are simply
 the Thom isomorphisms of Proposition \ref{P:orientablevb}.)
\end{lem}

%\begin{proof}
%\end{proof}

Finally, we prove a lemma about compatibility of Thom isomorphism
with excision.

\begin{lem}{\label{L:Thomex}} Let $X$ be a  manifold, and
  let $E \to X$ and $N\to X$ be  vector bundles of rank $n$. Assume
  that $E$ is oriented. Let $i \: N \to E$ be an open embedding
  which sends the zero section of $N$ to the zero section of $E$.
  (Note that $N$ is naturally isomorphic to $E$, hence oriented,
  via the isomorphisms
  $TX\oplus N\cong TE\cong TX\oplus E$.) Then, the
  following diagram commutes:
         $$\xymatrix@C=0pt@R=10pt@M=8pt{H^{\scriptstyle \bullet+n}(N,N-X) \ar[rr]^{\text{excision}}_{\cong} &&
              H^{\scriptstyle \bullet+n}(E,E-X) \\
                 & H^{\scriptstyle \bullet}(X)\ar[ul]_(0.4){\cong}^{\tau_{N}}
                                \ar[ur]^(0.4){\cong}_{\tau_{E}} & }$$
\end{lem}

\section{Loop stacks}

\subsection{Mapping stacks and the free loop
  stack}\label{mappingstack}

Let $\XX$ and $\YY$ be   stacks over $\Top$. We define the stack
$\bfhom(\YY,\XX)$, called the {\bf mapping stack} from $\YY$ to
$\XX$, by the rule
  $$ T\in \Top \ \ \   \mapsto  \ \ \ \Hom(T\times \YY,\XX)\,,$$
where $\Hom$ denotes the groupoid of stack morphisms.
This is easily seen to be a stack. It follows from the exponential
law for mapping spaces (\cite{Whitehead})
that when $X$ and $Y$ are spaces, with $Y$ Hausdorff,
then $\bfhom(Y,X)$ is representable by the usual mapping space from
$Y$ to $X$ (endowed with the compact-open topology).

\smallskip

The mapping stacks $\bfhom(\YY,\XX)$ are  functorial in $\XX$ and $\YY$.

\begin{prop}{\label{P:MapSt}}
 Let $\XX$ be a topological stack and $A$ a compact %Hausdorff
 topological space. Then $\bfhom(A,\XX)$ is a topological stack.
\end{prop}

\begin{pf}
This follows from Theorem 1.1 of \cite{Mapping}.
\end{pf}

 Let $\XX$ be a topological stack. Then
 $\LXX=\bfhom(S^1,\XX)$ is  also  a topological stack.
 It is called the {\bf loop stack} of $\XX$. By functoriality of
mapping stacks, for every $t \in S^1$ we have the corresponding
evaluation map $\ev_t : \LXX \to \XX$. In particular, denoting by
 $0\in S^1$ the standard choice of a base point, there is an evaluation map
 \begin{eqnarray}  \label{eq:ev0} \ev_0 : \LXX \to \XX.\end{eqnarray}
 Similarly,
 the {\bf path stack} of $\XX$, which is defined to be $\bfhom(I,\XX)$,
 is a topological stack.

\medskip

For the next result, we need to assume that $\XX$ is a Hurewicz
topological stack.

\begin{lem}{\label{L:glue}}
  Let $A$, $Y$, and $Z$ be as in  Proposition \ref{P:glue}.
  Let $\X$ be a Hurewicz topological stack. Then the  diagram
     $$\xymatrix{
       \bfhom(Z\vee_A Y, \X) \ar[r] \ar[d] &  \bfhom(Y, \X) \ar[d] \\
       \bfhom(Z, \X) \ar[r]   & \bfhom(A, \X)   }$$
  is a 2-cartesian diagram   of topological stacks.
\end{lem}

\begin{proof}
  We have to verify that for every topological space $T$ the
  $T$-points of
  the above mapping stacks form a
  2-cartesian diagram of groupoids. This follows from Proposition
  \ref{P:glue} applied to $A\times T$, $Y\times T$, and $Z\times T$.
\end{proof}

We denote by `$8$' the wedge $S^1\vee S^1$ of two circles.
\begin{cor}{\label{L:8}}
  Let $\X$ be a Hurewicz topological stack, and let $L\X$ be its loop stack.
  Then, the diagram
    $$\xymatrix@=12pt@M=10pt@M=10pt@=12pt{
       \Map(8, \X) \ar[r] \ar[d] &  L\X \ar[d] \\
       L\X  \ar[r] &  \X   }$$
  is 2-cartesian.
\end{cor}

\subsection{Groupoid presentation}
\label{S:freeloopgpoid}

Let us now describe a particular groupoid presentation of the loop
stack.  For this we will assume that $\XX$ is a {\em Hausdorff
Hurewicz  topological stack}.  Thus $\XX$ admits a groupoid
presentation $\gm\:\gm_1\toto \gm_0$, where $\Gamma_0$ and
$\Gamma_1$ are Hausdorff topological spaces, $\Gamma_1\to
\Gamma_0\times \Gamma_0$ is proper, and source and target maps are
local fibrations. We will fix the groupoid $\Gamma$.

We will construct a groupoid $\Lo\Gamma\:\Lo_1\Gamma\toto \Lo_0\Gamma$
out of $\Gamma$ which represents $\Lo\X$. This groupoid presentation
is useful in computations (see Section~\ref{frobenius}). Our
construction resembles the construction of the fundamental groupoid of
a groupoid~\cite{MM}.

Let $M\gm=[M_1\gm\toto
M_0\gm]$ be the morphism groupoid of $\gm$. Its object set is
$M_0\Gamma=\Gamma_1$ and its morphism set $M_1\gm$ is the set of commutative
squares in the underlying category of $\Gamma$:
 \begin{equation}\label{eq:MGamma}
\vcenter{ \xymatrix{
      t(h) & t(k) \lto_{g} \\
      s(h) \uto^{h} & s(k) \uto_{k} \lto_{h^{-1}gk}}}
\end{equation}
 The source and
target maps are the horizontal arrows in
square~(\ref{eq:MGamma}). The groupoid multiplication is by (vertical)
superposition of such squares.  Thus  we have
$M_1\gm\cong
\gm_3=\gm_1\times_{\gm_0}\gm_1\times_{\gm_0}\gm_1 $. The groupoid
$M\Gamma$ is another presentation of the stack $\X$ and is Morita
equivalent to $\Gamma$.

Let $P\subset S^1$ be
a finite subset of $S^1$ which contains the base point of $S^1$. The
points of $P$ are labeled according to increasing angle as
$P_0,P_1,\ldots,P_n$ in such a way that $P_0=P_n$ is the base point of
$S^1$.  Write $I_i$ for the closed interval $[P_{i-1},P_i]$.  Let
$S^P_0$ be the disjoint union $S^P_0=\coprod_{i=1}^n I_i$. There is a
canonical map $S_0^P\to S^1$. Let $S_1^P$ be the fiber product
$S^P_1=S^P_0\times_{S^1}S^P_0$. There is an obvious topological
groupoid structure $S^P_1\toto S^P_0$.  The compact-open topology
induces a topological groupoid structure on $\Lo^P\gm: \Lo_1^P\gm\toto
\Lo_0^P\gm$, where $\Lo_0^P\gm$ is the set of continuous strict
groupoid morphisms $[S_1^P\toto S^P_0]\to[\gm_1\toto \gm_0]$ and
$\Lo_1^P\gm$ is the set of strict continuous groupoid morphisms
$[S_1^P\toto S_0^P]\to[M_1\gm\toto M_0\gm]$.

The finite subsets of $S^1$ including the base point are ordered by
inclusion.  The ordering is {\em directed}.  For $P\leq Q$ there is a
canonical morphism of groupoids $\Lo^P\Gamma\longrightarrow \Lo^Q\Gamma$.
Using the fact that $\Gamma_0$ and $\Gamma_1$ are Hausdorff, it is not
difficult to prove that $\Lo^P\Gamma\to \Lo^Q\Gamma$ is an isomorphism
onto an open subgroupoid. Define the topological groupoid
$$\Lo\Gamma=\varinjlim_{P\subset S^1}\Lo^P\Gamma=\bigcup_{P\subset
  S^1}\Lo^P\Gamma\,.$$

\begin{prop}
\label{freeloopgpoid}
The groupoid $\Lo\Gamma$ presents the  loop stack $\Lo\X$.
\end{prop}
\begin{pf}
First, we need to construct a morphism $\Lo_0^P\Gamma\to\Lo\XX$, for
every $P$. The presentation $\Lo_0\Gamma\to\Lo\XX$ will then be
obtained by gluing these morphisms using the stack property of
$\Lo\X$ and the fact the $\Lo_0^P\Gamma$ form an open covering of
the topological space $\Lo_0\Gamma$.

The structure map $\Lo_0^P\times S_0^P\to\Gamma_0$ gives rise to a
morphism $\Lo_0^P\times S_0^P\to \XX$.  This morphism descends to
$\Lo_0^P\times S^1\to \XX$, by Proposition~\ref{P:glue}, because $S^1$
is obtained from $S_0^P$ as a pushout covered by that proposition.
By adjunction, we obtain the required morphism $\Lo_0^P\to\Lo\XX$.

The fact that $\bigcup_P \Lo_0^P\Gamma\to\Lo\XX$ is an epimorphism of
stacks, follows as in Proposition~\ref{P:MapSt}.

The fact that $\Lo_1\Gamma$ is the fibered product of $\Lo_0\Gamma$
with itself over $\Lo\X$ reduces immediately to the case of
$\Lo_1^P\Gamma$, for which it is immediate.
\end{pf}

It is easy to represent evaluation map and functorial properties
of the free loop stack at the groupoid level with this model.

\begin{numrmk}
In particular,  there is an equivalence of the underlying categories
between $\Lo\gm$ and the groupoid whose objects are the set of
generalized morphisms from the space $S^1$ to $\gm$ and has
equivalences of such as arrows.
\end{numrmk}
\begin{cor}
If $\X$ is a differentiable stack then $\Lo\X$ is regular \Lind.
\end{cor}

\subsubsection{Target connected groupoid}

Assume the groupoid $\Gamma$ is target connected.  This means that if
$T$ is a topological space, and $\phi\:T\to\Gamma_1$ a continuous map, then
for every point of $T$ there exist an open neighborhood $T'\subset
T$ and a homotopy $\Phi\:T'\times I\to \Gamma_1$, such that
$\Phi_0=\phi$ and $\Phi_1=t\circ\phi$, where $t:\Gamma_1\to\Gamma_0$
is the target map.  For example, any transformation groupoid with
connected Lie group is target connected.

For every  finite subset $P\subset S^1$
  and $x\in \Lo^P_0\gm$, there
  are  arrows $g_i \in \gm_1$ with $t(g_i)=P_i\in I_i$ and
  $s(g_i)=P_{i}\in I_{i+1}$ (or $P_0\in I_1$ if $i=n$).
These arrows can
  be continuously deformed to the identity
  point $P_n\in I_n$. Thus there is an element $\tilde{x}\in
  \bfhom(S^1,\gm_0)\subset \Lo^{\{0\}}\gm_0$
  and an arrow $\gamma \in \Lo^P_1\gm$ with $s(\gamma)=\tilde{x}$ and
  $t(\gamma)=x$.  From this observation, we deduce:

\begin{prop}\label{prop:targetconnected}
If $\Gamma$ is target connected, then the groupoid
  $\Lo\gm_1\toto \Lo\gm_0$ with pointwise source map, target map and
  multiplication presents the loop stack $\Lo\X$. Here $\Lo\gm_i$
  is the usual  free loop space of
  $\gm_i$ endowed with the compact-open topology.

In particular, $\LG$ is Morita equivalent to the groupoid
$\Lo\gm_1\toto \Lo\gm_0$.
\end{prop}
\begin{numex}
If $G$ is a connected Lie group acting on a manifold $M$, then Proposition~\ref{prop:targetconnected} implies that $\Lo [M/G] \cong [\Lo M/ \Lo G]$.
\end{numex}

% -------------------------------
\subsubsection{Discrete group action}

 To the contrary, if $G$ is a discrete group acting on a space $M$ one
  can form the global quotient $[M/G]$ which is represented by the
  transformation groupoid $\gm\:M\times G\toto M$. For any $x\in
  \Lo^P\gm_0$ one can easily find an arrow $\gamma\in \Lo^P\gm_1$ such
  that $s(\gamma)=x$ and $t(\gamma)\in \Lo^{\{0\}}\gm_0$. Furthermore,
  since $G$ is discrete, an element of $\Lo^P\gm_1$ is described by its
  source and one element $g_i\in G$ for $i=0,\ldots,|P|$. From these
  two observations one proves easily:

\begin{prop}\label{loopdiscrete}
Let $G$ be a discrete group acting on a space $M$. Then ${\rm L}[M/G]$
is presented by  the transformation groupoid
$$\left(\coprod_{g\in G}\mathcal{P}_g M\right) \times G \toto
\coprod_{g\in G}\mathcal{P}_g M$$ where $\mathcal{P}_g M=\{ f\:[0,1]\to
M \mbox{ such that } f(0)=f(1).g\}$ and $G$ acts by pointwise
conjugation.
\end{prop}
Note that if $G$ is finite, one recovers the loop orbifold
  of~\cite{LU/loop}.

% ----------------------------------------------------------------
\section{Bounded proper morphisms of topological stacks}
{\label{S:Proper}}

\begin{defn}{\label{D:proper}}
Let  $f \: \X \to \Y$ be a morphisms of topological stacks and $\E$
a metrizable vector bundle over $\Y$. A lifting $i \: \X \to \E$ of
$f$,
 $$\xymatrix@=16pt@M=8pt{  & \E \ar[d] \\
                          \X\ar@{^(..>}[ru]^i \ar[r]_{f} & \Y     }$$
is called {\bf bounded} if there is a choice of metric on $\E$ such
that $i$ factors through the unit disk bundle of $\E$. A morphism $f
\: \X \to \Y$ of topological stacks is called {\bf \proper} if there
exists a metrizable orientable vector bundle $\E$ on $\Y$ and a
bounded lifting $i$ as above such that $i$ is a closed embedding.
\end{defn}

\begin{defn}{\label{D:superproper}}
   A \proper morphism $f \: \X \to \Y$ is called {\bf \superproper}
   \comment{Until we decide on a less stupid name!}
   if every orientable metrizable vector bundle $\E$ on $\X$
   is a direct summand  of $f^*(\E')$ for some orientable metrizable vector
   bundle $\E'$ on $\Y$. (Note that, possibly
after multiplying by a positive $\bbR$-valued function on $\Y$, we
can arrange the inclusion $b^*\E\hookrightarrow (gb)^*(\E')$ to be
contractive, i.e., have norm at most one.)
\end{defn}

\begin{ex}{\label{E:proper}}
\begin{itemize}
  \item[$\mathbf{1.}$] Every \proper map $f \: X \to Y$ of a
    topological spaces with $Y$ compact is \superproperpt. In that case, one can use the fact that every vector bundle on
    a compact space is a subbundle of a trivial bundle.

  \item[$\mathbf{2.}$] Let $\X$ be a topological stack such that
  $\Delta \: \X \to \X\times\X$ is \properpt. Then $\Delta$ is
  \superproperpt. This follows from the fact that every vector bundle on
  $\X$ can be naturally extended to $\X\times\X$. Similarly, the iterated diagonal $\Delta^{(n)}\:\X \to \X^{n}$ is \superproperpt.
\item[$\mathbf{3.}$] Let $X,Y$ be compact $G$-manifolds (with $G$  compact) and $f\:X\to Y$ be a $G$-equivariant map. Then the induced map of stacks  $[f/G]\: [X/G]\to [Y/G]$ is \superproperpt.
\end{itemize}

\end{ex}

It does not seem to be true in general that two \proper maps compose
to a \proper map, but we have the following.

\begin{lem}{\label{L:superproper}}
  Let $f \: \X \to \Y$ and $g \: \Y \to \Z$ be \superproper morphisms.
  Then $g\circ f \: \X \to \Z$ is \superproperpt.
\end{lem}

\begin{proof}
  It is trivial that every orientable metrizable bundle on $\X$
  is a direct summand of one coming from $\Z$. Let us now prove that
  $g\circ f$ is proper.
  Suppose given factorizations
     $$\xymatrix@=16pt@M=8pt{  & \E \ar[d] \\
                          \X\ar@{^(..>}[ru]^i \ar[r]_{f} & \Y} \ \ \ \
      \xymatrix@=16pt@M=8pt{  & \F \ar[d] \\
                          \Y\ar@{^(..>}[ru]^j \ar[r]_{g} & \Z     }$$
  for $f$ and $g$. By enlarging $\E$, and using that $g$ is \superproper,
  we may assume that $\E=g^*(\E')$, for some oriented metrized
  vector bundle $\E'$ on $\Z$. Let $i' \: \X \to \E'$ be the
  composition $\pr\circ i$ where $\pr \: \E \to \E'$ is the
  projection map. The following diagram shows that $g\circ f$ is
  proper:
       $$\xymatrix@=16pt@M=8pt{  & \E'\oplus\F \ar[d] \\
            \X\ar@{^(..>}[ru]^{(i',jf)} \ar[r]_{g\circ f} & \Z}$$
\end{proof}

% -----------------------------------------------------------------
\subsection{Some technical lemmas}{\label{S:Techinical}}

In this section we prove a few technical lemmas that will be needed
in Section \ref{S:Bivariant} to define bivariant groups.

Let $f \: \X \to \Y$ be a morphism of topological stack that admits
a factorization
    $$\xymatrix@=16pt@M=8pt{  & \E \ar[d] \\
                          \X\ar@{^(..>}[ru]^i \ar[r]_{f} & \Y     }$$
For example, every \proper $f$ has this property (Definition
\ref{D:proper}). The following series lemmas investigate certain
properties of the relative cohomology groups
$H^{\scriptstyle \bullet}(\E,\E-\X)$.
% In none of these lemmas the assumption that $i \: \X \to \E$ factors
%through the unit disk bundle of $\E$ is needed.

\begin{lem}{\label{L:isom}}
  Let $f \: \X \to \Y$ be a morphism of topological stacks,
  and assume we are given two different factorizations $(i,\E)$
  and $(i',\E')$ for it. Then, there is a canonical isomorphism
   $H^{\scriptstyle \bullet+\rk\E}(\E,\E-\X) \cong H^{\scriptstyle \bullet+\rk\E'}(\E',\E'-\X)$.
\end{lem}

\begin{proof}
    Embed $\X$ in $\E\oplus\E'$ via $(i,i') \: \X \to \E\oplus\E'$.
    Consider the diagram
    $$(\E',\E'-\X)  \la (\E'\oplus\E,\E'\oplus\E-\X) \to (\E,\E-\X)   $$
   of pairs  of topological stacks. It follows from
   Proposition
   \ref{P:orientablevb}
   that we have natural isomorphisms
      $$H^{\scriptstyle \bullet+\rk\E'}(\E',\E'-\X) \lisom
   H^{\scriptstyle \bullet+\rk\F+\rk\E}(\E'\oplus\E,\E'\oplus\E-\X) \risom
   H^{\scriptstyle \bullet+\rk\E}(\E,\E-\X).$$
  We can now apply Lemma \ref{L:section}.
\end{proof}

Using a triple direct sum argument, it can be shown that given three
factorizations $(i,\E)$, $(i',\E')$,
  and $(i'',\E'')$ for $f$, the corresponding isomorphisms defined in the above
  lemma are compatible. Also, if we switch the order of $(i,\E)$ and  $(i',\E')$
we get the inverse isomorphism. Finally, when $(i,\E)$ and
$(i',\E')$ are equal we get the identity isomorphism. Therefore, the
group $H^{\scriptstyle \bullet}(\E,\E-\X)$ only depends on the morphism $f$.

\begin{lem}{\label{L:invariance}}
  Let $f \: \X \to \Y$ be a morphism of topological stacks, and
  let $\varphi:=f\circ\pr \: \X\times I \to \Y$, where $I$ is the unit
  interval and  $\pr$ stands for projection. Suppose we are give
  a factorization
          $$\xymatrix@=16pt@M=8pt{  & \E \ar[d] \\
                          \X\times I\ar@{^(..>}[ru]^{\iota} \ar[r]_{\varphi} & \Y     }$$
  for $\varphi$. Let $0\leq a\leq 1$, and define $\iota_a \: \X \to \E$ to be the restriction of $\iota$ to $\X=\X\times\{a\}$.
  Then, the natural map $\phi_a \: H^{\scriptstyle \bullet}(\E,\E-\iota_a(\X)) \to H^{\scriptstyle \bullet}(\E,\E-\iota(\X\times I))$
  induced by the map of pairs $(\E,\E-\iota(\X\times I)) \to (\E,\E-\iota_a(\X))$ is
  an isomorphism and it is independent of $a$.
\end{lem}

\begin{proof}
  We may assume that the image of $\iota$ does not intersect the zero section of $\E$.
  (For example, we  lift everything to $\E\oplus\bbR$ via  $(\iota,1)\: \X \to \E\oplus\bbR$
  and apply  Proposition \ref{L:section} to the vector bundle $\E\oplus\bbR \to
  \E$.).

 Let $\E'=\E\oplus\bbR$ and define $\beta \: \X \times I \hookrightarrow \E'$ by
 $\beta(x,t)=(\iota(x,0),t)$. This is a closed embedding, so by Lemma
   \ref{L:isom}, we have a commutative diagram
     $$\xymatrix@=16pt@M=8pt{ H^{\scriptstyle \bullet}(\E,\E-\iota(\X\times I)) \ar[r]^{\cong} & H^{\scriptstyle \bullet}(\E',\E'-\beta(\X\times I))  \\
                          H^{\scriptstyle \bullet}(\E,\E-\iota_a(\X)) \ar[u]^{\phi_a} \ar[r]^{\cong} & H^{\scriptstyle \bullet}(\E',\E'-\beta_a(\X))   \ar[u]_{\phi'_a}  }$$
    This reduces the problem to the case where our map is $\beta$
    instead of $\iota$, in which case the result is obvious.

%   Consider another copy of $\E$, which we denote by $\E'$ to avoid
%   confusion, and define $\beta \: \X \times I \to \E'$ by
%   $\beta(x,t)=t\iota(x,0)$. This is a closed embedding, so by Lemma
%   \ref{L:isom}, we have a commutative diagram
%     $$\xymatrix@=16pt@M=8pt{ H^{\scriptstyle \bullet}(\E,\E-\X\times I) \ar[r]^{\cong} & H^{\scriptstyle \bullet}(\E',\E'-\X\times I)  \\
%                          H^{\scriptstyle \bullet}(\E,\E-\X) \ar[u]^{\phi} \ar[r]^{\cong} & H^{\scriptstyle \bullet}(\E',\E'-\X)   \ar[u]_{\phi'}  }$$
%    So, to show that $\phi$ is an isomorphism, it s enough to show that $\phi'$ is an isomorphism. We prove this using a
%    deformation retraction argument.........
\end{proof}

% ----------------------------------------------------------------
\section{Bivariant theory for topological stacks} {\label{S:Bivariant}}

We define a bivariant cohomology theory \cite{FuMac} on the category
of topological stacks whose associated covariant and contravariant
theories are singular homology and cohomology, respectively. Our
bivariant theory satisfies  weaker axioms than those of \cite{FuMac}
in that  products are not always defined. We show, however, that
there are enough products  to enable us to define Gysin morphisms as
in \cite{FuMac}.

The underlying category of our bivariant theory is the category
$\TopSt$ of topological stacks. \comment{we ignore 2-isomorphisms}
The confined morphisms are all maps and independent squares are
2-cartesian diagrams.

% ----------------------------------------------------------------
\subsection{Bivariant groups}
\label{SS:bivariant}

To a morphism $f \: \X \to \Y$  of topological stacks, we associate
a category $\sfC(f)$ as follows. The objects of $\sfC(f)$ are
morphism $a \: \K \to \X$ such that $fa \: \K \to \Y$ is \proper
(Definition \ref{D:proper}). A morphism in $\sfC(f)$ between $a\: \K
\to \X$ and $b \: \LL \to \X$ is a homotopy class (relative to $\X$)
of morphisms $g \: \K \to \LL$ over $\X$.  \comment{Remark that,
when working with topological stacks, homotopy is {\em not}
necessarily an equivalence relation. So we have to consider the
equivalence relation generated by it}

\begin{lem}{\label{L:cofiltered}}
  The category $\sfC(f)$ is cofiltered.
\end{lem}

%\begin{proof}
%\end{proof}

Once and for all, we choose, for each object $a \: \K \to \X$,  a
vector bundle $\E \to \Y$ through which $fa$ factors, as in
Definition \ref{D:superproper}.% We will make this part of the notation.
%That is, we denote the object $a \: \K \to \X$ as a pair $(\K,\E)$.

We define the {\bf bivariant singular homology} of an arbitrary
morphism $f \: \X \to \Y$ to be  the $\bbZ$-graded abelian group
  $$H^{\scriptstyle \bullet}(\X\llra{f}\Y)=\varinjlim_{\sfC(f)}H^{\scriptstyle \bullet+\rk\E}(\E,\E-\K).$$

The homomorphisms in this direct limit are defined as follows.
Consider a morphism $\varphi \: \K \to \K'$ in $\sfC(f)$. From this
we will construct a natural graded pushforward homomorphism
$\varphi_* \: H^{\scriptstyle \bullet+m}(\E,\E-\K) \to
H^{\scriptstyle \bullet+n}(\E',\E'-\K')$, where   $m=\rk\E$ and
$n=\rk\E'$.

 Let $\F=\E\oplus\E'$ with the sum orientation.
Let $p \: \E' \to \Y$ be the projection map. Then, $p^*(\E)$ is an
oriented vector bundle over $\E'$. Note that the projection map $\pi
\: p^*(\E) \to \E'$ is naturally isomorphic to the second projection
map $\F=\E\oplus\E' \to \E'$; this allows us to view $\F$ as an
oriented vector bundle of rank $m$ over $\E'$.  Let $\D \subseteq
\F$ be the unit disc bundle. It follows from the assumptions that
$\K \subseteq \D$, hence also $\K \subseteq
\LL:=\pi^{-1}(\K')\cap\D$. The restriction homomorphism
  $$\varphi_* \: H^{\scriptstyle \bullet+m+n}(\F,\F-\K) \to H^{\scriptstyle \bullet+m+n}(\F,\F-\LL)\cong
  H^{\scriptstyle \bullet+n}(\E',\E'-\K'),$$
 induced by the  inclusion of
pairs  $(\F,\F-\LL) \to (\F,\F-\K)$  is the desired pushforward
homomorphism; here, we have used the isomorphism of Proposition
\ref{P:metrizedorientablevb}.

Next we have to show that the $\varphi_*$ is independent of the
homotopy class of $\varphi$. Consider $a\circ\pr\:\K\times I \to
\X$, and let $\rho_0, \rho_1 \: \K \to \K\times I$ be the times 0
and time 1 maps. Note that $a\circ\pr\:\K\times I \to \X$ is an
object of $\sfC(f)$.
 Since every homotopy (relative to $\X$)  between maps with domain $\K$
  factors through $\K
\times I$, it is enough to show that $\rho_{0,*}=\rho_{1,*}$. This
follows from Lemma \ref{L:invariance}.

\begin{rem}
  Let $\K \to \Y$ be a \proper morphism. It follows from Lemma
  \ref{L:isom}, that the cohomology $H^{\scriptstyle \bullet}(\E,\E-\K)$ is
  independent of choice of the vector bundle $\E$ and the embedding $i \:
  \K \hookrightarrow \E$, up to a canonical isomorphism.
  Furthermore, the pushforward maps constructed above are
  compatible with these canonical isomorphisms. So, $H^{\scriptstyle \bullet}(f)$ is independent of  all choices involved in its definition.
\end{rem}

\begin{lem}{\label{L:proper}}
 Let $f \: \X \to \Y$ be a \proper morphism and $\X \torel{i} \E \to
 \Y$ a factorization for $f$, where $i$ is a closed embedding (but $\E$
 is not necessarily metrizable).
 Then we have a natural isomorphism
   $$H^{\scriptstyle \bullet}(\X\llra{f}\Y)\cong H^{\scriptstyle \bullet+\rk\E}(\E,\E-\X).$$
\end{lem}

\begin{proof}
 Follows from Lemma \ref{L:isom}.
\end{proof}

% ----------------------------------------------------------------
\subsection{Independent pullbacks}{\label{SS:independent}}

Consider a cartesian diagram
    $$\xymatrix@=16pt@M=8pt{ \X' \ar[r]^{f'} \ar[d]& \Y' \ar[d]^h  \\
                          \X  \ar[r]_{f} & \Y     }$$
We define the pullback $h^* \: H(\X\llra{f}\Y) \to
H(\X'\llra{f'}\Y')$ as follows.

 Pullback along $h$ induces a functor $h^* \: \sfC(f) \to \sfC(f')$,
 $\K\mapsto h^*\K:=\X'\times_{\X}\K$.
 Furthermore,  we have a
 natural homomorphism
  $$H^{\scriptstyle \bullet+\rk\E}(\E,\E-\K) \to   H^{\scriptstyle \bullet+\rk
 \E}(h^*\E,h^*\E-h^*\K)$$
induced by the map of pairs $(h^*\E,h^*\E-h^*\K) \to (\E,\E-\K)$.
Using Lemma \ref{L:isom}, this induces the desired homomorphism of
colimits
  $$h^* \: \varinjlim_{\sfC(f)}H^{\scriptstyle \bullet+\rk\E}(\E,\E-\K) \to
                   \varinjlim_{\sfC(f')}H^{\scriptstyle \bullet+\rk\E'}(\E',\E'-\K').$$

% ----------------------------------------------------------------
\subsection{Confined pushforwards}{\label{SS:confined}}

Let $h \: \X \to \Y$ be a  morphism of topological stacks
(Definition \ref{D:proper}) fitting in a commutative triangle
$$\xymatrix@C=0pt@R=10pt@M=8pt{ \X \ar[rd]_{f} \ar[rr]^h & & \Y \ar[ld]^{g} \\
    & \Z & }$$
  We define the pushforward homomorphism $h_* \: H(\X\llra{f}\Z) \to
  H(\Y\llra{g}\Z)$ as follows.

  There is a natural functor $\sfC(f) \to \sfC(g)$, which sends $a \: \K \to \X$\
  to $ha \: \K \to \Y$. A factorization for $fa$ gives a
  factorization for $gha$ in a trivial manner:
     $$\xymatrix@R=-2pt@C=16pt@M=8pt{ \K \ar@{^(->}[r]^i \ar[dd]_a & \E \ar[dd]
                    &  & \K \ar@{^(->}[r]^i \ar[dd]_{ha} & \E \ar[dd] \\
                                         && \mapsto &&                    \\
                   \X  \ar[r]_{f} & \Z & &   \Y  \ar[r]_{g} & \Z     }$$
  Using Lemma \ref{L:isom}, this induces the desired homomorphism
   $$h_* \: \varinjlim_{\sfC(f)}H^{\scriptstyle \bullet+\rk\E}(\E,\E-\K) \to
                   \varinjlim_{\sfC(g)}H^{\scriptstyle \bullet+\rk\E}(\E,\E-\K).$$

% ----------------------------------------------------------------
\subsection{Products}{\label{SS:products}} Unfortunately, we are
not able to define product for arbitrary pairs of composable
morphisms $f$ and $g$. However, under an extra assumption on $g$
this will be possible.

\begin{defn}{\label{D:super}}
  A morphism $f \: \X \to \Y$ of topological stacks is called
  {\bf adequate}\comment{Until further notice!} if in the cofiltered category
  $\sfC(f)$ the subcategory consisting of $a \: \K \to \X$
  such that $fa \: \K \to \Y$ is \superproper is cofinal.
\end{defn}

\begin{ex}{\label{E:super}}
\begin{itemize}
\item[$\mathbf{1.}$] Every \superproper morphism is adequate.
  (Because in this case $\sfC(f)$ has a final object that is
  \superproper over $\Y$.)

\item[$\mathbf{2.}$] A morphism $f \: \X \to Y$ in which $Y$
is a paracompact topological space is adequate. (In this case every
object in $\sfC(f)$ is \superproper over $Y$; see Example
\ref{E:proper})
\end{itemize}
\end{ex}

Let $f \: \X \to \Y$ and $g \: \Y \to \Z$ be morphisms of
topological stacks, and assume $g$ is adequate. Then we can define
products of any two classes $\alpha\in H(f)$ and $\beta \in H(g)$.
The construction of the product is as follows. Consider objects
$(\K,a) \in \sfC(f)$ and $(\LL,b) \in \sfC(g)$, and choose
factorizations
    $$ \xymatrix@=16pt@M=8pt{ \K \ar@{^(->}[r]^i \ar[d]_a & \E \ar[d] \\
                          \X  \ar[r]_{f} & \Y     }  \ \ \ \
    \xymatrix@=16pt@M=8pt{ \LL \ar@{^(->}[r]^j \ar[d]_b & \F \ar[d] \\
                         \Y  \ar[r]_{g} & \Z     } $$
We may assume $gb \: \LL \to \Z$ is \superproperpt. There exists a
metrizable oriented vector bundle $\E'$ over $\Z$ such that $b^*\E$
is isomorphic to a subbundle of $(gb)^*(\E')$ as vector bundles over
$\LL$. Note that, possibly after multiplying by a positive
$\bbR$-valued function on $\Z$, we can arrange the inclusion
$b^*\E\hookrightarrow (gb)^*(\E')$ to be contractive (i.e., have
norm at most one). Let us denote $b^*\E$ by $\E_0$, $(gb)^*(\E')$ by
$\E_1$, and the codimension of $\E_0$ in $\E_1$ by $c$.

We  define the product
  $$H^{r}(\E,\E-\K)\otimes H^{s}(\F,\F-\LL) \to
  H^{r+s+c}(\E'\oplus\F,\E'\oplus\F-\K\times_{\Y}\LL).$$
 as follows. (Note that $(\K\times_{\Y}\LL,a\circ\pr)$ belongs to $C(g\circ f)$ and
 we have a factorization
     $$\xymatrix@R=16pt@M=8pt{ \K\times_{\Y}\LL \ar@{^(->}[r]^{(i,j)}
                                          \ar[d] & \E'\oplus\F \ar[d] \\
                          \X  \ar[r]_{g\circ f} & \Z    }$$
 for it. We explain this in more detail shortly.) By pulling back the map $i$ along
  $\varpi \: \E_0 \to \E$, we obtain a closed embedding $\K\times_{\Y}\LL
  \hookrightarrow \E_0$. On the other hand, we have a closed
  embedding $\E_1 \hookrightarrow \E'\oplus\F$; this is
  simply
  the pullback of $j$ along the projection map $\pi\: \E'\oplus\F \to \F$.
  Using the inclusion $\E_0\hookrightarrow \E_1$, we find
  a factorization
    $$\xymatrix@=16pt@M=12pt{\K\times_{\Y}\LL \ar@{^(->}[r] \ar@/^1pc/@<1ex> [rrr]^{(i,j)}
       & \E_0 \ar@{^(->}[r]&  \E_1 \ar@{^(->}[r] & \E'\oplus\F.}$$
  Now, let $\alpha \in H^{r}(\E,\E-\K)$ and $\beta \in H^{s}(\F,\F-\LL)$
  be two cohomology classes. We define $\alpha\cdot\beta \in
  H^{r+s+c}(\E'\oplus\F,\E'\oplus\F-\K\times_{\Y}\LL)$ to be
  $\tau(\varpi^*(\alpha))\cdot\pi^*(\beta)$, where the latter $\cdot$ is the product of Proposition
  \ref{P:dot}. In more detail, we have $\pi^*(\beta) \in
  H^s\big(\E'\oplus\F,\E'\oplus\F-\E_1\big)$, $\varpi^*(\alpha) \in
  H^r(\E_0,\E_0-\K\times_{\Y}\LL)$, and $\tau \: H^r(\E_0,\E_0-\K\times_{\Y}\LL)
  \to H^{r+c}\big(\E_1,\E_1-\K\times_{\Y}\LL\big)$ is
  the Thom isomorphism of Proposition \ref{P:orientablevb} for the
  vector bundle $\E_1$  over $\E_0$; to obtain this Thom isomorphism, we have used that,
  since the bundles are metrizable, $\E_0$ is a direct summand of
  $\E_1$ and its complement is oriented (Lemma \ref{L:summand}). Finally, our $\cdot$ is
  the one of  Proposition
  \ref{P:dot} with the inclusions $\K\times_{\Y}\LL \hookrightarrow \E_1 \hookrightarrow
  \E'\oplus\F$, $n=r+c$ and $m=s$.

%\begin{prop}{\label{P:fully2}}
%  Let $f \: \X \to \Y$ and $g\: \Y \to \Z$ be proper morphisms of
%   topological stacks, and assume that $g\circ f$ is also proper. Then, we have a fully
%   defined product
%      $$H^{r}(\X\llra{f}\Y)\otimes H^{s}(\Y\llra{g}\Z)
%      \to
%      H^{r+s}(\X\llra{g\circ f}\Z).$$
%\end{prop}
%
%\begin{proof} We will use the notation of the previous paragraphs.
%  We will also assume $\K=\X$ and $\LL=\Y$.
%
% Since $g\circ f$ is proper, we can find a vector bundle
%   $\E'$  over $\Z$ through which $g\circ f$ factors via a bounded map. Set $\E=g^*\E'$.
%   Thus, we obtain a factorization for $f$ through $\E$ via a bounded morphism.
%   By Lemma   \ref{L:proper}, we have
%     $$H^{r+\rk\E}(\E,\E-\X) \cong
%     H^r(\X\llra{f}\Y),  \ \ H^{s+\rk\F}(\F,\F-\Y)
%     \cong H^s(\Y\llra{g}\Z), \ \ \text{and}$$
%       $$H^{r+s+\rk\F+\rk E'}(\F\oplus\E',\F\oplus\E'-\X)
%     \cong H^{r+s}(\X\llra{g\circ f}\Z).$$
%    Thanks to the existence of $\E'$,  the product
%        $$H^{r+\rk\E}(\E,\E-\X)\otimes H^{s+\rk\F}(\F,\F-\Y) \to
%        H^{r+s+\rk\F+\rk E'}(\F\oplus\E',\F\oplus\E'-\X)$$
%   is defined. This is what we wanted to prove.
%\end{proof}

% ----------------------------------------------------------------
\subsection{K\"unneth formula}{\label{S:Kunneth}}

For bivariant singular cohomology with field coefficients we have the following K\"unneth formula.

\begin{prop}[K\"unneth formula]{\label{P:bivariantKunneth}}
  Let $f \: \X \to \Y$ and $f' \: \X' \to \Y'$ be morphisms of topological stacks. Then, we have
a natural isomorphism of graded groups
  $$H^{\bullet}(f\times f') \cong H^{\bullet}(f)\otimes H^{\bullet}(f'),$$
  where $f\times f' \: \X\times \X' \to \Y\times\Y'$ is the product map. 
\end{prop}

\begin{proof}
%    In the case where $f$ and $f'$ are closed embeddings this follows from the K\"unneth formula
%    for singular cohomology (Proposition \ref{P:Kunneth}). 
    First, suppose that $f$ and $f'$ are bounded proper,
    and choose factorizations as in Definition \ref{D:superproper}. We obtain  a factorization
       $$\xymatrix@=16pt@M=8pt{  & \E\boxtimes\E' \ar[d] \\
                          \X\times\X' \ar@{^(..>}[ru]^{(i ,i')}\ar[r]_{(f,f')} & \Y\times\Y'     }$$   
   for $f\times f'$. Note that the total space of the vector bundle  $\E\boxtimes\E'$ is $\E\times\E'$.
   Let $\U=\E-i(\X)$ and $\U'=\E'-i(\X')$.  Then,  
   $\E\times\E'-(i,i')(\X\times\X')=\E\times\U'\cup\U\times\E'$. So, by  Proposition \ref{P:Kunneth},
   we have
    {\small \begin{equation}\label{eq:bivariantKunneth}  H^{\bullet}(f\times f')\cong H^{\bullet+n+n'}(\E\times\E'-(i,i')(\X\times\X'))\cong   
      H^{\bullet+n}(\E,\U)\otimes H^{\bullet+n'}(\E',\U'),\end{equation}  }      
   and the latter is equal to $\cong H^{\bullet}(f)\otimes H^{\bullet}(f')$.    

   To prove the isomorphism for general $f$ and $f'$, consider the functor
   $P\: \sfC(f)\times\sfC(f') \to \sfC(f\times f')$ which sends a pair $(a,a')  \in \sfC(f)\times\sfC(f')$, with
   $a \: \K \to \X$ and $a' \: \K' \to \X'$, 
   to $a\times a' \: \K\times\K' \to \X\times\X'$. 
   Since we know the result for bounded proper morphisms, to prove the K\"unneth isomorphism for $f$
   and $f'$ we observe that, in general, for every  directed system indexed by  $\sfC(f\times f')$, the induced 
   directed system (via $P$) indexed by  $\sfC(f)\times\sfC(f')$ 
   has the same colimit. This is due to the  fact that $P$ has a left adjoint
   $Q\: \sfC(f\times f')  \to \sfC(f)\times\sfC(f')$,  defined by sending $a \: \K \to \X\times\X'$ to 
   $(\pr_1\circ a,\pr_2\circ a) \in \sfC(f)\times\sfC(f')$.
\end{proof}

When the coefficient is only a ring, the cohomology cross-product (see Remark~\ref{R:cross}) yields a bivariant cross-product.
\begin{prop}\label{P:bivariantcross}
 Let $f \: \X \to \Y$ and $f' \: \X' \to \Y'$ be morphisms of topological stacks. Then, we have
a natural homomorphism of graded groups
  $$  H^{\bullet}(f)\otimes H^{\bullet}(f') \to H^{\bullet}(f\times f')$$
  where $f\times f' \: \X\times \X' \to \Y\times\Y'$ is the product map.
\end{prop}
\begin{proof}
 The proof of Proposition~\ref{P:bivariantKunneth} applies with the only difference that the last isomorphism in the sequence~\eqref{eq:bivariantKunneth} of isomorphisms is replaced by the cross product $H^{\bullet+n+n'}(\E\times\E'-(i,i')(\X\times\X'))\leftarrow  
      H^{\bullet+n}(\E,\U)\otimes H^{\bullet+n'}(\E',\U')$ homomorphism (Remark~\ref{R:cross}).
\end{proof}

% ----------------------------------------------------------------
\subsection{Associated covariant and contravariant theories}{\label{SS:covariant}}

By definition, the $n^{th}$ graded piece of the contravariant theory
associated to the bivariant theory $H$ is give by
  $$H^n(\XX)=H^n(\X\llra{\id}\X)=\varinjlim_{\sfC(\id_{\X})}H^{n+\rk\E}(\E,\E-\K).$$
The category $\sfC(\id_{\X})$ has a final object $(\X,\X)$, so the
above colimit is isomorphic to $H^n(\X,\X-\X)=H^n(\X)$, the usual
singular cohomology.

The $n^{th}$ graded piece of the covariant theory associated to $H$
is defined to be
  $$H_n(\XX)=H^{-n}(\X \to pt)=\varinjlim_{\sfC(\X)}H^{e-n}(E,E-K)\cong
  \varinjlim_{K \to \X}H_n(K).$$
Here, $\sfC(\X)$ is the category whose object are pairs $(E,K)$
where $E$ is a Euclidean space of dimension $e$ and $K$ is a compact
subspace of $E$ together with a map $K \to \X$. In the latter
colimit, we have used the Spanier-Whitehead duality $H_n(K)\cong
H^{e-n}(E,E-K)$, and the limit is taken over the category of all
maps $K \to \X$ with $K$ a compact topological space that is
embeddable in some Euclidean space.  By the following proposition,
the latter colimit is, indeed, isomorphic to the singular cohomology
$H_n(\X)$.

\begin{prop}
  Let $\X$ be a topological stack. Then, we have a natural
  isomorphism
    $$\varinjlim_{K \to \X}H_n(K)\cong H_n(\X),$$
  where the limit is taken over the
  category of all maps $K \to \X$ with $K$ a compact topological space
  that is embeddable in some Euclidean space.
\end{prop}

%\begin{proof}
%\end{proof}

It is possible to generalize the axiomatic framework for (skew-symmetric) bivariant theories~\cite{FuMac} to include the present case, where products are only defined for a composition $X\llra{f} Y\llra{g} Z$ if $Y\llra{g} Z$ belongs to a subclass of morphisms called adequate. \notsure{See  Appendix~\ref{S:Generalizedbivariant} for the axioms.} Details will appear elsewhere.

% ----------------------------------------------------------------

\section{\Regular embeddings, submersions, and normally nonsingular morphisms}{\label{S:nnsre}}

% ------------------------------------
\subsection{Submersions}{\label{SS:submersions}}

\begin{defn}{\label{D:submersion}}
  Let $p \: \X \to \Y$ be a morphism of differentiable stacks. For $p$ representable, we
  say it is a {\bf submersion} if its base extension along any differentiable map
  $Y \to \Y$ from a manifold $Y$ is a submersion of manifolds.
  (It is enough to check this for one atlas $Y \to \Y$.) If $p$ is not necessarily representable,
  we say that $p$ is a submersion if for some (hence every) atlas $q \: X \to \X$, the composition
  $p\circ q \: X \to \Y$ is a submersion.
\end{defn}

\begin{ex}{\label{E:submersion}}
     The following are some simple examples of submersions.
     \begin{itemize}
       \item[\bf{1.}] A differentiable map $p \: X \to Y$ of manifolds is a submersion 
         in the above sense if  and only if it is a submersion in the usual sense.
       \item[\bf{2.}] Any projection $\X \times \Y \to \X$ is a submersion.
       \item[\bf{3.}] Let $\E$ be a vector bundle over $\X$. Then, the base map $p \: \E \to \X$
         is a submersion. More generally, if $p$ is an affine bundle (for example, a surjection
         $\E \to \F$ of vector bundles over a base $\X$), then $p$ is a submersion.
       \item[\bf{4.}]    Let $p \: X \to \X$ be an atlas for the differentiable stack $\X$. Then $p$ is a submersion.
           In other words, if $\X=[X_0/X_1]$ is the quotient stack of a differentiable groupoid $[X_1 \toto X_0]$,
            then the quotient map $p \: X_0 \to \X$ is a submersion.        
     \end{itemize} 
\end{ex}

\begin{lem}{\label{L:submersion}}
   Let $p \: \X \to \Y$ and $q \: \Y \to \Z$ be  submersions. Then, we have the following.
  \begin{itemize}
      \item[1.]  The composition $q\circ p \: \X \to \Z$ is a submersion.
      \item[2.]  For  an arbitrary morphism  $\Y' \to \Y$ of differentiable stacks, 
         the base extension $p' \: \X' \to \Y'$ of $p$ is a submersion.
  \end{itemize}   
\end{lem}

\begin{lem}{\label{L:submersiontransversal}}
   Consider the 2-cartesian diagram of differentiable stacks
            $$\xymatrix@=16pt@M=8pt{ \X' \ar@{^(->}[r]^{j} \ar[d]_{p} & \Y' \ar[d]^{q} \\
                          \X   \ar@{^(->} [r]_{i} & \Y    }$$
   in which the horizontal  morphisms are embeddings and $q$ is a submersion.
   Then, $q$ is transversal to $i$. That is,  $p^*\N_{\XX/\YY}=\N_{\XX'/\YY'}$
   (equivalently, the excess bundle $\E$ is trivial).              
\end{lem}

\begin{proof}
  Precomposing $q$ with an  atlas $Y' \to \Y'$, we are reduced to the case where $q$ is representable.
  By making a base change along an atlas $Y \to \Y$, we reduce further to the case where
  we have a diagram of smooth manifolds, in which case the result is clear.
\end{proof}

% ----------------------------------------------------
\subsection{\Regular embeddings}{\label{SS:regular}}
 
In differential topology existence of a tubular neighborhood for a submanifold is a strong tool 
which allows one to linearize the situation by passing to the normal
bundle of the submanifold. Unfortunately, this tool is not always available in the world of 
differentiable stacks, as a substack may not necessarily admit a tubular neighborhood.
To our knowledge, the only situation  where existence of tubular neighborhoods is guaranteed is when
the ambient differentiable 
stack is the quotient stack of a compact Lie group action on a smooth manifold (Example \ref{E:Gtubular}).

In this section, we introduce a class of
embeddings of differentiable stacks, called {\em \regular embeddings}, 
which behave as if they have tubular neighborhoods.
We begin with a preliminary definition.

\begin{defn}{\label{D:tubular}}
   We say that an embedding $i \: \X \hookrightarrow \Y$ of 
   topological stacks admits a {\bf tubular neighborhood} if there is a
   vector bundle $\N$ over $\X$ and a factorization
    $$\X  \stackrel{s}{\hookrightarrow} \N  \stackrel{j}{\hookrightarrow}   \Y$$
 for $i$, where $s$ is the zero section of  $\N$ and $j$ is an open embedding.
 (Note that the vector bundle $\N$ is 
canonically isomorphic to the normal bundle $N_{\XX/\YY}$.)
\end{defn}

\begin{ex}
  Let $\X$ be a topological stack and $\E$ a vector bundle over $\X$. Let
  $s \: \X \to \E$ be the zero section. Then $s$ admits
  a tubular neighborhood. The normal bundle and
  the tubular neighborhood of $\X$ in $\E$ are both $\E$ itself. 
\end{ex}

\begin{ex}{\label{E:Gtubular}}
  Let $\Y=[Y/G]$ be the quotient stack
  $[Y/G]$ of a  topological group $G$ action on  a topological space $Y$. Let $\X \subseteq \Y$ be a
  closed substack of $\Y$, and let $X \subseteq Y$ be the corresponding invariant subspace. Then,
  tubular neighborhoods of $\X$ in $\Y$ are in bijection with $G$-equivariant tubular neighborhoods
  of $X$ in $Y$. In particular, if $\Y$ is a differentiable stack which is isomorphic to the quotient
  stack $[Y/G]$ of a compact Lie group $G$ action on a smooth manifold $Y$, then any differentiable substack $\X$
  of $\Y$ admits a tubular neighborhood. This follows from 
  the $G$-equivariant tubular neighborhood theorem; see  (\cite{Bredon}, $\S$\;VI,
  Theorem 2.2) and also the proof of   Proposition \ref{P:nns}. 
\end{ex}
 
Embeddings which admit tubular neighborhoods have the expected nice properties, but they are not
flexible enough for our purposes. For instance, composition of two such embeddings does not appear 
to admit a tubular neighborhood in general. Also, pullback of such an embedding $\XX \hookrightarrow \YY$
along a submersion onto $\YY$ (even along the base map of vector bundle $\EE \to \YY$) does not
seem to always admit a tubular neighborhood. Definition \ref{D:regular} 
is devised to fix these deficiencies.

\begin{defn}{\label{D:compatible}}
     Let $\YY$ be a differentiable stack. Let  $i \: \XX \hookrightarrow \YY$ a differentiable substack with normal bundle
     $\N=N_{\XX/\YY}$.
     Let $c \in H^k(\N,\N-\XX)$ be a cohomology class.
     We say that a class $\bar{c} \in H^k(\YY,\YY-\XX)$ is {\em compatible} with $c$ if for every
     differentiable atlas $q \: Y \to \YY$, the class $q^*(c) \in H^k(N,N-X)$ corresponds to the class
     $q^*(\bar{c} ) \in H^k(Y,Y-X)$ under the isomorphism
     $H^k(N,N-X)\cong H^k(Y,Y-X)$ obtained by identifying $N:=p^*\N\cong N_{X/Y}$ with a tubular 
     neighborhood of $X:=p^{-1}\XX$ in $Y$ (and applying excision). 
               $$\xymatrix@=16pt@M=8pt{N \ar@{-}[r] & X \ar@{^(->}[r] \ar[d]_{p} & Y \ar[d]^{q} \\
                        \N \ar@{-}[r] & \X   \ar@{^(->} [r]_{i} & \Y    }$$
\end{defn}

Since tubular neighborhoods of submanifolds are unique up to isotopy,
the isomorphism $H^k(N,N-X)\cong H^k(Y,Y-X)$ in the above definition is 
independent of the choice of the tubular neighborhood.

Lemma \ref{L:compatible} justifies the above definition. Before proving it we quote
a useful Lemma from \cite{Zhu}.

\begin{lem}{\label{L:Zhu}}
  Let $\X$ be a differentiable stack and $M$ a smooth manifold. 
  Let $f \: M \to \X$ be a continuous map (i.e., a morphism of underlying topological stacks).
  Then, there exits a differentiable atlas $q \: X \to\X$ such that $f$ lifts to a continuous
  map $\tilde{f} \: M \to Y$. If $f$ is differentiable (i.e., a morphism of differentiable stacks),
  then  $\tilde{f}$ can also be taken to be differentiable.
\end{lem}

\begin{proof}
  The case where $f$ is differentiable is Lemma 3.10 of \cite{Zhu}. 
  The case where $f$ is only continuous is proved using the same argument given 
  in {\em loc. cit}.
\end{proof}

\begin{lem}{\label{L:compatible}}
  Notation being as in Definition \ref{D:compatible}, the class $\bar{c} \in H^k(\YY,\YY-\XX)$ is
  unique (if it exists).
\end{lem}

\begin{proof}
   The statement is true when $\YY$ is a manifold. So, it follows that if $\bar{c}_1,  \bar{c}_2 \in H^k(\YY,\YY-\XX)$ 
   are two such classes, the difference $d := \bar{c}_1-\bar{c}_2$ has the property that 
   $q^*(d) \in H^k(Y,Y-X)$ is zero for every differentiable atlas $q \: Y \to \YY$. We claim that this can only happen
   if $d=0$. 
   
   Suppose that $d \in H^k(\YY,\YY-\XX)$ is a nonzero cohomology class. 
   Choose a classifying space $\varphi \: Y_0 \to \YY$ for $\YY$ (see Section \ref{SS:classifying}).
   Then $\varphi^*(d)  \in H^k(Y_0,Y_0-X_0)$ is nonzero, where $X_0=\varphi^{-1}(\X) \subseteq Y_0$. We can find a finite simplicial
   complex $K$ and a map $f \: K \to Y_0$ such that $f^*\varphi^*(d) \in H^k(K,K-L)$ is nonzero, where $L=f^{-1}\varphi^{-1}(\X)
   \subseteq K$. By embedding $K$ in some Euclidean space and choosing a small tubular neighborhood $M$ of $K$  which retracts to $K$,
   we obtain a manifold $M$ and a continuous map $g \: M \to Y_0$ such that $g^*\varphi^*(d)$ is nonzero. Therefore, we have succeeded
   in finding a manifold $M$ and a morphism (of topological stacks) $\varphi \circ g \: M \to \YY$ which sees $d$. 
   It follows from Lemma
   \ref{L:Zhu} that the map $\varphi \circ g \: M \to \YY$ factors through a differentiable atlas $q \: Y \to \YY$. So, 
   $q$ also sees $d$, that is $q^*(d) \in H^k(Y,Y-X)$ is nonzero, which is what we wanted to prove. 
\end{proof}

\begin{defn}{\label{D:regular}}
    We say that an embedding $\X \hookrightarrow \Y$ of 
   differentiable stacks is a {\bf \regular embedding} if for every
   orientation class $\mu \in H^n(\N,\N-\XX)$  (i.e., a Thom class as in Definition \ref{D:orientablevb}),
   there is a class $\bar{\mu} \in H^n(\YY,\YY-\XX)$  compatible with it in the sense
   of Definition \ref{D:compatible}.
   Here, $\N=\N_{\XX/\YY}$ is the normal bundle to $\XX$ in $\YY$ and $n$ is its rank.
\end{defn}

\begin{lem}{\label{L:regular}}
    \Regular embeddings enjoy the following properties:
    \begin{itemize}
       \item[1.] Any embedding of smooth manifolds is a \regular embedding.
       \item[2.] Any embedding $i \: \X \hookrightarrow \Y$ which admits
         a tubular neighborhood (Definition \ref{D:tubular}) is a \regular embedding. 
     \end{itemize}    
\end{lem}
  
\begin{proof}
    By excision, $H^{\bullet}(\N,\N-\XX)\cong H^{\bullet}(\YY,\YY-\XX)$. This implies (2), and (2) implies (1). 
\end{proof}

\begin{lem}[Composition]{\label{L:regularcompose}}
        If $i \: \X \hookrightarrow \Y$ and $j \: \Y \hookrightarrow \Z$ are \regular embeddings, 
       then $j\circ i \: \X \hookrightarrow \Z$ is also a \regular embedding. Moreover, if $\mu_i \in 
       H^m(\N_{\XX/\YY},\N_{\XX/\YY}-\XX)$ and  $\mu_j\in 
       H^n(\N_{\YY/\ZZ},\N_{\YY/\ZZ}-\YY)$ are orientations  (Definition \ref{D:orientablevb}) and
       $\overline{\mu_i} \in 
       H^m(\YY,\YY-\XX)$ and   $\overline{\mu_j}\in 
       H^n(\ZZ,\ZZ-\YY)$ are the corresponding compatible classes, then the induced orientation
       $\mu_{j\circ i}\in 
       H^{m+n}(\N_{\XX/\ZZ},\N_{\XX/\ZZ}-\XX)$ (see Lemma \ref{L:summand2} and 
       Example \ref{E:normalbundle2})
       is compatible with $\overline{\mu_i} \cdot \overline{\mu_j}$, where the latter product is the one
       of Proposition \ref{P:dot}. In other words, $\overline{\mu_{j\circ i}}=\overline{\mu_i} \cdot \overline{\mu_j}$.
\end{lem}

\begin{proof}
      As we saw in Example \ref{E:normalbundle2}, we have a short exact sequence
         $$ 0 \to  N_{\XX/\YY}  \to N_{\XX/\ZZ} \to i^*N_{\YY/\ZZ} \to 0$$
      of vector bundles over $\X$.   By Lemma \ref{L:summand2}, the orientation classes
      $\mu_i$ and $i^*\mu_j$ induce an orientation class 
      $\mu_{j\circ i}:=\mu_i\cdot p^*i^*(\mu_j)$ on $N_{\XX/\ZZ}$,
      where $p$ stands for the map of pairs $(N_{\XX/\ZZ},N_{\XX/\ZZ}-N_{\XX/\YY}) \to
      (i^*N_{\YY/\ZZ},i^*N_{\YY/\ZZ} -\XX)$. To prove the lemma, we have to show that
      $\overline{\mu_i} \cdot \overline{\mu_j}$ is compatible with $\mu_{j\circ i}$. 
      By definition of compatibility, it is enough to verify this after pulling back everything
      along an atlas $Z \to \ZZ$. That is, it is enough to prove the lemma in the case of
      smooth manifolds. By choosing a suitable tubular neighborhood
      for the embedding $X \hookrightarrow Z$, we can further reduce
      to the case where the given embeddings are zero sections of vector bundles,
      in which case the result is trivial.
\end{proof}

\begin{lem}[Pullback]{\label{L:regularpullback}}
    Consider the 2-cartesian diagram 
                        $$\xymatrix@=16pt@M=8pt{ \X' \ar@{^(->}[r]^{i'} \ar[d]_{p} & \Y' \ar[d]^{q} \\
                              \X   \ar@{^(->} [r]_{i} & \Y    }$$
    Suppose that $i$ is a \regular embedding and $q$ is a submersion (Definition
    \ref{D:submersion}). Then, $i'$ is a \regular embedding.    
     Furthermore, if $\mu \in H^n(\N_{\XX/\YY},\N_{\XX/\YY}-\XX)$
     is an orientation class and $\bar{\mu} \in H^n(\YY,\YY-\XX)$ is compatible with $\mu$, then
     $q^*(\bar{\mu})$ is compatible with the pullback orientation on $\N_{\XX'/\YY'}=p^*\N_{\XX/\YY}$. 
     Here, $q^*\: H^n(\YY,\YY-\XX) \to H^n(\YY',\YY'-\XX')$
     is the induced map on relative cohomology.
\end{lem}

\begin{proof}
      It follows from
      the definition of compatibility (Definition \ref{D:compatible}) that an embedding
      $i \: \X \hookrightarrow \Y$ is a \regular embedding if and only if its base extension along every
      atlas $Y \to \Y$ is. This, combined with Lemmas \ref{L:submersion}.1 and
      \ref{L:submersiontransversal}, proves the lemma.
\end{proof}

In the above lemma, the case we are particularly interested in is where $q$ is an affine bundle (e.g,  when $q$
is the base map  of a vector bundle, see Example \ref{E:submersion}.3).

% ----------------------------------------
\subsection{Normally nonsingular morphisms of stacks and oriented stacks}{\label{S:nns}}

\begin{defn}{\label{D:nns}}
 We say that a representable morphism $f \: \X \to \Y$ of stacks is
{\bf normally nonsingular}, ({\bf nns} for short), if there exist vector bundles $\N$ and  
$\E$ over the stacks $\X$ and $\Y$, respectively, and a commutative
diagram
 $$\xymatrix@=16pt@M=8pt{ \N \ar@{^(->}^i [r] & \E \ar[d]^p \\
                          \X\ar[u]^s \ar[r]_{f} & \Y     }$$
where $s$ is the zero section of  $\N$,  $i$ is an open embedding,
and $\E$ is oriented.  When $\N$ is also oriented, we say that the
diagram is {\bf oriented}. (For the definition of on orientation on
a morphism $f$ see Definition \ref{D:oriented} below.)
 The integer $c = \rk \N - \rk \E$ depends only on $f$ and
is called the {\bf codimension} of $f$. 
(Note that in the case where $\E$ is
of rank zero this coincides with Definition \ref{D:tubular}, so $f$ admits a tubular neighborhood.)
\end{defn}

A diagram as above is called a {\em normally nonsingular diagram}
for $f$. The vector bundle $\N$ is sometimes referred to as the {\em
normal bundle} of $\X$ in $\E$, and $i(\N)$ as a {\em tubular
neighborhood} of $\X$ in $\E$. 
%As a matter of convention, whenever we
%talk about a regular embedding $f \: \X \hookrightarrow \Y$, we implicitly
%fix a tubular neighborhood for it as well (that is, an nns diagram with
%$\E$ of rank zero).

The following are two extreme examples of nns morphisms.

\begin{ex}{\label{E:nns1}}
  Let $\X$ be a topological stack and $\E$ a vector bundle over $\X$. Let
  $s \: \X \to \E$ be the zero section. 
    Then, the diagram
     $$\xymatrix@=16pt@M=8pt{ \E \ar@{^(->}^{\id} [r] & \E \ar[d]^{\id} \\
                          \X \ar[u]^s \ar[r]_{s} & \E    }$$
    is an nns diagram for $s$. Here, we are regarding $\E$ as a rank zero
    vector bundle over itself. The diagram is oriented if and only if $\E$ is. The
    codimension of $s$ is equal to $\rk \E$.  (For future use, let
    us also record the fact that $s$ is a strongly proper morphism.)               
\end{ex}

\begin{ex}{\label{E:nns2}}
  Let $\X$ be a topological stack and $\E$ an oriented vector bundle over $\X$. Let
  $p \: \E \to \X$ be the base map. 
  Then, the diagram
     $$\xymatrix@=16pt@M=8pt{ \E \ar@{^(->}^{\id} [r] & \E \ar[d]^{p} \\
                          \E \ar[u]^{\id} \ar[r]_{p} & \X    }$$
    is an oriented nns diagram for $p$. 
    Here, we are regarding $\E$ as a rank zero
    vector bundle over itself. The codimension of $p$ is
    equal to $-\rk \E$. The normal bundle and the tubular neighborhood
    of $\E$ in $\E$ are both $\E$ itself. 
  \end{ex}

\begin{prop}{\label{P:nns}} Let $G$ be a compact Lie group, and $X$
  and $Y$ smooth $G$-manifolds, with $\X=[X/G]$ and $\Y=[Y/G]$ the corresponding
  quotient stacks. Assume further that $X$ is of finite orbit type.
  Then, for every $G$-equivariant smooth map $X \to Y$, the induced
  morphism $f \: \X \to \Y$ of quotient stacks is normally
  nonsingular. %If both $X$ and $Y$ are oriented, and the
%  orientations are preserved under the $G$-action, then $f$ is
%  oriented.
\end{prop}

\begin{proof}
  First we claim that,
   there is a vector bundle $\V \to BG$ and a smooth embedding $j \:
   \X
   \to \V$, as in the following commutative diagram:
      $$\xymatrix@=16pt@M=8pt{ & & \V \ar[d]^p \\
         \X \ar@{..>}^j [urr] \ar[r]_f & \Y \ar[r]_{\pi_{\Y}}  & BG     }$$
   This statement is equivalent to the fact that every
   $G$-manifold $X$ of finite orbit type embeds $G$-equivariantly into a linear
   $G$-representation $V$ (\cite{Bredon}, $\S$\;II, Theorem 10.1).
   We can arrange for the $G$-action on $V$ to be orientation
   preserving by simply replacing $V$ with $V\oplus V$.

   Let $\E:=\Y\times_{BG}\V$ be the pullback of $\V$ over $\Y$.
   We obtain the following commutative diagram
     $$\xymatrix@=16pt@M=8pt{ &  \E \ar[d]^p \\
                         \X \ar@{..>}^{(f,j)} [ur] \ar[r]_f & \Y    }$$
   Observe that $(f,j)$ is a smooth closed embedding (this can be checked by pulling
   back the whole picture along a chart, say $* \to BG$, for $BG$).
   Let $\N$ be the normal bundle of $(f,j)(\X)$ in $\E$. By
   the existence of $G$-equivariant tubular neighborhoods (\cite{Bredon}, $\S$\;VI,
   Theorem 2.2), we find a vector bundle $\N$ over $\X$ and an open
   embedding $i \: \N \to \E$ making the following diagram
   commutative:
     $$\xymatrix@=24pt@M=8pt{ \N \ar@{^(->}^i [r] & \E \ar[d]^p \\
                          \X\ar[u]^s  \ar@{..>}|{(f,j)} [ur] \ar[r]_{f} & \Y     }$$
 This is exactly what we were looking for.
\end{proof}

\begin{ex}
   The action of a finite group on a manifold has finite orbit
   type. More interestingly,
    the action of a compact Lie group on a
   manifold  whose $\bbZ$-coefficient homology groups are finitely generated
   has finite orbit type. This is Mann's Theorem, see \cite{Bredon},
   $\S$\;IV.10.
\end{ex}

\begin{defn}{\label{D:strong}}
     Let $f \: \X \to \Y$ be a \superproper morphism.
  A bivariant class $\theta \in H(\X\llra{f}\Y)$, not necessarily homogenous, is called a {\bf
  strong orientation} if for every $g \: \Z \to \X$, multiplication
  by $\theta$ is an isomorphism $H(\Z\llra{g\circ f}\Y)\risom
  H(\X\llra{f}\Y)$.\footnote{Note that this definition can be made
  for every adequate morphism in a (generalized) bivariant theory. However we do not need this generality}
\end{defn}

\begin{defn}{\label{D:oriented}}
 A \superproper morphism $f \: \X \to \Y$ of topological stacks
 is called {\bf strongly oriented}, if it is  normally
 nonsingular and it is endowed with a strong orientation $\theta_f
 \in H^c(f)$, where $c=\codim f$;   see Definition \ref{D:nns}.
 A topological stack $\X$ is
  called  strongly oriented if the diagonal $\Delta \: \X \to \X\times \X$ is
  strongly oriented. In this case, we define $\dim\X:=\codim\Delta$.
\end{defn}

\begin{rem}
As we have avoided the discussion of  2-vector bundles
in this paper, we will not  give an intrinsic definition of orientation in terms of the tangent
2-vector bundle of a differentiable stack $\XX$. However, we point out that an orientation
for $\XX$ (in the bivariant sense) amounts to an orientation for 
the ``tangent complex'' of $\XX$, by which we (rather imprecisely) mean the
anchor map $E \stackrel{\rho}\to TX$ of the Lie algebroid  associated to a Lie groupoid presentation 
$X_1\toto X$  for $\XX$; see Example~\ref{E:tgtbundle}. (By an orientation on a complex $V_1 \to V_0$ of vector bundles 
on a manifold we simply mean an orientation on $V_1\oplus V_0$. Note that this sloppy
definition works because we are using singular (co)homology.) Conversely,
assuming that $\Delta\: \XX \to \XX\times\XX$ has a normally nonsingular
diagram as in Definition \ref{D:nns}, an orientation
for  $E \stackrel{\rho}\to TX$ gives rise to an orientation for $\Delta$ (hence, by definition,
for $\XX$).
\end{rem}

\begin{lem}{\label{L:supercompose}}
  Let $f \: \X \to \Y$ and $g \: \Y \to \Z$ be \superproper
  morphisms, 
  and let $\theta \in H(f)$ and $\psi \in H(g)$ be strong
  orientation classes. 
  Then, $\theta\cdot\psi$ is a strong
  orientation class  for $g\circ f \: \X \to \Z$. (Note that
  $g\circ f$ is \superproper by Lemma \ref{L:superproper}.)
\end{lem}

%\begin{proof}
 % Trivial.
%\end{proof}

\begin{lem}{\label{L:unit}}
  Let $f \: \X \to \Y$ be a \superproper map and  $\theta \in H(f)$
  a strong orientation class for it. Then multiplication by $\theta$
  induces an isomorphism $H(\X) \risom H(f)$.
  If $\theta' \in H(f)$ is another orientation class for $f$,
  then there is a unique unit $u \in H(\X)$ such
  that $\theta'=u\cdot\theta$.
\end{lem}

%\begin{proof}
 % Trivial.
%\end{proof}

The following result states that an oriented normally nonsingular
diagram gives rise a canonical strong orientation.

\begin{prop}{\label{P:class}}
 Let $f \: \X \to \Y$ be a \superproper  morphism
 of topological stacks equipped with  an oriented normally nonsingular diagram.
 Then, $f$ has a canonical strong orientation class $\theta_f \in  H^c(f)$,
 where $c=\codim f$.\comment{This is true for every adequate
 morphism. Do we need that?}
\end{prop}

\begin{proof}
  The proof is essentially the same as the one given in \cite{FuMac}. (The `strongly proper' assumption  is a 
  technical condition we need to impose on $f$ in order to be able to multiply bivariant classes. This does not come      
  up  in \cite{FuMac} as they only use trivial vector bundles when defining bivariant classes and the decent condition of 
  Definition \ref{D:superproper} is automatic
  in this case.)
\end{proof}

\begin{ex}[Euler class]{\label{E:Euler}}
  Let $\X$ be a topological stack and $\E$ an oriented vector bundle over $\X$.
  Let
  $s \: \X \to \E$ be the zero section. As we saw in Example \ref{E:nns1},
  $s$ is a strongly proper morphism equipped with a natural nns diagram.
  It follows from Proposition \ref{P:class} that $s$ has a canonical strong
  orientation class $\theta \in H^n(s)$, where $n=\rk \E$. Consider the
  following 2-cartesian diagram:
     $$\xymatrix@=16pt@M=8pt{ \X \ar[r]^{\id} \ar[d]_{\id} & \X \ar[d]^{s} \\
                          \X   \ar[r]_{s} & \E    }$$
  The pullback $s^*(\theta) \in H^n(\id_{\X})=H^n(\X)$ is the {\em Euler class}
  of $\E$.                        
\end{ex}

\begin{lem}{\label{L:nnsembedding}}
  Let $i \: \X \hookrightarrow \Y$ be an  embedding of codimension $n$ of differentiable stacks. 
  If $i$ is nns then it is  a \regular  embedding (Definition \ref{D:regular}).
  In fact, for any choice of orientation $\mu \in H^n(\N,\N-\X)$ for the normal bundle $\N=\N_{\XX/\YY}$, 
  the canonical strong orientation class $\theta_i\in H^n(i \: \X \hookrightarrow \Y)
  \cong H^n(\Y,\Y-\X)$ (Proposition \ref{P:class})
  coincides with the compatible class $\bar{\mu} \in H^n(\Y,\Y-\X)$ (Definition \ref{D:compatible}).
\end{lem}

\begin{proof}
   Pick an nns diagram 
      $$\xymatrix@=16pt@M=8pt{ \N' \ar@{^(->}^j [r] & \E \ar[d]^p \\
                          \X\ar[u]^s \ar[r]_{i} & \Y     }$$
   for $i$.  We have  a natural isomorphism $\N'\cong \N_{\X/\E}$.  
   This gives rise to a split  short exact sequence  
     $$0 \to \N \to \N' \to \E|_{\X} \to 0$$
    of vector bundles over $\X$.
   Recall that, by definition, $\E$ is oriented. By regarding
   $\N'$ as a vector bundle over $\N$ which is the pullback of $ \E|_{\X}$ 
   along the base map $\N \to \X$, we obtain  isomorphisms
             $$H^{\bullet}(\Y,\Y-\X) \cong H^{\bullet+r}(\E,\E-\X) \cong H^{\bullet+r}(\E,\E-js(\X)),$$          
  where $r=\rk E$.  In the last equality  we have used Lemma \ref{L:section}.
  Also, we have isomorphisms
     $$H^{\bullet}(\N, \N-\XX) \cong H^{\bullet+r}(\N',\N'-\X)\cong H^{\bullet+r}(\E,\E-js(\X)),$$
  where in the last equality  we have used excision because $\N'$ can be identified with
  an open substack of $\E$ via $j$.  Under the above identifications, the orientation class
  $\mu \in H^n(\N, \N-\XX)$ corresponds to its compatible class  $\bar{\mu} \in H^n(\Y,\Y-\X)$
  and the strong orientation class $\theta_i\in H^n(i \: \X \hookrightarrow \Y)=H^{n+r}(\E,\E-js(\X))$.  
  (This latter equality is the very definition of bivariant cohomology.)
%  (We have  implicitly      used Lemma \ref{L:summand2} as well.)      
\end{proof}

%\begin{lem}{\label{L:supercompose1}}
%  Let $f \: \X \hookrightarrow \Y$ and $g \: \Y \hookrightarrow \Z$ be  
%  oriented regular embeddings, 
%  and let $\theta \in H(f)$ and $\psi \in H(g)$ be the canonical strong
%  orientation classes for them.
%  Suppose that $g\circ f \: \X \hookrightarrow \Z$ is also a 
%  regular embedding. Then $g\circ f$ is naturally oriented and 
%  $\theta\cdot\psi$ is the canonical strong
%  orientation class for it.
%\end{lem}

%\begin{proof}
%   By taking tubular neighborhoods we are reduced to the case where
%   $f$ and $g$ are zero sections of oriented vector bundles, in which case
%   the result follows from  Lemma \ref{L:summand} (also see 
%   Example \ref{E:excess}.
%\end{proof}

The following proposition shows that any morphism between strongly
oriented topological stacks has a natural strong orientation.
Proposition \ref{P:multiplicative} shows that this class is
multiplicative.

\begin{prop}{\label{P:canonicalorientation1}}
   Let $f \: \X \to \Y$ be a \superproper normally nonsingular morphism of
   topological stacks, and assume that $\X$ and $\Y$ are both
   strongly
   oriented (Definition \ref{D:oriented}). Let  $d=\dim \X$ and
    $c=\dim \Y-\dim \X$. Then, there is a unique strong orientation class
    $\theta_f \in H^c(f)$ which satisfies the equality
   $\theta_f\cdot\theta_{\Y}=(-1)^{cd}\theta_{\X}\cdot(\theta_f\times\theta_f)$,
   as in the diagram
        $$\xymatrix@=26pt@M=8pt{ \X \ar[r]_f^*+<4pt>[o][F-]{\theta_f}
             \ar[d]^{\Delta}_*+<4pt>[o][F-]{\theta_{\X}} & \Y \ar[d]_{\Delta}^*+<4pt>[o][F-]{\theta_{\Y}} \\
                          \X\times\X  \ar[r]^{f\times f}_*+<4pt>[o][F-]{\theta_f\times\theta_f} & \Y\times\Y     }$$
\end{prop}

\begin{proof}
  By Proposition \ref{P:class}, there exists a strong orientation
  $\theta$ for $f$. It is easy to see that $\theta\times\theta$
  is a strong orientation for $f\times f \: \X\times\X \to
  \Y\times\Y$. By Lemma \ref{L:supercompose}, $\theta\cdot\theta_{\Y}$ and
  $\theta_{\X}\cdot(\theta\times\theta)$ are both strong orientation
  classes for $\X \to \Y\times\Y$. Therefore, by Lemma \ref{L:unit},
  there is a unit $u \in H^0(\X)$ such that
  $\theta\cdot\theta_{\Y}=u\cdot\theta_{\X}\cdot(\theta\times\theta)$.
  It follows that $\theta_f\:=(-1)^{cd}u\cdot\theta$ has the desired property;
  see Lemma \ref{L:diagonal} below.
\end{proof}

\begin{lem}{\label{L:diagonal}}
  Let $\X$ be a topological stack and  $\theta \in
  H(\X\llra{\Delta}\X \times\X)$. Let $u,v \in H^0(\X)$,
  and let $u\times v \in H^0(\X)\times H^0(\X)$ be their exterior
  product. Then, $\theta\cdot(u\times v)=u\cdot v\cdot \theta$,
  as classes in $H(\Delta)$.
\end{lem}

\begin{proof}
   Since $u\times v=(u\times 1)\cdot(1\times v)$, it is enough
   to prove the statement in the case where $v=1$.
   Recall that $u\times1$ is defined via independent pullback in
   the righthand square in the diagram
       $$\xymatrix@=28pt@M=8pt{\X \ar[r]_{\Delta}^*+<4pt>[o][F-]{\theta} \ar[d]^{\id}_*+<4pt>[o][F-]{u} &
          \X\times\X \ar[r]_{\pr_2} \ar[d]^{\id}_*+<4pt>[o][F-]{u\times 1} & \X \ar[d]_{\id}^*+<4pt>[o][F-]{u} \\
        \X \ar[r]^{\Delta}_*+<4pt>[o][F-]{\theta} & \X\times\X \ar[r]^{\pr_2} & \X }$$
The equality follows from the skew commutativity of the bivariant theory
applied to the lefthand square.
\end{proof}

\begin{prop}{\label{P:multiplicative}}
  Assume $f \: \X \to \Y$ and $g \: \Y \to \Z$ are  \superproper
 normally nonsingular
  morphisms of strongly oriented topological stacks. Let  $\theta_f \in
  H^c(f)$,  $c=\codim f$,
  and $\theta_g \in H^d(g)$, $d=\codim g$, be the strong orientations
  constructed in Proposition \ref{P:canonicalorientation1}.
 Then, $g\circ f$ is a \superproper normally nonsingular. Furthermore,
  $\theta_f\cdot\theta_g=\theta_{g\circ f}$.
\end{prop}

\begin{proof}
  By Lemma \ref{L:superproper}, $g\circ f$ is \superproperpt.
   Consider the normally nonsingular
diagrams for $f$ and $g$
$$\xymatrix@=16pt@M=8pt{ \N \ar@{^(->}^i [r] & \E \ar[d]  \\
                          \X\ar[u] \ar[r]_{f} & \Y    } \ \
                          \xymatrix@=16pt@M=8pt{ \M \ar@{^(->}^j  [r] & \F \ar[d]  \\
                          \Y\ar[u] \ar[r]_{g} & \Z    }$$
   By adding a vector bundle to $\E$, we may assume that $\E=g^*(\E')$ for some
   orientable vector bundle $\E'$ over $\Z$.
   The following is a normally nonsingular diagram for $g\circ f$
  $$\xymatrix@=16pt@M=8pt{ f^*\M \oplus  \N \ar@{^(->}^k [r] & \F\oplus\E' \ar[d]  \\
                          \X\ar[u] \ar[r]_{g\circ f} & \Z    }$$
  where $k$ is the composite
    $$\xymatrix@=16pt@M=8pt@C=28pt{ f^*\M \oplus  \N \ar@{^(->}^{(\pr,i)} [r]
                   & \M\oplus\E  \ar@{^(->}^{(j,\pr)} [r] &
                   \F\oplus\E'. }$$
  This proves that $g\circ f$ is normally nonsingular.

    The equality $\theta_f\cdot\theta_g=\theta_{g\circ f}$ follows from
    the identity in
    Proposition \ref{P:canonicalorientation1}.
\end{proof}

If $\X$ is strongly oriented  (Definition~\ref{D:oriented}), its iterated diagonals $\Delta^{(n)}\:\X \to \X^{n}$ are \superproper (Example~\ref{E:proper}).
\begin{cor}\label{C:orienteddiagonal}
Let $\X$ be a oriented stack. Then the diagonals $\Delta^{(n)}\:\X \to \X^{n}$ are canonically strongly oriented.
\end{cor}

\begin{prop}{\label{P:Gorinetation}}
   Notation being as in Proposition \ref{P:nns}, assume further that
   $X$ and $Y$ are oriented and that the $G$-actions are orientation
   preserving.
  Then, every normally nonsingular diagram for $f \: \X \to \Y$
   is naturally  oriented. In particular,
   when $f$ is \superproperpt, we have a strong  orientation class
    $\theta_f \in  H^c(f)$,  $c=\dim Y-\dim X$. Furthermore, this
    class is independent of the choice of the normally nonsingular diagram.
\end{prop}

\begin{proof} Let us first fix a notation: given a manifold $X$ with
an action of $G$, we denote $[TX/G]$ by $T\X$. (So, $T\X$ does
depend on $X$, and not just on $\X$. Since in what follows all
stacks are quotients of a $G$-action on a given manifold, this
should not cause confusion.)

Consider a normally nonsingular diagram
    $$\xymatrix@=16pt@M=8pt{ \N \ar@{^(->}^i [r] & \E \ar[d]^p \\
                          \X\ar[u]^s \ar[r]_{f} & \Y     }$$
as  in the proof of Proposition \ref{P:nns}.
 We show that $\N$ is naturally oriented.
 By Lemma \ref{L:summand2}, there is a natural orientation on
$T\E$, because it fits in the following short exact sequence
  $$0 \to p^*\E \to T\E \to p^*T\X\to 0.$$
In particular, we have an orientation on $f^*(T\E)$. We have an
isomorphism of vector bundles over $\X$
   $$T\X\oplus\N \cong f^*(T\E).$$
It now follows from Lemma \ref{L:summand} that $\N$  also carries a
natural orientation.
  This proves the first part of the proposition. In
 particular, when $f$ is proper, we obtain a class
 $\theta_f \in H^{c}(\X\llra{f}\Y)$
 as in Proposition \ref{P:class}.

 Now, we show that the class $\theta_f$ is independent of the
 normally nonsingular diagram above.
Consider another oriented normally
 nonsingular diagram for $f$
   $$\xymatrix@=16pt@M=8pt{ \M \ar@{^(->}^j [r] & \F \ar[d]^q \\
                          \X\ar[u]^t \ar[r]_{f} & \Y     }$$
We have to show that the following diagram commutes
  $$\xymatrix@=16pt@M=8pt{H^{\scriptstyle \bullet}(\X)  \ar[r]^(0.34){\cong} \ar[d]_{=}
  &
       H^{\scriptstyle \bullet+\rk\N}(\E,\E-\X) \ar[d]^{\cong} \\
                          H^{\scriptstyle \bullet}(\X) \ar[r]^(0.34){\cong} &
       H^{\scriptstyle \bullet+\rk\M}(\F,\F-\X)
      }$$
 where the horizontal isomorphisms are the one of Proposition
 \ref{P:class}, and the vertical isomorphism is the one of Lemma
 \ref{L:isom}.
 First we prove a special case.

\vspace{0.1in}

\noindent{\em Special case.}
  Assume $\E=\F$, and $is=tj$. In this case, we can
 choose a third vector bundle $\LL \to \X$ and an open embedding $k
 \: \LL \hookrightarrow \E$ that factors through both $\N$ and $\M$.
 The two orientations induced on $\LL$ from $\M$ and $\N$, as in
 Lemma \ref{L:Thomex}, are the same (because they are equal to the
 orientation induced from $\E$, as  described
 above). The claim now follows from the commutative diagram of Lemma
 \ref{L:Thomex} (applied once to the open embedding $\LL
 \hookrightarrow \N$
 and once to the open embedding $\LL \hookrightarrow \M$).

\vspace{0.1in}

\noindent{\em General case.}
  To prove the general case, we make use of the
 following auxiliary oriented nonsingular diagrams:
   $$\xymatrix@C=30pt@R=16pt@M=8pt{
       \N\oplus f^*\F \ar@{^(->}^(0.52){(i,\pr)} [r] & \E\oplus\F  \ar[d] \\
                          \X\ar[u]^{(s,jt)} \ar[r]_{f} & \Y     } \ \ \ \
   \xymatrix@C=30pt@R=16pt@M=8pt{
       \M\oplus f^*\E \ar@{^(->}^(0.52){(j,\pr)} [r] & \F\oplus\E \ar[d] \\
                          \X\ar[u]^{(t,is)} \ar[r]_{f} & \Y     }$$
Here, the two maps $\pr$ stand for the projection maps
$f^*\F=\X\times_{\Y}\F \to \F$ and $f^*\E=\X\times_{\Y}\E \to \E$.
 Let us denotes the ranks of $\E$, $\F$, $\N$, and
$\M$ by $e$, $f$, $n$, and $m$. (Hopefully, presence of two
different $f$ in the notation will not cause confusion!)
 The
first normally nonsingular diagram gives rise to the following
commutative diagram of isomorphisms: {\small
$$\xymatrix@=16pt@M=8pt{
   H^{\scriptstyle \bullet}(\X)\ar[d]_{=}\ar[r]^(0.25){\cong} \ar@/^1pc/@<2ex>[rr]^{\varphi} &
          H^{\scriptstyle \bullet+n+f}(\N\oplus f^*\F,\N\oplus f^*\F-\X)
             \ar[d]_{\cong}\ar[r]^(0.54){\cong}         &
          H^{\scriptstyle \bullet+n+f}(\E\oplus \F,\E\oplus\F-\X) \ar[d]^{\cong} \\
   H^{\scriptstyle \bullet}(\X) \ar[r]^(0.4){\cong}            &
          H^{\scriptstyle \bullet+n}(\N,\N-\X) \ar[r]^{\cong}  &
          H^{\scriptstyle \bullet+n}(\E,\E-\X)}$$
}
 The commutativity of the left square is   because of Lemma
\ref{L:section2}, and the commutativity of the right square is
because  Thom isomorphism (vertical) commutes with  excision
(horizontal).

Similarly, the second normally nonsingular diagram gives rise to the
following commutative diagram of isomorphisms {\small
$$\xymatrix@=16pt@M=8pt{
   H^{\scriptstyle \bullet}(\X)\ar[d]_{=}\ar[r]^(0.25){\cong}  \ar@/^1pc/@<2ex>[rr]^{\psi} &
          H^{\scriptstyle \bullet+m+e}(\M\oplus f^*\E,\M\oplus f^*\E-\X)
             \ar[d]_{\cong}\ar[r]^(0.54){\cong}         &
          H^{\scriptstyle \bullet+m+e}(\F\oplus \E,\F\oplus\E-\X) \ar[d]^{\cong} \\
   H^{\scriptstyle \bullet}(\X) \ar[r]^(0.4){\cong}            &
          H^{\scriptstyle \bullet+m}(\M,\M-\X) \ar[r]^{\cong}  &
          H^{\scriptstyle \bullet+m}(\F,\F-\X) }$$
          }

On the other hand, using the special case that we just proved, the
two normally nonsingular diagrams give rise to the following
commutative diagram:
  $$\xymatrix@C=20pt@R=16pt@M=8pt{
     H^{\scriptstyle \bullet}(\X)\ar[d]_{=}\ar[r]_(0.25){\cong}^(0.25){\varphi}    &
          H^{\scriptstyle \bullet+n+f}(\E\oplus \F,\E\oplus\F-\X) \ar[d]^{=} \\
     H^{\scriptstyle \bullet}(\X)\ar[r]_(0.25){\cong}^(0.25){\psi}    &
          H^{\scriptstyle \bullet+m+e}(\F\oplus \E,\F\oplus\E-\X)
    }$$
The general case now follows from combining this diagram with (the
other rectangles) of the previous two diagrams.
\end{proof}

\begin{cor}{\label{C:nns}}
Let $\X$ be a stack that is equivalent to the quotient stack $[X/G]$
of smooth orientation preserving action of a compact Lie group $G$
on a smooth oriented  manifold $X$ having finitely generated
homology groups. Then, the diagonal $\X \to \X \times \X$ is
naturally oriented. In particular, the diagonal of the classifying
stack $BG$ of a compact Lie group $G$ is naturally oriented.
\end{cor}

\begin{rem}
  Let $\X$, $\Y$ and $f$ be as in Proposition \ref{P:Gorinetation}.
  There are two ways of giving a strong orientation to $f$. Either
  we can use Proposition \ref{P:Gorinetation} directly, or we
  first apply Corollary \ref{C:nns} to endow $\X$ and $\Y$ with
  a strong orientation, and then apply Proposition
  \ref{P:canonicalorientation1}. The  orientations we get are the same for
  $f$. We denote $\theta_f$ this strong orientation.
\end{rem}

\begin{prop}\label{C:orbifoldnnns}
Let $\XX$ be a paracompact  orbifold whose tangent bundle (Example \ref{E:tgtbundle}) is oriented. 
Then the diagonal $\XX\to \XX \times \XX$ is strongly oriented and in particular, $\XX$ is naturally oriented.
\end{prop}
\begin{proof}
Locally, we can find a tubular neighborhood for the diagonal. The result follows using partition of unity.
\end{proof}

% ----------------------------------------------------------------

\section{Gysin maps}

As in \cite{Ch, CoJo, CoVo}, the main step in our construction of 
the \BV-structure on the homology of the loop stack is the systematic
development of Gysin maps for oriented morphisms of stacks. To this end, we
use (a slightly generalized version of) Fulton-MacPherson's bivariant theory. This
is in spirit very close to Chataur's bordism approach which relies on Jakob's
bivariant theory for differentiable manifolds \cite{Ch}, although bivariant
theories are not explicitly mentioned in \cite{Ch}. 

% ----------------------------------------------------------------
\subsection{Construction of the Gysin maps}
{\label{S:Gysin}}

We recall the construction of Gysin homomorphisms associated to a
bivariant class \cite{FuMac}.

Fix an element $\theta \in H^i(\X\llra{f}\Y)$. Let $u:\Y'\to \Y$ be
an arbitrary morphism of topological stacks  and $\X'=\X \times_{\Y}
\Y'$ the base change given by the cartesian square:
  \begin{equation}\label{basechange}
     \xymatrix@C=8pt@R=10pt@M=8pt{\XX'\dto
      \rto^{f'} &  \YY' \dto^{u}\\
      \XX \rto^{f} & \YY   .}
   \end{equation} 
Then $\theta$ determines {\bf Gysin homomorphisms}
   $$\theta^!\: H_j(\Y') \to H_{j-i}(\X')$$
and
   $$\theta_!\: H^{j}(\X') \to H^{j+i}(\Y').$$
For the cohomology Gysin map, we need to assume that $f'$ is adequate.
 These homomorphisms are defined by
   $$\theta^!(a)=\big(u^*(\theta)\big)\cdot a, \ \ \ \text{for} \ \ a \in
   H_{j}(\Y')=H^{-j}(\Y' \to pt),$$
and
  $$\theta_!(b)=f_*'\big(b\cdot u^*(\theta)\big), \ \ \ \text{for} \ \ b \in
   H^{j}(\X')=H^{j}(\X'\llra{\id}\X').$$
The homology Gysin map is defined because the map $\X' \to *$ is
adequate (see Example \ref{E:super}).

\smallskip

%\noindent{\bf Notation:} When $\X \llra{f} \Y$ is strongly oriented and $\theta=\theta_f$ is the strong orientation
%(Proposition \ref{P:class}), we denote  the Gysin map $({\theta_f})^!$ by $\gy{f}$.

% -------------------------------------------------
\subsection{Standard Properties of Gysin maps}
\label{PropertiesGysin}

By Proposition \ref{P:class}, when the map $f \: \X \to \Y$ in Diagram~\eqref{basechange}
is strongly oriented, it has a canonical strong orientation $\theta_f$. In this case, we have
a canonical Gysin morphism 
        $$\gy{f} :=({\theta_f})^!\: H_{\bullet}(\Y') \to H_{\bullet-c}(\X'),$$
where $c$ is the codimension of $f$. In this subsection we collect some of the
standard properties of these Gysin morphisms.

\begin{enumerate}
  \item{\bf Functoriality.}  Assume given  a commutative diagram of cartesian squares
     \eq
           \xymatrix{
           \XX'\dto \rto &  \YY' \dto \rto & \ZZ' \dto\\
            \XX \rto^{f} & \YY \rto^{g} & \ZZ}
     \eneq 
with $f\: \XX \to \YY$ and $g\: \YY \to \ZZ$  strongly oriented of codimensions $c$ and $d$, 
respectively. 
Then, the induced Gysin morphisms $\gy{f} \: H_{\bullet}(\Y') \to H_{\bullet-c}(\X')$ and
$\gy{g} \: H_{\bullet}(\Z') \to H_{\bullet-d}(\Y')$
satisfy the functoriality identity
            $$\gy{(g\circ f)}=\gy{f} \circ \gy{g}.$$

\item{\bf Naturality.} Assume given a a commutative diagram of cartesian squares
    \eq\label{naturality}
       \xymatrix{\XX'' \dto_{v} \rto & \YY'' \dto^{u}  \\
       \XX'\dto \rto &  \YY' \dto \\
       \XX \rto^{f} & \YY }
    \eneq
with $f$ strongly oriented. Then, the induced Gysin morphisms satisfy
             $$v_* \circ \gy{f}=\gy{f} \circ u_* .$$

\item{\bf Commutation with cross product.} Given two cartesian squares 
          $$\xymatrix{\XX'_1\dto \rto &  \YY'_1 \dto^{u_1}\\
            \XX_1 \rto^{f_1} & \YY_1} \quad \quad 
            \xymatrix{\XX'_2\dto \rto &  \YY'_2 \dto^{u_2}\\
            \XX_2 \rto^{f_2} & \YY_2}$$
consider the induced product square
        \eq        
            \xymatrix{\XX'_1\times\XX'_2\dto \rto &  \YY'_1\times\YY'_2 \dto^{u_1\times u_2}\\
            \XX_1\times\XX_2 \rto^{f_1\times f_2} & \YY_1\times\YY_2} 
        \eneq
If  $f_1$ and $f_2$ are strongly oriented, then so is $f_1\times f_2$. Moreover, the  three 
Gysin morphisms satisfy the equation 
       $$\gy{(f_1 \times f_2)}(-\times -)=\gy{f_1}(-)\times \gy{f_2}(-).$$

\item{\bf Commutation with pullback.} Given a cartesian square
     \eq 
        \xymatrix{\XX'\dto \rto^{f'} &  \YY' \dto^{u}\\
        \XX \rto^{f} & \YY  }
      \eneq
 with $f$ strongly oriented and $y'\in H^*(\YY')$, we have
 $$\gy{f}(y'\cap-)  = (-1)^{\deg(y')\codim(f)} f'^*
 (y') \cap \gy{f}(-).$$
\end{enumerate}

\noindent {\bf Proof of Properties 1,2,3 and 4}

Everything follows from the axioms of a bivariant theory. By
Proposition~\ref{P:multiplicative}, the products of two strongly
oriented maps is canonically strongly oriented. Thus Property 1
follows from Axiom~{\bf A13}. Property 2 is Axiom~{\bf A3} followed
by Axiom~{\bf A123}. Taking direct products of vector bundles shows
that \superproper and normally nonsingular morphisms are stable by
products. Hence Property 3 follows from the definition and naturality of the cross
product (Proposition~\ref{P:bivariantcross}). Property 4 is a consequence of the skew-commutativity.

\begin{numrmk}[Cohomology Gysin maps]
When  in  Diagram~\eqref{basechange} $f'$ is adequate, there is an induced cohomology
Gysin map $\gy{f}\:H\com(\XX') \to H^{\scriptstyle \bullet+c}(\YY')$. Properties 1,2,3
and 4 above have obvious analogs in cohomology  when all the relevant 
maps involved are adequate.  Recall
that a strongly oriented map is \superproper hence adequate.
\end{numrmk}

\begin{numrmk}
We have emphasized the case of strongly oriented maps for simplicity
and because it is sufficient for our purpose. Nevertheless, by
pullback axiom, any bivariant class $\theta \in
H^r(f)$ yields a bivariant class
$u^*(\theta)\:H^{r}(f')$ and thus  a Gysin map
$H_{\bullet}(\YY') \to H_{\bullet - r}(\XX')$. Properties 1,2,3 and 4
above will hold true in this more general setting.
\end{numrmk}

% --------------------------
\subsection{A special case: $G$-equivariant Gysin maps}

Let  $M$, $N$ be oriented compact manifolds and $G$  a Lie group
  acting on $M$, $N$ by orientation preserving diffeomorphisms.   By
  Proposition~\ref{P:Gorinetation}, if
  $f\:M \to  N$ is a $G$-equivariant map, then $f$ is canonically
  strongly oriented.  Gysin maps
  for equivariant (co)homology were already considered, for example,
  by Atiyah and Bott~\cite{AB}.
  
\begin{prop}\label{prop:AB} The Gysin maps $f_!, f^!$ associated to $f$ in
  (co)homology coincide with the equivariant Gysin maps in the sense
  of Atiyah and Bott~\cite{AB}.
\end{prop}

\begin{pf}
Gysin map in~\cite{AB} are obtained by the use of fiber integration
and Thom classes over the spaces $M_G=M\times_G EG$ and $N_G=N\times_G
EG$. These spaces are respectively classifying spaces of the stacks
$[M/G]$ and $[N/G]$ and thus are respectively the pullbacks
$[M/G]\times_{[*/G]}BG$, $[N/G]\times_{[*/G]}BG$. The pullback of the
normally nonsingular diagram of Proposition~\ref{P:nns} along the natural maps
$[M/G]\times_{[*/G]}BG\to [M/G] $   and $[N/G]\times_{[*/G]}BG\to
[N/G]$ yields a bundle $ \N_G=\N\times_{[*/G]}BG$ over $M_G$ and a
bundle $\E_G=\E\times_{[*/G]}BG$ over $N_G$. This defines a nonsingular
diagram for the induced map $f_G\: M_G\to N_G$. Unfolding the
definition of bivariant classes, it is straightforward
to check that the Gysin map associated to the strong orientation class
of Proposition~\ref{P:Gorinetation} is induced by the Thom isomorphism
associated to the bundle $ \N_G$ over $M_G$.
\end{pf}

Let $G$ be a subgroup of a finite (discrete) group $H$. Let
$Y$ be a manifold endowed with a (right)
$H$-action (and thus a $G$-action).
Consider the quotient stacks $[Y/H]$ and $[Y/H]$.   There  are well known ``transfer maps''
${\rm tr}_H^G\:H_*^H(Y)\to H_*^G(Y)$ (see~\cite{Ben})

\begin{lem}\label{Transfer}
When $G$ is a finite group, the Gysin map associated to the cartesian square
    $$\xymatrix{[Y/G]\dto \rto & [Y/H]  \dto\\
               [*/G] \rto & [*/H]   }$$
where the lower map is induced by the inclusion $ G\hookrightarrow H$,  is the
usual ``transfer map'' $H^{H}_{*}(Y)\to H^{G}_{*}(Y)$ in equivariant homology.
\end{lem}

\begin{pf}
The space $Y\times H$ is endowed with a natural right $H$-action given
by $(y,h).k=(y.k,k^{-1}h)$ as well as a right $G$-action
$(y,h).g=(y,hg)$. These two actions commutes hence we can form the
quotient stack $[Y\times H/H\times G]\cong [Y\times (H/G)
/H]$. Clearly the map $(y,h)\mapsto yh$ is equivariant with respect to
the $G$ action on the target and $H\times G$-action on the
source. One easily checks that this map induces an equivalence
$[Y\times (H/G) /H] \cong [Y/G]$. We are thus left to study the Gysin
map of an equivariant covering with fibers the set $H/G$. The argument
of  Proposition~\ref{prop:AB} easily shows that it coincides with the
usual transfer maps for coverings by a finite group and thus with the
transfer.
\end{pf}

 Assuming we take coefficient
in a field of characteristic coprime with $|H|$ for the singular homology, we have
$$H\lcom([Y/H])\cong H_*^{H}(Y)\cong
\big(H_*(Y)\big)_{H}.$$ In that case the map ${\rm
tr}_H^G\:\big(H_*(Y)\big)_{H} \to \big(H_*(Y)\big)_{G}$ is explicitly
given by \eq \label{eq:transfer} {\rm tr}_H^G(x)&=&\sum_{h\in H{/G}}
h.x. \eneq

% -------------------------------------------------
\subsection{The excess formula}
\label{SS:excess}

The main result of this subsection is the following.

\begin{prop}[Excess formula] \label{prop:Excess} \label{pullbacksubmsersion} 
Consider the 2-cartesian diagram
        $$\xymatrix@=16pt@M=8pt{\X'\ar[d]_{p} \ar@{^(->} [r]^{j} &  \Y' \ar[d]^{q}\\
                 \X \ar@{^(->} [r]_{i} & \Y   }$$ 
in which $i$ and $j$ are \regular
embeddings with normal bundles $\N$, $\N'$, respectively. 
Let $\E=p^*(\N)/\N'$ be the excess bundle (see Section \ref{SS:tangent}). Fix orientations on $\N$ 
and $\N'$ and endow $\E$ with the induced orientation as in Lemma \ref{L:summand2}. (In the case where 
$\N$ and $\N'$
have equal ranks the orientation on $\N'$ is uniquely determined by the one on $\N$.) 
Let $\theta_i \in H^n(\Y,\Y-\X)=H^n(i)$ and $\theta_j \in H^n(\Y',\Y'-\X')=H^{n'}(j)$
be the classes compatible with the orientations on $\N$ and $\N'$, respectively (Definition \ref{D:compatible}).
Then,
\eq \label{eq:excessformula}
q^*(\theta_i) = e(\E)\cdot \theta_j,
\eneq
where $e(\E)\in H\com(\XX')$ is the Euler class of $\E$ (see Example \ref{E:Euler}).
\end{prop}

\begin{proof}
  In the case where $q$ is a submersion the proposition follows from 
  Lemmas \ref{L:submersiontransversal} and \ref{L:regularpullback}. We use this to reduce the problem
  to the case of manifolds.

  By Lemmas \ref{L:submersiontransversal}, \ref{L:compatible},  and \ref{L:regularpullback}, 
  and the fact that bivariant product commutes with pullback, 
  it is enough  to prove the formula after passing to an arbitrary atlas $Y' \to \Y'$. So, we may assume that 
  $\Y'=:Y'$ and $\X'=:X'$ are smooth manifolds. By Lemma \ref{L:Zhu}, we can find  an atlas
  $Y \to \Y$ through which $q \: Y' \to \Y$ factors. Since the atlas $Y \to \Y$ is a
  submersion and the proposition is true for sumbersions, we may assume, after pulling back
  everything along $Y \to \Y$,  that $\Y=:Y$ and $\X=:X$ are also smooth manifolds.
  From now on, we use the notation $N$, $N'$ and $E$ instead of $\N$, $\N'$ and $\E$.

  We are reduced to proving the result in the case of manifolds.
  Since $q \: Y'\to Y $ factors as the composition 
    $$Y'\llra{(1,q)}  Y' \times Y \lra Y$$ 
  of an embedding and a submersion, it is enough, by
  functoriality of pullbacks, to consider the cases where
  $q$ is a submersion and $q$ is an embedding separately. 
  The former case is easy, as $\E$ is the zero bundle  
  and $q^*(\theta_i)=\theta_j$ by Lemma \ref{L:regularpullback}.
  
  It remains to prove the proposition in the case where $q$ is an embedding. By choosing appropriate
  tubular neighborhoods, we reduce to the case where $Y=F$ is a vector bundle over $X$
  and $i \: X \to F$ is the zero section. We may also assume that $X'$ is a submanifold of $X$
  and $Y'=N'$ is a vector bundle over $X'$ which is a
  subbundle of $F|_{X'}$, having zero section $j \: X' \to N'$. 
  Moreover, after choosing a metric on $F$, we may write
  $F|_{X'}=E\oplus N'$. Finally,  replacing $F$ with $F|_{X'}$, we may assume that $X=X'$.
  Summarizing all the reduction we have made, we are in a situation where we 
  have a manifold $X$, with vector bundles $N'$ and $E$ on it, so that the 2-cartesian square
  of the proposition has the form
            $$\xymatrix@=16pt@M=8pt{ X \ar[d]_{p=\id} \ar@{^(->} [r]^{j} &  N' \ar[d]^{q}\\
                 X \ar@{^(->} [r]_(0.35){i} &  E\oplus N'  }$$ 
  Here, the horizontal maps are zero sections and $q$ is the inclusion of the summand $N'$.       
   We expand this square to the 2-cartesian diagram
                    $$\xymatrix@=20pt@M=8pt{ X \ar@{^(->} [r]^{\id} \ar[d]_{\id} & X \ar[d]_s \ar@{^(->} [r]^{j} &  N' \ar[d]^{q}\\
                         X \ar@{^(->} [r]^s \ar@/_1pc/ [rr]_i &  E \ar@{^(->} [r] &  E\oplus N'  }$$ 
  where $s$ stands for the zero section.  Let $\theta_{E}$ be the strong
  orientation class of $X \hookrightarrow E$, and $\theta'$ the strong orientation class of $E \hookrightarrow E\oplus N'$.
  By Lemma \ref{L:regularcompose}, we have $\theta_i=\theta_{E}\cdot\theta'$.
  Then, since pullback respects bivariant product, and $s^*(\theta_E)=e(E)$ (Example \ref{E:Euler}), we find 
    $$q^*(\theta_i) = e(E)\cdot q^*(\theta').$$
  Making the rightmost square in the above diagram upside down and using the obvious projection maps,
  as in the diagram
               $$\xymatrix@=16pt@M=8pt{  E \ar@{^(->} [r] \ar[d] &  E\oplus N' \ar[d]^{\pi}\\
                                X  \ar@{^(->} [r]^{j} &  N'  }$$ 
  we see that $E \hookrightarrow E\oplus N'$ has $\pi^*(\theta_j)$ as its canonical strong orientation,
  that is $\theta'=\pi^*(\theta_j)$. Hence, $q^*(\theta')=\theta_j$, and the above displayed
  formula becomes $q^*(\theta_i) = e(E)\cdot \theta_j$, which is the desired excess formula.
\end{proof}

It is worth  noticing that when $i$ and $j$ are nns, then, by Lemma \ref{L:nnsembedding}, the classes $\theta_i \in H^n(i\: \X \to \Y)$
and $\theta_j \in H^n(j \: \X' \to \Y')$ are precisely the canonical strong orientations constructed in Proposition \ref{P:class}.

The following  immediate corollary of Proposition \ref{prop:Excess} is useful in computing Gysin maps.

\begin{cor}{\label{C:excess}}
  Consider the 2-cartesian diagram
        $$\xymatrix{\XX''\dto_u \rto &  \YY'' \dto \\
         \XX'\dto_{p} \rto^{j} &  \YY' \dto^{q}\\
        \XX \rto_{i} & \YY   }$$ 
in which the lower square is as in Proposition \ref{prop:Excess}.
Let $n$ and $n'$ be the ranks of $\N$ and $\N'$ respectively.
Let 
      $$i^!\: H_{\bullet}(\YY'') \to H_{\bullet-n}(\XX''),$$  
      $$j^! \: H_{\bullet}(\YY'') \to H_{\bullet-n'}(\XX'')$$ 
be the corresponding Gysin maps. Then, for any $c \in  H_{\bullet}(\YY'')$ we have the equality
       $$i^!(c)=u^*e(\E)\cdot j^!(c).$$ 
In particular, if $q$ is transversal to $i$ (e.g., when $q$ is a submersion), then
$i^!=j^!$.
\end{cor}

% ----------------------------------------------------------------
\section{The loop product}\label{Loopproduct}
In this section we consider (Hurewicz) {\em strongly oriented} stacks
(Definition~\ref{D:oriented}). We  obtain a loop product on the
homology of the free loop stack of an oriented stack which
generalizes Chas-Sullivan product for  the homology of a loop
manifolds~\cite{CS}. Recall that  a stack $\XX$ is called strongly
oriented if the diagonal $\Delta\:\XX \to \XX \times \XX$ has a
strong orientation class (Definition~\ref{D:strong}). For instance,
oriented manifolds and oriented orbifolds are oriented stacks.
 More generally,  the quotient stack of a compact Lie group
acting by orientation preserving automorphisms on an oriented
manifold is an oriented stack.
 
Note that it is possible to have two different group actions, 
say a Lie group $G$ acting on a manifold $X$ and 
another Lie group $H$ acting on another manifold $Y$, 
which give rise to the same quotient stacks, i.e.,
$[X/G]\cong[Y/H]$.  By definition, our notion of
orientation, as well as our construction of the loop product (and 
all other string operations that we construct), 
are independent of the choice of the presentation 
and only depend on the resulting quotient stack.
Put differently, and slightly more generally, what we do is 
we use a   Morita invariant  notion of orientation for  Lie groupoids,
and for such oriented Lie groupoids we construct Morita invariant
string operations.

\subsection{Construction of the loop product}\label{Construction}
Let $\XX$ be a Hurewicz oriented stack of finite dimension $d$. The
construction of the loop product $$H\lcom(\LXX)\otimes
H\lcom(\LXX)\to  H\lcom(\LXX)$$ is divided into  3 steps.
\begin{description}
\item[Step 1] There is a well-known external product (called the ``cross product'')
$$H_p(\LXX)\otimes H_q(\LXX) \stackrel{S}\longrightarrow H_{p+q}(\LXX). $$

\item[Step 2] The diagonal $\Delta\: \XX\to \XX\times \XX$ and the
  evaluation map $\ev_0:\LXX\to \XX$~\eqref{eq:ev0} yield
  the cartesian square
  \eq \label{eq:cartesianloop}\xymatrix@=16pt@M=8pt{ \LXX\times_{\XX}\LXX \ar[r]\ar[d] &
                          \LXX\times \LXX \ar[d]^{(\ev_0,\ev_0)} \\
                          \XX \ar[r]^(0.4){\Delta} & \XX\times\XX
                          .}\eneq
We will usually denote by $e:\LXX\times \LXX \to \XX\times \XX$ the
                          map $(\ev_0,\ev_0)$. Since $\XX$ is Hurewicz, Corollary~\ref{L:8} implies  that there is a natural equivalence of stacks
    $$\LXX\times_{\XX}\LXX\cong\map(8,\XX),$$ where the figure ``8'' stands for the topological stack
    associated to the topological space $S^1\vee S^1$. The wedge
    $S^1\vee S^1$  is taken with respect to the basepoint $0$ of $S^1$.
Since $\XX$ is oriented, its diagonal $\Delta\:\XX\to \XX\times \XX$ is oriented normally nonsingular and
    according to    Section~\ref{S:Gysin}, there is
                          a Gysin map
$$\gy{\Delta}\:H_{\scriptstyle \bullet}(\LXX\times \LXX) \to
H_{\scriptstyle \bullet-d}(\LXX\times_{\XX}\LXX)\cong H_{\scriptstyle \bullet-d}(\map(8,\XX)). $$

\item[Step 3]
The map $S^1 \to S^1\vee S^1$ that pinches $\frac{1}{2}$ to $0$,
induces a natural map of stacks $m: \map(8,\XX) \to \LXX$, called the {\it
  Pontrjagin multiplication}.  Hence we have an induced  map
on homology
$$m_*\:H_{\scriptstyle \bullet}(\map(8,\XX) \to  H_{\scriptstyle \bullet} (\LXX ).$$
\end{description}

 We define the {\em loop product} to be
the following composition

\begin{equation}\label{defloopproduct}
 H_p(\LXX)\otimes H_q(\LXX)\stackrel{S}\to H_{p+q}(\LXX \times
 \LXX)\stackrel{\gy{\Delta}}\to
 H_{p+q-d}(\map(8,\XX))\stackrel{m_*}\to  H_{p+q-d}(\LXX).
\end{equation}

\begin{them}\label{th:Loop}Let $\XX$ be an oriented (Hurewicz\footnote{recall that every differentiable stack is Hurewicz.}) stack of dimension $d$. The
loop product induces a structure of associative and graded commutative
algebra for the shifted homology
$\BH\lcom(\LXX):=H_{\scriptstyle \bullet+d}(\LXX)$.
\end{them}

%\begin{defn}\label{defnloopproduct}
%The Loop product $\star: H_p({\rm L}\Gamma)\otimes H_q({\rm L}\Gamma)
%\stackrel{}\to H_{p+q-d}({\rm L}\Gamma )$
%  is the map~\eqref{deflooproduct} {\it i.e.} the composition  $m_*\circ G_{\Delta}^{e}\circ S$.
%\end{defn}

The loop product is of degree $d=\dim(\XX)$ because the Gysin map
involved in Step 2 is of degree $d$. If  we denote
$\sH\lcom(\LXX):=H_{\scriptstyle \bullet +\dim(\XX)}(\LXX)$ the shifted homology
groups, then the loop product induces a degree 0 multiplication
$\sH\lcom(\LXX)\otimes \sH\lcom(\LXX) \to \sH\lcom(\LXX)$.

\medskip

  Indeed one can introduce  a ``twisted'' version of  loop product.
 Let $\alpha$ be a
class in $\bigoplus_{r\geq 0} H^r(\LXX\times_{\Gammaa}\LXX)$. The {\em twisted
loop product} $\star_{\alpha}\:H\lcom(\LXX)\otimes H\lcom(\LXX)
\to H_{\scriptstyle \bullet}(\LXX) $ is defined, for all $x,y\in
H\lcom(\LXX)$,
$$x\star_\alpha y = m_*\big(\gy{\Delta}(x\times y)\cap \alpha\big).$$

\begin{numrmk}
The twisted product $\star_\alpha$ is not graded since we do not
assume $\alpha$ to be homogeneous. However, if $\alpha \in
H^r(\LXX\times_\XX\LXX)$ is homogeneous of degree  $r$,  then
$\star_{e}\: H\lcom(\LXX)\otimes H\lcom(\LXX) \to H_{\scriptstyle
\bullet-d-r}(\LXX)$  is of degree $r+\dim(\XX)$.
\end{numrmk}

Let us introduce some notations. We denote, respectively, $p_{12}, p_{23}\:\LXX\times_{\Gammaa} \LXX
 \times_{\Gammaa} \LXX \to \LXX\times_{\Gammaa} \LXX$   the projections on the
   first two and the  last two
 factors. Also let $(m\times 1)\: \LXX\times_{\Gammaa}
\LXX\times_{\Gammaa} \LXX \to \LXX\times_{\Gammaa} \LXX$ and $(1\times m)\:\LXX\times_{\Gammaa}
\LXX\times_{\Gammaa} \LXX \to \LXX\times_{\Gammaa} \LXX$         be   the Pontrjagin
 multiplication of the two first factors and two last factors respectively. Furthermore, there are flip maps $\sigma\: \LXX\times
\LXX \to \LXX\times
\LXX$, $\widetilde{\sigma}\: \LXX\times_{\XX}
\LXX \to \LXX\times_{\XX}
\LXX$ permuting the two factors of $\LXX\times
\LXX$.
\begin{them}\label{Associativitycocycle} Let $\alpha$ be a class in
 $\bigoplus_{r\geq 0}H^r(\LXX\times_{\Gammaa}\LXX)$.\begin{itemize}
 \item  If $\alpha$ satisfies the $2$-cocycle 
condition
\eq \label{eq:Associativitycocycle} p_{12}^*(x)\cup (m\times 1)^*(\alpha)&= & p_{23}^*(\alpha)\cup (1\times m)^*(\alpha)\eneq
in $H\com(\LXX\times_{\Gammaa}\LXX\times_{\Gammaa} \LXX )$, then
$\star_e:H\lcom(\LXX)\otimes H\lcom(\LXX) \to H_{\scriptstyle \bullet}(\LXX)$ is
associative.
\item If $\alpha$ satisfies the flip condition $\widetilde{\sigma}^*(\alpha)= \alpha$, then the twisted Loop product $\star_\alpha\:\sH(\LXX)\otimes \sH(\LXX)\to
 \sH(\LXX)$ is graded commutative.
\end{itemize}
\end{them}

\begin{example} If $E$ is an oriented vector  bundle over a stack $\XX$ it has  a Euler class $e(E)$.  Note that
 the rank may vary on different connected components of $\XX$.
In particular, any vector bundle $E$ over $\LXX\times_{\XX}\LXX$
 defines a twisted
loop product $\star_E:=\star_{e(E)}\:H(\LXX)\otimes H(\LXX) \to
H(\LXX)$. Moreover,
   $\widetilde{\sigma}^*( e(E))\cong e(E)$ whenever $\widetilde{\sigma}^* E \cong E $. Since  identities between Euler classes
 are equivalent to  identities in $K$-theory we have:

\begin{cor}\label{cor:twistedloop}
Let $\XX$ be an oriented (Hurewicz) stack and $E$  a vector
 bundle over $\LXX\times_{\XX} \LXX$.\begin{itemize}
\item  If $E$ satisfies the cocycle condition $$p_{12}^*(E)+(m\times 1)^*(E)= p_{23}^*(E)+(1\times m)^*(E)$$ in $K$-theory, then $\star_E$ is associative.
\item If
 $\widetilde{\sigma}^* E \cong E$, then the twisted Loop product
 $\star_E:H(\LXX)\otimes H(\LXX)\to H(\LXX)$ is
 graded commutative. \end{itemize}
\end{cor}
\end{example}

\begin{numrmk} Let $M$ be an oriented manifold and $G$ a finite group acting on $M$ by orientation preserving diffeomorphisms and $\XX=[M/G]$ be the associated global quotient orbifold.
Using Proposition~\ref{loopdiscrete}, Proposition~\ref{prop:AB} and the argument  of the proof of
Proposition~\ref{loopsphere} below to identify evaluation maps and
Pontrjagin map, it is straightforward to prove
 the Loop product $\star:
  \sH\lcom(\LXX)\otimes
  \sH\lcom(\LXX)\to \sH_{\scriptstyle \bullet}(\LXX)$ coincides with the
  one introduced in \cite{LUX}.
\end{numrmk}

\subsection{Proof of Theorems}\label{Associativity}

The Pontrjagin multiplication $m:\map(8,\XX) \to \LXX$ is induced by
the pinch map $S^1 \to S^1\vee S^1$. The latter is  homotopy
coassociative, thus there is a chain homotopy equivalence between
$$m(m\times \id)\:C\lcom(\LXX\times_{\XX} \LXX\times_{\XX} \LXX)\to
C\lcom(\LXX)$$ and $m(\id \times m)$. 
This proves the next lemma:
\begin{lem}\label{homas}
The Pontrjagin multiplication  satisfies $$m_*\big(({\rm id}\times m)_*\big)=m_*\big((m\times {\rm id})_*\big). $$
\end{lem}

\begin{prop}\label{Associativityloop}
 The loop product $H\lcom(\LXX) \otimes H\lcom(\LXX)\stackrel{\scriptstyle \bullet}\to H_{\scriptstyle \bullet-d}(\LXX)$ is associative.
\end{prop}
\begin{pf}
It is well known that the cross product is associative so
that $$S^{(2)}\:H\lcom(\LXX)\otimes H\lcom(\LXX) \otimes  H\lcom(\LXX)\to H\lcom(\LXX
\times \LXX \times \LXX)$$ is equal to both  $S(S\times 1)$ and $S(1\times S)$.
 We write $m^{(2)}$ for the iterated  map
$$m_*(m\times 1)_*=m_*(1\times m)_*$$ as in Lemma~\ref{homas} and
$\Delta^{(2)}$ the iterated diagonal  $$\Delta(\Delta \times
1)=\Delta(1\times \Delta)\:\XX\to \XX^{\times 3}.$$ Also let
$e^{(2)}\:\LXX^{\times 3} \to \XX^{\times 3}$ denote the product
$\ev_0\times \ev_0\times \ev_0$ of the evaluation map on each component.  It is enough to
prove that, for all  $x,y,z \in H\lcom(\LXX)$, $$(x\bullet
y)\bullet z= m^{(2)}\big(\gy{\Delta^{(2)}} (x\times y\times z)\big)=x\bullet (y\bullet z).$$
The first equality  is given by the commutativity of the following diagram:
{\footnotesize 
\begin{equation}\label{AsD} \hspace{-0.2cm}
\vcenter{\xymatrix@C=8pt@R=10pt@M=7pt{H(\LXX)\otimes H(\LXX) \otimes H(\LXX)\dto_{S \otimes 1} \drto^{S^{(2)}} &&&\\
\ar @{} [dr] |{(3)} H(\LXX \times \LXX)\otimes H(\LXX) \rto^{S}
\dto_{\gy{\Delta}\otimes 1} & \ar @{} [dr] |{(1)} H(\LXX\times \LXX
\times \LXX) \dto^{\gy{\Delta\times 1}}
\rto^{\gy{\Delta^{(2)}}} &H(\LXX \times_{\XX} \LXX \times_{\XX} \LXX)\ar @{=}[d] & \\
H(\LXX \times_{\XX} \LXX)\otimes H(\LXX) \rto^{S} \dto_{m_*\otimes 1} & \ar @{} [dr] |{(2)} H(\LXX \times_{\XX}\LXX \times \LXX) \dto_{(m\times 1)_*} \rto^{\gy{\Delta}} & H(\LXX \times_{\XX} \LXX \times_{\XX} \LXX) \dto^{\widetilde{m\times 1}} \drto^{p^{(2)}}& \\
  H(\LXX)\otimes  H(\LXX) \rto^{S} & H(\LXX \times \LXX) \rto^{\gy{\Delta}} & H(\LXX \times_{\XX} \LXX)\rto^{m_*} & H(\LXX).  }}
\end{equation}}
The commutativity of bottom left square follows from the naturality of the cross
product, and the bottom right  triangle  from the associativity of $m_*$
 according to Lemma~\ref{homas}.  The three remaining squares commutes thanks to  the following reasons:
\begin{description}
\item[(Square 1)] There is a diagram of  cartesian squares
$$\xymatrix{\LXX \times_{\XX} \LXX \times_{\XX} \LXX \rto  \dto &\LXX \times_{\XX} \LXX \times \LXX \dto^{\widetilde{\ev_0}\times \ev_0} \rto & \LXX\times \LXX \times \LXX\dto^{e^{(2)}}  \\
\XX \rto^{\Delta} & \XX \times \XX \rto^{\Delta \times 1} & \XX \times \XX \times \XX}. $$ Thus the commutativity follows from the functoriality of Gysin maps.
\item[(Square 2)] Note that the map $\widetilde{\ev_0}$ in square (1)
   is equal to $\ev_0\circ m$. The commutativity follows, by naturality of Gysin maps,  from the tower of  cartesian diagrams:
$$\xymatrix{ \LXX \times_{\XX} \LXX \times_{\XX} \LXX \dto_{\widetilde{m\times 1}} \rto & \LXX \times_{\XX} \LXX \times \LXX \dto^{m\times 1} \\
\LXX \times_{\XX} \LXX \dto_{\widetilde{\ev_0}} \rto &\LXX \times \LXX \dto^{e} \\
\XX \rto^{\Delta} &\XX \times \XX.}$$
\item[(Square 3)]   It is commutative by compatibility of Gysin maps with the cross product.
\end{description}

Hence it follows that, for all $x,y,z \in H(\LXX)$, one has $(x\bullet y)\bullet z= m^{(2)}(\gy{\Delta^{(2)}}(x\times y\times z))$.

One proves in a similar way the identity  $m^{(2)}(\gy{\Delta^{(2)}}(x\times y\times z))=x\bullet (y\bullet z) $ from which the equation  $(x\bullet y)\bullet z=\bullet (y\bullet z)$ follows.
\end{pf}

\begin{prop}\label{Commutativityloop}
The loop product $\sH\lcom(\LXX) \otimes
\sH\lcom(\LXX)\stackrel{\star}\to \sH_{*}(\LXX) $ is graded commutative.
\end{prop}
\begin{pf} Essentially, this result follows from the homotopy commutativity of the
  Pontrjagin map $m:\LXX\times_{\XX}\LXX=\map(8,\XX)\to \XX$.
   More precisely we need to prove that, for $x\in \sH_p(\LXX),y \in \sH_q  (\LXX)$, we have $$m_*\big(\gy{\Delta}(x\times y)\big)=(-1)^{pq}\big(m_*\big(\gy{\Delta}(y\times x)\big).$$
The tower of pullback squares
$$\xymatrix{\LXX \times_\XX \LXX \rto \dto_{\widetilde{\sigma}}&  \LXX\times \LXX \dto^\sigma\\
\LXX \times_\XX \LXX \rto^{{\rm id}} \dto & \LXX\times \LXX \dto^{e}\\
\XX \rto^{\Delta} & \XX \times \XX  } $$
implies that
 $${\rm \widetilde{\sigma}}_*\circ \gy{\Delta}(x\times y)=
(-1)^{pq}\gy{\Delta}(y\times x).$$
 Here $\widetilde{\sigma}\:\LXX \times_\XX \LXX\to \LXX
 \times_\XX \LXX$ and $\sigma\: \LXX\times \LXX\to \LXX\times \LXX$ are  flip maps.
Hence the result  follows from $m_* \circ \widetilde{\sigma}_*
=m_*$ in homology. The latter is an immediate consequence of the existence of a homotopy
 between the pinch map $p\:S^1 \to S^1\vee S^1$ and $\sigma\circ p\:
 S^1 \to S^1\vee S^1$ obtained by making the base point $0\in S^1$
 goes to $\frac{1}{2}\in S^1$. Passing to the mapping stack functor
 $\map(-,\XX) $ yields a homotopy equivalence  between $m \circ
 \widetilde{\sigma}$ and $m$.
\end{pf}
\begin{rmk}
Note that the homotopy between the two pinch maps does not
 preserve the canonical basepoints. Hence it is crucial to work with
 the free loop stack (in other words with nonpointed mapping stack
 functors) in this proof.
\end{rmk}

\begin{prop}\label{pr:twistedassociativity}
If $\alpha\in H\com(\LXX\times_{\XX} \LXX)$ satisfies the cocycle
equation~\eqref{eq:Associativitycocycle}, then the twisted loop
product $\star_\alpha\: H(\LXX)\otimes H(\LXX)\to
 H(\LXX)$ is associative.
\end{prop}
\begin{pf}
We write $f_1:\LXX\times_{\XX}\LXX \times \LXX
\to\LXX\times_{\XX}\LXX $  and $f_3:\LXX\times
\LXX\times_{\XX}\LXX\to \LXX\times_{\XX}\LXX$ for the canonical
projections. Also we denote $j_3:\LXX \times_{\XX}\LXX \times_{\XX}
\LXX\hookrightarrow\LXX \times_{\XX}\LXX \times \LXX$ and  $j_1:
\LXX \times_{\XX}\LXX \times_{\XX} \LXX \hookrightarrow  \LXX \times
\LXX \times_{\XX} \LXX$ the canonical maps. Using the naturality of
cup product and cross product, we can write an associativity diagram
similar to (\ref{AsD}) for $\star_{\alpha}$, for which the only non
obviously commuting square is the one labeled by (1) which becomes :
$$\xymatrix{ H(\LXX\times \LXX \times \LXX) \dto^{\gy{\Delta\times 1}} \rto^{\gy{1\times\Delta}} &H(\LXX \times_{\XX} \LXX \times_{\XX} \LXX)\rto^{\cap f_3^*(\alpha)} & H(\LXX \times_{\XX} \LXX \times_{\XX} \LXX)\dto^{\gy{\Delta}} \\
 H(\LXX \times_{\XX}\LXX \times \LXX) \dto_{\cap f_1^*(\alpha)}  & &H(\LXX \times_{\XX} \LXX \times_{\XX} \LXX)\dto^{\cap p_{23}^*(\alpha)} \\
H(\LXX \times_{\XX}\LXX \times \LXX)\rto^{\gy{\Delta}} & H(\LXX \times_{\XX}\LXX \times_{\XX} \LXX)\rto^{\cap p_{12}^*(\alpha)} & H(\LXX \times_{\XX} \LXX \times_{\XX} \LXX).}$$
 Since Gysin maps commute with pullback, for any $y\in  H\lcom(\LXX \times_{\XX}\LXX \times \LXX)$,
\begin{eqnarray*}
\gy{\Delta}(y\cap f_1^*(\alpha)) &=& \gy{\Delta}(y)\cap ( f_1\circ j_3)^*(\alpha)\\
&=& \gy{\Delta}(y)\cap (m\times 1)^*(\alpha).
\end{eqnarray*}
Similarly, $\gy{\Delta}(y\cap f_3^*(\alpha))
=\gy{\Delta}(y)\cap (1\times m)^*(\alpha)$.
From  square (1) of diagram~\ref{AsD} we deduce that the commutativity of the square is equivalent to the identity
$$\big( \gy{\Delta}\circ \gy{\Delta\times
  1}\big)\cap (m\times 1)^*(\alpha)\cap p_{12}^*(\alpha) =
\big( \gy{\Delta}\circ \gy{1\times\Delta}\big)\cap (1\times
m)^*(\alpha)\cap p_{23}^*(\alpha)$$ $$\Longleftrightarrow
\gy{\Delta^{(2)}} \cap \big((m\times 1)^*(\alpha)\cup
p_{12}^*(\alpha)\big) = \gy{\Delta^{(2)}}\cap \big((1\times
m)^*(\alpha)\cup p_{23}^*(\alpha)\big) .$$  The last equality
follows immediately from the 2-cocycle
condition~\eqref{eq:Associativitycocycle}.
\end{pf}

\begin{prop}\label{pr:twistedcommutativity}
If $\widetilde{\sigma}^*(\alpha)=\alpha$, then the twisted loop product $\star_\alpha\: \sH(\LXX)\otimes \sH(\LXX)\to
 \sH(\LXX)$ is commutative.
\end{prop}
\begin{pf}
The proof of Propositions~\ref{Commutativityloop} applies verbatim.
\end{pf}

\section{Hidden loop product for family of groups over a stack} \label{S:hidden}
The Chas-Sullivan product  generalizes the intersection
product for a manifold $M$. Indeed, the  embedding of $M$ as the space of
constant loop in $LM$ makes  $H\lcom(M)$ a subalgebra of the loop
homology and the restriction of the loop product to this subalgebra is
the intersection product~\cite{CS}.

 In the context of stacks, there
are more interesting "constant" loops, that is loops which are
constant on the coarse space. From a mathematical physics point of
view, these spaces of loops are sometimes called {\em ghost loops}. The
canonical ghost loop stack is the inertia stack.

\subsection{Hidden loop product}

\label{stringproduct} In this section we construct a \hidden product
for the inertia stack. From the categorical point of view the {\bf
inertia stack} $\IXX$ of a stack $\XX$ is the stack of pairs
$(X,\varphi)$ where $X$ is an object of $\XX$ and $\varphi$ an
automorphism of $X$. If $\XX$ is a Hurewicz topological stack then
so is $\IXX$. However, if $\XX$ is differentiable, $\IXX$ is not
necessarily differentiable.
  Let $\gm$ be    a topological  groupoid representing $\XX$.
 Let $S\gm=\{g\in\gm_1\vert\; s(g)=t(g)\}$ be the space
of closed loops.
There is a natural action of $\gm$ on $S\gm$ by conjugation.
 Thus one  forms the transformation
groupoid $\Lambda\gm\:S\gm\rtimes \gm \toto S\gm$, which  is
always a topological groupoid, called the {\em inertia groupoid}.
It is a presentation of the inertia stack.
 and denoted $\IXX$.
 Indeed
 one obtains the following morphism of groupoids
\begin{equation}\label{doublegroupoid}\xymatrix{\Lambda \Gamma \ar@<-0.5ex>[r]
\ar@<0.5ex>[r]\ar[d] &  S\gm  \ar[d]  \ar[d]\\
\gm_1 \ar@<-0.5ex>[r]  \ar@<0.5ex>[r] & \Gamma_0 .} \end{equation}
There is a morphism of topological groupoids  $\ev_0: \IG \to
\gm$ given, for $(x,g)\in S\gm\rtimes \gm$, by
$\ev_0((x,g))=g\in \gm_1$. This groupoid morphism $\ev_0:\IG\to \gm$
induces the evaluation map \eq \label{eq:ev0string}\ev_0:\IXX\to \XX\eneq
on the corresponding stacks.

\smallskip

 The construction of the \hidden product can be made in 3 steps.
\begin{description}
\item[Step 1] The external product induces  a map:

$$ H_p(\IXX)\otimes H_q(\IXX)\stackrel{S}\to H_{p+q}(\IXX \times \IXX).$$

\item[Step 2] We can form the pullback of the evaluation map  $\ev_0:\IXX\to \XX$
along  the diagonal
  $\Delta\: \XX\to \XX\times \XX$, thus obtaining the cartesian square
  \eq\label{eq:cartesianghost}&&\xymatrix@=16pt@M=8pt{ \IXX\times_{\XX}\IXX \ar[r]\ar[d] &
                          \IXX\times \IXX \ar[d]^{(\ev_0,\ev_0)} \\
                          \XX \ar[r]^(0.4){\Delta} & \XX\times\XX
                          .}\eneq
Again we  denote by $e:\IXX\times \IXX \to \XX\times \XX$ the
                          map $(\ev_0,\ev_0)$. Since $\XX$ is
                         strongly oriented, so is its diagonal
 $\Delta\:\XX\to \XX\times \XX$. Hence we have
                          a Gysin map:
$$\gy{\Delta}\:H_{\scriptstyle \bullet}(\IXX\times \IXX) \to
H_{\scriptstyle \bullet-d}(\IXX\times_{\XX}\IXX). $$

\item[Step 3]  The stack $\IXX\times_{\XX} \IXX$ is known as the {\it double inertia stack}. Its objects are triples $(X,\varphi,\psi)$ where $X$ is an object of $\XX$ and $\varphi, \psi$ are automorphisms of $X$.
On the groupoid level the stack $\IXX \times_{\XX}\IXX $ is presented
 by the transformation groupoid
 $\big(S{\Gamma}\times_{\Gamma_0} S{\Gamma}\big) \rtimes \Gamma_1\toto S{\Gamma}
 \times_{\Gamma_0} S{\Gamma}$ where $\gm_1$ acts on
 $S{\Gamma}\times_{\Gamma_0} S{\Gamma}$ by conjugation diagonally.
 The double inertia stack is endowed with a ``Pontrjagin'' multiplication map $ m:\IXX\times_{\XX} \IXX\to \IXX$ given by $m(X,\varphi,\psi)=(X,\varphi\psi)$.
  It induces a morphism
on homology
$$ m_*\:H\lcom(\IXX\times_{\XX}\IXX) \to H\lcom(\IXX). $$
\end{description}

Composing the three maps in the above
 steps one obtains  a  product
$$\star: H_p(\IXX)\otimes H_q(\IXX)
\stackrel{}\to H_{p+q-d}(\IXX ), $$
 called the \hidden product:
    \begin{equation}\label{defstringproduct}
 H_p(\IXX)\otimes H_q(\IXX)\stackrel{S}\to H_{p+q}(\IXX \times
 \IXX)\stackrel{\gy{\Delta}}\to
 H_{p+q-d}(\IXX\times_{\XX}\IXX)\stackrel{m_*}\to  H_{p+q-d}(\IXX).
\end{equation}

As for the loop product, the \hidden product is a degree $0$
multiplication on the shifted homology
groups: $\sH\lcom(\IXX)=H_{\scriptstyle \bullet+d}(\IXX)$.

\begin{them}
\label{th:stringproduct}
Let $\XX$ be an oriented stack of dimension $d$.
The shifted  homology $\sH\lcom(\IXX)$ of the inertia stack
 is an associative graded
commutative algebra.
\end{them}
Before proving Theorem~\ref{th:stringproduct}, let us remark that the
``Pontrjagin'' map $m:\IXX\times_{\XX} \IXX\to \IXX$ corresponds to
the multiplication  is \underline{associative}. Thus, passing to
homology one has the lemma:
\begin{lem}\label{pisassociative}
$m_*\:H\lcom(\IXX\times_{\XX}\IXX) \to H\lcom(\IXX)$
satisfies the associativity condition:
$$m_*\big(({\rm id}\times m)_*\big)=m_*\big((m\times {\rm id})_*\big). $$
\end{lem}
Less obvious is that it is also \underline{commutative}: indeed there is a $2$-arrow $\alpha$
\eq\label{eq:inertiacommutative}
\xymatrix@C=6pt@R=10pt@M=6pt{ \IXX\times_{\XX} \IXX \ar[rr]^{\quad m} \ar[dd]_{\rm flip}   \ar@{=>}[dr]^{\alpha}&& \IXX \\ &&\\
\IXX\times_{\XX} \IXX \ar[rruu]_{m} & & }
\eneq
 which associates to $(X,\varphi,\psi)$ in the double inertia the isomorphism $\varphi^{-1}$ :
 \eq \label{eq:phiinertia}
 \xymatrix{X \ar[r]^{\varphi \psi} \ar[d]_{\varphi^{-1}} \ar @{} [dr] |{\cong}&  X \ar[d]^{\varphi^{-1}} \\
 X\ar[r]^{\psi \varphi}  & X   }
 \eneq
 between $(X,\varphi\psi)$ and $(X,\psi\varphi)$ in $\IXX$.

\begin{trivlist}\item[]{\sc Proof of Theorem~\ref{th:stringproduct}.}
Associativity follows {\em mutatis mutandis} from the proof of

Theorem~\ref{Associativityloop}, substituting $\LXX$ with $\IXX$ in the
argument. Similarly, the proof of Theorem~\ref{Commutativityloop}
leaves us to proving that the induced map $m_*\circ
\widetilde{\sigma}\:H\lcom(\IXX\times_{\XX}\IXX) \to H\lcom(\IXX)$
in homology
is equal to $m_*$. Here again $\widetilde{\sigma}$ is the  flip
map. Passing to any groupoid $\gm$ representing $\XX$ and denoting $\IG\times_{\gm}\IG= \big(S\gm\times_{\gm_0}S\gm\big)\rtimes \gm_1\toto\big(S\gm\times_{\gm_0}S\gm\big)$, it is enough to
check that  the induced map $$m_*\circ
\widetilde{\sigma}_*\:H\lcom\left(\IG\times_{\gm}\IG\right) \to
H\lcom(\IG)$$ in groupoid homology is equal to $m_*$.
At level of groupoids, the $2$-arrow $\alpha$ of diagram~\eqref{eq:inertiacommutative} yields the identity
\eqn
m_*(\sigma(n_1,n_2)) &=& \mu(n_2,n_1)\\
&=& \big(\mu(n_1,n_2)\big)^{n_2^{-1}}
\eneqn
for all $x=(n_1,n_2,\gamma)\in \big(S\gm\times_{\gm_0}S\gm\big)\rtimes \gm_1$.
 Here $\mu: S\gm\times_{\gm_0}S\gm\to S\gm$ is the restriction of the groupoid multiplication of $\gm$. Thus $m_*(\sigma(n_1,n_2))$ is canonically conjugate to $m_*(n_1,n_2)$
and in a equivariant way. It follows
that after passing to  groupoid homology, one has $m_*=m_*\circ
\widetilde{\sigma}$. An explicit homotopy
$h:C_n(\IG\times_{\gm}\IG)\to C_{n+1}(\IG)$ between $m_*$ and $m_*\circ
\widetilde{\sigma}$ at the chain level is given by $h: \sum_{i=0}^n
(-1)^ih_i$ where \eqn h_i((n_1,n_2),g_1,\dots,g_n)&=&
\left(\mu(n_1,n_2)\big)^{n_2^{-1}}, g_1,\dots\right.\\
&&\left. \dots,g_{i},(g_1\dots
  g_i)^{-1}n_2(g_1\dots g_i), g_{i+1},\dots,g_n \right)\eneqn for $i>0$
and $h_0((n_1,n_2),g_1,\dots,g_n)=
\left(\mu(n_1,n_2)\big)^{n_2^{-1}}, n_2,g_1,\dots,g_{n} \right) $.
\nolinebreak $\Box$ \end{trivlist}

\medskip

If $\alpha$ is  a cohomology  class in
 $\bigoplus_{r\geq 0}H^r(\IXX\times_{\XX} \IXX)$, one defines
the twisted \hidden product $$\star_{\alpha}\:H\lcom(\IXX)\otimes H\lcom(\IXX) \to H_{\scriptstyle \bullet}(\IXX) $$ as follows. For any  $x,y\in H\lcom(\IXX)$,
$$x\star_\alpha y = m_*\big(\gy{\Delta}(x\times y)\cap
\alpha\big).$$

We use similar notations as for Theorem~\ref{Associativitycocycle}:
 denote  $p_{12}, p_{23}\:\IXX\times_{\Gammaa} \IXX
 \times_{\Gammaa} \IXX \to \IXX\times_{\Gammaa} \IXX$   the projections on the
   first two and the  last two
 factors

\begin{prop}\label{twistedstring}
 Let $\alpha$ be a class in $\bigoplus_{r\geq 0} H^r(\IXX\times_{\XX}\IXX)$.
\begin{itemize}
\item If $\alpha $ satisfies the cocycle condition:
\eq \label{eq:stringcocycle}p_{12}^*(\alpha )\cup (m\times 1)^*(\alpha )&=& p_{23}^*(\alpha)
\cup (1\times m)^*(\alpha )\eneq
 in $H\com(\IXX\times_{\XX}\IXX\times_{\XX} \IXX )$,
then $\star_\alpha\:H(\IXX)\otimes H(\IXX) \to H(\IXX)$ is
associative.
\item If $\alpha$ satisfies the flip condition $\widetilde{\sigma}^*(\alpha)= \alpha$, then the twisted \hidden product $\star_\alpha\:\sH(\IXX)\otimes \sH(\IXX)\to
 \sH(\IXX)$ is graded commutative.
\end{itemize}
\end{prop}
\begin{pf}
The argument of Proposition~\ref{pr:twistedassociativity} and
Proposition~\ref{pr:twistedcommutativity} applies.
\end{pf}
Corollary~\ref{cor:twistedloop} has an obvious counterpart for inertia stack.
\begin{cor}\label{cor:twistedstring}
Let $\XX$ be an oriented stack and $E$  a vector
 bundle over $\IXX\times_{\XX} \IXX$.\begin{itemize}
\item  If $E$ satisfies the cocycle condition $$p_{12}^*(E)+(m\times 1)^*(E)= p_{23}^*(E)+(1\times m)^*(E)$$ in $K$-theory, then $\star_E$ is associative.
\item If
 $\widetilde{\sigma}^* E \cong E$, then the twisted Loop product
 $\star_E:H(\IXX)\otimes H(\IXX)\to H(\IXX)$ is
 graded commutative. \end{itemize}
\end{cor}

%%%%%%%%%%%%%%%%%%%%%%%%%%%%%%%%%%%%%%%%%%%%%%%%%%%%%%%%%%%%%%%%%%%%%%%%%%%%%%%%%%%%%%%%%%%%%%%%%%%%%%%

\subsection{Family of commutative groups and crossed modules}\label{crossedmodules}

The \hidden product can be defined for more general  "ghost loops"
stacks than the mere inertia stack. In fact, we can replace the commutative family  $\IXX \to \XX$ by an arbitrary commutative family of groups.

A {\bf family of groups} over a (topological) stack $\XX$ is a
(topological) stack $\GG$ together with a morphism of (topological)
stacks $\ev:\GG \to \XX$ and an associative multiplication
$m:\GG\times_{\XX} \GG \to \GG$. A family of groups $\GG\to \XX$
(over $\XX$) is said to be a {\bf commutative family of groups}
(over $\XX$) if  there exists an invertible $2$-arrow $\alpha$
making the following diagram \eq\label{eq:ghostcommutative}
\xymatrix@C=6pt@R=10pt@M=6pt{ \GG\times_{\XX} \GG \ar[rr]^{\quad m} \ar[dd]_{\rm flip}   \ar@{=>}[dr]^{\alpha}&& \GG \\ &&\\
\GG\times_{\XX} \GG \ar[rruu]_{m} & & }
\eneq commutative. Clearly, the inertia stack is a commutative family of groups (see Equation~\ref{eq:phiinertia}).

\smallskip

In the groupoid language, a commutative family of groups can be represented as follows.  A {\em crossed module}  of (topological) groupoids is a morphism of groupoids
\eqn
\xymatrix{N_1\rto^{i} \ar@<0.5ex>[d] \ar@<-0.5ex>[d]& \gm_1 \ar@<0.5ex>[d] \ar@<-0.5ex>[d]\\
N_0 \rto^{=} & \gm_0 }
\eneqn
which is the identity on the base spaces (in particular $N_0=\gm_0$)
and where $N_1\toto N_0$ is a family of groups (i.e. source and target
are equal),
together with a right action $(\gamma,n)\to n^{\gamma}$ of $\gm$ on $N$ by automorphisms satisfying:
\begin{enumerate}
\item for all $(n,\gamma)\in N\rtimes\gm_1$,
  $i(n^{\gamma})=\gamma^{-1}i(n)\gamma$;
\item for all $(x,y)\in N\times_{\gm_0} N$, $x^{i(y)}=y^{-1}xy$.
\end{enumerate}
Note that the equalities in (1) and (2) make sense because $N$ is a
family of groups. We use the short notation
$[N\stackrel{i}\longrightarrow \gm]$ for a crossed module.
\begin{numrmk} In the literature, groupoids for which source equals target are sometimes called bundle of groups. Since we do not assume the source to be locally trivial, we prefer the terminology family of groups.
\end{numrmk}

Since a crossed module $[N\stackrel{i}\longrightarrow \gm]$
comes with an action of $\gm$ on $N$,  one can form the transformation
groupoid $\ING\:N_1\rtimes\gm_1 \toto N_1$, which is  a
topological groupoid. Furthermore, the projection  $ N_1\rtimes_{\gm_0}\gm_1
\to \gm_1$ on the second factor induces a (topological) groupoid morphism $\ev:\ING
\to \gm$. Let $\GG$ and $\XX$ be the quotient stack $[N_1/N_1\rtimes\gm_1]$ and $[\gm_0/\gm_1]$ respectively. Then $\ev: \GG\to \XX$ is a commutative family of groups over  $\XX$.
Clearly, the inertia stack $\IXX$ corresponds to the crossed module $[S\gm \hookrightarrow \gm]$ for any groupoid presentation $\gm$ of $\XX$.
Obviously $\Lambda [S\gm \hookrightarrow \gm]$
is the inertia groupoid $\IG$.
The inertia stack  is universal among commutative family of groups over $\XX$:
\begin{lem}\label{lm:ghosttoinertia}
Let $\ev: \GG\to \XX$ be a commutative family of groups over $\XX$. There exists a unique factorization
{\small \eqn\xymatrix@R=14pt{\GG \rto^{\quad e} \ar[rd]_{ \ev} & \IXX \dto^{\ev_0} \\
& \XX.}\eneqn }
 In fact, for any crossed module $[N\stackrel{i}\longrightarrow \gm]$, there is a
unique map $e: \ING\to \IG$ making the following diagram
commutative:
{\small \eqn \xymatrix@R=14pt{\ING \rto^{\quad e} \ar[rd]_{ \ev} & \IG \dto^{\ev_0} \\
& \gm.}\eneqn}
\end{lem}

\begin{example}
Let $\XX$ be an abelian orbifold, that is an orbifold which can be
locally represented by quotients $[X/G]$ where $G$ is (finite)
abelian. Then the $k$-twisted sectors of~\cite{CR} carries a natural
crossed module structure $[S^k_\gm\stackrel{\mu}\longrightarrow \gm]$
where $\mu$ is the $k-1$-fold multiplication $S^k_\gm\to S_{\gm}$
followed by the inclusion $\iota$. Of course, for $k=1$, it is
well-known that the induced stack is the inertia stack and that the
abelian hypothesis can be droped. The associated commutative family of groups is $\Lambda_k \XX \to \XX$ where $\Lambda_k \XX=\IXX\times_{\XX} \cdots \times_{\XX} \IXX$ is the $k^{\rm th}$-inertia stack.
\end{example}

Let $\GG\to \XX$ be a commutative family of groups over a  stack $\XX$. If $\XX$ is strongly oriented, Section~\ref{S:Gysin} yields a canonical Gysin map
$$\gy\Delta\:H\lcom(\GG\times \GG)\to H_{\scriptstyle \bullet-d}(\GG\times_\XX \GG).$$
Thus one can form the composition
 \begin{equation}\label{ghostproduct}
 \star: H_p(\GG)\otimes H_q(\GG)\stackrel{S}\to H_{p+q}(\GG \times
 \GG)\stackrel{\gy{\Delta}}\to
 H_{p+q-d}(\GG\times_{\XX}\GG)\stackrel{m_*}\to  H_{p+q-d}(\GG)
\end{equation}
 Since $m:\GG\times_{\XX} \GG\to \GG$ is associative and commutative as for the inertia stack in Section~\ref{stringproduct}, Step~3,
the argument of Theorem~\ref{th:stringproduct} yields easily
\begin{prop}\label{pro:ghoststring}
Let $\GG$ be a commutative family of groups over an oriented stack  $\XX$ (with $\dim(\XX)=d$). The multiplication $\star$ (see Equation~\eqref{ghostproduct}) endows
the shifted  homology  groups $\sH\lcom(\GG)\cong H_{\scriptstyle \bullet+d}(\GG)$ with a structure of associative, graded
commutative algebra.
\end{prop}
\begin{numrmk}
It is easy to define twisted ring structures on $\sH\lcom(\GG)$
along the lines of Theorem~\ref{twistedstring}. Details are left to the reader.
\end{numrmk}
\begin{numrmk} If $\XX$ is a oriented stack
and if $\GG \to \XX$ is a family of groups which is not supposed to be commutative, the product $\star$ (Equation~\eqref{ghostproduct}) is still defined. Moreover the proof of Theorem~\ref{th:stringproduct} shows that $(\sH\lcom(\GG),\star)$ is an associative algebra. The $k^{\rm th}$-inertia stack $\Lambda_k \XX=\IXX \times_{\XX} \cdots\times_{\XX} \IXX$ is an example of non (necessarily) commutative family of groups.
\end{numrmk}

\begin{numrmk}
Unlike for free loop stacks in Section~\ref{Loopproduct}, we do not need to assume $\XX$ to be Hurewicz in this Section. However, we do not know any interesting example in which it is not the case.
\end{numrmk}

%%%%%%%%%%%%%%%%%%%%%%%%%%%%%%%%%%%%%%%%%%%

\section{Frobenius algebra structures}\label{frobenius}
The loop homology (with coefficients in a field) of a  manifold
carries a rich algebraic structure besides the loop product.
It is known~\cite{CoGo} that there exists also a coproduct, which makes it a   Frobenius algebra (without counit).

It is natural to
expect  that  such a structure also exists
  on $H\lcom(\LXX)$ for an oriented stack $\XX$. In Section~\ref{Frobeniusforloops} we show that this is
  indeed the case. We
  also prove a similar  statement for the homology of inertia stacks.

In this section we assume that our coefficient ring $\kor$ is a field, since
  we will use the  K{\"u}nneth formula $H\lcom(\XX\otimes
  \YY)\stackrel{\sim}\longrightarrow H\lcom(\XX) \otimes H\lcom(\YY)$ (Proposition~\ref{P:Kunneth}).
\subsection{Quick review on Frobenius algebras}\label{reviewFrobenius}
 Let $\kor$ be a field and $A$ a $\kor$-vector space. Recall that $A$ is said
 to be a {\it
 Frobenius algebra} if there is an associative commutative  multiplication
$\mu: A^{\otimes 2}\to A$ and a coassociative cocommutative comultiplication $\delta\: A\to
A^{\otimes 2}$  satisfying  the following compatibility condition
\eq \label{eq:frobenius}
\delta \circ \mu = (\mu\otimes 1)\circ (1\otimes \delta)= (1\otimes \mu)\circ (\delta\circ 1)\eneq
  in ${\rm Hom}(A^{\otimes 2},A^{\otimes 2})$.
 Here we  do not require the existence of
a unit nor a counit. Also we allow  $A$ to be graded and the maps
$\mu$ and $\delta$ to be graded as well. When both maps are of the
same degree $d$, we say that $A$ is a Frobenius algebra of {\em
  degree $d$}.
 The tensor
 product of two Frobenius algebras $A$ and $B$ is naturally  a Frobenius
 algebra with the   multiplication $(\mu \otimes \mu) \circ
 (\tau_{23})$ and comultiplication $\tau_{23}^{-1}\circ (\delta\otimes
 \delta)$ where
 $\tau_{23}\:A\otimes B\otimes A\otimes B \to A^{\otimes 2}\otimes
 B^{\otimes 2}$ is the map permuting the two middle components.

 \begin{warning}
We need a few words of caution concerning  our definition of
Frobenius algebras. In the literature, one often encounters
(commutative) Frobenius
algebras which are both unital and counital such that, if $c: A\to \kor$ is the counit, then $c\circ \mu: A\otimes A\to \kor$ is a nondegenerate  pairing.
\end{warning}

\begin{rmk} It is well-known~\cite{Seg, Abr} that a structure of $1+1$-dimensional Topological Quantum
Field Theory on $A$ is equivalent to a structure of unital and counital Frobenius algebra on $A$ such that the pairing $c\circ \mu: A\otimes A\to
\kor$, where $c$ is the counit and $\mu$ the multiplication, is non-degenerate.  Theorem~\ref{frobeniusloop} and Theorem~\ref{frobeniusstring} below imply that $H\lcom(\LXX)$ and $H\lcom(\IXX)$ have the structure of
$1+1$-positive boundary TQFT in the sense of~\cite{CoGo}. Positive boundary TQFT are obtained by considering only cobordism $\Sigma$ with boundary $\partial \Sigma= -S_1\coprod S_2$ such that $S_1,S_2\neq \emptyset$ (see~\cite{CoGo} for details).

Further, in the case of the free loop stack, we will see in section~\ref{S:HCFT} that the TQFT structure on $H\lcom(\LXX)$ can be extended to a whole homological conformal field theory with positive closed boundary (Theorem~\ref{T:HCFT}). 
\end{rmk}

\subsection{Frobenius algebra structure for loop stacks}\label{Frobeniusforloops}
In this subsection we prove the existence of a Frobenius algebra
structure on the homology of the free loop stack of an oriented (Hurewicz)
stack. Let $\ev_0, \ev_{1/2}\: \LXX \to \XX$ be the evaluation maps
defined as in  Equation~\eqref{eq:ev0}, where  $\XX$ is a Hurewicz
topological stack. To simplify the notations, let $\eb$ be the
evaluation map $(\ev_{0},\ev_{1/2})\:\LXX\to \XX\times \XX$.
\begin{lem} \label{PullbackFrobenius}The stack $\LXX\times_{\XX} \LXX$ fits into a cartesian
 square
\begin{equation}\label{eq:evcartesian}{\xymatrix{\LXX\times_{\XX} \LXX \ar[r]^{m} \ar[d] & \LXX \ar[d]^{\eb} \\
\XX \ar[r]^{\Delta} & \XX\times \XX  }} \end{equation} where $\Delta\:\XX\to
\XX\times \XX$ is the diagonal.
\end{lem}
\begin{pf}
Since $S^1$ is compact and $\XX$ is a Hurewicz topological stack,
Lemma~\ref{L:glue} ensures that the pushout diagram of topological
spaces \eqn
\xymatrix@R=18pt{{\rm pt} \coprod {\rm pt} \rto^{0\coprod \frac{1}{2}} \dto & S^1 \dto \\
{\rm pt} \rto & S^1\vee S^1  }
\eneqn
becomes a pullback diagram after applying the mapping stack  functor
$\map(-,\XX)$. This is precisely diagram~\eqref{eq:evcartesian}.
\end{pf}
\begin{numrmk}\label{rmk:iteratedev}
The argument of Lemma~\ref{PullbackFrobenius} can be applied to iterated diagonals
as well. In particular, $\LXX \times_{\XX} \cdots
\times_{\XX} \LXX$ (with $n$-terms) is the mapping stack
$\map(S^1\vee \dots \vee S^1,\XX)$ (with $n$ copies of $S^1$) and moreover there
is a cartesian square
\eq
{\xymatrix{\LXX\times_{\XX}\cdots\times_{\XX} \LXX \ar[r] \ar[d] & \LXX \ar[d]^{(\ev_0,\ev_{1/n},\dots, \ev_{n-1/n})} \\
\XX \ar[r]^{\Delta} & \XX\times \cdots \times \XX  .}}
\eneq
\end{numrmk}

\medskip

 Now assume further that $\XX$ is oriented of dimension $d$. According to  Section~\ref{S:Gysin},
 the cartesian square~\eqref{eq:evcartesian} yields a Gysin map
$$\gy{\Delta}\:H\lcom(\LXX)\longrightarrow H_{\scriptstyle \bullet -d}(\LXX \times_{\XX}
\LXX).$$ By diagram (\ref{eq:cartesianloop}), there is a canonical map
 $\map (8, \XX)\cong \LXX\times_\XX \LXX \stackrel{j}{\to}
\LXX\times \LXX$. Thus we
obtain a degree $d$ map
{\small \eq  \label{defncoproduct}\delta\: H\lcom (\LXX)\stackrel{\gy{\Delta}}
\rightarrow H_{\scriptstyle \bullet -d}(\LXX\times_{\XX}\LXX)\stackrel{j_* }
\rightarrow H_{\scriptstyle \bullet-d}(\LXX\times \LXX)\cong \hspace{-0.1cm} \bigoplus_{i+j=\bullet-d} \hspace{-0.1cm}
H_{i}(\LXX)\otimes H_{j}(\LXX). \eneq }
\begin{them}\label{frobeniusloop}
Let $\XX$ be an oriented  Hurewicz stack  of dimension $d$. Then
$\left( H\lcom(\LXX),\star,\delta\right)$ is a
 Frobenius algebra, where  both operations $\star$ and $\delta$
are of degree $d$.
\end{them}
\begin{pf} It remains to prove the coassociativity,
  cocommutativity of the coproduct and the Frobenius compatibility relation.
Denote by $\delta_S:H\lcom(\XX\times \YY)\to H\lcom(\XX)\otimes
H\lcom(\YY)$ the inverse of the cross product induced by the
K{\"u}nneth isomorphism and $\delta_S^{(n)}$ for its iteration.

\smallskip

 {\bf i) Coassociativity}  Let $\breve{e}^{(2)}\:\LXX\to
\XX\times \XX\times \XX$ be the iterated evaluation map
$(\ev_0,\ev_{1/3},\ev_{2/3})$. According to Corollary~\ref{C:orienteddiagonal}, the
iterated diagonal $\Delta^{(2)}\:\XX\to \XX\times \XX\times \XX$ is naturally
normally nonsingular oriented. Thus,  Remark~\ref{rmk:iteratedev} implies that  there
is a Gysin map $\gy{\Delta^{(2)}}$. Similarly there
is a canonical map $j^{(2)}\:\LXX\times_\XX \LXX \times_\XX \LXX \cong
\map(S^1\vee S^1\vee S^1,\XX)\to \LXX \times \LXX\times \LXX$.    The argument of
the proof of Theorem~\ref{Associativityloop} shows that it is sufficient
to prove that the following diagram is commutative (which is, in a certain sense, the  dual
of diagram~\eqref{AsD}).
{\footnotesize \begin{equation}\label{coAsD}
\vcenter{\xymatrix@C=7pt@R=10pt@M=6pt{H(\LXX)\otimes H(\LXX) \otimes H(\LXX)  &&&\\
 H(\LXX \times \LXX)\otimes H(\LXX) \uto^{\delta_S \otimes 1}
 &  H(\LXX\times \LXX
\times \LXX) \lto_{\delta_S} \ulto_{\delta_S^{(2)}}
 &H(\LXX \times_{\XX} \LXX \times_{\XX} \LXX)\lto_{j_*^{(2)}}& \\
H(\LXX \times_{\XX} \LXX)\otimes H(\LXX)
\uto^{j_*\otimes 1} &  H((\LXX \times_{\XX}\LXX)
\times \LXX) \ar @{} [dl] |{(3)}  \ar @{} [ul] |{(5)}\lto_{1\times \delta_S}
  \uto_{(j\times 1)_*} & H(\LXX \times_{\XX} \LXX \times_{\XX} \LXX)\ar @{} [ul] |{(1)}
  \lto_{(1\times j)_*} \ar @{=}[u]& \\
  H(\LXX)\otimes  H(\LXX) \uto^{\gy{\Delta}\otimes 1} & H(\LXX \times \LXX)
  \lto^{\delta_S}\uto^{\gy{\Delta} }
   & H(\LXX \times_{\XX} \LXX) \lto^{j_*}
   \uto^{\gy{\Delta}}\ar @{} [ul] |{(2)} & H(\LXX)\lto_{(4)}^{\gy{\Delta}}  \ulto_{\gy{\Delta^{(2)}}}}}
\end{equation}}
where $p$ and $\tilde{p}$ denote, respectively,  the projections
   $\LXX\times \LXX \to \LXX$ and $\LXX\times_\XX \LXX\to \LXX$ on the first factor.
   Square (5) is commutative by naturality of the cross
coproduct $\delta_S$ and the upper left triangle by its coassociativity.  We are left to study
the three remaining squares (1), (2), (3) and  triangle (4).
\begin{description}
\item[Square (1)] commutes in view of the  identity
  $j^{(2)}=(j\times 1) \circ (1\times j)$ which follows from
  the natural isomorphism $(\LXX \times_{\XX}\LXX)
\times \LXX \cong \LXX\times_{\XX}(\LXX\times \LXX) $. Here  the map
  $\LXX\times \LXX\to \XX $ is the composition $\ev_0 \circ p$. In the
  sequel, we use this isomorphism without further notice.

\item[Square (2)] Since  $\breve{e}\circ \tilde{p} =
  \big(\breve{e}\circ p
\big) \circ j:\LXX\times_\XX \LXX \to \XX\times \XX$, the  commutativity of square (2) follows immediately, by naturailty of Gysin maps,
  from  the tower of
  cartesian diagrams
$$\xymatrix@C=8pt@R=10pt@M=8pt{ \LXX \times_{\XX} \LXX \times_{\XX} \LXX \dto_{1\times j} \rto & \LXX \times_{\XX} \LXX \dto^{j} \\
\LXX \times_{\XX} \LXX \times \LXX \dto \rto^{\quad m\times 1} &\LXX \times \LXX
\dto^{\breve{e}\circ p} \\
\XX \rto^{\Delta} &\XX \times \XX.}$$
\item[Square (3)]   commutes by the same argument as for  square (3)
   in diagram~\eqref{AsD}.
\item[Triangle (4)]  The sequence of cartesian diagrams
\eq \label{eq:Delta2}&&
\xymatrix{\LXX\times_{\XX}\LXX\times_{\XX}\LXX \rto^{\quad \quad m\times 1} \dto
  &\LXX\times_{\XX}\LXX \rto^{m} \dto^{\breve{e}\circ \tilde{p}} &\LXX \dto^{(\ev_0,\ev_\half,\ev_{\frac{1}{4}})} \\
\XX \rto^{\Delta} & \XX\times \XX \rto^{\Delta\times 1 \quad  } & \XX\times
\XX \times \XX .}
\eneq  implies, by naturality of Gysin maps, that
\eq \label{eq:coass1}
\gy{\Delta} \circ \gy{ (\Delta\times 1)} &=& \gy{\Delta^{(2)}}.
\eneq
There is an homeomorphism $h:S^1\stackrel{\sim}\to S^1$ which,
 together with the flip map $\sigma$,
induces a commutative
diagram
\eq\label{eq:homeo}
&&\xymatrix{
 \LXX\times_{\XX}\LXX\times_{\XX} \LXX \dto \rto \ar@/^/[rr] & \LXX \rto_{h^*} \dto_{\breve{e}^{(2)}} & \LXX
 \dto^{(\ev_0,\ev_\half,\ev_{\frac{1}{4}})} \\
\XX \rto^{\Delta^{(2)}} & \XX\times \XX \times \XX \ar[r]^{1\times \sigma} &
 \XX\times \XX\times \XX }
\eneq
 As $h^*= \map(-,\XX)(h)$ is
a homeomorphism and $(1\times \sigma) \circ \Delta^{(2)}
 =\Delta^{(2)}$, diagram~\eqref{eq:homeo} identifies  $
\gy{\Delta^{(2)}}$ with the Gysin map (denoted $\gy{\Delta^{(2)}}$ by abuse of notation) associated to Diagram~\eqref{eq:Delta2}. Since  $ \gy{(\Delta\times
 1)}=\gy{\Delta}$ the
 commutativity of Triangle~(4)  follows  from
Equation~\eqref{eq:coass1}.
\end{description}

{\bf ii)} Let's turn  to the point of cocommutativity. It is sufficient to
prove that  \eq\label{eq:cocom}
\gy{\Delta} &=& \widetilde{\sigma}_* \circ
\gy{\Delta},
\eneq
where $\widetilde{\sigma}\: \LXX\times_{\XX}\LXX \to
\LXX\times_{\XX}\LXX$ is the flip map.  There is a natural homotopy $F\:
I\times \LXX \times_{\XX}\LXX \to \LXX$
 between
$m\circ \widetilde{\sigma}$ and
$m$ (see the proof of
Theorem~\ref{Commutativityloop}). Equation~\eqref{eq:cocom} follows easily by naturality of Gysin
 maps applied to the cartesian squares below (where $t\in I$)
\eqn
\xymatrix@C=16pt@R=10pt@M=8pt{ \LXX\times_{\XX} \LXX \ar[r]^{\quad F(t,-)} \ar[d]_{(t,1)} &
  \LXX \ar[d]^{(t,1)} \\I\times \LXX\times_{\XX} \LXX \ar[r]^{\quad
    \quad (1\circ p,F)} \ar[d] & I\times \LXX \ar[d]^{\eb} \\
\XX \ar[r]^{\Delta} & \XX\times \XX  .}
\eneqn
The map $(t,1)\:\LXX\to  I\times \LXX$ is the
map $\LXX\stackrel{\sim}\to \{t\}\times \LXX \to I\times \LXX$.
 The left upper vertical map is similar.
The maps $(t,1)$ are homotopy equivalences inverting
the canonical projections $I\times \LXX\to \LXX$, $I\times
\LXX\times_{\XX} \LXX\to   \LXX\times_{\XX} \LXX$.
\smallskip

{\bf iii)}  It remains to prove the Frobenius relation~\eqref{eq:frobenius}.
 To avoid confusion between different Gysin maps,  we
 now denote $m^!\:=\gy{\Delta}\:H\lcom(\LXX)\to
 H_{\scriptstyle \bullet-d}(\LXX\times_{\XX} \LXX )$ and $j^!:=
 \gy{\Delta}\:H\lcom(\LXX\times \LXX)\to
 H_{\scriptstyle \bullet-d}(\LXX\times_{\XX} \LXX )$ the Gysin maps inducing the
 product and coproduct. The cartesian squares
{\footnotesize \eqn
\xymatrix{\LXX\times_{\XX} \LXX \times_{\XX} \LXX \ar[r]^{\;\;\;
 1\times m}\ar[d] & \LXX\times_{\XX} \LXX
 \ar[d]^{\breve{e}\circ \tilde{p}} \\
\XX \ar[r]^{\Delta} & \XX\times \XX   ,} && \xymatrix{\LXX\times_{\XX} \LXX \times_{\XX} \LXX \ar[r]^{\tilde{j}}\ar[d] & \LXX \times \LXX\times_{\XX} \LXX
 \ar[d]^{(\ev_0\times \ev_0)\circ (1\times \tilde{p})} \\
\XX \ar[r]^{\Delta} & \XX\times \XX   }
\eneqn} give rise to  Gysin maps (see Section~\ref{S:Gysin}) \eqn (1\times m)^!\: H_{\scriptstyle \bullet}(\LXX\times_\XX \LXX) \to
 H_{\scriptstyle \bullet-d}(\LXX\times_{\XX} \LXX \times_{\XX} \LXX) \quad \mbox{ and }\\
 \tilde{j}^!\:  H_{\scriptstyle \bullet}(\LXX
 \times \LXX\times_{\XX} \LXX)\to H_{\scriptstyle \bullet-d}(\LXX\times_{\XX} \LXX \times_{\XX}
 \LXX).\eneqn There is a canonical map
 $\tilde{\tilde{j}}\:\map(S^1\vee S^1\vee S^1,\XX)\cong \LXX\times_\XX
 \LXX\times_\XX \LXX\to \LXX\times_\XX \LXX \times \LXX$ sitting in
 the pullback diagram
\eqn
\xymatrix{\LXX\times_\XX
 \LXX\times_\XX \LXX \ar[d]\ar[r]^{\tilde{\tilde{j}}}& (\LXX\times_\XX \LXX) \times \LXX\ar[d]^{\tilde{\ev_0}
  \times \ev_0}\\
\XX \ar[r]^{\Delta} & \XX\times \XX}
\eneqn
Consider the following diagram
{\footnotesize \begin{equation}\label{frobD}\xymatrix{H\lcom(\LXX \times \LXX) \ar[d]_{j^!} \ar @{} [dr] |{(a)}  \ar[rrd]^{m_{23}^!}  \ar[rr]^{(1\times m)^!}& & H_{\scriptstyle \bullet-d}(\LXX\times \LXX\times_{\XX}\LXX)  \ar @{} [dl] |{(a')}\ar[r]^{(1\times j)_*} \ar[d]^{\tilde{j}^!}& H_{\scriptstyle \bullet-d}(\LXX\times \LXX\times\LXX)\ar[d]^{(j\times 1)^!}\\
 H_{\scriptstyle \bullet-d}(\LXX\times_{\XX} \LXX)  \ar[d]_{m_*}\ar @{} [drr] |{(b)}\ar[rr]_{ 1\times_{\XX} m^! }&& H_{\scriptstyle \bullet-2d}(\LXX\times_{\XX} \LXX \times_{\XX} \LXX) \ar @{} [dr] |{(c)}\ar @{} [ur] |{(b')}\ar[d]^{ m_*\times_{\XX}1} \ar[r]^{\tilde{\tilde{j}}_*}&  H_{\scriptstyle \bullet-2d}(\LXX\times_{\XX} \LXX\times \LXX )   \ar[d]^{m_*\times 1} \\
 H_{\scriptstyle \bullet-d}(\LXX) \ar[rr]_{m^!} && H_{\scriptstyle \bullet-2d}(\LXX\times_{\XX} \LXX) \ar[r]_{j_*} & H_{\scriptstyle \bullet-2d}(\LXX\times \LXX)
}\end{equation}} where $m_{23}^!\:H\lcom(\LXX \times \LXX)\to
 H_{\scriptstyle \bullet-2d}(\LXX\times_{\XX} \LXX \times_{\XX} \LXX)$ is the Gysin
 map determined by the cartesian square (applying Corollary~\ref{C:orienteddiagonal})
$$\xymatrix{\LXX\times_{\XX} \LXX \times_{\XX} \LXX  \ar[d] \ar[rr]^{j\circ (1\times_{\XX} m)}  && \LXX \times \LXX \ar[d]^{\ev_0\times \eb}\\
\XX \ar[rr]^{  \Delta^{(2)}} && \XX\times \XX\times \XX
.}$$
The  triangle  $(a)$ in diagram~\eqref{frobD} is  commutative because we have a
 sequence of
cartesian squares
{\footnotesize \eq \xymatrix{\LXX\times_{\XX} \LXX \times_{\XX} \LXX \ar[r]^{\quad
    1\times_{\XX} m} \ar[d]&  \LXX\times_{\XX} \LXX \ar[r]^{j}\ar[d]^{\eb\circ \tilde{p}}&
  \LXX \times \LXX \ar[d]^{\ev_0\times \eb}\\
\XX \ar[r]^{\Delta}& \XX\times \XX \ar[r]^{ \Delta\times 1} &
\XX\times \XX\times \XX .} \eneq } Similarly,  triangle $(a')$ is
commutative, i.e., $\tilde{j}^! \circ (1\times m)^!=m_{23}^! .$ By
naturality of Gysin maps,  the towers  of cartesian squares
{\footnotesize \eqn \xymatrix{\LXX\times_{\XX} \LXX \times_{\XX} \LXX \ar[r]^{\;\;\; 1\times_{\XX} m}\ar[d]_{1 \times_{\XX} m} & \LXX\times_{\XX} \LXX \ar[d]^{m} \\
\LXX\times_{\XX} \LXX \ar[d]\ar[r]^{m} & \LXX \ar[d] \\
\XX \ar[r]^{\Delta} & \XX\times \XX \;  ,}&& \xymatrix{\LXX\times_{\XX}
  \LXX \times_{\XX} \LXX \ar[r]\ar[d]_{\tilde{\tilde{j}}} & \LXX\times
  \LXX\times_{\XX} \LXX \ar[d]^{1\times j} \\
\LXX\times_{\XX} \LXX \times \LXX \ar[d]_{\tilde{\ev_0}\circ (1\times p)}\ar[r]^{j\times 1} &
\LXX\times \LXX\times \LXX
\ar[d]^{(\ev_0\times \ev_0)\circ (1\times p)} \\
\XX \ar[r]^{\Delta} & \XX\times \XX   }\eneqn} give the commutativity of squares
 (b) and (b') in diagram~\eqref{frobD}. The commutativity of Square (c) is trivial.
 Thus diagram~\eqref{frobD} is commutative. Up to the identification
 $H\lcom(\LXX \times \LXX)\cong H\lcom(\LXX) \otimes H\lcom(\LXX)$,
 the composition of the bottom horizontal map and the left vertical
 one in diagram~\eqref{frobD} is the composition $\delta ( -\star
 -)$. The composition of the right vertical map with the upper arrow
 is $  (a\star b^{(1)})\otimes b^{(2)}$. Finally commutation of the Gysin maps with  the cross product yields the identity
$$ \delta(a\star b)=a\star b^{(1)}\otimes b^{(2)}.$$ The proof of  identity
 $ \delta(a\star b)=a^{(1)}\otimes a^{(2)}\star b$  is similar.
\end{pf}

\subsection{Frobenius algebra structure for inertia stacks}\label{frobeniusstring}
 In this section we show that the homology of the inertia stack  is also a Frobenius algebra, similarly to Theorem~\ref{frobeniusloop}.

\smallskip

 Let $\XX$ be a 
topological stack of dimension $d$ and $\gm$  a topological groupoid
representing $\XX$. Thus its inertia stack $\IXX$  is the stack
associated to the inertia groupoid $\IG\:=S\gm\rtimes \gm \toto
S\gm$, where  $S\gm$ is the space of closed loops.  Any loop $S^1\to
X$ on a topological space $X$ can be evaluated in $0$ but also in
$1/2$. It is folklore to think of $\IXX$ as a  ghost loop stack.
Hence evaluation map at $0$ and $1/2$ should make sense as well. We
first
 construct these evaluation maps for the inertia
stack which leads to the construction of the Frobenius structure on
$H\lcom(\IXX)$ when $\XX$ is oriented.

\smallskip

First of all, let us   introduce another groupoid $\IWG$ which is
Morita equivalent to $\IG$.
Objects of $\IWG$ consist of all diagrams
\begin{equation}\label{arrows}\xymatrix{ x & y \ar@/_/[l]_{g_1}   & x\ar@/_/[l]_{g_2}  } \end{equation} in $\gm$. Note  that the composition $g_1g_2$ is a loop over $x$.  Arrows of $\IWG$ consist  of commutative diagrams
{\footnotesize $$\xymatrix@R=14pt{ x & y \ar@/_/[l]_{g_1}   & x\ar@/_/[l]_{g_2} \\
x'\ar[u]^{h_0}  & y' \ar[u]^{h_{1/2}} \ar@/^/[l] & x'
\ar[u]_{h_{0}}  \ar@/^/[l] .} $$} Note that the left and right
vertical arrows are the same.
  The target map is the top row
$$\xymatrix{ x & y \ar@/_/[l]_{g_1}   & x\ar@/_/[l]_{g_2}  } $$

while the source map is the bottom row $$\xymatrix{ x' & y'
\ar@/_/[l]_{h_0^{-1}g_1h_{1/2}}   &
x'\ar@/_/[l]_{h_{1/2}^{-1}g_2h_0}  }.$$ The unit map is  obtained
by taking identities as vertical arrows. The composition is
obtained by superposing two diagrams and deleting the middle row of
the diagram, i.e.
{\footnotesize $$\xymatrix@R=14pt{ x & y \ar@/_/[l]_{g_1}   & x\ar@/_/[l]_{g_2} \\
x' \ar[u]^{h_0}  & y' \ar[u]^{h_{1/2}} \ar@/^/[l] & x'
\ar[u]_{h_{0}}  \ar@/^/[l] }
*
\xymatrix@R=14pt{ x'& y' \ar@/_/[l]   & x'\ar@/_/[l] \\
x" \ar[u]^{h'_0}  & y"  \ar[u]^{h'_{1/2}} \ar@/^/[l] & x"
\ar[u]_{h'_0}  \ar@/^/[l] }$$} is mapped to {\footnotesize
$$\xymatrix@R=14pt{
 x & y \ar@/_/[l]_{g_1}   & x\ar@/_/[l]_{g_2} \\
x" \ar[u]^{h_0h'_0}  & y"  \ar[u]^{h_{1/2}h'_{1/2}} \ar@/^/[l] & x" \, .
\ar[u]_{h_0h'_{0}}  \ar@/^/[l] }$$} In other words, $\IWG$ is the
transformation groupoid $\widetilde{S\gm}\rtimes_{\Gamma_0\times
\Gamma_0}\big(\Gamma_1\times \Gamma_1\big)$, where
$\widetilde{S\gm}=\{(g_1,g_2)\in \gm_2 \, |\,
t(g_1)=s(g_2)\}$, the momentum map $\widetilde{S\gm}\to
\Gamma_0\times \Gamma_0$ is $(t,t)$,  and the action is given, for
all compatible $(h_0,h_{1/2})\in \Gamma_1\times \Gamma_1$,
$(g_1,g_2)\in \widetilde{S\gm}$, by
  $$(g_1,g_2)\cdot (h_0, h_{1/2})= (h_0^{-1}g_1h_{1/2}, h_{1/2}^{-1}g_2h_0).$$

One defines evaluation maps taking  by the vertical arrows of
$\widetilde{S\gm}$, i.e.  $\forall (g_1,g_2,h_0,h_{1/2})\in \IWG_1$
define
$$ \ev_0:(g_1,g_2,h_0,h_{1/2})\mapsto h_0, \quad
\ev_{1/2}\:(g_1,g_2,h_0,h_{1/2})\mapsto h_{1/2}.$$ It is simple to prove
\begin{lem}\label{evaluations} Both evaluation maps
$\ev_0:\IWG\to \gm$ and $\ev_{1/2}\:\IWG\to \gm$ are groupoid
morphisms.
\end{lem}

There is a  map \eq \label{eq:p} p\:\IWG \to \IG \eneq obtained by   sending a
diagram in $\IWG_1$ to the composition of the horizontal arrows,
i.e.,
{\footnotesize \eqn\xymatrix@R=14pt{ x & y \ar@/_/[l]_{g_1}   & x\ar@/_/[l]_{g_2} \\
x' \ar[u]^{h_0}  & y' \ar[u]^{h_{1/2}} \ar@/^/[l] & x' \ar[u]_{h_{0}}
\ar@/^/[l] } &\mbox{is mapped to} & \xymatrix@R=14pt{ x & x\ar@/_/[l]_{g_1g_2} \\
x'\ar[u]^{h_0}   & x' \ar[u]_{h_{0}}  \ar@/^/[l] .} \eneqn} In other words
$ p(g_1,g_2,h_0,h_{1/2}) = (g_1g_2, h_0)$.

\begin{lem}\label{Moritaequiv}
The map $p\:\IWG\to \IG$ is a Morita morphism.
\end{lem}
\begin{pf}
The map $p_0:\IWG_0\to \IG_0$ is a surjective submersion with a  section given by $g\mapsto (g,1_{s(g)})$ for $g\in S\gm$. Let $g,g'\in S\gm$. Assume given $(g_1,g_2) \in \gm_2$ with $g_1g_2=g$ and  $(g_1',g_2') \in \gm_2$ with $g_1'g_2'=g'$. Then any arrow in $\IWG$ from $\xymatrix{ x & y \ar@/_/[l]_{g_1}   & x\ar@/_/[l]_{g_2}}$ to $\xymatrix{ x & y \ar@/_/[l]_{g_1'}   & x\ar@/_/[l]_{g_2'}}$ is uniquely determined by  $h_0\in \gm_1$ satisfying $h_0^{-1}g_1g_2h_0=g_1'g_2'$. Indeed, $h_{1/2}$ is given by $h_{1/2}=g_2h_0{g_2'}^{-1}$.
\end{pf}
As a consequence the groupoid $\IWG$ also represents the inertia stack $\IXX$, and
Lemma~\ref{evaluations} implies that there are two stack maps $\ev_0,
\ev_{1/2}\:\IXX\to \XX$.

\medskip

We now proceed to construct the \hidden coproduct.
As in Section~\ref{stringproduct} above,   $\IG \times_{\gm}\IG$ is  the transformation groupoid
 $\big(S{\Gamma}\times_{\Gamma_0} S{\Gamma}\big) \rtimes \Gamma\toto
S{\Gamma}\times_{\Gamma_0} S{\Gamma}$, where $\Gamma$ acts on
$S{\Gamma}\times_{\Gamma_0} S{\Gamma}$ by conjugations diagonally.
 Its corresponding stack is $\IXX\times_{\XX} \IXX$.

\begin{lem} \label{PullbackIWG}The stack $\IXX\times_{\XX} \IXX$ fits into the cartesian square
{\footnotesize $${\xymatrix@R=18pt{\IXX\times_{\XX} \IXX \ar[r]^{m} \ar[d] & \IXX \ar[d]^{(\ev_0, \ev_{1/2})} \\
\XX \ar[r]^{\Delta} & \XX\times \XX  }} .$$}
\end{lem}
As in Section~\ref{frobeniusloop}, we denote $\eb:=
(\ev_0,\ev_{1/2})\:\IXX\to \XX\times \XX$ the right vertical map in the
diagram of lemma~\ref{PullbackIWG}.

\begin{pf}We use $\IWG$ as a groupoid representative of $\IXX$.
By the definition of the evaluation maps, the fiber product \eqn
\xymatrix@R=18pt{\gm\times_{\gm\times \gm} \IWG \ar[r] \ar[d] & \IWG \ar[d]^{(ev_0, ev_{1/2})} \\
\gm \ar[r]^{\Delta} & \gm\times \gm }
\eneqn can be identified with the subgroupoid of $\IWG$,  which
consists of $(g_1,g_2,h_0,h_{1/2})$ such that $h_0=h_{1/2}$. The
latter is simply  the transformation groupoid
$${S\gm\times_{\gm_0}S\gm \rtimes \gm_1}\toto
S\gm\times_{\gm_0}S\gm$$ which is precisely $\IG\times_{\gm}
\IG$. Moreover the composition  $$\IG\times_{\gm}
\IG\cong \gm\times_{\gm\times \gm} \IWG
\to  \IWG \stackrel{p}\to \IG, $$ where $p$ is defined by equation~\eqref{eq:p}, is precisely the ``Pontrjagin'' map
$m:\IG\times_{\gm}\IG \to \IG $ in
Section~\ref{stringproduct}.
\end{pf}
\begin{numrmk}
It is not hard to generalize the above  construction to any finite
number  of evaluation maps and obtain the following cartesian square
(see the proof of Theorem~\ref{frobeniusstring} below)
{\footnotesize $$\xymatrix@R=18pt{\IXX\times_{\XX} \cdots
\times_{\XX} \IXX \rto \dto & \IXX
  \dto \\
\XX \rto & \XX\times \cdots \times \XX. } $$}
\end{numrmk}

\medskip

If $\XX$ is oriented of dimension $d$,  the cartesian square of
Lemma~\ref{PullbackIWG} yields a Gysin map (Section~\ref{S:Gysin})
$$\gy{\Delta}\:H\lcom(\IXX)\longrightarrow H_{\scriptstyle \bullet -d}(\IXX
\times_{\XX}\IXX).$$ As shown in Section~\ref{stringproduct}, there
is also a canonical map   $j:\IXX\times_{\XX} \IXX\to \IXX\times \IXX$.
\begin{them}\label{th:frobeniusstring}
Assume $\XX$ is an oriented stack of dimension $d$. The composition $$
H_n(\IXX)\stackrel{\gy{\Delta}}\longrightarrow H_{n
-d}(\IXX \times_{\XX}\IXX)\stackrel{j}\longrightarrow
H_{n-d}(\IXX\times \IXX)\cong \bigoplus_{i+j=n-d}H_{i}(\IXX)\otimes
H_{j}(\IXX)$$ yields a coproduct $\delta\:H\lcom(\IXX)\to
\bigoplus_{i+j=\bullet -d} H_{i}(\IXX)\otimes
H_{j}(\IXX)$ which is  a coassociative and graded cocommutative
coproduct on
the shifted homology $\sH\lcom(\IXX):=H_{\scriptstyle \bullet+d}(\IXX)$,
called the \hidden coproduct of $\IXX$.
\end{them}
\begin{pf}
The proof is very similar to that  of
Theorem~\ref{frobeniusloop}. We only explain the difference.

{\bf i)} First we introduce a third evaluation map $\ev_{2/3}\:\IXX
\to \XX$ similar to $\ev_{1/2}$. Taking a representative $\gm$
of $\XX$,  the idea is to replace $\IG$ by another groupoid
$\widetilde{\IWG}$, where
$\widetilde{\IWG}_1$consists of
commutative diagrams:
{\footnotesize $$\xymatrix@R=14pt{ x & y \ar@/_/[l]_{g_1}   & z\ar@/_/[l]_{g_2} &x\ar@/_/[l]_{g_3} \\
x'\ar[u]^{h_0}  & y' \ar[u]^{h_{1/2}} \ar@/^/[l] & z'
\ar[u]^{h_{2/3}}  \ar@/^/[l] & x'\ar[u]_{h_0}  \ar@/^/[l] .} $$}
The source and target maps are, respectively, given by the bottom and
upper lines. The multiplication is by superposition of diagrams. There are evaluation maps
$\ev_0,\ev_{1/2},\ev_{2/3}\:\widetilde{\IWG}\to\gm$, respectively,
given by $h_0,h_{1/2},g_{2/3}$. A proof similar to those
of Lemmas~\ref{Moritaequiv} and Lemmas \ref{PullbackIWG} gives the following facts :
\begin{enumerate}
\item the groupoid $\widetilde{\IWG}$ is Morita equivalent to
  $\IG$. Hence it also  represents the stack $\IXX$.
\item The evaluation maps induce a cartesian square
{\footnotesize $$\xymatrix{\IXX\times_{\XX} \IXX\times_{\XX}\IXX \rto \dto & \IXX
  \dto^{\breve{e}^{(2)}} \\
\XX \rto^{\Delta^{(2)}} & \XX \times \XX \times \XX} $$} which yields a Gysin
map $\gy{\Delta^{(2)}}\:H\lcom(\IXX) \to H_{\scriptstyle \bullet-2d}(\IXX\times_{\XX} \IXX\times_{\XX}\IXX)$.
\end{enumerate}
It follows that one can form a diagram similar to~\eqref{coAsD} for $\IXX$ and
prove that all its squares (1),(2), (3), (5) are commutative {\em
  mutatis mutandis}.
The proof of the commutativity of triangle (4) is even easier: it
follows immediately
from the sequence of cartesian square
\eqn &&
\xymatrix{\IXX\times_{\XX} \IXX\times_{\XX}\IXX \ar[r] \ar[d] & \IXX\times_{\XX} \IXX
  \ar[r]^{m}\ar[d]^{\eb\circ \tilde{p}} & \IXX \ar[d]^{\breve{e}^{(2)}}\\ \XX  \ar[r]^{\Delta}& \XX
  \times \XX \ar[r]^{\Delta\times 1} & \XX\times \XX \times \XX   .}
\eneqn

{\bf ii)} Since $p\circ \widetilde{\sigma}$ is conjugate to $p$, the proof of the cocommutativity of $\delta$ is similar to the proof
of Proposition~\ref{frobeniusloop} and of Proposition~\ref{th:stringproduct}.
\end{pf}
\begin{them}
The homology groups $(H_{\scriptstyle \bullet}(\IXX),\bullet,\delta)$ form a (non unital, non counital) Frobenius algebra of degree $d$.
\end{them}
\begin{pf}
According to Theorems~\ref{th:stringproduct}, \ref{frobeniusstring} it
suffices to
prove the compatibility condition between the \hidden product and
\hidden coproduct.  The argument of the proof  of Theorem~\ref{frobeniusloop}.iii) applies.
\end{pf}
\begin{numrmk}\label{universalcoefficient}
If $\XX$ has finitely generated homology groups in each degree, then
by universal coefficient theorem, $H\com(\IXX)$ inherits a
Frobenius coalgebra structure which is unital iff
$(H_{\scriptstyle \bullet}(\IXX),\delta)$ is counital.
\end{numrmk}

%%%%%%%%%%%%%%%%%%%%%%%%%%%%%%%%%%%%%%%%%%%%%%%%%%%%%%%%%%%%%%%%%%%%%

\subsection{The canonical morphism $\IXX\to \LXX$}\label{S:Frobmap}
There is a   morphism of stacks $\Phi\:\Lambda \XX \to \LXX$ generalizing
the canonical inclusion of a space into its loop space (as a constant
loop).
\begin{numrmk}\label{rm:Phi}
Objects of $\IXX$ are pairs $(X,\varphi)$ where $X$ is an object of $\XX$ and $\varphi$ an automorphism of $X$. The morphism $\Phi$ may be thought to maps $(X,\varphi)\in \IXX$ to the isotrivial family over $S^1$, which is obtained from the constant family $X_I$ over the interval by identifying the two endpoints via $\varphi$.
\end{numrmk}
We  show in this Section that $\Phi$ induces a morphism of
Frobenius algebras in  homology.

\medskip

Let $\gm$
be a groupoid representing the oriented stack $\XX$ (of dimension $d$) and $\Lambda \gm$ its inertia
groupoid representing $\IXX$. Proposition~\ref{freeloopgpoid} gives a
groupoid $\LG$  representing  the free loop stack $\LXX$. We use the
notations of Section~\ref{S:freeloopgpoid}. Recall that the topological groupoid $\LG$ is a limit of topological groupoids $L^P\gm$ where $P$ is a finite subset of $S^1$. We take $P=\{1,1\}\subset S^1$ the trivial subset of $S^1$. We will construct a morphism of groupoids $\IG\to L^P\gm$ inducing the map $\IXX\to \LXX$.

\smallskip

 Any $(g,h)\in S\gm \rtimes \gm= \IG_1$
(i.e. $g\in \gm_1$ with $s(g)=t(g)$) determines a commutative diagram $\Phi(g,h)$ in the  underlying category of the groupoid
$\gm$~:
 \begin{equation}\label{D(g,h)} \xymatrix@=16pt@M=8pt@C=22pt{
      t(h) & t(h) \lto_{g} \\
      s(h) \uto^{h} & s(h) \uto_{h} \lto_{h^{-1}gh}\;. }\end{equation}
      The square $\Phi(g,h)$ (defined by diagram~\eqref{D(g,h)}) being
      commutative, it is an element of $M_1\gm$. Since $P$ is a
      trivial subset of $S^1$, a morphism $[S_1^P\toto S_0^P]\to [M_1\gm\toto M_0\gm]$ is given by a path $f\:[0,1] \to \gm_1$ and  elements $k,k'\in \gm_1$ such that the diagram
\eqn  \xymatrix@=16pt@M=8pt@C=22pt{
      t(f(0)) & t(f(1)) \lto_{k} \\
      s(f(0)) \uto^{f(0)} & s(f(1))) \uto_{f(1)} \lto_{k'}}
\eneqn
commutes.
 In particular, the diagram $\Phi(g,h) \in M_1\gm$ yields  a (constant) groupoid
morphism $[S_1^P\toto S_0^P]\to [M_1\gm\toto M_0\gm]$ defined by $t\mapsto f(t)=h$. The
map $(g,h)\mapsto \Phi(g,h)$ is easily seen to be a groupoid morphism.
 We denote by  $\Phi\: \Lambda \gm \to
       \Lo\gm$ its composition with the inclusion $\Lo^P\gm\to
       \Lo\gm$. It is still a morphism of groupoids. Hence we have the
      following
\begin{lem}\label{lem:frobeniusmorphism}
The map $\Phi\:\Lambda \gm \to \Lo\gm$ induces a functorial map of stacks $\Lambda \XX \to \LXX$.
\end{lem}
In particular there is an induced map $\Phi_*\:H\lcom(\Lambda \XX) \to
H\lcom(\LXX)$ in  homology.

\begin{them}\label{Frobeniusmorphism}
Let $\XX$ be an oriented Hurewicz stack. The map $\Phi_*\:(H\lcom(\Lambda \XX),\bullet,\delta) \to (H\lcom(\LXX),\star,\delta)$ is a
morphism  of Frobenius algebras.
\end{them}
\begin{pf}
Let $\gm$ be a groupoid representing $\XX$.   For any
$(g,h)\in S\gm\rtimes \gm_1=\IG_1$, one has
\eqn
\ev_0\big(\Phi(g,h)\big) &=& h =\ev_0(g,h)
\eneqn
where $\ev_0$ stands for both evaluation maps \comment{I
  use the gpoid representation LG at this level implicitly}  $\LG\to \gm$,  $\IG\to \gm$.  Thus the
cartesian square of Step (2) in the construction of the \hidden product
factors through the one of the loop product and we have a tower of
cartesian squares:
\eq \label{eq:frobeniusmorphism} \xymatrix@R=17pt{\IXX\times_{\XX}\IXX \rto \dto_{\widetilde{\Phi}}  &
  \IXX\times \IXX\dto^{\Phi\times \Phi}  \\
\LXX\times_{\XX} \LXX \dto \rto & \LXX \times \LXX\dto^{\ev_0\times \ev_0}\\
\XX \rto^{\Delta}& \XX \times \XX}\eneq where $\tilde{\Phi}$ is
induced by $\Phi\times \Phi$. The square~\eqref{eq:frobeniusmorphism} shows that\eq \label{eq:frobeniusGysin}
\gy{\Delta}\circ \big(\Phi\times \Phi\big)_* &=& \widetilde{\Phi}_*\circ \gy{\Delta}.
\eneq

\smallskip

Since $\LG$ is a presentation of $\LXX$, the cartesian square
$\LG\times_{\gm}\LG$ represents
the stack $\LXX\times_{\XX} \LXX$. Given any $(g_1,g_2,h)$  in $\big(S\gm\times_{\gm_0}
S\gm \big)\rtimes\gm_1=\IG\times_{\gm}\IG$, one can form a
commutative diagram  $\tilde{\Phi}(g_1,g_2,h)$:
\begin{equation*} \label{eq:D2}\xymatrix@=14pt@M=8pt@C=20pt{
      t(h) & t(h) \ar@/_/[l]_{g_1} & t(h) \ar@/_/[l]_{g_2} \\
      s(h) \uto^{h} & s(h) \uto_{h} \ar@/^/[l]^{h^{-1}g_1h} &  s(h) \uto_{h} \ar@/^/[l]^{h^{-1}g_2h}, }\end{equation*}
      which induces canonically an arrow of $\LG\times_{\gm}\LG$ as in the construction of $\Phi$. The
      map $(g_1,g_2,h)\mapsto  \tilde{\Phi}(g_1,g_2,h)$ represents the
      stack morphism $\widetilde{\Phi}$. Since 
\eqn m\big(\widetilde{\Phi}(g_1,g_2,h)\big) &=& \Phi(g_1g_2,h)\eneqn
 the diagram
\eq \label{eq:frobeniusmorphism2}
\xymatrix{\IXX\times_{\XX}\IXX \rto^{m} \dto_{\widetilde{\Phi}} & \LXX\times_{\XX}\LXX \dto^{m} \\
\IXX \rto^{\Phi} & \LXX }
\eneq is commutative. Hence, diagram~\eqref{eq:frobeniusmorphism2} and
      Equation~\eqref{eq:frobeniusGysin} implies  that $\Phi_*$ is an algebra morphism.
       Similarly $\Phi_*$ is a coalgebra morphism since the diagram
\eqn
&&\xymatrix@R=16pt{\IXX\times_{\XX}\IXX \rto^{m} \dto_{\widetilde{\Phi}} & \IXX \dto^{\Phi} \\
\LXX\times_{\XX} \LXX \rto^{m} \dto & \LXX \dto^{(\ev_0,\ev_{1/2})}\\
X\rto^{\Delta} & \XX\times \XX.  }
\eneqn
is commutative.
\end{pf}
\begin{numrmk}
If the stack  $\XX$ is actually a manifold $X$, then its inertia
stack is $X$ itself and $\LXX=LX$ the free loop space of $X$. It is
clear that the map $\Phi$ becomes the usual inclusion
$X\hookrightarrow LX$ identifying $X$ with constant loops. For
manifolds, the map $\Phi_*$ is injective but not surjective (except
in trivial cases). However, for general stacks, $\Phi_*$ is not
necessary  injective nor surjective. See Section~\ref{Liegroups}.
\end{numrmk}

%%%%%%%%%%%%%%%%%%%%%%%%%%%%%%%%%%%%%%
%%%%%%%%%%%%%%%%%%%%%%%%%%%%%%%%%%%%%%%%%

\section{The \BV-algebra  on the homology of free loop stack}\label{BVstructure}

\subsection{\BV-structure}
In this section we construct a \BV-algebra structure on the homology of the
loop stack. First we recall the definition of a \BV-algebra.

A {\bf Batalin-Vilkovisky algebra} ({\bf \BV-algebra} for short) is a graded commutative unital associative   
algebra with a degree $1$ operator $D$  such that $D(1)=0$, $D^2=0$, and the  following identity
is satisfied:
     \begin{multline} \label{eq:BVidentity} D(abc)-D(ab)c-(-1)^{|a|}aD(bc)-(-1)^{(|a|+1)|b|} bD(ac)+
      \\ +D(a)bc+(-1)^{|a|}aD(b)c+
      (-1)^{|a|+|b|}abD(c)=0.\end{multline}
In other words, $D$ is a second-order differential operator. 

Now, let $\XX$ be a topological stack and $\LXX$ its loop stack. 
The circle  $S^1$ acts on itself by left multiplication. By
functoriality  of the mapping stack, this $S^1$-action  
confers an $S^1$-action to $\LXX$ for any topological stack $\XX$.
 This action
endows $H\lcom(\LXX)$ with a degree one operator $D$ as follows. Let
$[S^1]\in H_1(S^1)$ be the fundamental class. Then a linear map $D:
H\lcom(\LXX)\to H_{\scriptstyle \bullet+1}(\LXX)$  is defined by the
composition
$$H\lcom(\LXX)\stackrel{\times [S^1]}\longrightarrow  H_{\scriptstyle \bullet+1}(\LXX\times S^1)
\stackrel{\rho_*}\longrightarrow H_{\scriptstyle \bullet+1}(\LXX),$$ where  the last arrow is induced by the action $\rho:S^1\times \LXX \to \LXX$.

\begin{lem} \label{L:Dsquare}
The operator $D$ satisfies $D^2=0$, {\it i.e.} is a differential.
\end{lem}
\begin{pf}
Write $m:S^1\times S^1\to S^1$ for the group multiplication on
$S^1$. The  naturality of the cross product  implies, for any $x\in
\sH\lcom(\LXX)$, that
$$D^2(x)=\rho_* \big(m_*([S^1]\times [S^1])\times x) \big)=0$$
since $m_*([S^1]\times [S^1])\in H_2(S^1)=0$.
\end{pf}

\begin{them}\label{BV} Let $\XX$ be an oriented (Hurewicz\footnote{Recall that any differentiable stack is Hurewicz}) stack of dimension $d$.  Then the shifted homology
$\sH\lcom(\LXX)=H_{\scriptstyle \bullet +d}(\LXX)$ admits a \BV-algebra structure
given by the  loop product $\star :\sH\lcom(\LXX)\otimes \sH\lcom(\LXX)\to
\sH\lcom(\LXX)$ and the operator  $D:\sH\lcom(\LXX)\to
\sH_{\scriptstyle \bullet+1}(\LXX)$.
\end{them}
\begin{rmk}
We will get a proof of Theorem~\ref{BV} using conformal field theory in Section~\ref{S:HCFT}. However, here we wish to give a direct proof, hoping this proof can also be applied to some family of commutative groups as introduced in Section~\ref{S:hidden}.
\end{rmk}

\subsection{Gerstenhaber bracket and proof of Theorem~\ref{BV}}
We start by some well-known facts on \BV-algebras. Let $(A, \cdot, D)$ be a \BV-algebra. We can define a degree 1 binary operator $\{\,;\,\}$ by the following formula:
\begin{eqnarray}\label{eq:Gerstenhaberbracket}
 \{a;b\} &=& (-1)^{|a|}D(a\cdot b) -(-1)^{|a|}D(a)\cdot b -a \cdot D(b) 
\end{eqnarray}
The \BV-identity~\eqref{eq:BVidentity} and commutativity of the product implies that $\{\,;\,\}$ is a derivation of each variable (and anti-symmetric with respect to the degree shifted down by 1). Further the relation $D^2=0$ implies the (graded) Jacobi identity for $\{\,;\,\}$. In other words, $(A,\cdot, \{\,;\,\})$ is a \emph{Gerstenhaber} algebra, that is a commutative gradeed algebra equipped with a bracket $\{\,;\,\}$ that makes $A[1]$ a graded Lie algebra and satisfying a graded Leibniz rule~\cite{CoJo}. 

Indeed it is standard (see~\cite{Get}) that a graded commutative algebra $(A,\cdot)$ equipped with a degree 1 operator $D$, such that $D^2=0$, is a \BV-algebra if and only if the operator $\{\,;\,\}$ defined by the formula~\eqref{eq:Gerstenhaberbracket} is a derivation of the second variable, that is 
\begin{eqnarray}\label{eq:BVidentity2}
 \{a;bc\}= \{a;b\}\cdot c + (-1)^{|b|(|a|+1)}b\cdot\{a;c\}.
\end{eqnarray}

By Theorem~\ref{th:Loop} and Lemma~\ref{L:Dsquare}, the shifted homology $(H\lcom(\LXX),\star ,D)$, equipped with the loop product and operator $D$ induced by the circle action on $\LXX$, is a  graded commutative algebra and $D^2=0$. In order to prove Theorem~\ref{BV}, we will thus prove the identity~\eqref{eq:BVidentity2}. First, we identify the bracket $\{\,;\,\}$ given by formula~\eqref{eq:Gerstenhaberbracket}. We need to introduce some notations to do so.

Let $\ev_\rho: S^1\times \LXX \times \LXX \to \XX\times \XX$ be the (twisted by $\rho$) evaluation map defined by $\ev_\rho(t,\gamma,\beta)= \ev_0(\gamma) \times \ev_{2t}(\beta)$ for $0\leq t\leq 1/2$ and $\ev_\rho(t,\gamma,\beta) = \ev_{2t-1}(\gamma)\times \ev_{1}(\beta)$ for $1/2\leq t\leq 1$ (where $\ev_x: \LXX\to \XX$ is the evaluation map defined in Section~\ref{mappingstack}). We let $L_\rho(\XX\times \XX)$ be the pullback stack of the diagonal along $\ev_\rho$:
\begin{eqnarray}\label{eq:Lrho}
 \xymatrix{L_\rho(\XX\times \XX) \ar[r]^{i_\rho} \ar[d] & S^1\times \LXX\times \LXX \ar[d]^{\ev_\rho}\\
 \XX \ar[r]^{\Delta} & \XX\times \XX }
\end{eqnarray}
Note that $L_\rho(\XX\times \XX)\cong S^1\times \Map(S^1\vee S^1,\XX)$ (by Lemma~\ref{L:glue}) and that, under this identification, $i_\rho$ becomes the map $(t,\gamma)\mapsto \big(t,  \gamma^{(1)}, \rho(2t)(\gamma^{(2)})\big)$ for $0\leq t\leq 1/2$ and $(t,\gamma)\mapsto \big(t,\rho(2t-1)(\gamma^{(1)}),\gamma^{(2)} \big)$ for $1/2\leq t\leq 1$ (here $\gamma\mapsto (\gamma^{(1)},\gamma^{(2)})$ is the map $\Map(S^1\vee S^1,\XX)\cong \LXX\times_{\XX}\LXX \to \LXX\times \LXX$).

For any $0\leq t\leq 1/2$, we have a pinching map $p_t:S^1\to S^1$ defined by $p_t(u)=(0,2u-2t)\in S^1\vee S^1\subset S^1\times S^1$ for $0\leq u\leq t$, $p_t(u)=(2u-2t,0)$  for $t\leq u\leq t+1/2$ and $p_t(u)=(0,2u-1-2t)$ for $t+1/2\leq u$. For $ 1/2\leq t\leq 1$, we define similarly $p_t(u)=(2u-2t+1,0)$ for $u\leq t-1/2$, $p_t(u)=(0, 2u-2t+1)$ for $t-1/2\leq u\leq t$ and $p_t(u)=(2u-2t,0)$ for $u\geq t$. Note that $p_0=p_1$ is the pinching map of Section~\ref{Loopproduct}. We let $m_\rho: L_\rho(\XX\times \XX)\to \LXX$ be the map $(t,\gamma)\mapsto \gamma \circ p_t$ induced by the above pinching maps. Informally, the map $m_{\rho}$ can be described as follows:  an element in $L_\rho(\XX\times \XX)$ is given, for each $1/2\geq t\in S^1$, by two loops $a,b\in \LXX$ such that $\ev_0(a)=\ev_{2t}(b)$. Then $m_\rho(t,a,b)$ is the loop starting at $\ev_0(b)$, describing $b$ until it reaches $\ev_{2t}(b)=\ev_0(a)$ where it follows the loop $a$ and then follows $b$ back to $\ev_0(b)$. There is a similar picture for $t\geq 1/2$. Note that, in the proof of Theorem~\ref{BV}, we will use several times an informal description similar to the one of $m_\rho$ to descrive various maps that are rigorously defined using a parametrized pinching procedure as above.    

\smallskip

When $\XX$ is an oriented stack of dimension $d$, then the pullback diagram~\eqref{eq:Lrho} induces a Gysin map $\Delta_\rho^!:H_\bullet\big(S^1\times \LXX\times \LXX\big) \to  H_{\bullet+d}\big(L_\rho(\XX\times \XX)\big)$.  

\begin{lem}\label{L:Gerstenhaberbracket} If $\XX$ is an oriented stack and $a,b\in \BH\lcom(\LXX)$, one has
 $$\{a;b\} = {m_{\rho}}_*\circ \Delta_\rho^!\Big([S^1]\times a\times b\Big). $$ 
\end{lem}
\begin{proof} The proof is the stack analogue of a result of Tamanoi~\cite[Theorem 5.4]{Ta} and~\cite[Definition 3.1]{Ta}. To apply the proof of~\cite{Ta}, we only need to use the evaluation maps as we did to define $L_\rho(\XX\times \XX)$ and Gysin maps induced by the pullback along the diagonal $\Delta : \XX\to \XX\times \XX$ which is strongly oriented by assumption. Then all identities involving Gysin maps in~\cite{Ta} follow using  the Gysin map given by the bivariant theory (see Section~\ref{S:Gysin}, and the technics of the proofs of Theorem~\ref{th:Loop} and Theorem~\ref{frobeniusloop}) so that the proof of the Lemma   for manifolds~\cite{Ta} goes through the category of oriented stacks.  
\end{proof}

\noindent{{\sc Proof of Theorem~\ref{BV}}.
We already have proved that  $(\BH\lcom(\LXX),\star)$ is a graded commutative algebra (see Theorem~\ref{th:Loop}) and that the operator $D\:\BH\lcom(\LXX)\to \BH_{\bullet+1}(\LXX)$ squares to zero: $D^2=0$ (Lemma~\ref{L:Dsquare}). Thus, we only need to prove   identity~\eqref{eq:BVidentity2} in order to prove Theorem~\ref{BV}.

By Lemma~\ref{L:Gerstenhaberbracket}, given $a,b,c\in H(\LXX)$, the left hand side of identity~\eqref{eq:BVidentity2} is $f\big([S^1]\times a\times b\times c\big)$ where $f:H(S^1\times \big(\LXX\big)^3) \to H (\LXX)$ is the composition
\begin{multline*}
 H(S^1\times \big(\LXX\big)^3) \stackrel{\id\times \Delta^!}\longrightarrow H(S^1\times \LXX\times \LXX\times_\XX\LXX) \\ \stackrel{\id\times m_*}\longrightarrow H(S^1\times \LXX\times \LXX)
 \stackrel{\Delta_\rho^!}\longrightarrow H(L_\rho(\XX\times \XX)) \stackrel{{m_\rho}_*}\longrightarrow H(\LXX).
\end{multline*}
We denote $L_{\rho,1}(\XX\times \XX\times \XX)$ the pullback 
\begin{eqnarray*}
\xymatrix{L_{\rho,1}(\XX\times \XX\times \XX) \ar[r] \ar[d]_{\mu} & S^1\times \LXX\times \LXX\times_\XX\LXX \ar[d]^{\id\times m} \\  
 L_\rho(\XX\times \XX) \ar[r]^{i_\rho} & S^1\times \LXX\times \LXX } 
\end{eqnarray*}
of $S^1\times \LXX\times \LXX\times_\XX\LXX$ along $i_\rho$. We also denote $m_{\rho,1}:L_{\rho,1}(\XX\times \XX\times \XX)\to \LXX$ the composition $m_\rho \circ \mu$.
Consider the following tower of pullback squares:
\begin{eqnarray}\label{eq:Lrhodef}
 \xymatrix{L_{\rho,1}(\XX\times \XX\times \XX) \ar[rr] \ar[d]_{\mu} & & S^1\times \LXX\times \LXX\times_\XX\LXX \ar[d]^{\id\times m}   \\ L_{\rho}(\XX\times \XX) \ar[rr]^{i_\rho} \ar[d] & & S^1\times \LXX\times \LXX \ar[d]^{\ev_{\rho}} \\ \XX \ar[rr]^{\Delta} && \XX\times \XX }
\end{eqnarray}
The naturality of Gysin maps with respect to this tower yields that 
\begin{eqnarray}\label{eq:abstarc}\{a;b\star c\} & = & {m_{\rho,1}}_* \circ \Delta_{\rho,1}^! \circ (\id\times \Delta^!) \big( [S^1]\times a\times b\times c\big)\end{eqnarray} where $\Delta_{\rho,1}^!:H(S^1\times \LXX\times \LXX\times_\XX\LXX)\to H(L_{\rho,1}(\XX\times \XX\times \XX))$ is the Gysin map obtained by pulling-back the Gysin map of the strongly oriented diagonal $\Delta: \XX\to \XX\times \XX$ along    $\ev_{\rho}\circ (\id\times m)$. 

Now the main point is to analyze the composition $\Delta_{\rho,1}^! \circ (\id\times \Delta^!)$. By Section~\ref{S:Gysin} each Gysin map is obtained by pulling back the  the normally non-singular diagram of the diagonal $\Delta:\XX\to \XX \times \XX$ (see Definition~\ref{D:oriented}) along, respectively, $\ev_{\rho}\circ (\id\times m)$ and $\id \times \ev_0\times \ev_0$ and taking the Thom class of the associated diagram. Recall that each nns diagram yieds a tubular neighborhood in the sense of Definition~\ref{D:tubular} (after replacing the target by a fiber bundle) and similarly after taking pullbacks. Then the composition $\Delta_{\rho,1}^! \circ (\id\times \Delta^!)$ is the product in the bivariant theory (see Section~\ref{S:Bivariant}) of these Thom classes.  It is essentially obtained by considering a fiber products (over $S^1\times \LXX\times \LXX\times_{\XX}\LXX$) of the above pulled-back normally non singular diagrams; more precisely by taking suitable pullbacks of the (pulled-back along $\ev_{\rho}\circ (\id\times m)$ and $\id \times \ev_0\times \ev_0$) tubular neighborhoods (as in Section~\ref{SS:products}) and composing them in a way similar to the proof of Lemma~\ref{L:regularcompose}. As in Definition~\ref{D:oriented}, we let $\theta_\Delta \in H^d(\XX\to \XX\times \XX)$ be the strong orientation class of $\XX$ and $\theta_{\Delta^{(2)}}\in H^{2d}(\XX\to \XX\times \XX\times \XX)$ be the strong orientation class of the iterated diagonal (see Corollary~\ref{C:orienteddiagonal}).

In order to carry on the analysis, we divide the circle $S^1$ into the joint of 4-intervals $I_i$ ($i=1,\dots,4$) corresponding to $[0,1/4]$, $[1/4, 1/2]$, $[1/2, 3/4]$ and $[3/4,1]$ (here we identify $S^1=[0,1] \slash (0\sim 1)$) with the obvious identifications. Note that the regular embeddings (induced by the nns diagram of the diagonals) inducing the Thom classes can be obtained by gluing together the regular embeddings restricted over each $I_i\times \LXX\times \LXX\times_{\XX}\LXX$, that is by taking the fiber product (over the regular embeddings obtained by restricting to $\{i/4\}\times  \LXX\times \LXX\times_{\XX}\LXX$)   of the regular embeddings over each  $I_i\times \LXX\times \LXX\times_{\XX}\LXX$. We first consider a restriction of $S^1\times (\LXX)^3$ to $[1/2,1]\times  (\LXX)^3$. It yields a commutative diagram of pullback squares:
\begin{eqnarray}\label{eq:halfto1}
 \xymatrix{ P_{\rho,1}(\XX\times \XX\times \XX) \ar[r] \ar@{^{(}->}[d] &  
 [\frac{1}{2},1]\times \LXX\times \LXX\times_\XX\LXX \ar@{^{(}->}[d] \ar[r]&  
 [\frac{1}{2},1]\times (\LXX)^3\ar[dd]^{\ev_{\rho\times \rho}}
\\ L_{\rho,1}(\XX\times \XX\times \XX) \ar[r]^{i_\rho} \ar[d] &  S^1\times \LXX\times \LXX\times_\XX\LXX \ar[d]^{\ev_{\rho}^{(2)}} &  \\ \XX^2 \ar[r]_{\Delta\times \id} & \XX^3 \ar[r]_{\tau_{34}\circ (\id \times \Delta\times\id)} & \XX^4 }
\end{eqnarray}
where $\tau_{3,4}:\XX^4\to \XX^4$ switches the last two factors. The vertical map $\ev_{\rho}^{(2)}$  is the composition  $\Big(\ev_\rho\circ (\id\times m), \pi_1\circ\ev_\rho \Big) $ (where $\pi_1$ is the projection on the first component), that is the map defined by $\ev_\rho^{(2)}(t,\gamma,\beta,\eta)= \big(\ev_0(\gamma), \ev_{2t}(\beta),\ev_0(\gamma)$ for $0\leq t\leq 1/2$ and $\ev_\rho^{(2)}(t,\gamma,\beta,\eta) = (\ev_{2t-1}(\gamma), \ev_{1}(\beta),\ev_{2t-1}(\gamma)) $ for $1/2\leq t\leq 1$. Similarly, the map $\ev_{\rho\times \rho}$ is the map defined by $\ev_{\rho\times \rho}(t,\gamma,\beta,\eta)=  (\ev_{2t-1}(\gamma), \ev_{1}(\beta),\ev_{2t-1}(\gamma), \ev_{1}(\eta))$ (for $1/2\leq t\leq 1$). By the above diagram~\eqref{eq:halfto1}, the (restriction to $[1/2,1]\times \big(\LXX\big)^3$ of the) product  ${\ev_{\rho}^{(2)}}^*(\theta_\Delta)\cdot \big(\id\ev_0\times _ev_0\big)^*(\theta_{\Delta})$ is computed by   
${\ev_{\rho\times \rho}}^*(\theta_{\Delta^{(2)}})$.

\smallskip

We let $\Delta_2$ be the standard two dimensional simplex $\Delta_2:=\{(u,s),\, 0\leq u\leq s\leq 1\}$ and consider the following diagram (in which $P_{H}(\XX\times \XX\times \XX)$ is defined so that the square is cartesian)
\begin{eqnarray*}
\xymatrix{\LXX &  P_{H}(\XX\times \XX\times \XX) \ar[l]_{m_H\hspace{11pt}} \ar[r] \ar[d] & \Delta_2 \times \LXX\times \LXX\times \LXX \ar[d]^{\ev_{H}}\\ 
& \XX\times \XX \ar[r]^{\Delta \times \Delta} & \XX^4}
\end{eqnarray*}
Here the vertical (evaluation) map $\ev_H$ is defined by
$\ev_H\big((u,s),a,b,c\big) =\big(\ev_u(a), \ev_0(b),\ev_s(a), \ev_0(c)\big)$ and $m_H:P_{H}(\XX\times \XX\times \XX)\to \LXX$ is defined similarly to the map $m_\rho: L_\rho(\XX\to \XX) \to \LXX$; that is, $m_H\big((u,s),a,b,c\big)$ (with $u,s,a,b,c$ satisfying the pullback relation) is the loop  starting at $\ev_0(a)$, describing $a$ until it reaches $\ev_u(a)$, then describing $b$ back to $\ev_0(b)=\ev_u(a)$, following $a$ again until it reaches $\ev_s(a)$, then following $c$ until reaches $\ev_1(c)=\ev_s(a)$ again and finishes to follow $a$ back to $\ev_1(a)=\ev_0(a)$. Here we use an informal description of the map $m_H$ that can be obtained by a pinching procedure as explained before Lemma~\ref{L:Gerstenhaberbracket} above.   

Let $\iota: [1/2,1] \to \Delta_2$ be the map $\iota(t)=(2t,2t)$. The diagram
\begin{eqnarray}\label{eq:halfto1homotopy}
\xymatrix{
 & P_{\rho,1}(\XX\times \XX\times \XX) \ar[ld]_{m_{\rho,1}} \ar[d] \ar[r] &    
 [\frac{1}{2},1]\times (\LXX)^3\ar[d]_{\iota\times \id} \ar@/^3pc/[dd]^{\ev_{\rho\times \rho}}
\\
\LXX &  P_{H}(\XX\times \XX\times \XX) \ar[l]_{m_H\hspace{11pt}} \ar[r] \ar[d] & \Delta_2 \times (\LXX)^3 \ar[d]_{\ev_{H}}\\ 
& \XX\times \XX \ar[r]^{\Delta \times \Delta} &  \XX^4}
\end{eqnarray}
is commutative. Since the subset $\{(u,s),\, u=s\}\subset \Delta_2$ is a boundary component of the dimension 2 simplex, it follows from Lemma~\ref{L:invariance} and diagram~\eqref{eq:halfto1homotopy} that, to compute the Thom class of the composition of the top lines in diagram~\eqref{eq:halfto1}, we can replace the top line of diagram~\eqref{eq:halfto1} by the restriction to the subset $\{u=0\mbox{ or } s=1\}\subset \Delta_2$ (that is the complementary boundary of $\Delta_2$)  of the top line of diagram~\eqref{eq:halfto1homotopy}.

Now we analyze the restriction  of diagram~\eqref{eq:halfto1homotopy} to $I_4\times (\LXX)^3=[3/4,1]\times (\LXX)^3$ which reduces to the analysis of the pullback square 
\begin{eqnarray}\label{eq:quarterto1}
\xymatrix{ {P_4}_{\rho}(\XX\times \XX\times \XX) \ar[r] \ar[d] &  
 \big[\frac{3}{4},1\big]\times (\LXX\times \LXX)\times \LXX   \ar[d]^{\ev_\rho\times \ev_0} 
   \\ \XX^2 \ar[r]^{\Delta\times \Delta} & \XX^4  }
\end{eqnarray} and its induced tubular neighborhood (given by the pullback along $\ev_\rho\times \ev_0$ of the nns diagram of the diagonal). 
Here the vertical map $\ev_\rho\times \ev_0$ is  the map $\ev_\rho\times \ev_0(t,a,b,c)=\big(\ev_{4t-3}(a), \ev_0(b), \ev_0(a),\ev_0(c)\big)$ (and ${P_4}_{\rho}(\XX\times \XX\times \XX)$ is defined by the pullback property). Furthermore, the restriction of the map $m_H\:P_H(\XX\times \XX\times \XX) \to \LXX$ (considered above) to this boundary component identifies with the map $m_{4,\rho}: {P_4}_{\rho}(\XX\times \XX\times \XX)\to \LXX$ which maps $(t,a,b,c)$ to the loop starting at $\ev_0(a)$, following $a$ until it reaches $\ev_{4t-3}(a)=\ev_0(b)$ then follows $b$ until it gets back to $\ev_1(b)=\ev_0(b)$, follows $a$ back to $\ev_1(a)=\ev_0(c)$ and then goes through $c$. Similarly, the restriction of diagram~\eqref{eq:halfto1homotopy} to $I_3\times (\LXX)^3 =[1/2, 3/4]\times (\LXX)^3$  yields a  pullback ${P_3}_{\rho}(\XX\times \XX\times \XX)$ similar to the pullback~\eqref{eq:quarterto1} as well as a map $m_{3,\rho}\:{P_3}_{\rho}(\XX\times \XX\times \XX)\to \LXX$. 

Furthermore, restricting diagram~\eqref{eq:Lrhodef} to $[0,1/4]\times (\LXX)^3$  yields a  cartesian square
\begin{eqnarray}\label{eq:0toquarter}
\xymatrix{ {P_1}_{\rho}(\XX\times \XX\times \XX) \ar[r] \ar[d] &  
 [0,\frac{1}{4}]\times \LXX\times \LXX\times \LXX   \ar[d]^{\ev_\rho\times \ev_0} 
   \\ \XX^2 \ar[r]^{\Delta\times \Delta} & \XX^4  }
\end{eqnarray}
where $\ev_\rho\times \ev_0(t,a,b,c)=\big(\ev_0(a),\ev_{4t}(b),\ev_0(b),\ev_0(c)\big)$ and,  similary, restricting $[1/4,1/2]\times (\LXX)^3$, a cartesian diagram
\begin{eqnarray}\label{eq:quartertohalf}
\xymatrix{ {P_2}_{\rho}(\XX\times \XX\times \XX) \ar[r] \ar[d] &  
 [\frac{1}{4},\frac{1}{2}]\times \LXX\times \LXX\times \LXX   \ar[d]^{\ev_\rho\times \ev_0} 
   \\ \XX^2 \ar[r]^{\Delta\times \Delta} & \XX^4  }
\end{eqnarray}
where $\ev_\rho\times \ev_0(t,a,b,c)=\big(\ev_0(a),\ev_{4t-1}(c),\ev_0(b),\ev_0(c)\big)$. 
Note also that restricting the map  $m_{\rho,1}:L_{\rho,1}(\XX\times \XX\times \XX)\to \LXX $ to $I_1=[0,1/4]$ and $I_2=[1/4,1/2]$ gives two maps $m_{1,\rho}:{P_1}_{\rho}(\XX\times \XX\times \XX)\to \LXX$ and  $m_{2,\rho}:{P_2}_{\rho}(\XX\times \XX\times \XX)\to \LXX$.

Now, to compute  ${m_{\rho,1}}_* \circ \Delta_{\rho,1}^! \circ (\id\times \Delta^!)$, \emph{i.e.}, the Thom classes $(\ev_\rho\times \id)\circ (\id\times m)^*(\theta_\Delta) \cdot (\id\times \ev_0\times \ev_0)(\theta_\Delta)$ we are left to study the pullbacks of the tubular neighborhood induced by the normally non-singular diagrams of the iterated diagonal $\Delta\times \Delta$ along the 4 various $\ev_\rho\times \ev_0$-maps (corresponding to the spaces $P_{i,\rho}(\XX\times \XX \times \XX)$) and to join them (as noted above). 
First, we remark that the restrictions to the points $\{1/4\}$ and $\{3/4\}$ of the cartesian squares~\eqref{eq:0toquarter}, \eqref{eq:quarterto1}, \eqref{eq:quartertohalf} and of the maps $m_{i,\rho}$ are identical. Thus it is enough to study first the joint $S^1_+:=[3/4,1]\cup[0,1/4]\slash (3/4\sim 1/4)$ of $I_1$ and $I_4$ with the boundary points identified and the Thom class induced by the cartesian squares~\eqref{eq:quarterto1},~\eqref{eq:0toquarter}.
Note that diagram~\eqref{eq:quarterto1} factors into the following diagram whose vertical squares  and top horizontal square are cartesian:
{\small\begin{eqnarray}\label{eq:quarterto1factorized}
\xymatrix@=8pt@R=8pt@M=8pt{ \LXX & \LXX\times_{\XX} \LXX \ar[l]_{m_\rho} \ar@^{^{(}->}[rr] \ar@{-}[d]& & \LXX\times \LXX \ar@{-}[d]&  \\
{P_4}_{\rho}(\XX\times \XX\times \XX) \ar[u]^{m_{4,\rho}} \ar[ur]^{m_{4,\rho}} \ar@^{^{(}->}[rr] \ar[dd]& \ar[d] & P_{4\rho}(\XX\times \XX)\times \LXX \ar@^{^{(}->}[rr] \ar[dd]\ar[ur]^{m^+_{\rho}\times \id}& \ar[d]^{\ev_0\times \ev_0} & [\frac{3}{4},1]\times (\LXX)^3 \ar[dd]^{\ev_{\rho}\times\ev_0}\\
 & \XX \ar@{-}[r]_{\Delta\hspace{12pt}}& \ar[r] & \XX\times \XX & \\
\XX^2 \ar[rr]_{\id \times \Delta} \ar[ru]_{\pi_2} & &\XX\times \XX^2 \ar[rr]_{ \Delta  \times \id} \ar[ru]_{\pi_{2,3}}& &\XX^4
 }
\end{eqnarray}}
Here $\pi_2$ and $\pi_{2,3}$ denote the projection on the last factors.  By functoriality and naturality of the construction of Gysin maps, we get that $${m_{4,\rho}}_*\big((\ev_\rho\times \ev_0)^*(\theta_{\Delta\times \Delta}) \big)= {m_\rho}_*\big((\ev_0\times \ev_0)^*(\theta_{\Delta})\big)\cdot \,  {(m_{\rho}^+\times \id)}_*\big((\ev_\rho\times \ev_0)^*(\theta_{\Delta}\big)$$ There is a  diagram similar to~\eqref{eq:quarterto1factorized} associated to the cartesian square~\eqref{eq:0toquarter}. Joining these two diagrams, we get:
{\small\begin{eqnarray}\label{eq:S1plusfactorized}
\xymatrix{ \LXX & \LXX\times_{\XX} \LXX \ar[l]_{m_\rho\hspace{5pt}} \ar@^{^{(}->}[rr] \ar[dd]& & \LXX\times \LXX \ar@{-}[d]&  \\
&  & L_{1,\rho}(\XX\times \XX)\times \LXX \ar@^{^{(}->}[rr] \ar[dd]\ar[ur]^{m_{\rho}^+\times \id}& \ar[d]^{\ev_0} & S^1_+\times (\LXX)^3 \ar[dd]^{\ev_{\rho}}\\
 & \XX \ar@{-}[r]_{\Delta\hspace{12pt}}& \ar[r] & \XX\times \XX & \\
 & &\XX  \ar[rr]_{ \Delta  } & &\XX^2
 }
\end{eqnarray}}
which exhibits the restriction to $S^1_+\times \big(\LXX\big)^3$ of  ${m_{1,\rho}}_*\big(\Delta_{\rho,1}^! \circ (\id\times \Delta^!)\big)$ as ${m_\rho}_*\Big( \ev_0^*(\theta_{\Delta})\cdot \, ({m_{\rho}^+\times \id})_*\big(\ev_0^*(\theta_\Delta)\big) \Big)$. 
Note that the above diagram~\eqref{eq:S1plusfactorized} is precisely the diagram defining $\{a,b\}\star c$; in other words, 
\begin{multline*}{m_\rho}_*\Big( \ev_0^*(\theta_{\Delta})\cdot \, ({m_{\rho}^+\times \id})_*\big(\ev_0^*(\theta_\Delta)\big) \Big)\big( [S^1_+]\times a\times b\times c\big) \\ =\; {m_\rho}_*\circ\Delta^!\circ ({m_{\rho}^+\times \id})_*\circ \Delta_\rho^! \big([S^1_+]\times a\times b\times c\big)\; = \; \{a,b\}\star c\end{multline*}

The above arguments for $S^1_+$  apply  similarly to study the joint  $S^1_{-}:=[3/4,1]\cup[0,1/4]\slash (3/4\sim 1/4)$ of $I_1$ and $I_4$ with the boundary points identified. It yields that 
$${m_\rho}_*\Big( \ev_0^*(\theta_{\Delta})\cdot \, ({m_{\rho}^-\times \id})_*\big(\ev_0^*(\theta_\Delta)\big) \Big)\big( [S^1_+]\times a\times b\times c\big) \; = \, (-1)^{|b| |c|} \{a,c\}\star b$$
where the sign comes from the fact that ones has to exchange $b$ and $c$ in that case (with a transposition similar to the one appearing in the bottom line of diagram~\eqref{eq:halfto1}). Recall from above that the restrictions to the points $\{1/4\}$ and $\{3/4\}$ of the cartesian squares~\eqref{eq:0toquarter}, \eqref{eq:quarterto1}, \eqref{eq:quartertohalf} and of the maps $m_{i,\rho}$ are identical. It follows that the computation of the Gysin maps for $S^1\times \big(\LXX\big)^3$ factors through the one of $(S^1_+\vee S^1_{-})\times \big(\LXX\big)^3$. Thus,  we deduce from the above computations for $S^1_+$ and $S^1_{-}$ and identity~\eqref{eq:abstarc} that
$$ \{a,b\star c\} \;= \; \{a,b\}\star c+ (-1)^{|b| |c|} \{a,c\}\star b$$
that is identity~\eqref{eq:BVidentity2} holds, by graded commutativity of the loop product. 
 \nolinebreak $\Box$ }

\section{Homological conformal field theory  and free loop stacks}\label{S:HCFTbig}
In this section, we extend the \BV-structure of Theorem~\ref{BV} and the Frobenius structure of 
Theorem~\ref{frobeniusloop} into the whole structure of an \emph{homological conformal field theory}
 (with positive boundaries) following ideas of Cohen-Godin~\cite{CoGo, Go} for manifolds and 
Chataur-Menichi~\cite{CM} for classifying spaces of groups. As in Section~\ref{frobenius} 
(and for the same reasons), we assume in this section  that our ground ring $k$ is a field.  Note that, unlike in Godin's paper~\cite{Go}, we will only allow \emph{closed} boundaries and a positive number of \emph{both} incoming and outgoing boundaries components.

\subsection{Quick review on Homological Conformal Field theory with positive boundaries }
We start by recalling some definitions of Homological Conformal Field theories. We will make strong restrictions on the type of boundary we consider (which simplify greatly the theory). We follow~\cite{Cos, Cos2, Go, Seg04}.

\smallskip

A (closed)\footnote{unless otherwise stated, we only consider closed cobordism in this paper} homological conformal field theory is an algebra over the PROP of the homology of the stack 
(or moduli space) of compact oriented Riemann surfaces or, equivalently a symmetric monoidal functor 
from the (homology of the) Segal category of Riemann surfaces~\cite{Seg04} to the 
category of graded vector spaces.  
 Let us start by giving more details on what
 this definition means, following~\cite{CoVo, Cos, Cos2, Get}. 
We first recall that, a \emph{complex cobordism} from a family $\coprod_{i=1}^n S^1$ of circles 
to another family $\coprod_{i=1}^m S^1$ of circles is a closed (non-necessarily connected) Riemann 
surface $\Sigma$ equipped with two holomorphic embeddings (with disjoint images) 
$\rho_{in}\: \coprod_{i=1}^n D^2 \hookrightarrow \Sigma$ and $\rho_{out}\: \coprod_{i=1}^M D^2$ 
of closed disks. The image of $\rho_{in}$ is called the \emph{incoming} boundary and the image of
 $\rho_{out}$ the outgoing boundary. Two complex cobordism $\Sigma_1$ and $\Sigma_2$ 
(from $\coprod_{i=1}^n S^1$ to $\coprod_{i=1}^m S^1$) are equivalent if there exists a 
biholomorphism $h\: \Sigma_1\stackrel{\sim}\to \Sigma_2$ which fixes the boundary (\emph{i.e.} 
commutes with $\rho_{in}$ and $\rho_{out}$). We denote $\MM_{n,m}$ the moduli space of equivalences 
classes of complex cobordism from $\coprod_{i=1}^n S^1$ to $\coprod_{i=1}^m S^1$, that is the (coarse moduli space of the)
differentiable stack $[ \SS_{n,m}/Bihol]$ obtained as the quotient of the space $\SS_{n,m}$ of 
holomorphic embeddings of disks inside compact Riemann surfaces by the group of biholomorphism 
fixing the boundary. Note that, there is an isomorphism of stacks 
$[ \SS_{n,m}/Bihol]\cong \coprod_{[\Sigma]}[pt/ \Gamma_{n,m}(\Sigma)]$ where the union is over a set of 
representatives of the isomorphisms classes of cobordisms (with $n$ incoming and $m$ outgoing closed boundaries 
components) and $\Gamma_{n,m}(\Sigma)\:= \pi_0\Big(\mathop{Diff}^+_{n,m}(\Sigma)\Big)$ is the isotopy classes of
the group $\mathop{Diff}^+_{n,m}(\Sigma)$ of oriented diffeomorphisms preserving the boundaries pointwise of a 
surface $\Sigma$ with $n$ incoming and $m$ outgoing closed boundary components.

\smallskip

The disjoint union of surfaces yields a canonical morphism 
$\MM(n,m) \times \MM(n',m')\to \MM(n+n',m+m')$. Further, given $\Sigma_1 \in \MM(\ell,n)$ and
 $\Sigma_2 \in \MM(n,m)$, using the embeddings of disks 
$ \Sigma_1 \hookleftarrow \coprod_{i=1}^n D^2\hookrightarrow \Sigma_2$, we can glue $\Sigma_2$ on $\Sigma_1$ 
along their common boundary. We denote $\Sigma_2 \circ \Sigma_1 \in \MM_{\ell,m}$ the Rieman surface 
thus obtained. Applying the singular homology functor to the above operations yields linear map 
$$H\lcom(\MM_{n,m)}) \otimes H\lcom(\MM_{n',m'})\stackrel{H\lcom(\coprod)}\to H\lcom(\MM_{n+n',m+m'}) $$ and 
$$H\lcom(\MM_{\ell,n})\otimes H\lcom (\MM_{n,m)}) \stackrel{H\lcom(\circ)}\longrightarrow H\lcom(\MM_{\ell,m})  
$$ that satisfies natural associativity and compatibility relations.
 It follows that the collection $\big(H_\bullet(\MM_{n,m})\big)_{n,m\geq 0}$ are 
the morphisms of a graded linear symmetric monoidal category $\mathcal{C}_{\MM}$ whose objects are 
the nonnegative integers $n\in \mathbb{N}$ and the monoidal structure is induced by $k\otimes \ell= k+\ell$ 
on the objects and disjoint union of surfaces on morphisms. 
An \emph{homological conformal field theory} is a symmetric monoidal functor from the category 
$\mathcal{C}_\MM$ to the symmetric monoidal category of graded vector spaces (equipped with the usual graded
 tensor product). 
  Informally, this definition simply means that an
 homological conformal field theory is a graded vector space $A$ with an operation 
$\mu(c)\: A^{\otimes n}\to A^{\otimes m}$  for any homology class $c \in H\lcom(\MM_{n,m})$ such 
that $\mu(c\circ d)=\mu(c) \circ \mu(d)$ and $\mu(c\coprod d) =\mu(c)\otimes \mu(d)$.

\smallskip

Unlike oriented closed manifolds, oriented stacks does not have unit for the loop product 
(nor counit) in general (see Section~\ref{Liegroups} for instance). 
This forces us to consider \emph{non-unital and non-counital} homological conformal field theory 
which are symmetric monoidal functor from the category $\mathcal{C}_{\MM}^{nu,nc}$ to the 
category of graded vector spaces , where $\mathcal{C}_{\MM}^{nu,nc}\subset \mathcal{C}_{\MM}$ 
is the (monoidal) subcategory  obtained 
by considering only cobordisms in $\MM_{n,m}$ for which every connected component has at least one 
ingoing 
\emph{and} one outgoing boundary component. We may also refer to such algebraic structure as an
 homological conformal field theory with positive (closed) boundaries.

\smallskip

We wish to make the homology $H\lcom(\LXX)$ of the free loop stack of an oriented stack an homological 
conformal field theory (without (co)units). However,
since the basic operations we consider  are non-trivially graded (for instance the loop product is of 
degree $\dim(\XX)$), we need to plug in a notion of dimension in the definition of conformal field 
theories to take care of this phenomenon and encode the sign issues. We follow the ideas and presentation 
of Costello~\cite{Cos2} and Godin~\cite{Go}, where
the grading is taken into account by a local coefficient system $\det^{\otimes d}$ on the moduli spaces 
$\MM_{n,m}$. The local coefficient system $\det$ is a graded invertible locally constant sheaf (\emph{i.e.} a 
graded $k$-linear locally constant sheaf  of dimension 1). 
 To a closed Riemann surface $\Sigma\in \MM_{n,m}$, we associate   a compact Riemann  surface $\Sigma^{bd}$ with
 boundary by removing from $\Sigma$ the interior of (the images of) the closed disks 
$ \rho_{in}\: \coprod_{i=1}^n D^2 \hookrightarrow \Sigma$ and $\rho_{out}\: \coprod_{i=1}^m D^2 
\hookrightarrow \Sigma$. Restricting $\rho_{in}$ to the boundary $\coprod_{i=1}^n S^1$ 
of the disks, we get a diffeomorphism from $\coprod_{i=1}^n S^1$ onto the incoming boundary of $\Sigma^{bd}$. 
 Following~\cite{Cos2} and~\cite[Section 4.1]{Go}, we define the fibre of the local coefficient system $\det$ at a
 surface $\Sigma\in \MM_{n,m}$ to be 
 $$\det(\Sigma)\:= \det\Big(H^0\big(\Sigma^{bd}, \rho_{in}\big(\coprod_{i=1}^n S^1\big)\big)\Big) 
\otimes \det\Big(H_1\big(\Sigma^{bd}, \rho_{in}\big(\coprod_{i=1}^n S^1\big)\big)\Big). $$ 
Here, given a finite dimensional $k$-vector space $V$, $\det$ denotes the determinant, that is  $\det(V)=\bigwedge^{\dim(V)}V $ 
is the top exterior power of $V$, and we consider the (relative) homology groups of a pair.  
 This defines the local coefficient system $\det$ on $(\MM_{n,m})_{n,m\geq 0}$ and similarly, for an integer 
$d$, the local coefficient system $\det^{\otimes d}$ obtained by tensoring $\det$ with itself $d$-times.
 
 It is proved in~\cite{Cos2, Go, CM} that the composition of surface induces a natural isomorphism
 $\det(\Sigma_2)\otimes \det(\Sigma_1)\to \det(\Sigma_2 \circ \Sigma_1)$ which is associative and compatible with
 the canonical isomorphism $\det(\Sigma'_1\coprod \Sigma'_2)\cong \det(\Sigma'_1)\otimes \det(\Sigma'_2)$.
  This allows to see the collection of homology groups
 $\Big(H\lcom\big(\MM_{n,m}; \det^{\otimes d}\big)\Big)_{n,m\geq 0}$ with value in the local coefficient
 $\det^{\otimes d}$ as the morphism of a graded linear symmetric monoidal category
 $\mathcal{C}_{\MM, \det^{\otimes d}}$, and as above we also get a graded linear symmetric monoidal category
 $\mathcal{C}^{nu, nc}_{\MM, \det^{\otimes d}}$ by restricting to cobordism with at least one incoming and one
 outgoing boundary on each connected component. According to~\cite{Cos2, Go}, a (non-unital, non-counital) 
\emph{$d$-dimensional homological conformal field theory} is a symmetric monoidal functor from the category
 $\mathcal{C}^{nu, nc}_{\MM, \det^{\otimes d}}$ to the category of graded vector spaces.

\subsection{The Homological Conformal Field Theory with positive closed boundaries associated to free loop stacks}\label{S:HCFT}

It is known (see~\cite{Get, CoVo, MSS}) that an (non-unital) homological conformal field theory carries 
the structure of a (non-unital) \BV-algebra as well as that of a Frobenius algebra (without (co)unit). 
For instance the associative and commutative operation of the {\BV} or Frobenius structure is induced by 
the pair of pants surface lying in $\MM_{2,1}$. The main result of this Section enriches the above structures already obtained for loop stacks into an HCFT over $\mathcal{C}^{nu, nc}_{\MM, \det^{\otimes d}}$.

\begin{thm}\label{T:HCFT} Let $\XX$ be an oriented (Hurewicz\footnote{Recall that any differentiable stack is Hurewicz}) stack of dimension $d$. There is a $d$-dimensional 
non-unital, non-counital homological conformal field theory on the homology $H\lcom(\LXX)$ of the free loop stack
 which induces 
 the \BV-algebra and Frobenius structure on the homology $H\lcom(\LXX)$ given by Theorem~\ref{BV} and 
Theorem~\ref{frobeniusloop}.
\end{thm}

\begin{rmk}
 The proof of the Theorem~\ref{T:HCFT} actually implies the ones of Theorems~\ref{BV} and Theorem~\ref{frobenius}
as well. However, this proof does \emph{not} apply to prove similar statements 
(for instance Theorem~\ref{frobeniusstring}) for inertia stack (and thus to define the intersection pairing as in Section~\ref{Orbifolds}) or any other family of groups over a stack 
considered in Section~\ref{S:hidden}. Further, it is not obvious that this proof will also applied to the twisted versions of the loop product studied in Section~\ref{Loopproduct} and aforementioned in Section~\ref{SS:OrbifoldsHidden}.
\end{rmk}
To prove the above Theorem~\ref{T:HCFT},
 we follow the approach of~\cite{CM, Go, Cos2}, using chord diagrams/ribbon graphs, but using a stack point of
 view (instead of a purely homotopical one) and the benefits of the bivariant theory of Section~\ref{S:Bivariant}. 
We will first determine the value on a particular cobordism $\Sigma_{g,n,m}$ of the HCFT, which will be given by the 
linear map~\eqref{eq:HCFTop} below.

First, we need to recall some preliminaries on Sullivan's chord diagrams and fat graphs which are taken 
from~\cite{I, CoGo, Cos}. By a graph,  we mean a pair $G=(V,H)$ consisting of a finite set of vertices $V$, 
of half-edges (which can be thought as oriented edges) $H$ equipped  with a map $s: H\to V$ and an involutive map
 with no-fixed points $ e\mapsto \overline{e}$ on the set of half-edges.   A fat graph is a graph equipped, at
 each vertex $v$, with a cyclic ordering of the half-edges emanating from $v$. The geometric realization of a fat 
graph is thus a 1 dimensional cell complex plus extra data. It is well-known that the classifying space of fat 
graphs is equivalent to the moduli space of Riemann surfaces (see~\cite{I, Cos, Cos2, Ti, Go} for much more
 precised statements). In particular, every (isomorphism class of a) Riemann surface $\Sigma\in \MM_{n,m}$ is a 
deformation retract of (the geometric realization of) a fat graph with $n$-incoming boundary cycles and m outgoing 
ones (we refer to~\cite{I, CoGo} for the definition of these boundary cycles).

A chord diagram (see~\cite{CoGo}) is a special kind of (geometric realization) of a fat graph.  
A \emph{chord diagram} of type $(g,n,m)$ is a  union of $n$ disjoint circles with a disjoint unions of trees whose
 endpoints are glued on the circles (on distinct points), and such that the induced cell complex is the   
(geometric realization) of a fat graph representing a surface of genus $g$ with $n+m$ boundary components. The
 first circles are refered to as the incoming circles and the set of (necessarily path-connected) trees will be 
denoted $\mathcal{T}(c_{g,n,m})$. The set of points in the circle in which the endpoints of the trees lies is 
denoted $\mathcal{V}(c_{g,n,m})$ (and called the set of circular vertices). 

We let $c_{g,n,m}$ denote a chord diagram of type $(g,n,m)$ and $\Sigma_{g,n,m}$ be the surface represented by 
$c_{g,n,m}$. Here $g$ is the genus of $\Sigma_{g,n,m}$, that is  the sum of the genera of its components. In
 particular, the Euler characteristic of $\Sigma_{g,n,m}$ is given by
 $\chi(\Sigma_{g,n,m})=2 \#\Sigma_{g,m,n} -2g-n-m$ where $\#\Sigma_{g,n,m}$ is the number of (arcwise) connected 
components of $\Sigma_{g,n,m}$. Then $c_{g,n,m}$ is a deformation retract of $\Sigma_{g,n,m}$. We let 
$r_{g,n,m}: \Sigma_{g,n,m}\to c_{g,n,m}$ be the retraction and $\iota_{g,n,m}:c_{g,n,m} \to \Sigma_{g,n,m}$  be 
the inclusion. Note that since every connected component is assumed to have positive incoming and outgoing boundary component,
 $\chi(\Sigma_{g,n,m})$ is always non-positive in our case.

Given a tree $t\in \mathcal{T}(c_{g,n,m})$, we can associate the subset $V(t)\subset \mathcal{V}(c_{g,n,m})$ of 
circular vertices given by the endpoints of $t$ and we get a canonical inclusion map 
$\coprod_{v \in \mathcal{V}(c_{g,n,m})} \{pt\} \to \coprod_{t\in \mathcal{T}(c_{g,n,m})} t$. Applying the mapping 
stack functor, we get a map 
$$d_{c_{g,n,m}}\:\prod_{t\in \mathcal{T}(c_{g,n,m})}\bfhom(t, \XX) \longrightarrow \XX^{\mathcal{V}(c_{g,n,m})}$$ 
which was already considered in~\cite{CM}.
\begin{lem}\label{L:orienttree}
 Let $\XX$ be an oriented stack of dimension $d$. Then there is an orientation class 
$$\theta_{d_{c_{g,n,m}}} \in H^{-d\chi(\Sigma_{g,n,m})}\Big(\prod_{t\in \mathcal{T}(c_{g,n,m})}\bfhom(t, \XX) \stackrel{d_{c_{g,n,m}}}\longrightarrow \XX^{\mathcal{V}(c_{g,n,m})}\Big).$$
\end{lem}

\begin{proof}
Each tree is a deformation retract onto any of its vertex, hence we have deformation retract 
$t \stackrel{\iota_t}{\underset{r_t}{\leftrightarrows}} pt $ for each $t\in \mathcal{T}(c_{g,n,m})$ and a
 factorisation
$$\xymatrix{\XX^{\mathcal{T}(c_{g,n,m})} \ar[rr]^{\hspace{-4pt}\prod r_t^*} 
\ar@/_1pc/[rrrr]_{\Delta^{(\#\mathcal{V}(c_{g,n,m}) - \#\mathcal{T}(c_{g,n,m}) )}}&&
\prod_{t\in \mathcal{T}(c_{g,n,m})}\bfhom(t, \XX) \ar[rr]^{\hspace{7pt} d_{c_{g,n,m}}} & & 
\XX^{\mathcal{V}(c_{g,n,m})} } $$
The bottom line is an iterated diagonal, hence strongly oriented by Corollary~\ref{C:orienteddiagonal}. 
Then we take $\theta_{d_{c_{g,n,m}}}$ to be the    pushforward $\big(\prod r_t^*\big)_*(\theta)$ of the orientation
 class $\theta \in H^{d(\#\mathcal{V}(c_{g,n,m}) - \#\mathcal{T}(c_{g,n,m}))}(\XX^{\mathcal{T}(c_{g,n,m})}
\to \XX^{\mathcal{V}(c_{g,n,m})})$ of this iterated diagonal.
 Since the Euler characteristic of the chord diagram agrees with the one of the surface it represents, we get
 $\#\mathcal{V}(c_{g,n,m}) - \#\mathcal{T}(c_{g,n,m})=-\chi(\Sigma_{g,n,m})$ as in~\cite{CoGo} and the lemma follows. 
\end{proof}

\subsection{Construction of the operations}
We now define the operations associated to (the isomorphism class of) a Riemann surface
 $\Sigma_{g,n,m} \in \MM_{n,m}$ (not necessarily connected). We will first define some quotient stacks of maping stacks by diffeomorphism groups. Let $c_{g,n,m}$ be  a chord diagram representing
  $\Sigma_{g,n,m}$. We still denote  $r_{c_{g,n,m}}: \Sigma_{g,n,m}\to c_{g,n,m}$  the retraction and
 $\iota_{c_{g,n,m}}:c_{g,n,m} \to \Sigma_{g,n,m}$ the inclusion, which yield an homotopy equivalence
$r_{g,n,m}\:\bfhom({c_{g,n,m},\XX})\stackrel{\sim}\longrightarrow \bfhom(\Sigma_{g,n,m},\XX)$ (and $\iota_{g,n,m}\:\bfhom(\Sigma_{g,n,m},\XX)\stackrel{\sim}\longrightarrow \bfhom(c_{g,n,m},\XX)$). 

The circular vertices, which are the  points where the trees in $\mathcal{T}(c_{g,n,m})$ are glued 
to the $n$ disjoint circles,  yields pushout diagram 
$$ \xymatrix{\coprod_{v\in \mathcal{V}(c_{g,n,m})}\{v\} \ar[d] \ar[rr]^{i_{\mathcal{V}}}&& \coprod_{i=1}^n S^1 \ar[d]\\
\coprod_{t\in \mathcal{T}(c_{g,n,m})} t \ar[rr]&& c_{g,n,m}    }$$
which, by Lemma~\ref{L:glue}, induces a pullback of stacks
\begin{equation}\label{eq:chordpullback} \xymatrix{\bfhom(c_{g,n,m},\XX)\ar[d] \ar[rr]&& (\LXX)^{n} \ar[d]^{\ev_{\mathcal{V}}}\\
\prod_{t\in \mathcal{T}(c_{g,n,m})}\bfhom(t, \XX) \ar[rr]^{\hspace{3pc} d_{c_{g,n,m}}} &&
 \XX^{\mathcal{V}(c_{g,n,m})} } \end{equation}
where $\ev_{\mathcal{V}}\: \big(\LXX\big)^n \to \XX^{\mathcal{V}(c_{g,n,m})}$ is the evaluation map induced by the 
inclusions $i_{\mathcal{V}}\: \coprod_{v\in \mathcal{V}(c_{g,n,m})}\{v\} \to \coprod_{i=1}^n S^1$.
We denote
\begin{equation}
\label{eq:thetachord}
\theta_{g,n,m}^{\XX} = \ev_{\mathcal{V}}^*(\theta_{d_{c_{g,n,m}}}) \in H^{-d\chi(\Sigma_{g,n,m})}\Big( \bfhom(c_{g,n,m},\XX)\longrightarrow \big(\LXX\big)^n\Big) 
\end{equation}
the bivariant class induced by the pullback diagram~\eqref{eq:chordpullback} and the orientation class $\theta_{d_{c_{g,n,m}}}$ 
of Lemma~\ref{L:orienttree}.

The retraction $r_{g,n,m}\:\bfhom({c_{g,n,m},\XX})\stackrel{\sim}\to \bfhom(\Sigma_{g,n,m},\XX)$ sits inside a 
commutative diagram:
\begin{equation}\label{eq:pushSigma}
\xymatrix{
& \bfhom(\Sigma_{g,n,m},\XX) \ar@/_1pc/[ld]_{\iota_{g,n,m}\hspace{1pc}}\ar[rrd]^{\rho_{in}}& &  \\
\bfhom(c_{g,n,m},\XX)\ar[ru]_{\hspace{1pc}r_{{g,n,m}}} \ar[rrr]&& & (\LXX)^{n} } 
\end{equation}
where the map $\bfhom(\Sigma_{g,n,m},\XX)\to (\LXX)^{n}$ is induced by the inclusion of the incoming boundary in 
$\Sigma_{g,n,m}$. Applying the pushforward map along $r_{g,n,m}$ given by the above diagram~\eqref{eq:pushSigma}
 to the bivariant class $\theta^\XX_{g,n,m}$ (see~\eqref{eq:thetachord}) gives us a bivariant class
\begin{equation}\label{eq:rorientationclass}
 {r_{g,n,m}}_*\big( \theta^\XX_{g,n,m}\big ) \in  
H^{-d\chi(\Sigma_{g,n,m})}\Big(\bfhom(\Sigma_{g,n,m},\XX) \longrightarrow \big(\LXX\big)^n \Big).
\end{equation}

\paragraph{Quotienting by diffeomorphisms}
To shorten notations, we write
$G_{\Sigma}$ for the group $\mathop{Diff}^+_{n,m}(\Sigma_{g,n,m})$ of oriented diffeomorphisms of $\Sigma_{g,n,m}$
 preserving the boundaries pointwise. The group $G_{\Sigma}=\mathop{Diff}^+_{n,m}(\Sigma_{g,n,m})$ acts on $\Sigma_{g,n,m}$
and thus on $\bfhom(\Sigma_{g,n,m},\XX)$ (by functoriality, see Section~\ref{mappingstack}) and the restriction map
 $\rho_{in}\:\bfhom(\Sigma_{g,n,m},\XX)\longrightarrow (\LXX)^{n}$ is equivariant (where the action on 
$(\LXX)^{n}$ is trivial since the boundary is pointwise fixed). Similarly to the construction of the transformation groupoid of a topological group acting on a space (see Section~\ref{S:topological}), we can pass to the quotient of the above stacks by the $G_\Sigma$-action:
\begin{lem}\label{L:quotientbyGSigma}
 The action of $G_\Sigma$ on $\bfhom(\Sigma_{g,n,m},\XX)$  induces a quotient topological stack $\big[\bfhom(\Sigma_{g,n,m},\XX)/G_{\Sigma}\big] $ together with an natural topological stack epimorphism  $p_{G_\Sigma}\:\bfhom(\Sigma_{g,n,m},\XX)\to \big[\bfhom(\Sigma_{g,n,m},\XX)/G_{\Sigma}\big]$ and, similarly, a quotient stack $\big[(\LXX)^{n}/G_{\Sigma}\big]$ with an natural topological stack epimorphism $p_{G_\Sigma}\:(\LXX)^{n})\to \big[(\LXX)^{n}/G_{\Sigma}\big]$ such that 
 \begin{enumerate}
  \item there is a cartesian square
  $$\xymatrix{G_\Sigma \times  \bfhom(\Sigma_{g,n,m},\XX) \ar[r] \ar[d]_{pr_2} & 
  \bfhom(\Sigma_{g,n,m},\XX)   \ar[d]^{p_{G_\Sigma}} \\
  \bfhom(\Sigma_{g,n,m},\XX)   \ar[r]^{p_{G_\Sigma}} & \big[\bfhom(\Sigma_{g,n,m},\XX)/G_{\Sigma}\big] }$$ where the top arrow is given by the $G_\Sigma$-action as well as a similar cartesian square with $(\LXX)^{n}$ instead of $\bfhom(\Sigma_{g,n,m},\XX)$;
  \item The map $p_{G_\Sigma}\:\bfhom(\Sigma_{g,n,m},\XX)\to \big[\bfhom(\Sigma_{g,n,m},\XX)/G_{\Sigma}\big]$ makes $\bfhom(\Sigma_{g,n,m},\XX)$ a $G_\Sigma$-torsor  and similarly for $(\LXX)^n$;
  \item There is a topological stack isomorphism $\big[(\LXX)^{n}/G_{\Sigma}\big] \cong (\LXX)^n \times \big[*/G_\Sigma\big]$ and the following diagram is commutative $$\xymatrix{\big[(\LXX)^{n}/G_{\Sigma}\big]  \ar[r]^{\!\!\!\simeq} & (\LXX)^n \times \big[*/G_\Sigma\big] \\
  (\LXX)^n \ar[u]^{p_{G_\Sigma}} \ar[ru]_{\id \times {q}_{G_\Sigma}} & } $$
 where ${q}_{G_\Sigma}\: * \to \big[*/G_{\Sigma}\big]$ is the canonical map. 
 \end{enumerate}
\end{lem}
\begin{proof}
 We know from Section~\ref{mappingstack}, that the stack $\bfhom(\Sigma_{g,n,m},\XX)$ is the fibered groupoid over $\Top$ given by the rule $ T\in \Top  \mapsto  \  \Hom(T\times \Sigma_{g,n,m},\XX)\,,$ where $\Hom$ is the groupoid of (stack) morphisms. 
 Then we can define another fibered groupoid $T\mapsto \mathfrak{M}_{\XX,G_\Sigma}(T) $ where the set of objects of $\mathfrak{M}_{\XX,G_\Sigma}(T)$ is the set of stack morphisms $\mathop{hom}(T\times \Sigma_{g,n,m},\XX)$ (here $\mathop{hom}(\XX,\YY)$ denotes the set of objects of the groupoid of stack morphisms $\Hom(\XX,\YY)$).
 The morphisms of $\mathfrak{M}_{\XX,G_\Sigma}(T) $ are $$\mathop{Mor}(\mathfrak{M}_{\XX,G_\Sigma}(T))=\Big\{(g,f)\in G_\Sigma \times \mathop{Mor}(\Hom(T\times \Sigma_{g,n,m},\XX))\Big\}$$ with source and target maps given by $s(g,f)=s(f)$ and $t(g,f)=t(g\cdot f)$ (where $\cdot$ denotes the action of $G_\Sigma$).
 The rule $T\mapsto \mathfrak{M}_{\XX,G_\Sigma}(T) $ is easily seen to define a prestack. 
 We let  $\big[\bfhom(\Sigma_{g,n,m},\XX)/G_{\Sigma}\big]$ be the stackification 
 of  $\mathfrak{M}_{\XX,G_\Sigma}(-)$. We define in the same way the stack $\big[(\LXX)^{n}/G_{\Sigma}\big]$ as the stackification of a fibered groupoid $\mathfrak{L}_{\XX,G_\Sigma}(T)$.
 Since the action of $G_\Sigma$ on $(\LXX)^n$ is trivial, there is an isomorphism (of fibered groupoids) in between 
 $$\mathfrak{L}_{\XX,G_\Sigma}(T):=\Big\{(g,f)\in G_\Sigma \times \mathop{Mor}(\Hom(T\times \coprod_{i=1}^n S^1,\XX))\Big\} \toto  \mathop{hom}(T\times \coprod_{i=1}^n S^1,\XX))$$
 and  $ G_\Sigma \times \mathop{Mor}(\Hom(T\times \coprod_{i=1}^n S^1,\XX)) \toto \mathop{hom}(T\times \coprod_{i=1}^n S^1,\XX)).$
 Hence an isomorphism of stacks $\big[(\LXX)^{n}/G_{\Sigma}\big] \cong (\LXX)^n \times \big[*/G_\Sigma\big]$. 
 
 Choosing $g=1$, the unit of $G_{\Sigma}$ induces a prestack morphism $p\:\bfhom(\Sigma_{g,n,m},\XX)\to \mathfrak{M}_{\XX,G_\Sigma}$ which yields the canonical map $p_{G_\Sigma}\:\bfhom(\Sigma_{g,n,m},\XX)\to \big[\bfhom(\Sigma_{g,n,m},\XX)/G_{\Sigma}\big]$ after stackification. This map is shown to be an epimorphism and a $G_\Sigma$-torsor as in the usual case of transformation groupoids (see~\cite{BX, Noohi, Mapping}).  The case of $(\LXX)^n$ is similar and assertion 3. of the Lemma follows immediately. 
 
 In order to prove that these quotient stacks are topological stacks, we note that 
 if $X\to \bfhom(\Sigma_{g,n,m},\XX)$ is a chart for $\bfhom(\Sigma_{g,n,m},\XX)$ (that is a representable epimorphism from a topological space), then the composition $X\to \bfhom(\Sigma_{g,n,m},\XX) \to \big[\bfhom(\Sigma_{g,n,m},\XX)/G_{\Sigma}\big]$ is again a representable epimorphism. The existence of $X$ is given by Proposition~\ref{P:MapSt}.

 By construction, the fibered groupoid $\mathfrak{M}_{\XX,G_\Sigma}$ is defined so that the diagram $$\xymatrix{G_\Sigma \times  \Hom(T\times \Sigma_{g,n,m},\XX) \ar[r] \ar[d]_{pr_2} & \Hom(T\times \Sigma_{g,n,m},\XX)
    \ar[d]^{p} \\
 \Hom(T\times \Sigma_{g,n,m},\XX)   \ar[r]^{p} & \mathfrak{M}_{\XX,G_\Sigma}(T) }$$
is 2-cartesian. Since the stackification functors commutes with 2-fiber products, it induces the cartesian square asserted in the Lemma.\end{proof}

Since the map $\rho_{in}$ is equivariant, it passes to the
quotient to give a stack morphism 
$\big[\rho_{in}/G_{\Sigma}\big]\:\big[\bfhom(\Sigma_{g,n,m},\XX)/G_{\Sigma}\big] \longrightarrow \big[(\LXX)^{n}/G_{\Sigma}\big]$. Furthermore, we have the following lemma
\begin{lem} \label{L:pullbackquotientbyGSigma} Let $p_{/{G_\Sigma}}:\big(\LXX\big)^n \to \big[(\LXX)^{n}/G_{\Sigma}\big]$ be the quotient map (of stacks) given by Lemma~\ref{L:quotientbyGSigma} and similarly for $\bfhom(\Sigma_{g,n,m},\XX)$. 
 The following diagram
 \begin{equation*}
 \xymatrix{ \bfhom(\Sigma_{g,n,m},\XX) \ar[rr]^{\hspace{2pc}\rho_{in}} \ar[d] && \big(\LXX\big)^n \ar[d]^{p_{/{G_\Sigma}}}\\
 \big[\bfhom(\Sigma_{g,n,m},\XX)/G_{\Sigma}\big] \ar[rr]^{\hspace{2pc}\big[\rho_{in}/G_{\Sigma}\big]} &&  \big[(\LXX)^{n}/G_{\Sigma}\big].}
\end{equation*}
is a cartesian square.
\end{lem}
\begin{proof}
 Unfolding the definition of $\big[\bfhom(\Sigma_{g,n,m},\XX)/G_{\Sigma}\big] $ and $\big[(\LXX)^{n}/G_{\Sigma}\big]$ in the proof of Lemma~\ref{L:quotientbyGSigma}, we see that, since $\rho_{in}$ is $G_\Sigma$-equivariant it induces a map of fibered groupoids
 $\mathfrak{M}_{\XX,G_\Sigma}(T)\to \mathfrak{L}_{\XX,G_\Sigma}(T)$. After stackification, we get the map $\big[\rho_{in}/G_{\Sigma}\big]\:\big[\bfhom(\Sigma_{g,n,m},\XX)/G_{\Sigma}\big] \longrightarrow \big[(\LXX)^{n}/G_{\Sigma}\big]$. Further, the diagram pictured in the Lemma follows from the same diagram of fibered groupoids. In order to check that this digram is 2-cartesian, we note that both vertical arrows are $G_\Sigma$-torsors (by Lemma~\ref{L:quotientbyGSigma}) and the lemma follows from the usual interpretation of transformation fibered groupoid as groupoids of torsors recalled in Section~\ref{S:topological} (also see~\cite{BX, Noohi}). 
\end{proof}
\begin{rmk}
Lemma~\ref{L:quotientbyGSigma} and Lemma \ref{L:pullbackquotientbyGSigma} (as well as the constructions underlined there) basically follows because we are considering strict actions of a topological group on a stack (induced by the mapping stack construction). These statements are actually particular cases of more general statements about quotient of topological stacks by (topological) groups will be studied by the second and third author in a work in progress.
\end{rmk}

By Lemma~\ref{L:pullbackquotientbyGSigma}, any bivariant class in $H^\bullet(\big[\bfhom(\Sigma_{g,n,m},\XX)/G_{\Sigma}\big] \stackrel{\big[\rho_{in}/G_{\Sigma}\big]}\longrightarrow \big[(\LXX)^{n}/G_{\Sigma}\big])$ can be pulled-back along  $p_{/{G_\Sigma}}$.
\begin{lem}\label{L:Sigmaclass}
 There is a bivariant class $$\sigma^\XX_{g,n,m} \in H^{-d\chi(\Sigma_{g,n,m})}\Big( \big[\bfhom(\Sigma_{g,n,m},\XX)/G_{\Sigma}\big] 
\stackrel{\big[\rho_{in}/G_{\Sigma}\big]}\longrightarrow \big[(\LXX)^{n}/G_{\Sigma}\big]\Big)$$ 
whose pullback $p_{/{G_\Sigma}}^*\big(\sigma^\XX_{g,n,m}\big) $ is the class ${r_{g,n,m}}_*\big( \theta^\XX_{g,n,m}\big )$ (see~\eqref{eq:rorientationclass}).
\end{lem}
\begin{proof}
 An element $g\in G_\Sigma$ acts on $\Sigma_{g,n,m}$ and thus on $\iota_{g,n,m}(c_{g,n,m})$. The image 
$\sigma \cdot \iota_{g,n,m}(c_{g,n,m})$ is  a chord diagram diffeomorphic to $c_{g,n,m}$, by a diffeomorphism fixing 
the boundary circles of $c_{g,n,m}$. Further $\Sigma_{g,n,m}$ also retracts on 
$\sigma \cdot \iota_{g,n,m}(c_{g,n,m})$ and, indeed, $x\mapsto g\cdot \iota_{g,n,m}(r_{g,n,m}(g^{-1}\cdot x))$ is a 
retraction of $\Sigma_{g,n,m}$ on $c_{g,n,m}$.

A proof similar to the one of Lemma~\ref{L:quotientbyGSigma} shows that  we can form the quotient stack $\big[\bfhom{c_{g,n,m},\XX}/G_{\Sigma}\big]$ induced by  the above action on $c_{g,n,m}$. 
In particular the action of $G_{\Sigma}$ factors through an action of the orientation and circles preserving 
diffeomorphism $\mathop{Diff}^+_{\partial}(c_{g,n,m})$ group of $c_{g,n,m}$. Note that the $n$ disjoint circles are pointwise fixed
by  an element of $\mathop{Diff}^+_{\partial}(c_{g,n,m})$. It follows that this group actually acts on the disjoint 
union of the trees $\coprod_{\mathcal{T}(c_{g,n,m}} t$ (and preserving the circular vertices).

Thus, applying Lemma~\ref{L:pullbackquotientbyGSigma} and its proof we get the diagram of pullback squares
\begin{equation}\label{diag:quotientcSigma}
 \xymatrix{\bfhom(c_{g,n,m},\XX) \ar[r]^{r_{g,n,m}} \ar[d]& \bfhom(\Sigma_{g,n,m},\XX) \ar[rr]^{\hspace{2pc}\rho_{in}} \ar[d] 
&& \big(\LXX\big)^n \ar[d]^{p_{/{G_\Sigma}}}\\
 \big[\bfhom(c_{g,n,m},\XX)/G_{\Sigma}\big] \ar[r]^{[r/G_\Sigma]}& \big[\bfhom(\Sigma_{g,n,m},\XX)/G_{\Sigma}\big] \ar[rr]^{\hspace{2pc}\big[\rho_{in}/G_{\Sigma}\big]} 
&&  \big[(\LXX)^{n}/G_{\Sigma}\big].}
\end{equation}
For simplicity, we also denote $\rho_{in}$ the composition $\rho_{in}\circ r_{g,n,m}$. This pullback square will allow us to reduce the 
statement of the lemma to an analogue statement with $c_{g,n,m}$ instead of $\Sigma_{g,n,m}$. Indeed, assume we have a class 
$$\tilde{\sigma}^{\XX}_{g,n,m} \in H^{-d\chi(\Sigma_{g,n,m})}\Big( \big[\bfhom(c_{g,n,m},\XX)/G_{\Sigma}\big] 
\stackrel{\big[\rho_{in}/G_{\Sigma}\big]}\longrightarrow \big[(\LXX)^{n}/G_{\Sigma}\big]\Big)$$
such that $p_{/{G_\Sigma}}^*\big(\tilde{\sigma}^\XX_{g,n,m}\big) =\theta^{\XX}_{g,n,m}$. 
Then, since pushfoward and pullbacks commutes in a bivariant theory (see Axiom~A.13 in Appendix~\ref{Axioms}), we get that
$$p_{/{G_\Sigma}}^*\Big([r/G_\Sigma]_*\big(\tilde{\sigma}^\XX_{g,n,m}\big)\Big)={r_{g,n,m}}_*\big( \theta^\XX_{g,n,m}\big ). $$
Thus to finish the proof of the Lemma, it suffices to define the class $\tilde{\sigma}^{\XX}_{g,n,m}$ satisfying 
 $p_{/{G_\Sigma}}^*\big(\tilde{\sigma}^\XX_{g,n,m}\big) =\theta^{\XX}_{g,n,m}$ (where $\theta^{\XX}_{g,n,m}$ is the class defined by the the identity~\eqref{eq:thetachord}). To do so, we follow the proof of 
lemma~\ref{L:orienttree} to get an orientation class $\theta_{d_{/G_{\Sigma}}} \in H^{-d\chi(\Sigma_{g,n,m})}\big(\prod_{t\in \mathcal{T}(c_{g,n,m})}\big[\bfhom(t, \XX)/G_{\Sigma}\big]\to \big[\XX^{\mathcal{V}(c_{g,n,m})}/G_\Sigma\big]\big) $.
Each tree is a deformation retract onto any of its circular vertex, hence we have,  for each $t\in \mathcal{T}(c_{g,n,m})$,  a
 factorisation
$$\xymatrix{\big[\XX^{\mathcal{T}(c_{g,n,m})}/G_\Sigma\big] \ar[rr]^{\hspace{-2pc}\prod {r_t^*}_{/G_\Sigma}} 
\ar@/_1pc/[rrd]_{\hspace{-2pc}\Delta^{(\#\mathcal{V}(c_{g,n,m}) - \#\mathcal{T}(c_{g,n,m}) )}_{/G_\Sigma}}&&
\prod_{t\in \mathcal{T}(c_{g,n,m})}\big[\bfhom(t, \XX)/G_\Sigma\big] \ar[d]^{{d_{c_{g,n,m}}}_{/G_\Sigma}}\\ && 
\big[\XX^{\mathcal{V}(c_{g,n,m})}/G_\Sigma\big] } $$
 where the ${r_t^*}_{/G_\Sigma}$ are retractions.
Note that, similarly to the proof of Lemma~\ref{L:quotientbyGSigma}.3., there are  topological stacks isomorphisms $\big[\XX^{\mathcal{V}(c_{g,n,m})}/G_\Sigma\big]\cong \XX^{\mathcal{V}(c_{g,n,m})} \times \big[ pt/G_\Sigma\big]$ and $\big[\XX^{\mathcal{T}(c_{g,n,m})}/G_\Sigma\big]\cong \XX^{\mathcal{T}(c_{g,n,m})}\times \big[pt/G_\Sigma\big] $ since the action on $G_\Sigma$ on the circular vertices is trivial. 
Further, under this isomorphism, the map ${d_{c_{g,n,m}}}_{/G_\Sigma}$ identifies with $d_{c_{g,n,m}}\times \id_{/G_\Sigma}$. 
Thus, by Proposition~\ref{P:bivariantKunneth}, an orientation class $\theta$ (as the one given by Lemma~\ref{L:orienttree}) in $ H^{-d\chi(\Sigma_{g,n,m})}(\XX^{\mathcal{T}(c_{g,n,m})}\to \XX^{\mathcal{V}(c_{g,n,m})})$ determines an orientation 
$$ \theta \otimes 1 \in H^{-d\chi(\Sigma_{g,n,m})}\Big(\XX^{\mathcal{T}(c_{g,n,m})}\times \big[pt/G_{\Sigma}\big]\stackrel{d_{c_{g,n,m}}\times \id}\longrightarrow \XX^{\mathcal{V}(c_{g,n,m})}\times \big[pt/G_{\Sigma}\big]\Big).$$
Further, it follows from the proof of Proposition~\ref{P:bivariantKunneth} that the pullback $p_{/G_\Sigma}^*(\theta \otimes 1)$ is the original orientation class $\theta$. 
Taking, as in the proof of Lemma~\ref{L:orienttree}, 
\begin{equation}\label{eq:thetadG}\theta_{d_{/G_{\Sigma}}}=\big(\prod {r_t^*}_{/G_\Sigma}\big)_*(\theta\otimes 1)\end{equation}  to be the    pushforward 
of the orientation class $\theta\otimes 1$, we see that
\begin{equation}\label{eq:thetapulled}
p_{/G_\Sigma}^*\big(\theta_{d_{/G_{\Sigma}}}\big) \; =\; \theta_{d_{c_{g,n,m}}},
\end{equation}
 that is, the pullback
$p_{/G_\Sigma}^*\big(\theta_{d_{/G_{\Sigma}}}\big)$ is the class $\theta_{d_{c_{g,n,m}}}$ of Lemma~\ref{L:orienttree} (once again using that pullback and pushforward commute in a bivariant theory).

\medskip

Now, let us consider the following commutative diagram
{l \begin{equation}\label{diag:pullbackclass}
 \xymatrix@C=10pt@R=10pt@M=4pt{ & \bfhom(c_{g,n,m},\XX) \ar@{-}[d]\ar[rr]^{\hspace{2pc}\rho_{in}} \ar[dl] && \big(\LXX\big)^n \ar[dl]^{p_{/{G_\Sigma}}} 
\ar[dd]^{\ev_{\mathcal{V}}}  \\
\big[\bfhom(c_{g,n,m},\XX)/G_{\Sigma}\big] \ar[dd]\ar[rr]^{\hspace{8pc}\big[\rho_{in}/G_{\Sigma}\big]} &\ar[d]&  \big[(\LXX)^{n}/G_{\Sigma}\big] \ar[dd]\ar@{}[d]^{{\ev_{\mathcal{V}}}_{/G_\Sigma}} & \\
& \prod\limits_{t\in \mathcal{T}(c_{g,n,m})}\hspace{-1pc}\bfhom(t, \XX) \ar@{-}[r]^{ \hspace{4pc}d_{c_{g,n,m}}} \ar[dl]& \ar[r]& 
\XX^{\mathcal{V}(c_{g,n,m})}\ar[dl]^{p_{/G_\Sigma}} \\
\prod\limits_{t\in \mathcal{T}(c_{g,n,m})}\hspace{-1pc}\big[\bfhom(t, \XX)/G_{\Sigma}\big] \ar[rr]^{\hspace{3pc} {d_{c_{g,n,m}}}_{/G_\Sigma}} & & 
\big[\XX^{\mathcal{V}(c_{g,n,m})}/G_\Sigma\big] &  }
\end{equation}}
The bottom vertical and lower horizontal square are pullbacks and so is the top horizontal square. It follows that the front vertical square is a pullback too. We denote $$\tilde{\sigma}^{\XX}_{g,n,m}={\ev_{\mathcal{V}}}_{/G_\Sigma}^*(\theta_{d_{G_\Sigma}})\in H^{-d\chi(\Sigma_{g,n,m})}\Big(\big[\bfhom(c_{g,n,m},\XX)/G_{\Sigma}\big]\stackrel{\big[\rho_{in}/G_{\Sigma}\big]}\longrightarrow \big[(\LXX)^{n}/G_{\Sigma}\big]\Big)$$ the pullback of the class~\eqref{eq:thetadG} 
$$\theta_{d_{G_\Sigma}}\in H^{-d\chi(\Sigma_{g,n,m})}\Big( \prod_{t\in \mathcal{T}(c_{g,n,m})}\big[\bfhom(t, \XX)/G_{\Sigma}\big]\to \big[\XX^{\mathcal{V}(c_{g,n,m})}/G_\Sigma\big]\Big)$$ defined above.
By definition of $\theta^{\XX}_{g,n,m}$ (see identity~\eqref{eq:thetachord}), one has the identities
\begin{eqnarray*}
\theta^{\XX}_{g,n,m} &=& \ev_{\mathcal{V}}^*(\theta_{d_{c_{g,n,m}}}) \\
&=& \ev_{\mathcal{V}}^*\big( p^*_{/G_\Sigma}(\theta_{d_{/G_\Sigma}})\big) \mbox{ by relation~\eqref{eq:thetapulled} }\\
&=& p^*_{/G_\Sigma}\big({\ev_{\mathcal{V}}}_{/G_\Sigma}^*(\theta_{d_{/G_\Sigma}})  \big) \mbox{ by diagram~\eqref{diag:pullbackclass}} \\
&=& p^*_{/G_\Sigma}\big(\tilde{\sigma}^{\XX}_{g,n,m}\big)
\end{eqnarray*}
which prove that the class $\tilde{\sigma}^{\XX}_{g,n,m}$ satisfies the expected identity and thus finishes the proof of the lemma.
\end{proof}

\medskip

\paragraph{Defining the homology conformal field theory and Proof of Theorem~\ref{T:HCFT}}
The inclusion of the outgoing boundary components $\coprod_{i=1}^m S^1 \hookrightarrow \Sigma_{g,n,m}$ is also 
$G_{\Sigma}$-equivariant. Thus, similarly to the ingoing boundary case,  we get a stack morphism 
$\big[\rho_{out}/G_{\Sigma}\big]\:\big[\bfhom(\Sigma_{g,n,m},\XX)/G_{\Sigma}\big] \longrightarrow \big[(\LXX)^{m}/G_{\Sigma}\big]$.
Since $G_\Sigma$ acts trivially on $\big(\LXX\big)^m$, the terminal map $G_{\Sigma}\to \{1\}$ yields a stack 
morphism $[(\LXX)^{m}/G_{\Sigma}\big] \to [(\LXX)^{m}/\{1\}\big]\cong \big(\LXX\big)^m$.
 By Section~\ref{S:Gysin}, the bivariant class $\sigma^\XX_{g,n,m}$ given by Lemma~\ref{L:Sigmaclass}, yields a Gysin 
map $(\sigma^\XX_{g,n,m})^!\: H\big(\big[\big(\LXX\big)^n/G_{\Sigma}\big]\big) \to H\big(\big[\bfhom(\Sigma_{g,n,m},\XX)/G_{\Sigma}\big] \big)$.
Composing with the homology pushforward of the two preceeding stacks morphisms, we get, for every (isomorphism class of a) 
surface $\Sigma_{g,n,m}$, the following linear map:
\begin{multline}\label{eq:HCFTop}
\mu_{\Sigma_{g,n,m}}\: H\big(\big[\big(\LXX\big)^n/G_\Sigma\big]\big) \stackrel{(\sigma^\XX_{g,n,m})^!}\longrightarrow 
H\big(\big[\bfhom(\Sigma_{g,n,m},\XX)/G_{\Sigma}\big] \\
\stackrel{\big[\rho_{out}/G_{\Sigma}\big]_*}\longrightarrow 
H\big(\big[(\LXX)^{m}/G_{\Sigma}\big]\big) \to H\big(\big(\LXX\big)^{m}\big)
\end{multline}
This map indeed defines the (non-unital, non-counital) homological conformal field theory structure of $H\lcom(\LXX)$ 
asserted by Theorem~\ref{T:HCFT} as we will now prove (recall from Lemma~\ref{L:quotientbyGSigma} that $\big[\big(\LXX\big)^n/G_\Sigma\big]\cong \big(\LXX\big)^n \times \big[*/G_\Sigma \big]$).

\bigskip

{\noindent{\sc Proof of theorem~\ref{T:HCFT}.}}
We wish to define the $d$-dimensional homological conformal field theory structure by assiging to any positive integer $n$, 
the graded space $H\lcom\big(\big(\LXX\big)^n\big)$. 
First note that, by Lemmas~\ref{L:orienttree} and~\ref{L:Sigmaclass}, the map~\eqref{eq:HCFTop} above 
 $\mu_{\Sigma_{g,n,m}}\: H_{\bullet}\big(\big[\big(\LXX\big)^n/G_\Sigma\big]\big)\to 
H_{\bullet+d\chi(\Sigma_{g,n,m})}\big(\big(\LXX\big)^{m}\big)\big)$ is of degree $-d\chi(\Sigma_{g,n,m})$. 
Further, since $G_{\Sigma}=\mathop{Diff}^+_{n,m}(\Sigma_{g,n,m})$ acts trivially the incoming boundary of $\Sigma_{g,n,m}$ 
and thus on $(\LXX)^n$, there is (see Lemma~\ref{L:quotientbyGSigma}) a canonical isomorphism of topological stacks 
$\big[\big(\LXX\big)^n/G_\Sigma\big]\cong \big(\LXX\big)^n \times \big[*/G_{\Sigma}\big]$ and thus natural isomorphisms
$$ H_{\bullet}\big(\big[\big(\LXX\big)^n/G_\Sigma\big]\big) \cong H_{\bullet}\big(\big(\LXX\big)^n\big)\otimes H_\bullet(BG_{\Sigma})
\cong H_{\bullet}\big(\big(\LXX\big)^n\big)\otimes H_\bullet(B\Gamma_{n,m}(\Sigma)).$$
It follows that the maps $\mu_{\Sigma_{g,n,m}}$ induces, for any $\Sigma_{g,n,m}\in \MM_{n,m}$, a well defined map from 
$H\lcom\big(\big(\LXX\big)^n\big)$ to $H\lcom\big(\big(\LXX\big)^m\big)$, which we still denote $\mu_{\Sigma_{g,n,m}}$. This map has the correct degree shifting with respect to
the twisting local coefficient system $\det^{\otimes d}$. Further, since we are only considering closed boundaries,
the bundle corresponding to this local system is oriented by~\cite{Go} and locally trivial over the \emph{category} $\MM_{n,m}$ 
by~\cite{Cos, Cos2}. It follows that in order to check that the rule $n\mapsto H\lcom\big(\big(\LXX\big)^n\big)$ together with the 
maps $\mu_{\Sigma_{g,n,m}}$ defines a symmetric monoidal functor from $\mathcal{C}^{nu,nc}_{\MM, \det^{\otimes d}}$ to the category
 of graded vector spaces (\emph{i.e.} a non-unital non-counital $d$-dimensional homological conformal field theory), it suffices to
check the behavior of the maps  $\mu_{\Sigma_{g,n,m}}$ with respect to disjoint union and gluing of surfaces. This will be done below similarly
to the proof of associativity and coassociativity of the loop product in Theorems~\ref{th:Loop} and~\ref{frobeniusloop} as well as in 
Chataur-Menichi~\cite{CM}.  

We first deal with the gluing of surfaces. Let $\Sigma_{g,\ell,m}\in \MM_{\ell,n}$ and $\Sigma'_{g',n,m}\in \MM_{n,m}$ be two 
surfaces. We note $G_{\Sigma}=\mathop{Diff}^+_{\ell,n}(\Sigma_{g,\ell,n})$ and 
$G'_{\Sigma'}=\mathop{Diff}^+_{n,m}(\Sigma'_{g',n,m})$ the corresponding diffeomorphisms groups. Since these groups are fixing
the boundaries pointwise, it follows that we have an injective morphisms of topological groups 
$G_\Sigma \times G'_{\Sigma'}\hookrightarrow H_{\Sigma'\circ \Sigma}$ where 
$H_{\Sigma\circ \Sigma'}=\mathop{Diff}^+_{n,m}(\Sigma'_{g',n,m}\circ \Sigma_{g,\ell,n})$  is the group of oriented diffeomorphisms
fixing pointwise the boundary of the surface  $\Sigma'_{g',n,m}\circ \Sigma_{g,\ell,n}$ obtained by gluing $\Sigma_{g,\ell,n}$ and 
$\Sigma'_{g',n,m}$. For simplicity, henceforth, we denote
$\tilde{\Sigma}_{\circ} =\Sigma'_{g',n,m}\circ \Sigma_{g,\ell,n}$ this gluing. Since the boundary circles are fixed pointwise, both $G_\Sigma \times G'_{\Sigma'}$ and $H_{\Sigma'\circ \Sigma}$ acts trivially on $\big(\LXX\big)^{\ell}$ and $\big(\LXX\big)^{m}$ so that we have stacks morphisms $\big[\big(\LXX\big)^{\ell}/G_\Sigma \times G'_{\Sigma'}\big]\to \big(\LXX\big)^{\ell}$, $\big[\big(\LXX\big)^{\ell}/H_{\Sigma'\circ \Sigma}\big]\to \big(\LXX\big)^{\ell}$ and similarly with the outgoing circles (the proof that these quotient stacks are well defined and that we have these equivariant maps is the same as the ones of Lemma~\ref{L:quotientbyGSigma} and Lemma\ref{L:pullbackquotientbyGSigma}).   The morphism 
$G_\Sigma \times G'_{\Sigma'}\hookrightarrow H_{\Sigma'\circ \Sigma}$ induces a commutative diagram of stacks
$$\xymatrix{\big[pt/G_\Sigma\big] \times   \big[pt/G'_{\Sigma'}\big] \ar[rr]^{{\circ}} \ar@{^{(}->}[d] &&  
\big[pt/H_{\Sigma\circ \Sigma'}\big] \ar@{^{(}->}[d]\\ \MM_{\ell,n} \times \MM_{n,m} \ar[rr]^{\circ}&& \MM_{\ell,m}} $$
 and a similar diagram after passing to homology with twisted coefficient. 
We thus have to prove that $\mu_{\Sigma'_{g',n,m}} \circ \mu_{\Sigma_{g,\ell,n}}=\mu_{\tilde{\Sigma}_{\circ}}$. The above group morphism 
$G_\Sigma \times G'_{\Sigma'}\hookrightarrow H_{\Sigma'\circ \Sigma}$ induces an action of $G_\Sigma \times G'_{\Sigma'}$ on $\bfhom(\tilde{\Sigma}_{\circ},\XX)$ and a diagram of cartesian squares (applying, mutatis mutandis, the proof of Lemma~\ref{L:pullbackquotientbyGSigma}) of topological stacks
\begin{equation} \label{eq:HCFTpull1}
\xymatrix{\bfhom(\tilde{\Sigma}_{\circ},\XX) \ar[d]\ar[rr]^{\rho_{in}} && \big(\LXX\big)^{\ell} \ar[d] \\
\big[\bfhom(\tilde{\Sigma}_{\circ},\XX)/G_\Sigma \times G'_{\Sigma'}\big]\ar[rr]^{\hspace{1pc}\big[\rho_{in}/G_\Sigma \times G'_{\Sigma'}\big]} \ar[d] &&  \big[\big(\LXX\big)^{\ell}/G_\Sigma \times G'_{\Sigma'}\big]  \ar[d]\\
 \big[\bfhom(\tilde{\Sigma}_{\circ},\XX)/H_{\Sigma'\circ \Sigma}\big]\ar[rr]^{\big[\rho_{in}/H_{\Sigma'\circ \Sigma}\big]} && 
\big[\big(\LXX\big)^{\ell}/H_{\Sigma'\circ \Sigma}\big]\big]}
\end{equation}
Furthermore, the group $G_\Sigma \times G'_{\Sigma'}$ also acts on $\bfhom(\Sigma_{g,\ell,n},\XX)$ (through the above action of $G_\sigma$ and the trivial action of $G'_{\Sigma'}$) 
and similarly on $\bfhom(\Sigma'_{g',n,m},\XX)$. Since $\tilde{\Sigma}_{\circ}=\Sigma'_{g',n,m}\circ \Sigma_{g,\ell,n}$ is obtained by glueing $\Sigma_{g,\ell,n}$ on $\Sigma'_{g',n,m}$
 along the $n$-disjoint circles of their common boundaries, applying Lemma~\ref{L:glue} and (the proof of) Lemma~\ref{L:pullbackquotientbyGSigma}, we get another cartesian square
\begin{equation}\label{eq:HCFTpull2}
\xymatrix{\big[\bfhom(\tilde{\Sigma}_{\circ},\XX)/G_\Sigma \times G'_{\Sigma'}\big]\ar[rr]^{\big[\mathop{res}_{in}/G_\Sigma \times G'_{\Sigma'}\big]} \ar[d]_{\tilde{\mathop{res}}_{out}} 
&&  \big[\bfhom(\Sigma_{g,\ell,n},\XX)/G_\Sigma \times G'_{\Sigma'}\big] \ar[d] \\
  \big[\bfhom(\Sigma'_{g',n,m},\XX)/G_\Sigma \times G'_{\Sigma'}\big]\ar[rr]^{\hspace{1pc}\big[\mathop{\rho}_{in}/G_\Sigma \times G'_{\Sigma'}\big]} 
&&  \big[\big(\LXX\big)^n/G_\Sigma \times G'_{\Sigma'}\big]}
\end{equation}  
of stacks (also see~\cite{CM}), where $\mathop{res}_{in}$ is the restriction map induced by the inclusion $\Sigma_{g,\ell,n}\to \tilde{\Sigma}_{\circ}$.

\smallskip

The inclusion of the outgoing boundary of $\Sigma_{g,\ell,n}$ induces a stack morphism $$\tilde{\rho}_{mid}\: \big[\bfhom(\Sigma_{g,\ell,n},\XX)/G_\Sigma \times G'_{\Sigma'}\big]\to \big[\big(\LXX\big)^{n}/G_\Sigma \times G'_{\Sigma'}\big] \to  \big[\big(\LXX\big)^{n}/G'_{\Sigma'}\big]$$ since the group $G_\Sigma $ acts trivially on $\big(\LXX\big)^n$. Similarly we have a topological stacks morphisms
$$\tilde{\rho}_{out}\: \big[\bfhom(\Sigma'_{g',n,m},\XX)/ G'_{\Sigma'}\big] \to \big[\big(\LXX\big)^m/G'_{\Sigma'}\big]\to \big(\LXX\big)^m,$$
$$ \tilde{\tilde{\rho}}_{out}\: \big[\bfhom(\tilde{\Sigma}_{\circ},\XX)/ H_{\Sigma'\circ \Sigma}\big] \to \big[\big(\LXX\big)^m/H_{\Sigma'\circ \Sigma}\big]\to \big(\LXX\big)^m$$
and  a composition of morphisms of topological stacks
\begin{multline*}\tilde{\mathop{res}}_{out}\: \big[\bfhom(\tilde{\Sigma}_{\circ},\XX)/G_\Sigma \times G'_{\Sigma'}\big] \longrightarrow  \big[\bfhom(\Sigma'_{g',n,m},\XX)/G_\Sigma \times G'_{\Sigma'}\big]\\  \longrightarrow  \big[\bfhom(\Sigma'_{g',n,m},\XX)/G'_{\Sigma'}\big].\end{multline*}

\smallskip

By Lemma~\ref{L:Sigmaclass} applied to the surfaces $\Sigma_{g,\ell,n}$, $\Sigma'_{g',n,m}$ and $\tilde{\Sigma}_{\circ}$, there are bivariant classes  $\sigma^\XX_{g,\ell,n}$, ${\sigma'}^\XX_{g',n,m}$ and $\tilde{\sigma}^\XX_{\circ}$ inducing (as for the definition of $\mu_{\Sigma_{g,n,m}}$), respectively, the Gysin maps 
\begin{equation}\label{eq:GysHCFT1} (\sigma^\XX_{g,\ell,n})^!\: H\big(\big[\big(\LXX\big)^{\ell}/G_{\Sigma}\big]\big) \to H\big(\big[\bfhom(\Sigma_{g,\ell,n},\XX)/G_{\Sigma}\big] \big), \end{equation}
\begin{equation}\label{eq:GysHCFT2} ({\sigma'}^\XX_{g',n,m})^!\: H\big(\big[\big(\LXX\big)^{n}/G'_{\Sigma'}\big]\big) \to H\big(\big[\bfhom(\Sigma_{g',n,m},\XX)/G'_{\Sigma'}\big] \big),\end{equation}
 \begin{equation}\label{eq:GysHCFT3}(\tilde{\sigma}^\XX_{\circ})^!\: H\big(\big[\big(\LXX\big)^{\ell}/H_{\Sigma'\circ \Sigma}\big]\big) \to H\big(\big[\bfhom(\tilde{\Sigma}_{\circ},\XX)/H_{\Sigma'\circ\Sigma}\big] \big).
 \end{equation}
 The first  Gysin map above is equivariant with respect to the action of $G'_{\Sigma'}$ since this group acts trivially on $\Sigma_{g,\ell,n}$, thus passes to quotient stack to define a Gysin map
\begin{equation}\label{eq:GysHCFT4}(\sigma^\XX_{g,\ell,n})^!\: H\big(\big[\big(\LXX\big)^{\ell}/G_{\Sigma}\times G'_{\Sigma'}\big]\big) \to H\big(\big[\bfhom(\Sigma_{g,\ell,n},\XX)/G_{\Sigma}\times G'_{\Sigma'}\big] \big). \end{equation}
Further, the proof of Lemma~\ref{L:Sigmaclass} applied to $\Sigma_{g,\ell,n}$ with respect to the action of the group $G_{\Sigma}\times G'_{\Sigma'}$ (and not the full $H_{\Sigma'\circ \Sigma}$) also yields a bivariant class $\tilde{\sigma}^\XX_{G_{\Sigma}\times G'_{\Sigma'}}$ and an associated Gysin map
\begin{equation} \label{eq:GysHCFT5}(\tilde{\sigma}^\XX_{G_{\Sigma}\times G'_{\Sigma'}})^!\: H\big(\big[\big(\LXX\big)^{\ell}/G_{\Sigma}\times G'_{\Sigma'}\big]\big) \to H\big(\big[\bfhom(\tilde{\Sigma}_{\circ},\XX)/G_{\Sigma}\times G'_{\Sigma'}\big] \big).\end{equation}

\smallskip

The Gysin maps $(\tilde{\sigma}^\XX_{G_{\Sigma}\times G'_{\Sigma'}})^!$ (morphism~\eqref{eq:GysHCFT5}) and $(\sigma^\XX_{g,\ell,n})^!$ (morphism~\eqref{eq:GysHCFT4}) are related as follows. We choose $c_{g,\ell,n}$ to be a chord diagram associated to $\Sigma_{g,\ell,n}$, $c'_{g',n,m}$ to be 
associated to $\Sigma_{g',n,m}$ and $\tilde{c}_{g+g',\ell,m}=c'_{g',n,m} \circ c_{g,\ell,n}$ the one associated to $\tilde{\Sigma}_{\circ}$ obtained by gluing the two previous ones (see~\cite{CoGo} for the composition of chord diagrams). By construction (\emph{i.e.}, use of Lemma~\ref{L:Sigmaclass}), the bivariant class
$\sigma^\XX_{g,\ell,n}$ is obtained as a class whose pullback along the stack morphism $\big(\LXX\big)^\ell \to \big[\big(\LXX\big)^\ell/G_{\Sigma}\times G'_{\Sigma'}\big]$ is $r_{g,\ell,n,*}\big(\theta^\XX_{g,\ell,n}\big)$ where $\theta^\XX_{g,\ell,n}$ is given as in formula~\eqref{eq:thetachord} and Lemma~\ref{L:orienttree}. Similarly, the bivariant class
$\tilde{\sigma}^\XX_{G_{\Sigma}\times G'_{\Sigma'}}$ is obtained as a class whose pullback along the stack morphism $\big(\LXX\big)^\ell \to \big[\big(\LXX\big)^\ell/G_{\Sigma}\times G'_{\Sigma'}\big]$ is $r_{g+g',\ell,m,*}\big(\theta^\XX_{g+g',\ell,m}\big)$, where, again the class $\theta^\XX_{g+g',\ell,m}$  is given by formula~\eqref{eq:thetachord} and Lemma~\ref{L:orienttree}. 

\smallskip

We now relate the classes $\theta^\XX_{g+g',\ell,m}$ and $\theta^\XX_{g,\ell,n}$.
The chord diagram
$\tilde{c}_{g+g',\ell,m}$ has $\ell$-disjoint circles, a set of disjoint trees $\mathcal{T}(c_{g,\ell,n})$ (corresponding to to the trees of $c_{g,\ell,n}$) and an additional  set of disjoint trees $\mathcal{T}'(c'_{g',n,m})$ such that the union $\tilde{\mathcal{T}}(\tilde{c}_{g+g',\ell,m})=\mathcal{T}(c_{g,\ell,n})\coprod \mathcal{T}'(c'_{g',n,m})$ is the set of trees associated to $\tilde{c}_{g+g',\ell,m}$. We define similarly the set of circular vertices $\mathcal{V}(c_{g,\ell,n})$, $\mathcal{V}'(c'_{g',n,m}) $ and $ \tilde{\mathcal{V}}(\tilde{c}_{g+g',\ell,m})$.  We thus get a diagram of pullback squares (obtained as for the square~\eqref{eq:chordpullback}):
\begin{equation}\label{eq:chordpullback2} 
\xymatrix{\bfhom(\tilde{c}_{g+g',\ell,m},\XX) \ar[r]^{\mathop{res}_{in}^c}\ar[d] &\bfhom(c_{g,n,m},\XX)\ar[d]^{\ev_{\mathcal{V}'}} \ar[r]& (\LXX)^{n} \ar[d]^{\ev_{\tilde{\mathcal{V}}}}\\
\prod\limits_{t\in \tilde{\mathcal{T}}(\tilde{c}_{g+g',\ell,m})}\hspace{-1pc}\bfhom(t, \XX) \ar[r] &\XX^{\mathcal{V}'(c'_{g',n,m})}\times \prod\limits_{t\in \mathcal{T}(c_{g,n,m})}\hspace{-1pc}\bfhom(t, \XX) \ar[r]&
 \XX^{\tilde{\mathcal{V}}(\tilde{c}_{g+g',\ell,m})} } \end{equation}
The left cartesian square above~\eqref{eq:chordpullback2} and the proof of Lemma~\ref{L:orienttree} gives us a bivariant class 
$${\theta}^\XX_{\mathop{res}^c_{in}}=\ev_{ \mathcal{V}'}^*({\theta}^\XX_{d_{\mathop{res}^c_{in}}}) \in H\Big(\bfhom(\tilde{c}_{g+g',\ell,m},\XX) \stackrel{\mathop{res}^c_{in}}\longrightarrow \bfhom(c_{g,n,m},\XX)\Big) $$ 
Since the classes $\theta^\XX_{g+g',\ell,m}$ and $\theta^\XX_{g,\ell,n}$ are also  induced by pullbacks along $\ev_{\tilde{V}}$ of classes given by Lemma~\ref{L:orienttree} (corresponding to the bottom line of~\eqref{eq:chordpullback2} and the various chord diagrams involved here), it follows from the functoriality of pullbacks (Axiom A.3 in Appendix~\ref{Axioms}) that $\theta^\XX_{g+g',\ell,m}= {\theta}^\XX_{\mathop{res}^c_{in}}\cdot \theta^\XX_{g,\ell,n}$. 

Since the following diagram (induced by the various restriction maps)
$$\xymatrix{ \bfhom(\tilde{\Sigma}_{\circ},\XX) \ar[rr]^{\mathop{res}_{in}}  && \bfhom(\Sigma_{g,\ell,n},\XX) \\
\bfhom(\tilde{c}_{g+g',\ell,m},\XX)\ar[u]^{r_{g+g',\ell,m}} \ar[rr]^{\mathop{res}_{in}^c} && \bfhom(c_{g,\ell,n},\XX)\ar[u]_{r_{g,\ell,n}} } $$
is commutative, it follows that $$r_{g+g',\ell,m,*}\big(\theta^\XX_{g+g',\ell,m}\big)= r_{g+g',\ell,m,*}\big({\theta}^\XX_{\mathop{res}^c_{in}} \big)\cdot r_{g,\ell,n,*}\big(\theta^\XX_{g,\ell,n}\big).$$

Unfolding the argument of the proof of Lemma~\ref{L:Sigmaclass}, we get that the bivariant class $r_{g+g',\ell,m,*}\big({\theta}^\XX_{\mathop{res}^c_{in}} \big) $
is the pullback of a class $$\sigma^{\XX}_{\mathop{res}_{in}}\in H\Big (\big[\bfhom(\Sigma_{g,\ell,n},\XX)/G_\Sigma \times G'_{\Sigma'}\big] \stackrel{\big[\mathop{res}_{in}/G_{\Sigma}\times G'_{\Sigma'}\big]}\longrightarrow  \big[\bfhom(\tilde{\Sigma}_{\circ},\XX)/G_\Sigma \times G'_{\Sigma'}\big]\Big) $$
and further that the  Gysin maps induced by the classes $\tilde{\sigma}^\XX_{G_{\Sigma}\times G'_{\Sigma'}}$, $ \sigma^\XX_{g,\ell,n}$ and $\sigma^{\XX}_{\mathop{res}_{in}}$ satisfy the identity
\begin{eqnarray}\label{eq:rhoinGysin}
 (\tilde{\sigma}^\XX_{G_{\Sigma}\times G'_{\Sigma'}})^! &= & (\sigma^{\XX}_{\mathop{res}_{in}})^! \circ (\sigma^\XX_{g,\ell,n})^!.
\end{eqnarray}

\medskip

Now, using the above defined maps, we can consider the following diagram
{\small\begin{equation}\label{eq:bigdiagramHCFT}
\xymatrix{H\big(\big[\big(\LXX\big)^{\ell}/G_\Sigma \times G'_{\Sigma'}\big]\big) \ar[r]^{\hspace{-1pc}\!\!(\sigma^\XX_{g,\ell,n})^!}\ar[rd]_{\hspace{-1pc}\!\!(\tilde{\sigma}^\XX_{G_{\Sigma}\times G'_{\Sigma'}})^!}\ar[dd]  & H\big(\big[\bfhom(\Sigma_{g,\ell,n},\XX)/G_\Sigma \times G'_{\Sigma'}\big]\big) \ar[d]^{(\sigma^{\XX}_{\mathop{res}_{in}})^!}\ar[r]^{\hspace{2pc} \tilde{\rho}_{mid *}}  & H\big(\big[\big(\LXX\big)^n/G'_{\Sigma'}\big]\big) \ar[d]^{(\sigma^\XX_{g',n,m})^!}\\
\ar @{} [dr] |{(1)}&   H\big(\big[\bfhom(\tilde{\Sigma}_{\circ},\XX)/G_\Sigma \times G'_{\Sigma'}\big] \big)\ar[r]_{{\tilde{\mathop{res}}_{out *} } }\ar[d]\ar @{} [ur] |{(2)} \ar @{} [dr] |{(3)}&
 H\big(\big[\bfhom(\Sigma'_{g',n,m},\XX)/G'_{\Sigma'}\big]\big) \ar[d]^{{\tilde{\rho}_{out *} } } \\
H\big(\big[\big(\LXX\big)^{\ell}/H_{\Sigma'\circ \Sigma} \big]\big) \ar[r]_{(\tilde{\sigma}^\XX_{\circ})^!}& 
H\big(\big[\bfhom(\tilde{\Sigma}_{\circ},\XX)/H_{\Sigma'\circ \Sigma} \big]\big)  \ar[r]_{{\tilde{\tilde{\rho}}_{out *} } } & H\big(\big(\LXX\big)^m \big)}
\end{equation}}
The left vertical map ad the bottom line represents the map $\mu_{\Sigma'_{g',n,m}\circ \Sigma_{g,\ell,n}}$ while the top line and right vertical composition represents the composition $\mu_{\Sigma'_{g',n,m}} \circ \mu_{\Sigma_{g,\ell,n}}$. Thus to prove the glueing property, it is enough to show that the diagram~\eqref{eq:bigdiagramHCFT} is commutative. The commutativity of the upper left triangle is precisely identity~\eqref{eq:rhoinGysin} above.
The commutativity of the trapezoid labelled~(1) follows from the tower of pullback squares~\eqref{eq:HCFTpull1} and naturality of Gysin maps with respect to towers of pullback squares, since all the classes involved are obtained by pullback from a class induced by Lemma~\ref{L:Sigmaclass}. The commutativity of the square labelled~(3) follows from the fact that the restriction to the outgoing boundary of 
$\tilde{\Sigma}_{\circ}$ coincides with the restriction to the outgoing boundary of  $\Sigma'_{g',n,m}$. Finally the square labelled~(2) is commutative thanks to the naturality of Gysin maps applied to the cartesian square~\eqref{eq:HCFTpull2}.

\medskip

We are left to the case of disjoint union of surfaces, that is to prove that
$\mu_{\Sigma_{g,n,m}\coprod \Sigma'_{g',n',m'}} = \pm \,\mu_{\Sigma_{g,n,m}}\otimes \mu_{\Sigma'_{g',n',m'}}$ where the sign is induced by the Koszul rule and the local coefficient $\det^{\otimes d}$. The sign  follows as in~\cite{CM, Go, Cos, Cos2}.
Let $c(g,n,m)$ and $c'(g',n',m')$ be chord diagrams representing respectively $\Sigma_{g,n,m}$ and $\Sigma'_{g',n',m'}$. Then the disjoint union $c(g,n,m) \coprod c'(g',n',m')$ is a chord diagram representing $\Sigma_{g,n,m}\coprod \Sigma'_{g',n',m'}$ and further, the diffeomorphism group $\mathop{Diff}^+_{n+n',m+m'}(\Sigma_{g,n,m}\coprod \Sigma'_{g',n',m'})$ is the cartesian product of $G_{\Sigma}=\mathop{Diff}^+_{n,m}(\Sigma_{g,n,m})$ and $G_{\Sigma'}=\mathop{Diff}^+_{n',m'}(\Sigma'_{g',n',m'}) $. Let  $\theta_{d_{c_{g,n,m}}}$ and $\theta_{d_{c'_{g',n',m'}}}$ be the respective strong orientation classes given by Lemma~\ref{L:orienttree} applied to $\Sigma_{g,n,m}$ and $\Sigma'_{g',n',m'}$.  For simplicity we simply denote $\Sigma=\Sigma_{g,n,m}$, $\Sigma'=\Sigma'_{g',n',m'}$, $c=c_{g,n,m}$ and $c'=c'_{g',n',m'}$. Since we are working over a field, by Proposition~\ref{P:bivariantKunneth}, we get that the tensor product 
\begin{multline*}\theta_{d_{c_{g,n,m}}}\otimes \theta_{d_{c'_{g',n',m'}}}\hspace{-1pc} \in H^{-d\chi(\Sigma\coprod \Sigma')}\Big(\xymatrix{\prod\limits_{ \mathcal{T}(c\coprod c')}\bfhom(t, \XX) \ar[rr]^{\,\, d_c\times d_{c'}}&& \XX^{\mathcal{V}(c\coprod c')}} \Big) \end{multline*}
identifies with the class $\theta_{d_{ c_{g,n,m}\coprod c'_{g',n',m'} } }$ 
given by Lemma~\ref{L:orienttree} applied to the surface $\Sigma_{g,n,m}\coprod \Sigma_{g',n',m'}$. This is a consequence of the fact that   this class is induced by taking the products of the diagrams inducing $\theta_{d_{c_{g,n,m}}}$ and $\theta_{d_{c'_{g',n',m'}}}$. It follows that the same property holds for the classes obtained by applying Lemma~\ref{L:Sigmaclass}, namely $\sigma^\XX_{g+g',n+n',m+m'} =\sigma^\XX_{g,n,m}\otimes \sigma^\XX_{g,n,m}$. Further the restriction to the incoming and outgoing  boundary yields a commutative diagram
$$\xymatrix{ \big(\LXX\big)^{n+n'} \ar[d]_{\cong}
&\bfhom(\Sigma_{g,n,m}\coprod \Sigma'_{g',n',m'},\XX) \ar[d]_{\cong} \ar[r]^{\hspace{2pc}\rho_{out}} \ar[l]_{\rho_{in}\hspace{1pc}} & \big(\LXX\big)^{m+m'} \ar[d]^{\cong} \\
\big(\LXX\big)^{n}\times \big(\LXX\big)^{n'} 
&\bfhom(\Sigma_{g,n,m},\XX)\times \bfhom(\Sigma'_{g',n',m'},\XX)  \ar[r]_{\hspace{3pc}\rho_{out}\times \rho_{out}} \ar[l]^{\rho_{in}\times \rho_{in} \hspace{4pc}} & 
\big(\LXX\big)^{m}\times \big(\LXX\big)^{m'}.} $$
It now follows similarly to the case of the glueing of surfaces (and actually more easily) that the operation $\mu_{\Sigma_{g,n,m}\coprod \Sigma'_{g',n',m'}}$ is $\mu_{\Sigma_{g,n,m}}\otimes \mu_{\Sigma'_{g',n',m'}}$.

\medskip

Now we only need to identify the operations of the \BV-structure and Frobenius structure with  the one given by the homological conformal field theory we just defined. By~\cite{Get, Go, CM},
 we know that the \BV-operator $H_\bullet(\LXX) \to H_{\bullet+1}(\LXX)$ is induced by the generator of degree 1 in the homology $H_\bullet(\MM_{1,1})$ corresponding to the diffeomorphism of a cylinder given by
 the Dehn twist along a generator of the degree 1 homology of the cylinder $\Sigma_{0,1,1}$. In other words, this generator is induced by the fundamental class of $S^1$, and  passing to the quotient stack
 $\bfhom(\mathop{\Sigma_{0,1,1},\XX})\to \big[\bfhom(\mathop{\Sigma_{0,1,1},\XX})/S^1\big]$,
 we see that the action of this generator on $\rho_{out}\big(\bfhom(\mathop{\Sigma_{0,1,1},\XX})\big)=\LXX$ coincides 
with the operator $D$ of Theorem~\ref{BV}. Now the product and the coproduct are respectively given by pair of pants
 (with different incoming and utgoing boundaries). The first one correspond to the chord diagram with two circles and one edge connecting them while the second one correspond to the chord diagram with one circle and one diameter. 
They are given by degree $0$ homology classes in $H_\bullet(\MM_{2,1})$ and $H_{\bullet}(\MM_{1,2})$, thus to identify them,
 it is enough to consider the Gysin maps obtained via Lemma~\ref{L:orienttree} before passing to the quotient by the
 diffeomorphim groups. Unfolding the proof of Lemma~\ref{L:orienttree}, we see that these Gysin maps coincides
 with the ones defining the loop product in Section~\ref{Construction} and loop coproduct in Theorem~\ref{frobeniusloop}.
{\nolinebreak $\Box$ }

\section{Remarks on brane topology for stacks}\label{S:Brane}

Brane topology, as coined by Sullivan and Voronov~\cite{CoVo}, is an higher dimensional analogue of string topology defined for free mapping sphere spaces instead of free loop space. Many aspects of Brane topology for manifolds have been studied in~\cite{CoVo, Ch, KaSa, HKV} for instance. Roughly, Brane topology is
concerned with the algebraic structure of the homology of $M^{S^n} =
\map (S^n,M)$, where $M$ is an oriented manifold and $S^n$ the
standard $n$-dimensional sphere. In this section we sketch how to apply our general machinery to define Brane topological operations for oriented stacks.

By Proposition~\ref{P:MapSt}, for any topological stack $\XX$, the mapping stack $\bfhom(S^n,\XX)$ is a topological stack. Here $S^n$ is the $n$-dimensional sphere. Let $\ev_{0}\:\bfhom(S^n,\XX)\to \XX$ be the map induced by the evaluation in $0$, the based point of $S^n$. There is a standard pinching map $p_{S^n}\: S^n\to S^n\vee S^n$ obtained by collapsing an  equator to the based point. This map is homotopy coassociative. By Lemma~\ref{L:glue}, if $\XX$ is Hurewicz, there is a cartesian square
$$\xymatrix{\bfhom(S^n\vee S^n,\XX) \ar[r] \ar[d] & \bfhom(S^n,\XX)\times \bfhom(S^n,\XX)\ar[d]^{\ev_{0}\times \ev_{0}}\\
\XX \ar[r]^{\Delta} & \XX\times \XX } $$
and thus, if $\XX$ is oriented of dimension $d$, a Gysin map $\Delta^!\: H_\bullet\big(\bfhom(S^n,\XX)\times \bfhom(S^n,\XX)  \big) \to H_{\bullet-d}\big(\bfhom(S^n\vee S^n,\XX)\big) $. Composing $\Delta^!$ with the map induced by the pinching mpa $p_{S^n}$, yields the \emph{brane product} 
\begin{multline}\label{eq:braneproduct}
\star_{S^n}\: H\big(\bfhom(S^n,\XX)\big)^{\otimes 2} \cong H\big(\bfhom(S^n,\XX)\times \bfhom(S^n,\XX) \big)\\
 \stackrel{\Delta^!}\longrightarrow H\big(\bfhom(S^n\vee S^n,\XX)\big) 
\stackrel{(p_{S^n})_*}\longrightarrow H\big(\bfhom(S^n,\XX)\big).
\end{multline}

\begin{prop}\label{P:brane}Let $\XX$ be an oriented stack of dimension $d$.
The brane product makes the shifted homology $\mathbb{H}_{\bullet}\big( \bfhom(S^n,\XX)\big)=H_{\bullet+d}\big( \bfhom(S^n,\XX)\big)$ a graded commutative algebra.
\end{prop}
\begin{proof}
The argument is the same as the ones of the proof of Theorem~\ref{th:Loop} and Proposition~\ref{Commutativityloop}.
\end{proof}
\begin{numrmk} The brane product~\eqref{eq:braneproduct} can be ``twisted'' by a class $\alpha\in \bigoplus_{n\geq 0} H^n(\bfhom(S^n\vee S^n,\XX))$ similarly to the loop product as in Section~\ref{Loopproduct}. The analogue for brane product of Theorem~\ref{Associativitycocycle} hold true, the proof being the same. For instance, the twisted brane product is associative if the class $\alpha$ satisfies the 2-cocycle condition~\eqref{eq:Associativitycocycle}.
\end{numrmk}

Let $\phi^k\:S^n \to S^n$ be the \emph{$k^{th}$-iterated power map} defined by the composition \begin{equation}\label{eq:powermapbrane} \phi^k\: S^n \;\stackrel{p^{(k)}}\longrightarrow; \bigvee_{i=1}^k S^n \;\stackrel{\bigvee \id }\longrightarrow \; S^n\end{equation} 
where $p^{(k)}\:S^n\to \bigvee S^n$ is the $k^{th}$-iterated pinching map and $\bigvee \id: \bigvee S^n\to S^n$ is the identity on each sphere of the bouquet $\bigvee S^n$.
The map $\phi^k$ induce by precomposition maps $\lambda^k\:\bfhom(S^n,\XX)\to \bfhom(S^n,\XX)$, $f\mapsto \lambda^k(f)=f\circ \phi(k) $. We have 
\begin{eqnarray}\label{eq:lambdaring}
 \lambda^k \circ \lambda^\ell & = & \lambda^{k\ell}.
\end{eqnarray}

\begin{them}\label{T:branepower} Let $\XX$ be an oriented (Hurewicz) stack of dimension $d$ and assume $n\geq 2$.
Then the maps $\lambda^k\: \sH\lcom( \bfhom(S^n,\XX))\to \sH\lcom( \bfhom(S^n,\XX))$ are maps of algebras (for the brane product $\star_{S^n}$). 

Further, if the ground ring $k$ contains $\mathbb{Q}$, there is a decomposition
$$\sH\lcom( \bfhom(S^n,\XX)) \cong \prod_{i\geq 0} \sH\lcom^{(i)}(\XX)$$, where  $\sH\lcom^{(i)}(\XX)$ is the eigenspace of $\lambda^k$ (with eigenvalue $k^{i}$), that makes  $\Big(\sH\lcom( \bfhom(S^n,\XX)),\star_{S^n} \Big)$
  a bigraded commutative algebra (with respect to the shifted total degree and the grading induced by the decomposition).
\end{them}
\begin{proof}
 By identity~\eqref{eq:lambdaring}, the maps $\lambda^k\:\bfhom(S^n,\XX)\to \bfhom(S^n,\XX)$ makes $\sH\lcom( \bfhom(S^n,\XX))$ equipped with the \emph{null} multiplication (and not the brane product $\star_{S^n}$) a $\lambda$-ring with trivial multiplication as in~\cite{Lo}. The existence of the decomposition then follows from standard properties of $\lambda$-ring, see~\cite{AT, Lo}. Thus in order to prove the Lemma, we are left to prove that the maps $\lambda^k$ are maps of algebras with respect to the brane product.
 This is an easy consequence of the commutativity of the following diagram
{\scriptsize \begin{equation}\label{eq:lambdabrane}
  \xymatrix{  H\big(\bfhom(S^n\coprod S^n,\XX)\big) \ar[r]^{\Delta^!}\ar[d]_{(\vee \id)\coprod (\vee \id)}  &H\big(\bfhom(S^n\vee S^n,\XX)\big) \ar[r]^{p^{(2)}}\ar[d]^{(\vee \id)\vee (\vee \id)} & \bfhom(S^n,\XX) \ar[d]^{\vee \id} \\
  H\big(\bfhom(\bigvee_{i=1}^k S^n\coprod \bigvee_{i=1}^k S^n,\XX)\big) \ar[r]^{\Delta_{\bigvee_{i=1}^k S^n}^!}\ar[d]_{p^{(k)}\coprod p^{(k)}}  & H\big(\bfhom\big((\bigvee_{i=1}^k S^n) \vee (\bigvee_{i=1}^k S^n),\XX\big)\big) \ar[r]^{\hspace{4pc}\bigvee p^{(2)}}\ar[d]^{p^{(k)}\vee p^{(k)}} & \bfhom(\bigvee S^n,\XX) \ar[d]^{p^{(k)}}\\ 
  H\big(\bfhom(S^n\coprod S^n,\XX)\big) \ar[r]^{\Delta^!}  &H\big(\bfhom(S^n\vee S^n,\XX)\big) \ar[r]^{p^{(2)}} & \bfhom(S^n,\XX)  }
 \end{equation} }
where the Gysin map $\Delta_{\bigvee_{i=1}^k S^n}^!$ is obtained as the pullback by $\ev_0\times \ev_0\:  \bfhom\big((\bigvee_{i=1}^k S^n) \coprod (\bigvee_{i=1}^k S^n),\XX\big)\to \XX\times \XX$ (the evaluation at the base points of each component) of the diagonal $\Delta\: \XX\to \XX\times \XX$. The map 
$\bigvee p^{(2)}\:\bfhom\big((\bigvee_{i=1}^k S^n) \vee (\bigvee_{i=1}^k S^n),\XX\big)\to \bigvee_{i=1}^k S^n$ by  applying a permutation on the bouquet of spheres (so that the first and $k+1$-sphere are put next to each other, and then the second sphere with the $k+2$-sphere and so on)
and then applying $p^{(2)}$ $k$-times. 
It follows immediately from this definition that the top right square of diagram~\eqref{eq:lambdabrane} is commutative. The left squares are seen to be commutative by applying the naturality of Gysin maps as in the proof of Theorem~\ref{th:Loop}. 
The maps $ p^{(k)}\circ \vee p^{(2)}:\bfhom\big( (\bigvee_{i=1}^k S^n) \vee (\bigvee_{i=1}^k S^n),\XX\big) \to \bfhom(S^n,\XX)$ and $p^{(2)}\circ (p^{(k)}\vee p^{(k)})$ involved in the lower right square of diagram~\eqref{eq:lambdabrane} are not equal. However, since $n\geq 2$, they are homotopic to each other (the proof being similar to the commutativity of the higher homotopy groups). The result follows. 
\end{proof}

\begin{numrmk}
In view of Section~\ref{ChasSullivan}, Theorem~\ref{T:branepower} holds for oriented manifolds, and did not seem to be written in the litterature to the authors knowledge. 
\end{numrmk}

\medskip

The brane product described above shall belong to a bigger algebraic structure. Namely, we believe that the following
\begin{claim} Let $\XX$ be an oriented Hurewicz stack of dimension $d$. Then
$\sH_{\scriptstyle \bullet}\big(\bfhom(S^n,\XX)\big)$ is an algebra over the homology $H_\bullet(\cac^{(n)})$ of the $n$-dimensional cacti operad $\cac^{(n)}$ (see~\cite{CoVo, Vo/cacti} for the definition).
\end{claim} 
is true for stacks (the corresponding property for manifolds is due to Sullivan-Voronov~\cite{CoVo}).

The case $n=1$ follows from Theorem~\ref{T:HCFT} since the operad $H_\bullet(\cac^{(1)})$ is the \BV-operad (see~\cite{Get, CoVo, SaWa}). (Indeed, $1$-dimensional cacti can be seen as special kind of chord diagram).

We believe the methods introduced in Section~\ref{Loopproduct} and Section~\ref{S:HCFTbig} could be applied to prove the above claim provided one has a model for the  $n$-dimensional cacti operad in which cacti  are obtained by gluing 
 the various lobes using trees.

\begin{rmk}
According to a result of Sullivan and Voronov~\cite[Theorem 5.1.1]{CoVo}, there is an isomorphism of operads between $H\bullet(\cac^{(n)})$ and the homology operad $H_\bullet(\mathcal{D}^{fr}_n)$ of the framed little $n$-dimensional disks operad $\mathcal{D}^{fr}_n$ (studied in details in~\cite{SaWa}). Hence the claim, if proved, implies such a structure on $\sH_{\scriptstyle \bullet}\big(\bfhom(S^n,\XX)\big)$.
\end{rmk}

\begin{rmk}
As in~\cite{CoVo}, one can  prove that the claim follows from
a $n$-dimensional cactus algebra structure of the free sphere stack $\bfhom(S^n,\XX)$ in the
category of {\it correspondences of topological  stacks} (and not of
topological stacks). This category {\it Cor} has Hurewicz
topological stacks for objects and morphisms from $\XX$ to $\YY$
given by diagram $\XX \leftarrow \ZZ \rightarrow \YY$.  The
composition is defined by taking pullbacks. The proof
of~\cite{CoVo} applies {\it verbatim} to the framework   of
 stacks. However, applying this idea to pass to homology is rather subtle as one needs to be careful with Gysin maps and do not seem to be straightforward.
\end{rmk}

%%%%%%%%%%%%%%%%%%%%%%%%%%%%%%%%%%%%%%%%%%%%%%%%%%%%%%%%%%%%%%%%%%%%%%%%%%%%%%%%%%%%%%%%%%%%%%%%%%%%%%%%%%%%%%%%%%%%%%%%%%%%%%%%%%%%%%%%%%%%%%%%%%%%%%%%%%%%%

%%%%%%%%%%%%%%%%%%%%%%%%%%%%%%%%%%%%%%%%%%%

\section{Orbifold intersection pairing}\label{Orbifolds}
In this Section, unless otherwise stated, $\XX$ will be an almost complex orbifold. Then $\IXX$ is again an almost complex orbifold. In particular, $\XX$ and $\IXX$ are oriented orbifolds. Care has to be taken, because even if $\XX$ is connected and has constant dimension, $\IXX$ usually has many components of varying dimension (the so-called twisted sectors). Using the (almost) complex structure and our bivariant theory, we will define a refinement of the \hidden product, which is the (Poincar\'e) dual of the orbifold cup-product~\cite{CR}.

\begin{warning}In this section, all (co)homology groups are taken with coefficients in $\mathbb{C}$, the field of complex numbers. In particular this is true for singular homology $H\lcom(\XX)$, de Rham cohomology (denoted  $\DR{\scriptstyle \bullet}(\XX)$) and compactly supported de Rham cohomology (denoted  $\DRC{\scriptstyle \bullet}(\XX)$).
\end{warning}

\subsection{Poincar{\'e} duality and orbifolds}\label{Poincareduality}
For (any) oriented orbifolds, there is the {\bf Poincar\'e duality homomorphism}
$\PH: H_i(\XX) \to H^{d-i}(\XX)$~\cite{Behrend}. Here $\XX$ is an oriented orbifold which has constant (real) dimension $d=\dim (\XX)$. Let us recall briefly the definition of the Poincar\'e duality homomorphism, see~\cite{Behrend} for details. There is the canonical inclusion $H_i(\XX) \hookrightarrow \big(H^i(\XX) \big)^*$  which is an isomorphism if $H_i(\XX)$ is finite dimensional.   Since $\XX$ is of dimension $d$, there is the Poincar\'e duality isomorphism $\big(\DR{i}\big)^*\stackrel{\sim}\lra \DRC{d-i}$~\cite{Behrend}.  Let ${\rm inc}\: \DRC{\scriptstyle \bullet}(\XX) \to \DR{\scriptstyle \bullet}(\XX)$ be the canonical map. The Poincar\'e duality homomorphism $\PH$ is the composition
\begin{multline} \label{eq:defP}
H_i(\XX) \longrightarrow (H^i(\XX))^* \stackrel{\sim}\longrightarrow (\DR{i}(\XX))^* \stackrel{\sim}\longrightarrow \DRC{d-i}(\XX)\\
\stackrel{{\rm inc}}\longrightarrow \DR{d-i}(\XX) \stackrel{\sim}\longrightarrow H^{d-i}(\XX).
\end{multline}
If the orbifold $\XX$ is proper, then $\PH\: H\lcom(\XX) \to H^{d-\bullet}(\XX)$ is an isomorphism. 

\medskip

Recall that the inertia stack $\IXX$ has usually many components of varying dimension. The inverse map $I\:\IXX \to \IXX$ is the isomorphism defined for any object $(X,\varphi)$ in $\IXX$, where $X$ is an object of $\XX$ and $\varphi$ an automorphism of $X$, by $I(X,\varphi)=(X,\varphi^{-1})$. In the language of groupoids, if $\XX$ is presented by a Lie groupoid $\gm$, the map $I$ is presented by the map $(\gamma,\alpha)\mapsto (\gamma^{-1},\beta)$ for  any $(\gamma,\alpha) \in S{\gm}\times_{\gm_0} \gm_1$.

  The age is a locally constant function $\age :\IXX \to \mathbb{Q}$. If $\XX=[M/G]$ is a global quotient with $G$ a finite group, then $$\IXX=\left[\left(\coprod_{g\in G} M^g \right)/G \right]$$
and for $x\in M^g$, the age is equal to $\sum k_j$ if the eigenvalues of $g$ on $T_xM$ are $\exp(2i\pi k_j)$ with $0\leq k_j <1$. The age does not depend on which way $\XX$ is considered as a global quotient. So it is well-defined on $\IXX$ for any arbitrary almost complex orbifold, because any such $\XX$ can be locally written as a global quotient $[M/G]$.
Similarly, the dimension is a locally constant function $\dim\:\IXX \to \mathbb{Z}$. The age and the dimension are related by the formula (for instance see~\cite{CR, FG})
\eq \label{eq:age}
\dim &= & d -2\,\age -2\,  \age\circ I
\eneq
where $I\:\IXX \to \IXX$ is the inverse map (as above).
The {\bf orbifold homology} of $\XX$ is
\eqn
\Horb{\scriptstyle \bullet}(\XX) = H_{\scriptstyle \bullet -2\, \age \circ I}(\IXX) =\bigoplus_{q\in \mathbb{Q}} H_{\scriptstyle \bullet -2q} \big( [\IXX]_{\age \circ I =q}\big)
\eneqn
where $[\IXX]_{\age =n}$ is the component of $\IXX$ for which the age is equal to $n$. In plain english, we define the orbifold homology of $\XX$ to be the homology of $\IXX$ with a local degree shifting given, on a component of a certain fixed $\age$, by $-2\age \circ I$. According to formula~\eqref{eq:age}, the local degree shifting is also equal to $\dim-d+2\age$. 

\smallskip

 The orbifold cohomology is $H_{\rm orb}\com(\XX)=H^{\scriptstyle \bullet -2\, \age}(\IXX)$ (see~\cite{CR, FG}). Note that the shift of degrees are not the same, but rather are Poincar\'e dual. Indeed,
 since $\IXX$ is an oriented orbifold, there is the Poincar\'e duality homomorphism $\PH\:H\lcom (\IXX) \to H^{\scriptstyle \bullet}(\IXX)$ obtained as the composition~\eqref{eq:defP} above on every connected component of $\IXX$. Since $\IXX$ has in general several connected components of different dimensions, this is not a graded map with respect to the usual grading. However:
\begin{lem}\label{L:defP}
The Poincar\'e duality homomorphism
$$H\lcom (\IXX) \stackrel{\PH}\longrightarrow H^{\scriptstyle \bullet}(\IXX)$$
maps $\Horb{i}(\XX)$ into $H_{\rm orb}^{d-i}(\XX)$. We call it the {\bf orbifold Poincar\'e duality homomorphism} $\PHo\:\Horb{i}(\XX)\to H_{\rm orb}^{d-i}(\XX)$.
\end{lem}
\begin{pf}
It follows from formula~\eqref{eq:age}.
\end{pf}

\subsection{Orbifold intersection pairing and \hidden product} \label{SS:OrbifoldsHidden}

Recall that, if $\XX$ is a manifold, then the homology $H\lcom(\XX)$ has the intersection pairing and the cohomology $H\com(\XX)$ has the cup-product. The Poincar\'e duality homomorphism is an algebra map. However, if $\XX$ is not compact, the intersection ring and cohomology ring may be very different from each other (for instance, if $\XX$ is not compact, $H\lcom(\XX)$ has no unit).

Chen-Ruan~\cite{CR} defined the orbifold cup-product on the
cohomology $H_{\rm orb}\com(\XX)$ of an almost complex orbifold $\XX$, generalizing the cup-product for
manifolds. We will define the analogue of Chen-Ruan orbifold product
in homology. Our construction generalizes the intersection pairing
for manifolds. Note that we do not assume our orbifolds to be
compact.

\smallskip

Our definition of the  orbifold intersection pairing is as follows.
There are the canonical maps $j:\IXX\times_\XX \IXX \to \IXX \times \IXX$ and $m:\IXX \times_\XX \IXX\to \IXX$ (see Section~\ref{stringproduct}) and a Gysin homomorphism $j^!\:H(\IXX \times \IXX) \to H(\IXX \times_\XX \IXX)$ (in homology)  because  $j$ is strongly oriented. Note that this Gysin maps is \emph{not} the same as the one obtained by pulling back (as in Section~\ref{stringproduct}) the orientation class of the diagonal $\XX\to \XX\times \XX$ in general.

\smallskip

The main ingredient in the definition of Chen-Ruan orbifold cup-product is the so
called {\em obstruction bundle} whose construction is explained in
details in~\cite{CR} and~\cite{FG}. Another very nice reference for this is~\cite{JKK}. The obstruction bundle is a bundle over $\IXX\times_\XX \IXX$ denoted $\Ob_\XX$.  We denote $\mathfrak{e}_\XX =e(\Ob_\XX)$ the Euler class of $\Ob_\XX$.
The {\bf orbifold intersection pairing} is the composition:
\begin{multline} \label{eq:defint}
H(\IXX) \otimes H(\IXX) \stackrel{\times}\longrightarrow H(\IXX\times \IXX) \stackrel{j^!}\longrightarrow H(\IXX\times_\XX \IXX) \\ \stackrel{\cap \eu_\XX}\longrightarrow H(\IXX\times_\XX \IXX) \stackrel{m_*}\longrightarrow H(\XX).
\end{multline}
\begin{them}\label{th:intersection}
Suppose $\XX$ is an almost complex orbifold of (real) dimension $d$.
\begin{enumerate}
\item The orbifold intersection pairing defines a bilinear pairing
$$ \Horb{i}(\XX) \otimes \Horb{j}(\XX) \stackrel{\Cap}\longrightarrow \Horb{i+j-d}(\XX) .$$
\item The orbifold intersection pairing $\Cap$ is associative and graded commutative.
\item The orbifold Poincar\'e duality homomorphism $\PHo\:\Horb{\scriptstyle \bullet}(\XX) \lra H_{\rm orb}^{d-\bullet}(\XX)$ is a homomorphism of $\mathbb{C}$-algebras, where $H_{\rm orb}^{d-\bullet}(\XX)$ is equipped with the orbifold cup-product~\cite{CR}.
\end{enumerate}
\end{them}
Recall that graded commutative means that, for any $x\in H_{i}([\IXX]_{\age\circ I = k})\subset \Horb{i+2k}(\XX)$ and $y\in H_{j}([\IXX]_{\age\circ I =\ell})\subset \Horb{j+2\ell}(\XX)$, one has $$x\Cap y =(-1)^{(i+2k)(j+2\ell)} y\Cap x.$$

\begin{pf}

1. By Riemann-Roch, the obstruction bundle $\Ob_\XX$ satisfies the following well-known formula (see~\cite{FG} Lemma 1.12 and~\cite{CR} Lemma 4.2.2):
\eq \label{eq:rankO}
\rank(\Ob_\XX)= 2(\age \circ p_1 +\age \circ p_2 -\age \circ m) +\dim_2 -\dim \circ m
\eneq
where $p_1,p_2:\IXX\times_\XX \IXX \to \IXX$ are the projections on the first and second factor respectively, $\dim_2 :\IXX\times_\XX \IXX \to \mathbb{Z}$ is the dimension function of the orbifold $\IXX\times_\XX \IXX$ and $\rank: \Ob_\XX \to \mathbb{Z}$ is the rank function of the vector bundle $\Ob_\XX$ (as a real vector bundle).  Since $j:\IXX\times_\XX \IXX \to \IXX \times \IXX$ has codimension equal to $\dim \circ p_1 +\dim \circ p_2 -\dim_2$, the result follows from formula~\eqref{eq:rankO} and formula~\eqref{eq:age}.

\smallskip

2. Since ${\rm flip}(\Ob_\XX) \cong \Ob_\XX$ (for instance see~\cite{FG}) and $\eu_\XX$ is of even degrees (thus strictly commutes with any class), the commutativity follows  as  in the proof of~\ref{th:stringproduct}. It remains to prove the associativity. Consider the cartesian diagrams
{\small \eq \label{eq:Ex12}
\xymatrix@=10pt@M=6pt{ & \IXX\times_\XX \IXX \times \IXX \drto^{\; m_{12}} & \\
\IXX\times_\XX \IXX \times_\XX \IXX   \urto^{j_{(12)3}} \drto_{m_{12}}& &\IXX\times \IXX\\
& \IXX\times_\XX \IXX  \urto_{j} & }\\
\label{eq:Ex23}\xymatrix@=10pt@M=6pt{ & \IXX\times \IXX \times_\XX \IXX \drto^{\; m_{23}} & \\
\IXX\times_\XX \IXX \times_\XX \IXX   \urto^{j_{1(23)}} \drto_{ m_{23}}& &\IXX\times \IXX\\
& \IXX\times_\XX \IXX  \urto_{j} & }
\eneq}
The map $j_{(12)3}$, $j_{1(23)}$ are the canonical embeddings induced by $j:\IXX\times_\XX \IXX \to \IXX\times \IXX$ (applied, respectively, to the last two and first two factors). The maps $m_{ii+1}$ ($i=1,2$) are induced by multiplication of the components $i,i+1$. Similarly we denote $p_{ij}\:\IXX\times_{\XX}\IXX \times_\XX \IXX \to \IXX $ ($i\neq j$) the map $(p_i,p_j)$ induced by the projections on the component $i$ and $j$ and $$j_{12}=j\times \id: \IXX ^{\times 3}\to \IXX \times_\XX \IXX \times \IXX,$$  $$j_{23}=\id \times j:\IXX^{\times 3}\to \IXX \times \IXX\times_\XX \IXX$$ the respective embeddings.
The  excess bundle $\Ex_{12}$ associated to diagram~\eqref{eq:Ex12} is defined as follows .  There is a canonical map from the normal bundle $N_{j_{(12)3}}$ of $\IXX\times_\XX \IXX \times_\XX \IXX \stackrel{j_{(12)3}}\to \IXX\times_\XX \IXX \times \IXX$ to the restriction $m_{12}^* N_{j}$ of the normal bundle of $ \IXX\times_\XX \IXX\stackrel{j}\to \IXX$. By definition $\Ex_{12}=\coker(N_{j_{(12)3}}\hookrightarrow m_{12}^* N_{j})$. Similarly, there is the excess bundle $\Ex_{23}=\coker (N_{j_{(1(23)}}\hookrightarrow  m_{23}^* N_{j})$ associated to diagram~\eqref{eq:Ex23}.
The proof of Theorem~\ref{Associativityloop} together with the commutativity of $\eu_\XX$ with any class, shows that
\eqn (\alpha\Cap \beta)\Cap \gamma &=& m_*\left(j^!\big({m_{12}}_*\big(j_{12}^!(\alpha\times \beta\times \gamma) \cap p_{12}^*\eu_\XX\big)\big)\cap \eu_\XX\right)\\
&=& m_*\left({m_{12}}_*\big({j_{(12)3}}^!\big((j_{12}^!(\alpha\times \beta\times \gamma) \cap p_{12}^*\eu_\XX\big)\cap e(\Ex_{12})\big)\cap  \eu_\XX\right)\\
&=& m^{(2)}_*\left({j^{(2)}}^!(\alpha\times \beta\times \gamma) \cap p_{12}^*\eu_\XX \cap e(\Ex_{12}) \cap m_{12}^*\eu_\XX\right).\eneqn
The second line follows from the excess bundle formula (see Proposition~\ref{prop:Excess}) applied to diagram~\eqref{eq:Ex12}. Similarly,
\eqn
\alpha\Cap (\beta\Cap \gamma) &=&  m^{(2)}_*\left({j^{(2)}}^!(\alpha\times \beta\times \gamma) \cap p_{23}^*\eu_\XX \cap e(\Ex_{23}) \cap m_{23}^*\eu_\XX\right).
\eneqn
Hence we need to prove that the bundles $\Ob_\XX$ and $\Ex_{ij}$ satisfy the following identity
\eq\label{eq:ObIAs}
p_{12}^*(\Ob_\XX)+m_{12}^*(\Ob_\XX) + \Ex_{12} &= & p_{23}^*(\Ob_\XX)+m_{23}^*(\Ob_\XX) + \Ex_{23}
\eneq in the $K$-theory group of vector bundles over $\IXX\times_\XX \IXX \times_\XX \IXX$.

\smallskip

The main property of the obstruction bundle $\Ob_\XX$ is that it
precisely satisfies an "affine cocycle condition" see Equation~\eqref{eq:ObAs} below. In fact, there are two cartesian squares (for $i=1,2$), analogous to~\eqref{eq:Ex12}, \eqref{eq:Ex23} {\small \eq \label{eq:ExOb} &&
  \xymatrix@=10pt@M=6pt{ & \IXX\times_\XX \IXX \drto^{m} & \\
\IXX\times_\XX \IXX \times_\XX \IXX   \urto^{p_{ii+1}} \drto_{m_{ii+1}}& &\IXX\\
& \IXX\times_\XX \IXX  \urto_{p_i} & } \eneq }
Since $p_{12}=p_{12}\circ  j_{(12)3}$ and $p_1=p_1\circ j$, it is easy to check that the   "excess" bundles associated to diagram~\eqref{eq:ExOb} for $i=1,2$ coincide with $\Ex_{12}$ and $\Ex_{23}$ respectively. Indeed, there are the following identities
\begin{eqnarray} \label{eq:E12v}
\E_{12}\!\!\! &=& \!\!\!
p_{12}^*m^*T_{\IXX}+T_{\IXX\times_{\XX}\IXX\times_{\XX}\IXX}-p_{12}^*T_{\IXX\times_\XX
\IXX}
- m_{12}^*T_{\IXX\times_\XX \IXX},\\
\label{eq:E23v} \!\!\E_{23}\! \!\!\! &=& \!\!\!\!\!
p_{23}^*m^*T_{\IXX}+T_{\IXX\times_\XX \IXX\times_\XX
  \IXX}-p_{23}^*T_{\IXX\times_\XX \IXX}
- m_{23}^*T_{\IXX\times_\XX \IXX}.
\end{eqnarray}
in the $K$-theory group of vector bundles over $\IXX\times_\XX \IXX \times_\XX \IXX$.
 over $\IXX\times_\XX \IXX\times_\XX \IXX$ associated to the diagram~\eqref{eq:ExOb} defined by $\E_{ii+1}= \coker (N_{p_{ii+1}} \to m_{ii+1}^* N_{p_{ii+1}})$.
It follows from Lemma 4.3.2 and Proposition 4.3.4 in~\cite{CR} (also
see Lemma 1.20 and Proposition 1.25 of~\cite{FG} for more details)
that $\Ob_\XX$ satisfies the following "associativity" equation
\eq \label{eq:ObAs}
p_{12}^*(\Ob_\XX) +m_{12}^*(\Ob_\XX)+\E_{12}&=&p_{23}^*(\Ob_\XX) +m_{23}^*(\Ob_\XX)+\E_{23}
\eneq in the $K$-theory group of vector bundles over $\IXX\times_\XX
\IXX\times_\XX \IXX$. This is precisely identity~\eqref{eq:ObIAs}; the associativity of $\Cap$ follows.

\medskip

3. Since $\dim: \IXX\to \mathbb{Z}$ is always even, $\PH$ commutes with the cross product.   Using  general argument on the Poincar\'e duality homomorphism in~\cite{Behrend}, Proposition~\ref{prop:AB} and tubular neighborhood (see Section~\ref{S:nnsre}), it is straightforward that $\PH \circ f^!=f^* \PH$ for any strongly oriented map of orbifolds $f\:\XX\to \YY$. \comment{here we use that $\dim$ is always even (as a real manifold) thus $\PH$ commutes with $\times$ and we have no signs issues for  $\PH \circ f^!=f^* \PH$} Hence the following diagrams are commutative
{\eqn \xymatrix@M=3pt@C=16pt@R=16pt{  
H\com(\IXX) \otimes H\com(\IXX) \rto^{\; \times}&  H\com(\IXX\times \IXX) \rto^{j^*} &  H\com(\IXX\times_\XX \IXX)  \\
H\lcom(\IXX) \otimes H\lcom(\IXX) \rto^{\; \times}\uto^{\PH} & H\lcom(\IXX\times \IXX) \rto^{j^!}\uto^{\PH}& H\lcom(\IXX\times_\XX \IXX) \uto^{\PH},}
\eneqn}
{\eqn \xymatrix@M=3pt@C=16pt@R=16pt{  
  H\com(\IXX\times_\XX \IXX) \rto^{\cup \eu_\XX}  & H\com(\IXX\times_\XX \IXX) \rto^{\qquad m_!} & H\com(\XX) \\
 H\lcom(\IXX\times_\XX \IXX) \rto^{\cap \eu_\XX}\uto^{\PH} & H\lcom(\IXX\times_\XX \IXX) \rto^{\qquad m_*}\uto^{\PH} & H\lcom(\XX)\uto^{\PH}  .}
\eneqn}
Now the result follows from Lemma~\ref{pro:globalCR} below.
\end{pf}
\begin{numrmk}\label{rmk:compact}
If $\XX$ is compact, the orbifold Poincar\'e duality map is a linear isomorphism, thus an isomorphism of algebras according to Theorem~\ref{th:intersection}.3.
\end{numrmk}

\begin{lem}\label{pro:globalCR}
The Chen-Ruan orbifold cup-product~\cite{CR} is the composition
\begin{multline*}
H\com(\IXX)\otimes H\com(\IXX)\stackrel{\times}\longrightarrow
H\com(\IXX\times\IXX)\stackrel{i^*}\longrightarrow
H\com(\IXX\times_{\XX}\IXX)\\ \stackrel{\cup \eu_\XX}\longrightarrow
H\com(\IXX\times_{\XX}\IXX)\stackrel{m^!}\longrightarrow H\com(\IXX).
\end{multline*}
\end{lem}

\begin{pf}
The Chen-Ruan pairing in~\cite{CR} is defined, for compact orbifolds, by the formula
\eq \label{eq:CRdef}
 \langle \alpha\CR \beta,\gamma \rangle_{\rm orb}&=& \int_{\IXX\times_{\XX}\IXX}
  p_1^*(\alpha)\cup
 p_2^*(\beta) \cup m^*(I^*(\gamma)) \cup f.
\eneq
Until the end of this proof, let us write $\mu$ for the pairing given by the formula of
Proposition~\ref{pro:globalCR}. We  compute $\langle
\mu(\alpha,\beta),\gamma\rangle_{\rm orb}$. Denoting $\int_{\IXX}$ the orbifold integration map defined in~\cite{CR}, we find
\eqn
\langle
\mu(\alpha,\beta),\gamma\rangle_{\rm orb} &= & \int_{\IXX}
\mu(\alpha,\beta)\cup I^*(\gamma)\\
&=& \int_{\IXX} m^!\big(p_1^*(\alpha)\cup p_2^*(\beta)\cup \eu_\XX\big) \cup
I^*(\gamma)\\
&=& \int_{\IXX} m^!\big(p_1^*(\alpha)\cup p_2^*(\beta)\cup \eu_\XX\cup
m^*(I^*(\gamma))\big)\\
&=&  \int_{\IXX\times_{\XX}\IXX}
  p_1^*(\alpha)\cup
 p_2^*(\beta) \cup m^*(I^*(\gamma)) \cup \eu_\XX.
\eneqn By nondegeneracy of the orbifold pairing, we get $\alpha\CR
\beta =\mu(\alpha,\beta)$.
\end{pf}

Similarly to the twisted \hidden product~\ref{stringproduct}, we now introduce  orbifold intersection pairing twisted by a cohomology class.
\begin{defn} \label{def:twistOrb}\begin{enumerate}\item Let $\alpha\in H\com(\IXX\times_{\XX}\IXX)$ be a (not necessarily
  homogeneous) cohomology class. The {\em orbifold intersection pairing twisted by $\alpha$}, denoted $\Cap^\alpha$, is the composition
{\footnotesize \eqn
H(\IXX) \otimes H(\IXX) \stackrel{\times}\lra H(\IXX\times \IXX) \stackrel{j^!}\lra H(\IXX\times_\XX \IXX) \llra{\cap (\eu_\XX\cup \alpha)} H(\IXX\times_\XX \IXX) \stackrel{m_*}\lra H(\XX).
\eneqn }
\item Let $\EE$ be a vector bundle over $\IXX\times_\XX \IXX$. We call $\Cap^{e(\EE)}$ the orbifold intersection pairing twisted by $\EE$.
\end{enumerate}
\end{defn}

With similar notations as for Theorem~\ref{Associativitycocycle}, we prove
\begin{prop}\label{CRcocycle}
\begin{enumerate}
\item If $\alpha$ satisfies the cocycle condition:
$$p_{12}^*(\alpha )\cup (m\times 1)^*(\alpha )= p_{23}^*(\alpha)
\cup (1\times m)^*(\alpha )$$
 in $H\com(\IXX\times_{\XX}\IXX\times_{\XX} \IXX )$,
then $\Cap^\alpha\: H(\IXX)\otimes H(\IXX) \to H(\IXX)$ is
associative.
\item If  $\EE$ is a bundle over $\IXX\times_{\XX}\IXX$ which satisfies the cocycle condition
\begin{equation}
 \label{eq:cocycle}
p_{12}^*(\EE)+(m\times 1)^*(\EE)= p_{23}^*(\EE)+(1\times m)^*(\EE)
\end{equation}
 in the $K$-theory group of vector bundles over $\IXX\times_{\XX}\IXX \times_{\XX}\IXX$,  then
$\Cap^{e(\EE)}$ is  associative.
\end{enumerate}
\end{prop}
\begin{pf}
It follows as
 Theorem~\ref{twistedstring} and Theorem~\ref{th:intersection}.2.
\end{pf}

We will now explain the relationship between the \hidden product and the orbifold intersection product. Namely, the first one is obtained by twisting the second one by an explicit vector bundle over $\IXX\times_\XX \IXX$.

%%%%%%%%%%%%%%%%%%%%%%%%%%%%%%%%%%%%%%%%%%%%%%%%%%%%%%%%%%
%%%%%%%%%%%%%%%%%%%%%%%%%%%%%%%%%%%%

The inverse map $I\:\IXX \stackrel{\sim}\to \IXX$ induces  the "inverse" obstruction bundle  $\Ob^{-1}_\XX=(I\times_\XX I)^*(\Ob_\XX)$ which is also a bundle  over $\IXX\times_\XX \IXX$. \unsure{(one could also pull back the above bundle by ${\rm flip}$, since ${\rm flip}^* \Ob_\XX \cong \Ob_\XX$).} 
We let $\N_\XX$ be the normal bundle of the regular embedding $\IXX\times_{\XX} \IXX\stackrel{m}\to
\IXX$.
\begin{them}\label{them:globaltwisting}
For any almost complex orbifold $\XX$,
 the  \hidden product coincides with the orbifold intersection pairing twisted by $\Ob^{-1}_\XX\oplus \N_\XX$, {\it i.e.}, for  any
  $x \in H\lcom(\IXX)$, $$x\star y = x\Cap^{e(\Ob^{-1}_\XX \oplus \N_\XX)}\,y.$$
\end{them}
The proof will reduce to the following two lemmas.

\smallskip

The first lemma relates the \hidden product with the Gysin homomorphism $j^!\:H(\IXX \times \IXX) \to H(\IXX \times_\XX \IXX)$ and yet another canonical bundle on  $\IXX\times_\XX \IXX$. 
We define the \emph{full excess bundle}  $\F_\XX$ as the excess bundle associated to the cartesian diagram
{\small \eqn \xymatrix@=10pt@M=6pt{ & \IXX \drto^{\ev}& \\
\IXX_{\XX}\IXX \urto^{p_1} \drto_{p_2}& & \XX\\
& \IXX \urto_{\ev} & } \eneqn}
{\it i.e.}, $\F_\XX = \coker \big(N_{\IXX\times_\XX \IXX \stackrel{p_1}\to \IXX} \hookrightarrow N_{\IXX \stackrel{\ev}\to \XX}\big)$.

\begin{lem}
\label{globalloopasnap}Let $\XX$ be an almost complex orbifold.
The  \hidden product $\star:H(\IXX)\otimes H(\IXX)\to H(\IXX)$ is equal to the composition
{\begin{multline*}
H(\IXX) \otimes H(\IXX) \stackrel{\times}\lra H(\IXX\times \IXX) \stackrel{j^!}\lra H(\IXX\times_\XX \IXX)\\ \stackrel{\cap e(\F_\XX)}\lra H(\IXX\times_\XX \IXX) \stackrel{m_*}\lra H(\XX).
\end{multline*}}
\end{lem}
\begin{pf} Apply the excess formula (Proposition~\ref{prop:Excess}).
\end{pf}

The next Lemma relates the four bundles introduced above.
\begin{lem}\label{globaltwistCR}
The obstruction bundle satisfies the identity $$\Ob_\XX
+\N_\XX+\Ob_\XX^{-1}=\FF_\XX$$ in the $K$-theory group of vector
bundles over $\IXX\times_\XX \IXX$.
\end{lem}
\begin{pf}
Recall that $\Ob_\XX$ is solution of the equation~\eqref{eq:ObAs}:
 \eqn
p_{12}^*(\Ob_\XX) +m_{12}^*(\Ob_\XX)+\E_{12}&=&p_{23}^*(\Ob_\XX) +m_{23}^*(\Ob_\XX)+\E_{23}
\eneqn in the $K$-theory group of vector bundles over $\IXX\times_\XX
\IXX\times_\XX \IXX$.

\smallskip

For any permutation $\tau\in \Sigma_3$ of the set $\{1,2,3\}$, there is a map
$\mathcal{T}_\tau: \IXX\times_\XX \IXX \to \IXX\times_\XX \IXX$
defined as the composition { \small \eqn
 \xymatrix@M=2pt@C=30pt{\IXX\times_\XX \IXX \ar[r]^{\hspace{-0.4cm}(p_1,p_2,I\circ m)}  &
  \IXX\times_\XX \IXX\times_\XX \IXX \ar[r]^{\tilde{\tau}} &
  \IXX\times_\XX \IXX\times_\XX \IXX \ar[r]^{\quad p_{12}}& \IXX\times_\XX \IXX ,}
\eneqn}
where $\tilde{\tau}$ is the permutation of  factors induced by
$\tau$. It is well-known (see~\cite{CR}, \cite{FG} Lemma
1.10) that
$\mathcal{T}_\tau^*(\Ob_\XX)\cong \Ob_\XX$.

\smallskip

Let $r$ be the map
$(p_1,p_2, I\circ p_2)\:\IXX\times_\XX \IXX\to \IXX\times_\XX
\IXX\times_\XX \IXX$. Note that \eqn p_{12}\circ r=\id, \qquad
(m_{12}\circ r)= I \circ \mathcal{T}_{(13)}, \qquad   p_{23}\circ r =(p_2, I\circ p_2).\eneqn and furthermore  $r^* m_{23}^* T_{\IXX\times_\XX \IXX}\cong p_1^*(T_{\IXX})$. It follows (using Equation~\eqref{eq:E12v}, Equation~\eqref{eq:E23v} and $
\N_\XX = m^*T_{\IXX}-T_{\IXX\times_\XX \IXX}$)
 that the  pullback of
  Equation~\eqref{eq:ObAs}   along $r:\IXX\times_\XX \IXX\to \IXX\times_\XX
\IXX\times_\XX \IXX$, yields the identity
\eqn
\Ob_\XX +\Ob_\XX^{-1}+\NN_\XX&=&
{\ev}^*T_{\XX}-T_{\IXX\times_\XX \IXX}
-p_1^*T_{\IXX}- p_2^*T_{\IXX}\\
&& +r^*p_{23}^*(\Ob_\XX) +r^*m_{23}^*(\Ob_\XX)
\eneqn in the $K$-theory group of vector bundles over $\IXX\times_\XX
\IXX$. Since the right-hand side  of the first line is isomorphic to
$\FF_\XX$, it suffices to prove that $r^*p_{23}^*(\Ob_\XX)$ and
$r^*m_{23}^*(\Ob_\XX)$ have rank $0$.  It is an easy application of the Riemann-Roch formula~\eqref{eq:rankO}. 
\end{pf}
\begin{numrmk}
 One can easily check that the full excess bundle also satisfies the ``associativity'' equation~\eqref{eq:ObAs}. It thus follows from Lemma~\ref{globaltwistCR} that the twisting bundle $\N_\XX+\Ob_\XX^{-1}$ satisfies equation~\eqref{eq:cocycle}.
\end{numrmk}

\medskip

\noindent{{\sc Proof of Theorem~\ref{them:globaltwisting}}.
By Lemma~\ref{globalloopasnap}, it suffices to prove that $e(\FF_\XX)=\eu_\XX\cup e(\Ob^{-1}_\XX\oplus \NN_\XX)$ which is trivial by Lemma~\ref{globaltwistCR}.
 \nolinebreak $\Box$ }

\begin{numrmk}Let us sum up the philosophy underlying Theorem~\ref{them:globaltwisting}.
In Section~\ref{stringproduct}, we have defined an associative product on $H_{\bullet}(\IXX)$ by using the Gysin map $\Delta^!\:H_{\bullet}(\IXX\times \IXX)\to H_{\bullet-d}(\IXX\times_\XX \IXX)$ for general oriented stacks. In the case of orbifolds, one can directly use the regular embedding $j\:\IXX\times_\XX \IXX\to \IXX\times \IXX$ (more precisely the regular embbedings associated to the various connected components of $\IXX$) to define a Gysin map $j^!$. Due to the excess formula (Proposition~\ref{prop:Excess}), the map $j^!$ does not induce an associative product, so that to get such a product one needs to twist this map by a class satisfying the ``affine cocycle'' condition~\eqref{eq:ObIAs}. In the case of almost complex orbifold, the obstruction bundle is a small bundle satisfying this equation. The excess formula ensures that the  Gysin map $\Delta^!$ is the Gysin 
map $j^!$ twisted by the Euler class of a bundle (the full excess bundle by Lemma~\ref{globalloopasnap}).  Then Lemma~\ref{globaltwistCR} proves that the obstruction bundle is a (virtual) subbundle of the full excess bundle and indeed, measures the difference between these two bundles and hence between the two Gysin maps. 
\end{numrmk}

\begin{numrmk}
According to Theorem~\ref{them:globaltwisting}, Theorem~\ref{th:intersection}.3 and Remark~\ref{rmk:compact}, if $\XX$ is compact,  the orbifold Poincar\'e duality homomorphism $\PHo$ induces an isomorphism of algebras between the \hidden algebra $\big(\sH\lcom(\IXX),\star \big)$ and the orbifold cohomology equipped with Chen-Ruan orbifold cup-product twisted by the class $e(\Ob^{-1}_XX\oplus \N_\XX)$. A nice interpretation of this isomorphism has  recently been found by Gonz{\'a}lez, Lupercio, Segovia and Uribe~\cite{GLSU}. They proved that
the \hidden product of compact complex orbifolds is isomorphic to the Chen Ruan product of the
cotangent bundle $T^*\XX$ of $\XX$.
\end{numrmk}

\subsection{Examples of orbifold interesection pairing}
Let us consider now  a few examples. Note that the orbifold cup-product has been computed for many compact orbifolds. For instance one can refer to~\cite{CR, FG, CCLT, GHK}. Complex toric orbifolds were dealt on in~\cite{Ji, GH}. By Theorem~\ref{th:intersection}, the intersection pairing coincides with the orbifold cup-product in these cases.

\begin{ex}\label{ex:intersectiongroup}
 Let $G$ be a finite group. Then $[*/G]$ is a complex orbifold. It follows from Theorem~\ref{them:globaltwisting} that the intersection product is the same as the \hidden product in that case. Note that all the ages are trivial and thus $\Horb{i}(\XX)=H_i(\IXX)$ for all $i$'s. Thus the algebra $(\Horb{\bullet}, \Cap)$ is concentrated in degree $0$ and isomorphic to the center $Z(\mathbb{C}[G])$ of the group algebra $\mathbb{C}[G]$, see Example~\ref{prop:G}.
 
 \smallskip
 
 Now let $V$ be a $\mathbb{C}$-linear representation of $G$ of dimension $n>0$. Then the quotient stack  $[V/G]$ is a complex orbifold. Since all invariants subspaces $V^g$ ($g\in G$) are contractible, the homology $H_\bullet(\Lambda [V/G])$ is still concentrated in degree $0$ and isomorphic to $(\mathbb{C}[G])_G$ as a vector space. However, the ages are no longer trivial (and depends on the specific representation) and add an interesting combinatorial grading to $(\mathbb{C}[G])_G$. 
 It follows that $\Horb{\bullet}(G)=\bigoplus_{[g]\in C(G)} \mathbb{C}[2\age (g^{-1})]$ where  $C(G)$ is the set of conjugacy classes of $G$ and $\mathbb{C}[q]$ means $\mathbb{C}$ placed in degree $q$.
 Since the age is invariant by conjugation, $Z(\mathbb{C}[G])$ (which is isomorphic to $ \Horb{\bullet}(G)$ as a vector space) inherits a non-trivial grading induced by the age.  Since $[V/G]$ is of dimension $n>0$, the \hidden product is trivial. Similarly, the  pairing $[g]\Cap [h]=0$ if either $V^g$ or $V^h$ is  non-trivial, while if $V^g=V^h=\{0\}$, $[g]\Cap [h]$ is non-zero and behaves as in the case of $[*/G]$. 
 It follows that there is an isomorphism of graded algebra $(\Horb{\bullet}([V/G]),  \Cap)\cong Z_{V,G} \oplus Ann_{V,G} $ where $Z_{V,G}$ is the (graded by the age) subalgebra of $Z(\mathbb{C}[G])$ generated by the elements with trivial invariant subspace and $Ann_{V,G}$ is the complementary subspace equipped with the zero multiplication. Finally, the orbifold cohomology ring $H_{\rm orb}^\bullet(\XX)$ is isomorphic, as  a graded ring, to $Z(\mathbb{C}[G])$ equipped with the age grading. Thus for a generic representation none of the three rings are isomorphic as rings, though they are as  vector spaces.   
\end{ex}

\begin{ex}
We consider the Kummer surface orbifold $\KK=\big[(S^1)^4/\mathbb{Z}/2\mathbb{Z}\big]$ where  $\mathbb{Z}/2\mathbb{Z}$ acts diagonally by $z\mapsto 1/z$ (on each factor). Then $\KK$ is a complex orbifold with 16 orbifolds points $\big\{\big((-1)^\ep_1,\dots, (-1)^\ep_4\big)\big\}$. We identify the above sets with $({\mathbb{Z}/2\mathbb{Z}})^4$ (and we will denote $\ep=(\ep_1,\dots, \ep_4)$ an element in $({\mathbb{Z}/2\mathbb{Z}})^4$).  The age of the generator $\tau$ of $\mathbb{Z}/2\mathbb{Z}$ is one at every orbifold point.  Thus,  as a graded vector space, 
$\Horb{\bullet}(\KK)\cong H_{\bullet}(\KK) \oplus \big(\mathbb{C}[-2]\big)\big\langle{({\mathbb{Z}/2\mathbb{Z}})^4}\big\rangle $
(where $\mathbb{C}[-2]$) is the vector space $\mathbb{C}$ placed in homological degree 2. Since the coefficient ring is $\mathbb{C}$, the cohomology of $\KK$ is computed by the cohomology of invariant form $\Omega\big((S^1)^4\big)^{\mathbb{Z}/2\mathbb{Z}}$ and its homology by the homology of $\mathbb{Z}/2\mathbb{Z}$-coinvariant chains. Hence
$$H_\bullet(\KK)\cong \mathbb{C}\cdot 1 \oplus \mathbb{C}[-2]\langle a_{i,j} ,\, 1\leq i <j \leq 4\rangle \oplus \mathbb{C}[-4] \alpha_{\KK}$$
where $\alpha_{\KK}$ is the fundamental class of $\KK$ and $a_{i,j}$ is the class given by the cross product of the fundamental classes of the $i^{th}$ and $j^{th}$ circles in $\big(S^1\big)^4$. Thus $$\Horb{\bullet}(\KK)\cong \mathbb{C}\cdot 1 \oplus \mathbb{C}[-2]\langle a_{i,j} ,\, 1\leq i <j \leq 4\rangle \oplus \big(\mathbb{C}[-2]\big)\big\langle{({\mathbb{Z}/2\mathbb{Z}})^4}\big\rangle \oplus \mathbb{C}[-4] \alpha_{\KK}$$ 
 The computation of the intersection pairing is similar to the computations in Section~\ref{S:E:quotient}. We find that $\alpha_{\KK}$ is a unit for the orbifold intersection pairing and that $a_{i,j}\Cap a_{k,l}=1$ if the 4 indices are distinct and $0$ otherwise.
The (degree 2) generators $\ep \in {({\mathbb{Z}/2\mathbb{Z}})^4}$ comes from the degree $0$ homology of the corresponding orbifold point. Hence $\ep \Cap \ep= 1$ and $\ep \Cap \ep'=0$ if $\ep \neq \ep'$. A degree argument also shows that $\ep \Cap a_{i,j}=0$. 
This finishes the description of the whole intersection pairing of $\KK$. The \hidden product of $\KK$ is a little bit different  since $\ep \star \ep=0$ for the latter product.
\end{ex}

\begin{ex}{\label{E:noncompact}}
Here is a (rather trivial) non compact example of orbifold product. Let $d$ be a positive integer and
consider the weighted projective stack $\bbP(d,d,\cdots,d)=[V/\mathbb{C}^*]$,
where $V:=\mathbb{C}^{n+1}-\{(0,0,\cdots,0)\}$ and
the action of  $\lambda \in \mathbb{C}^{*}$ is multiplication by $(\lambda^d,\lambda^d,\cdots,\lambda^d)$. 
Let $\UU$ be an open substack
of $\bbP(d,d,\cdots,d)$. Denoting $\mu_d$ the group of $d^{th}$roots of unity, we have 
    $$\Lambda\UU=\underset{\xi\in \mu_d}\coprod \UU.$$
Note that the coarse moduli space of $\bbP(d,d,\cdots,d)$ is $\bbC\bbP^{n}$, and the coarse moduli space
of $\UU$ is an open $U$ in $\bbP^{n}$. Combining the results of example~\ref{ex:intersectiongroup} and remark~\ref{R:hiddenformanifolds} we see that the orbifold homology ring of $\UU$ is 
  $$\Horb{\bullet}(\UU)=\mathbb{C}[\mu_d]\otimes H_{\bullet}(U),$$
where $H_{\bullet}(U)$ is the usual homology of $U$ endowed with intersection product, and 
$\mathbb{C}[\mu_d]$ is the group ring of the group $\mu_d$ of $d^{th}$ roots of unity (sitting in degree zero).
\end{ex}

\begin{ex}{\label{E:wtdprojline}} We compute the orbifold homology of the weighted projective line
 $\bbP(m,n)$. Recall that, for $m$ and $n$ positive integers, $\bbP(m,n):=[V/\mathbb{C}^*]$,
 where $V:=\mathbb{C}^2-\{(0,0)\}$ and the action of  $\lambda \in \mathbb{C}^{*}$ 
 on $V$ is by multiplication by $(\lambda^m,\lambda^n)$. We do \emph{not} assume $n$ and $m$ to be relatively prime integers.
 
 Note that, except for $\bbC\mathbb{P}^1=\bbP(1,1)$, a weighted projective line is never a (orbifold) global quotient as it is simply 
 connected. We can, however, cover $\bbP(m,n)$ by two global quotient open substacks as follows.
 Denote the point $[0:1] \in \bbP(m,n)$ by $0$ and the point $[1:0] \in \bbP(m,n)$
 by $\infty$. Then
   $$\bbP(m,n)-\{\infty\}\cong [\mathbb{C}/\mu_{n}] \ \ \text{and} \ \ 
      \bbP(m,n)-\{0\}\cong[\mathbb{C}/\mu_{m}],$$
where $\mu_n\cong\mathbb{Z}/n\mathbb{Z}$ stands for the group of $n^{th}$ roots of unity. The action of $ \xi\in \mu_n$ on
$\mathbb{C}$ is given by  $x \mapsto \xi^dx$, where
$d=\gcd(m,n)$. It follows that the inertia group of the point  $0$ is $\mu_n$ and 
the inertia group of  $\infty$ is $\mu_m$. The inertia group of every other point is $\mu_d$. Note that $\mu_d$ sits inside $\mu_n$ as the cyclic subgroup generated by $e^{2i\pi  /d}=\big(e^{2i\pi /n}\big)^{n/d}$ and similarly, it also sits inside $\mu_m$ as the cyclic subgroup generated by $\big(e^{2i\pi /m}\big)^{m/d}$.

Let us describe the inertia stack $\Lambda \bbP(m,n)$. Over the point $0 \in \bbP(m,n)$
consider a disjoint union 
of $n$ copies of $[*/\mu_n]$ indexed by $\mu_n$. Similarly, over the point $\infty$ consider
a disjoint union of $m$ copies of $[*/\mu_m]$ indexed by $\mu_m$. For every $\psi \in \mu_d$, 
``join'' the  copies of $[*/\mu_n]$ and $[*/\mu_m]$ corresponding to $\psi \in \mu_n$ and 
$\psi\in\mu_m$, respectively,  by viewing them as the $0$ and the $\infty$
points of  a new copy of $\bbP(m,n)$ (indexed by $\psi$). Then, $\Lambda \bbP(m,n)$
is the union of all these, that is,
  $$\Lambda \bbP(m,n)=
  \big(\underset{\xi\in \mu_n}\coprod[*/\mu_n]\big)\cup_{\mu_d} 
  \big(\coprod_{\psi\in \mu_d}\bbP(m,n)\big)\cup_{\mu_d}  
  \big(\underset{\zeta\in \mu_m}\coprod[*/\mu_m]\big).$$
The map $\Lambda \bbP(m,n) \to \bbP(m,n)$ is the obvious one (it is the identity on each copy of $\bbP(m,n)$
in $\Lambda \bbP(m,n)$). 

Let us now describe the orbifold homology $\Horb{\bullet}(\bbP(m,n))$. First we determine the ages of the twisted sectors of
$\bbP(m,n)$. For every $\psi \in \mu_d$, the age of  $\bbP(m,n)_{\psi}$ 
is zero. For $\xi =e^{2i\pi  l/n} \in \mu_n$, the age of $[*/\mu_n]_{\xi}$  is $\{dl/n\}$. (For $r \in \mathbb{R}$,
we define $0\leq \{r\}< 1$ to be  the residue  of $r$ modulo 1.)
Similarly, the age of $[*/\mu_m]_{\zeta}$, $\zeta=e^{2i\pi  k/m}  \in \mu_m$, is $\{dk/m\}$. Observe that the ages of
two twisted sectors $[*/\mu_n]_{\xi}$ and $[*/\mu_m]_{\zeta}$ are distinct unless they lie on some  $\bbP(m,n)_{\psi}$
(that is, if $\xi=\zeta=\psi$, for some $\psi \in \mu_d$).

The above decomposition of $\Lambda \bbP(m,n)$ implies that $\Horb{\bullet}(\bbP(m,n))$ is a union
of three subrings $R_0$, $R$ and $R_{\infty}$ (given by the sectors/copies over the point $0$, $\infty$ and a generic point of $\bbP(n,m)$) . Only $R$ is a unital subring. 

The ring $R$ is the homology of $\coprod_{\psi\in \mu_d}\bbP(m,n)$
and is isomorphic, as a graded ring, to $\mathbb{C}[\mu_d]\otimes H_{\bullet}(\mathbb{CP}^1)$, where the elements of the group
ring $\mu_d$ are sitting in degree 0 (this computation is the one done in Example~\ref{E:noncompact}). More precisely, if $\beta \in H_0(\mathbb{CP}^1)=\mathbb{C}$ is the
homology class of a point and $\alpha \in H_2(\mathbb{CP}^1)=\mathbb{C}$ is the fundamental class, 
then $\psi\otimes \beta$ corresponds to the homology of a 
point in $H_0(\bbP(m,n)_{\psi})$ and $\psi\otimes \alpha$ corresponds to the fundamental class
of $\bbP(m,n)_{\psi}$. In particular, the unit element of $R$ (and also of the whole $\Horb{\bullet}(\bbP(m,n))$ is
$1\otimes \alpha$.

\smallskip

The ring $R_0$ is isomorphic, as a vector space, to $\mathbb{C}[\mu_n]$ (as in Example~\ref{ex:intersectiongroup}). For $\xi_l \in \mu_n$, 
the basis element $\xi_l \in \mathbb{C}[\mu_n]$
corresponds to the generator of $H_{\bullet}([*/\mu_n]_{\xi_l})=H_{0}(*)=\mathbb{C}$. Its degree is equal to $2-2\{-dl/n\}$,
where   $\xi_l=e^{2i\pi  l/n}$. Note that if $\psi=\xi_{\ell n/d}\in \mu_d\subset \mu_n$, then $[*/\mu_n]_{\psi}\subset \bbP(n,m)_\psi$ and the class $\xi_{\ell n/d}$ is identified with the degree $0$ generator of $H_0(\bbP(n,m_\psi)$.  Every class $\xi_i$ in $R_0$ is given by (a degree shifting of)  an ordinary degree 0 homology class of the point $0$ viewed as lying in the $\xi_i$-fixed point locus $0^{\xi_i}$ .
 The map $m\: \Lambda \bbP(m,n)\times_{\bbP(n,m)} \Lambda \bbP(n,m)\to \bbP(n,m)$ clearly maps $0^{\xi_i, \xi_j}$ (the intersection locus of the copies of $[*/\mu_n]_{\xi_i}$ and $[*/\mu]_{\xi_j}$  indexed by $\xi_i$ and $\xi_j$) to $0^{\xi_{i+j}}$ (the copy of $[*/\mu_n]$ indexed by $\xi_{i+j}$).  If  $\xi_i$  is in $\mu_d$ (that is if  is of the form $\ell n/d$), then $\xi_i\Cap \xi_j =0$ for degree reasons, since $\xi_i$ is a degree zero homology class. Thus to compute the intersection pairing, it now remains to analyze the obstruction bundle in the other cases. For an integer $i$, let us denote $0\leq \{i\}_n< n/d$ the residue of $i$ modulo $n/d$.  From  (the Riemann-Roch) formula~\eqref{eq:rankO}, we find that, for $\xi_i,\xi_j\notin \mu_d$, the obstruction bundle (over the copies indexed by $\xi_i$ and $\xi_j$) is of rank $0$ if  $$\{i\}_n +\{j\}_n\le n/d$$ and is of rank 2 otherwise.  Indeed, if $\{i\}_n +\{j\}_n< n/d$, then $\{d(i+j)/n\}=\{di/n\} +\{dj/n\}$ and $\dim_2=0=\dim\circ m$. If  $\{i\}_n +\{j\}_n> n/d$, then $\{d(i+j)/n\}=\{di/n\} +\{dj/n\}-1$  and $\dim_2=0=\dim\circ m$. If  $\{i\}_n +\{j\}_n= n/d$, then  $\{d(i+j)/n\}=\{di/n\} +\{dj/n\}-1$  and $\dim_2=0$ but now $\dim\circ m=2$ (the dimension of the copy  $ \bbP(n,m)_{\xi_{i+j}}$).

Thus, it follows that the orbifold intersection pairing $\xi_i\Cap \xi_j$ of basis elements $\xi_i$ and $\xi_j$ is given by 
\begin{equation}\label{eq:R0Cap}
\xi_i\Cap \xi_j = \left\{\begin{array}{lcc} \xi_{i+j} & \mbox{if } \{i\}_n +\{j\}_n\le n/d & \mbox{and } \xi_i, \xi_j \notin \mu_d\\ 
0 & \mbox{if }  \{i\}_n +\{j\}_n> n/d & \mbox{or } \xi_i \mbox{ or } \xi_j \in \mu_d \end{array} \right. 
\end{equation}
In particular, $R_0$ is generated, as a graded ring, by  the elements $\xi_1, \xi_{1+n/d}, \dots$. 
It is easy to describe the action of $R$ on $R_0$. Indeed, for degree reasons, $(\psi\otimes \beta) \Cap \xi_i=0$. However, since the fundamental class of $\bbP(n,m)_{\xi_{\ell n/d}}$ intersects $[*/\mu_n]_{\xi_i}$ for any $\xi_i \in \mu_n$ with a trivial obstruction bundle, we have $(\xi_{\ell n/d}\otimes \alpha)\Cap \xi_i =\xi_{i+\ell n/d}$ (in particular $1\otimes \beta$ acts as the unit). Summing up, we have
\begin{equation}\label{eq:RactsR0}
\left\{\begin{array}{rcl}
(\psi\otimes \beta) \Cap \xi_i & = & 0 \\
(\xi_{\ell n/d}\otimes \alpha)\Cap \xi_i & = & \xi_{i+\ell n/d}
\end{array}\right.
\end{equation}

\smallskip

The ring $R_\infty$ is isomorphic, as a vector space, to $\mathbb{C}[\mu_m]$. For $\zeta_k \in \mu_m$, 
the basis element $\zeta_k \in \mathbb{C}[\mu_m]$
corresponds to the generator of $H_{\bullet}([*/\mu_m]_{\zeta_k})=H_{0}(*)=\mathbb{C}$. Its degree is equal to $2-2\{-dk/m\}$,
where   $\zeta_k=e^{2i\pi  k/m}$. We denote  $0\leq \{i\}_m< m/d$ the residue of any integer $i$ modulo $m/d$. A computation similar to the one of $R_0$ shows that, the orbifold intersection pairing $\zeta_p\Cap \zeta_q $ of basis elements $\zeta_p, \zeta_q \in \mu_m$ is given by 
\begin{equation}\label{eq:RinftyCap}
\zeta_p\Cap \zeta_q = \left\{\begin{array}{lcc} \zeta_{p+q} & \mbox{if } \{p\}_m +\{q\}_m\le m/d & \mbox{and } \zeta_p, \zeta_q \notin \mu_d\\ 
0 & \mbox{if }  \{p\}_m +\{q\}_m> m/d & \mbox{or } \zeta_p \mbox{ or } \zeta_q \in \mu_d \end{array} \right. 
\end{equation}
and that the $R$ action on $R_\infty$ is given by
\begin{equation}\label{eq:RactsRinfty} 
\left\{\begin{array}{rcl} (\psi\otimes \beta)\Cap \zeta_p &=& 0 \\
(\zeta_{\ell m/d}\otimes \alpha) \Cap \zeta_p &= & \zeta_{p+\ell m/d}.\end{array}\right.
\end{equation}
Since the points $0$ and $\infty$ do not intersect, it is immediate that $R_0 \Cap R_\infty=\{0\}$.

\smallskip

From the above descriptions (identifying the sectors $[*/\mu_n]_{\xi_{\ell n/d}}$ with their image in $\bbP(n,m)_{\xi_{\ell n/d}}$), we find that the orbifold homology $\Horb{\bullet}(\bbP(n,m))$ is the graded ring 
$$ \Horb{\bullet} = R_0\oplus R\oplus R_\infty /(\xi_{\ell n/d}=\xi_{\ell n/d} \otimes \alpha =\zeta_{\ell m/d}\otimes \alpha = \zeta_{\ell m/d} )$$ where the product structure is given by formulas~\eqref{eq:R0Cap}, \eqref{eq:RactsR0}, \eqref{eq:RinftyCap} and~\eqref{eq:RactsRinfty}.

The reader can check that the orbifold product on $\Horb{\bullet}(\bbP(m,n))$  is Poincar\'e dual to the Chen-Ruan
orbifold cup product in Example 5.3 of \cite{CR} (note that {\em loc. cit.} only considers the relatively prime case). It is not hard to check that the \hidden product is different. Indeed, the grading is different (concentrated in degree 0 and 2) and the only non zero products of basis elements are those involving $\psi\otimes \beta$.
\end{ex}
%%%%%%%%%%%%%%%%%%%%%%%%%%%%%%
%%%%%%%%%%%%%%%%%%%%%%%%%%%
%%%%%%%%%%%%%%%%%%%%%%%%%%

\section{Examples}\label{S:Examples}
\subsection{The case of manifolds}\label{ChasSullivan}
Smooth manifolds form a special class of differentiable stacks with
normally non-singular diagonal (Definition~\ref{D:nns}). Denote by
the same letter $M$ a manifold and its associated (topological)
stack. The diagonal  $\Delta\: M\to M\times M$ is  strongly oriented
iff the manifold $M$ is oriented.
\begin{prop}\label{Manifold}Let $M$ be an oriented manifold. The \BV-algebra,   Frobenius algebra and (non-unital, non-counital) homological conformal field theories structures of $H\lcom(\LM)$ given by Theorem~\ref{BV}, Theorem~\ref{frobeniusloop} and Theorem~\ref{T:HCFT} coincide with Chas-Sullivan~\cite{CS}, Cohen-Jones~\cite{CoJo}, Cohen-Godin~\cite{CoGo} and Godin~\cite{Go} ones.
\end{prop}
\begin{pf}
By Proposition~\ref{prop:targetconnected},  the free loop
stack of $M$ is isomorphic to the free loop space $\LM$. It follows from Proposition~\ref{prop:AB} and Proposition~\ref{P:Gorinetation} (in the case $G=\{1\}$), that the Gysin maps of Sections~\ref{Loopproduct}, \ref{frobenius}, \ref{BVstructure} coincide with the Gysin maps in~\cite{CoJo} Section 1 (also see \cite{FTV} Section 3.1).
\end{pf}

\begin{numrmk}\label{R:hiddenformanifolds}
When $M$ is an oriented manifold, the \hidden product on $H\lcom(\Lambda M)\cong H\lcom(M)$ is simply the usual intersection pairing. It is also immediate that, if $M$ is an almost complex manifold, the orbifold intersection pairing of $M$  also coincides with the usual intersection pairing.
\end{numrmk}

\begin{numrmk}
 Similarly, the brane product for manifolds given by Proposition~\ref{P:brane} coincides with Sullivan-Voronov one~\cite{CoVo} and Chataur one~\cite{Ch}.
\end{numrmk}

%%%%%%%%%%%%%%%%%%%%%%%%%%%%%%%%%%%%%%%%%%%%%%%%%%%%%%%%%%%%%%%%%%%%%%%%%%%%%%%%%%%%%%%%%%%%%%%%%%%%%%%%%%%%%%%%%%%%%%%%%%%%%%%%%%%%%%%%%%%%%%%%%%%%%%%%%%%%%%%

\subsection{Hidden loop (co)product for global quotient by a finite group}\label{Orbifoldcomputations}
A  special important class of oriented orbifolds is  the
global quotient $[M/G]$, where $G$ is a finite group, $M$ is an oriented  manifold together with an action  of $G$ by orientation preserving diffeomorphisms. In this case, the homology of the inertia
stack $H([M/G])$  is well known.  Assume that our coefficient ring
$\kor$ is a field of
characteristic coprime with $|G|$ (or $0$). The inertia stack of $[M/G]$ is represented by the
transformation groupoid \eq &&\label{eq:inertiaGQ}\copd\limits_{g\in G} M^g\times G \rightrightarrows
\copd\limits_{g\in G} M^g\eneq
where the action of $h\in G$  moves $y\in M^g$  to
$ y\cdot h\in M^{h^{-1}gh}$. Furthermore, $\IXX\times_{\XX}
\IXX \cong \big[\coprod_{g,h\in G} M^{g,h}/G\big]$, where $M^{g,h}=M^g\cap M^h$, and the ``Pontrjagin'' map $m:
\IXX\times_{\XX} \IXX\to \IXX$ is induced by
 the embeddings $i_{g,h}\:M^{g,h}\hookrightarrow M^{gh}$. Since $|G|$ is coprime
 with ${\rm char}(\kor)$, the homology groups of the
 inertia stack $\Lambda[M/G]$  are
$$H\lcom(\Lambda [M/G])\cong H\lcom\left(\copd\limits_{g\in G} M^g\right)_{G}\cong
\left(\bigoplus_{g\in G} H\lcom(M^{g})\right)_G .$$

\smallskip

 The excess bundle $Ex(M,X,X')$ of the  diagram of embeddings
{\footnotesize $$ \xymatrix@C=5pt@R=7pt@M=5pt{ & X \drto & \\
Z=X\bigcap X' \urto \drto& & M\\
& X' \urto & }$$}
is the cokernel
of the bundle map $N_{Z\hookrightarrow X}\hookrightarrow (N_{X'\hookrightarrow M})_{/Z}$. Thus  $Ex(M,X,X')$ is the virtual bundle $T_M -T_X
-T_{X'}+T_{Z}$  (each component being restricted to
$Z$).
For $g,h\in G$, we denote
$Ex(g,h):=Ex(M,M^{g},M^{h})$. The bundles $Ex(g,h)$ induce a
 bundle $Ex$ on $\Lambda[M/G]\times_{[M/G]}\Lambda[M/G]$ whose Euler class is denoted  $e(Ex)$.
Since the diagonal
$G\to G\times G$ is a group monomorphism, there is    a
transfer map ${\rm tr}^G_{G\times G}\:(\bigoplus_{g,h\in G} H\lcom(M^{g})\otimes
H\lcom(M^h))_{G\times G} \to \big(\bigoplus_{g,h\in G}
H\lcom(M^{g})\otimes H\lcom(M^h)\big)_{G}$ explicitly given (see
Equation~\eqref{eq:transfer}) by
$${\rm tr}^G_{G\times G} (x)=\sum_{g\in G}x \cdot (g,1).$$    The maps $i_g\:M^{g,h}\hookrightarrow M^g$, $i_h:M^{g,h}\hookrightarrow M^h$ yield   Gysin morphims $(i_g\times i_h)^!\:H\lcom(M^g\times M^h)\to H\lcom(M^{g,h})$.

\begin{prop}\label{loopasnap}
The  \hidden product $\star:H(\Lambda [M/G])\otimes H(\Lambda [M/G])\to H(\Lambda [M/G])$ is
 the composition {\footnotesize \eqn
  \left(\bigoplus_{g\in G} H(M^g)\right)_G\otimes
  \left(\bigoplus_{h\in
  G}H(M^h)\right)_G\rightarrow \left(\bigoplus_{g,h\in
  G}H(M^g\times M^h)\right)_{G\times G}
  \hspace{-0.5cm}\stackrel{{\rm tr}^G_{G\times G}}\longrightarrow   \left(\bigoplus_{g,h\in
  G}H(M^g\times M^h)\right)_{G} \\
\stackrel{\bigoplus (i_g\times i_h)^!}\longrightarrow
  \left(\bigoplus_{g,h\in G} \, H(M^{g,h})\right)_{G}
\stackrel{\cap e(Ex)}\longrightarrow  \left(\bigoplus_{g,h\in G} \, H(M^{g,h})\right)_{G}
\stackrel{m_*}\longrightarrow\left(\bigoplus_{k\in G} H(M^k)\right)_G \qquad \eneqn}
\end{prop}
The proof of Proposition~\ref{loopasnap} relies on Lemma~\ref{identificationofnap} below,
which is  of independent interest.
 Note that there is a oriented stack morphism  \eq \label{eq:varphi} \varphi\: \IMG\times_{[M/G]}\IMG \cong \left[\copd_{g,h \in G}
 M^{g,h}/G\right]\to \IMG\times
 \IMG\eneq induced by the groupoid map $(x,g)\mapsto (i_g(x),g, i_h(x),g)$.
\begin{lem}\label{identificationofnap}
The Gysin map $\varphi^!$ is the composition
{\footnotesize \eqn \left(\bigoplus_{g,h\in G} \, H(M^g\times M^h)\right)_{G\times G}\hspace{-0.4cm}\stackrel{{\rm tr}^G_{G\times G}}\longrightarrow \left(\bigoplus_{g,h\in G} \, H(M^g\times M^h)\right)_G \stackrel{\bigoplus (i_g\times i_h)^!}\longrightarrow \left(\bigoplus_{g,h\in G} \, H(M^{g,h})\right)_{G}\hspace{-0.2cm}.\eneqn}
\end{lem}
\begin{pf} The $G$-equivariant map $M^{g,h}\to M^g \times M^h$, given by $x\mapsto (i_g(x),i_h(x))$, induces a oriented stack morphism  $\psi\: \IMG\times_{[M/G]} \IMG
 \to \IMG\times_{[*/G]} \IMG$. By Proposition~\ref{prop:AB},
 $\psi^!=\oplus (i_g\times i_h)^!$. Then, the result follows from the functoriality of Gysin maps and Lemma~\ref{Transfer}.
\end{pf}

\begin{trivlist}\item[]{\sc Proof of Proposition~\ref{loopasnap}.}
  We use the notations of Section~\ref{stringproduct}. The cartesian
  diagram~\eqref{eq:cartesianghost} (where $\XX=[M/G]$) and the excess formula~\ref{prop:Excess} shows that,
  $$\gy{\Delta}=\varphi^!(x)\cap
  e(Ex).$$ Thus the result   follows   from
  Lemma~\ref{identificationofnap}.
\nolinebreak $\Box$ \end{trivlist}

\smallskip

Similarly we compute
 the \hidden coproduct. For any $g\in G$, the unit $1_g\in H_0(M^g)$ induces a map $1_g\: H(M^g)\to H_0(M^g)\otimes H(M^g)\to H(M^g \times M^g)$.
\begin{prop}\label{stringcoproductorbifolds}
 The \hidden
coproduct is induced (after passing to $G$-invariant) by the composition
{\footnotesize \eqn \bigoplus_{g\in G} H(M^g)\stackrel{\oplus 1_g}\lra \bigoplus_{g\in G} H(M^g \times M^g)\stackrel{tr^G_{G\times G}}\lra
\bigoplus_{g,h\in G}   H(M^h\times M^g) \stackrel{\oplus i_{g,h}^!}\longrightarrow
\bigoplus_{g,h\in G}H(M^{g,h})\\
\stackrel{ \cap
  e(Ex)}\lra \bigoplus_{g,h\in G} H(M^{g,h})
\stackrel{(i_{g},i_h)_*}\lra \bigoplus_{g,h\in G}H(M^g\times M^h)\cong \bigoplus_{g,h\in G} H(M^g)\otimes H(M^h).
\eneqn }
\end{prop}
\begin{pf}
Let $\gm$ be the transformation groupoid $M\rtimes G\toto M$. Unfolding the definition
of the  groupoid $\IWG$ (see
Section~\ref{frobeniusstring}), one finds that $\IWG$ is the transformation groupoid $\left(G\times \coprod_{h\in G} M^h\right)\rtimes G^2 \toto G\times \coprod_{h\in G} M^h $, where  the action of $(h_0,h_{1/2})\in G^2$ on $(g,m)\in {g}\times M^h$ is  $(h_0^{-1}gh_{1/2}, m.h_0)$.
The Morita map $p\:\IWG \to \IMG$ (Equation~\eqref{eq:p})  has
 a section $\kappa$  defined, for $m\in M^h$ and $h_0\in G$, by  $\kappa(m,h_0)=(h,m,
 h_0,h_0)$. In particular  $\kappa$ induces an isomorphism in homology and commutes with Gysin maps.  Thus  the Gysin map $\gy{\Delta}$ of
 Section~\ref{frobeniusstring} is the composition of $\kappa_*$ with the Gysin map associated to the  sequence of cartesian diagrams
{\footnotesize\eq\label{eq:kappa} &&\xymatrix@R=12pt{[\coprod M^{g,h}/G] \rto \dto &
  [G\times \coprod M^h/G] \rto \dto & [G\times \coprod M^h/G\times G]\dto  \\
[M/G] \rto & [M\times M/G] \rto & [M/G]\times [M/G] & . }\eneq} By
Proposition~\ref{prop:Excess} and Lemma~\ref{Transfer}, the Gysin
maps   associated to the left square and the right square are,
respectively,
 $i_{g,h}^!(-)\cap e(Ex(g,h))$ and ${\rm Tr}_{G\times G}^G$. Since $\kappa_*=\oplus 1_g$, the result follows.
\end{pf}

\begin{numex}\label{prop:G}
Consider  $[*/G]$ where $G$ is a finite
group. By Proposition~\ref{loopdiscrete}, the stack morphism $\Phi\: \Lambda[*/G] \to {\rm L}[*/G]$ (see
Lemma~\ref{lem:frobeniusmorphism}) is an isomorphism.
Let $\kor$ be a field of characteristic coprime with
$|G|$. Then
$$H\lcom(\Lambda[*/G])= \big(\bigoplus_{g\in G}\kor
\big)_G\cong \big(\bigoplus_{g\in G}\kor
\big)^{G}\cong Z(\kor[G])$$
where $Z(\kor [G])$ is the center of the group algebra $\kor [G]$.
By Propositions~\ref{loopasnap}, the isomorphism $H\lcom(\Lambda[*/G])\cong Z(\kor[G])$ is an isomorphism of algebras. By Proposition~\ref{stringcoproductorbifolds}, the \hidden coproduct is given by $\delta([g])=\sum_{hk=g}[h]\otimes [k].$ Thus the Frobenius algebra structure coincides with the one given by Dijkgraaf-Witten~\cite{DiWi}.
\end{numex}

%%%%%%%%%%%%%%%%%%%%%%%%%%%%%%%%%%%%%%%%%%%%%%%%%%%%%%%%%%%%%%%%%%%%%%%%%%%%%%%%%%%%%%%%%%%%%%%%%%%%%%%%%%%%%%%%%%%%%%%%%%%%%%%%%%%%%%%%%%%%%%%%%%%%%%%%%%%%%%%%%%%%%%%%%

\subsection{String topology of $[S^{2n+1}/(\mathbb{Z}/2\mathbb{Z})^{n+1}]$}\label{S:E:quotient}

Let $S^{2n+1}$ be the euclidian sphere $\{
|z_0|^2+\cdots+|z_{n}|^2=1, z_i\in \mathbb{C}\}$ acted upon by $\sz$
identified with the group generated the reflections  across the
hyperplanes  $z_i=0$ ($0\leq i\leq n$). Let $\RR=[S^{2n+1}/\sz]$ be
the induced quotient stack  which is obviously an oriented orbifold
of dimension $2n+1$. We now   describe the Frobenius algebras
associated to $\IXX$ and $\LXX$. Until the end of this section we
denote $R=\sz$.

\smallskip

The \hidden product has a very simple combinatorial description.  Let $\Delta^n$
be a $n$-dimensional standard simplex. Denote $v_0,\dots,v_n$ its $n+1$-vertex and  $F_0,\dots, F_n$ its $n$-faces of dimension $n-1$. In other words $F_i=\Delta(v_0,\dots, \widehat{v_i},\dots, v_n)$ is the convex hull of all vertices but $v_i$. More generally we denote $F_{i_1 \dots i_k}:=F_{i_1}\cap \dots \cap  F_{i_k}$ the subface of dimension $n-k$ given by the convex hull of all vertices but $v_{i_1},\dots,v_{i_k}$. We assign the degree $2n -2k +1$ to a face $F_{i_1 \dots i_k}$ of dimension $n-k$.
\begin{prop}\label{stringsphere} Let $\kor$ be a ring with $1/2\in \kor$. Then $H\lcom(\Lambda \RR)$ is the free $\kor$-module with basis indexed by elements $r\in R-\{1\}$ of degree $0$ and all faces $F_{i_1 \dots i_k}$ in degree $2(n-k)+1$ (in particular $F_\emptyset=\Delta^n$ has degree $2n+1$), {\it i.e.},
\eqn
H\lcom(\Lambda \RR) &\cong&
 \kor^{|R|-1}\oplus \Big(\bigoplus_{\mbox{\scriptsize$\begin{array}{c} k=0\dots n \\0\leq i_1<\dots<i_k\leq n \end{array}$}
} \vspace{-1cm} \kor.F_{i_1\dots i_k}\Big).
\eneqn
  The
 \hidden product $\star$ is defined on the basis by the identities
 $$F_{i_1\dots i_k}\star F_{j_1 \dots j_l} = F_{i_1\dots
   i_k}\cap F_{j_1\dots j_l}$$ if the two subfaces have transversal
 intersection in $\Delta^n$,  and is $0$ otherwise. The element
   $\Delta^n=F_\emptyset$ is set to be the unit and
   all other products involving a generator of $\kor^{|R|-1}$ are trivial.
\end{prop}
In other words,  $H_0(\Lambda
\RR)=\kor^{|R|-1}$, and
$H_{2i+1}(\Lambda \RR)$ is the free module generated by the  subfaces of dimension $i$ of the
simplex  $\Delta^n$. The product is given by transverse intersection in $\Delta^n$.
\begin{pf}
Write $s_i$ (i=0\dots n) for the  reflection across the hyperplane $z_i=0$.
 Then, for $0\leq
k\leq n$, \eqn {(S^{2n+1})}^{s_{i_1}\dots s_{i_k}} &\cong& \left\{(z_0,\dots,z_n)\in \mathbb{C}^{n+1} \, / \sum_{j\neq
  i_1,\dots i_k}\! |z_j|^2 =1\right\}\cong S^{2n-2k+1}.\eneqn Thus  $H\lcom\left({(S^{2n+1})}^{s_{i_1}\dots s_{i_k}}\right)\cong \kor\, V'_{s_{i_1}\dots s_{i_k}} \oplus \kor\, F'_{i_1\dots i_k}[2(n-k)+1]$.  Since these
generators are $R$-invariant, $|R|$ is invertible in $\kor$ and
$(S^{2n+1})^{s_0\dots s_n}=\emptyset$, one has \eqn H\lcom(\Lambda \RR)\cong \bigoplus_{g\in R} H\lcom\big((S^{2n+1})^g\big)_R \cong
\bigoplus_{g-\{1\}\in R} H\lcom \big((S^{2n+1})^g\big)\eneqn
By Proposition~\ref{loopasnap}, the \hidden product is the composition of ${\rm tr}_{R\times R}^R$ with
{\footnotesize \eqn  H\big( (S^{2n+1})^g\times (S^{2n+1})^h\big)\stackrel{{(i_g\times i_h)}^!}\to H\big((S^{2n+1})^{g,h}\big)\stackrel{\cap e(Ex(g,h))}\lra H\big((S^{2n+1})^{g,h}\big)\stackrel{m_*}\to H((S^{2n+1})^{gh}) \eneqn}
Clearly ${\rm tr}^R_{R\times R}$ is multiplication by the order of $R$.
  Furthermore $F_{i_1 \dots i_k}$ and $F_{j_1 \dots j_l}$ are transversal iff    the sets $\{i_1,\dots,i_k\}$ and $\{j_1,\dots,
j_l\}$ are disjoint iff   the submanifolds $(S^{2n+1})^{s_{i_1}\dots
  s_{i_k}}$ and $(S^{2n+1})^{s_{j_1}\dots s_{j_l}}$ are transversal in $(S^{2n+1})$.  In particular, if $F_{i_1 \dots i_k}$ and $F_{j_1 \dots j_l}$ are transversal,  $(S^{2n+1})^{s_{i_1}\dots s_{i_k}, s_{j_1}\dots s_{j_l}}=(S^{2n+1})^{s_{i_1}\dots
  s_{i_k}.s_{j_1}\dots s_{j_l}}$, the excess bundle is of rank $0$, $m_*=\id$ and by
Poincar\'e duality,   \eqn (i_{s_{i_1}\dots
  s_{i_k}} \times i_{s_{j_1}\dots s_{j_l}})^!\, (F'_{i_1 \dots i_k} \times
F'_{j_1 \dots j_l})&=& F'_{i_1\dots i_k j_1\dots j_l}.\eneqn
If $F_{i_1 \dots i_k}$ and $F_{j_1 \dots j_l}$ are not transversal, one finds  \eqn (i_{s_{i_1}\dots
  s_{i_k}} \times i_{s_{j_1}\dots s_{j_l}})^!(F'_{i_1 \dots i_k} \times F'_{j_1 \dots j_l})&=& F'_{i_1\dots i_k}\cap F'_{ j_1\dots j_l}=F'_{\{i_1,\dots,i_k\}\cup \{j_1,\dots,j_l\}} \eneqn
 and $(S^{2n+1})^{s_{i_1}\dots s_{i_k} s_{j_1} \dots s_{j_l}}$ contains
  $(S^{2n+1})^{s_{i_1}\dots s_{i_k}, s_{j_1} \dots s_{j_l}}$ as a submanifold
  of codimension $>0$. It follows  that  $m_*\left(F'_{\{i_1,\dots,i_k\}\cup \{j_1,\dots,j_l\}}\cap e(Ex) \right) =0$ for degree reason. \unsure{$m_*$ is  null in positive degrees and the excess bundle is easily computed to be the subbundle of ${T_M}_{/M^{s_{i_1}\dots s_{i_k}, s_{j_1} \dots s_{j_l}}}$ generated by the coordinates $z_s, s\in A$.    It follows easily that $F_{\{i_1,\dots,i_k\}\cup \{j_1,\dots,j_l\}}\cap e(Ex)$ is not in degree $0$, thus vanish when composed with $m_*$.}
Similarly, $F'_{i_1,\dots,i_k}\star g=0$ for any $g\in R$. \unsure{
We now have to compute $F'_{i_1,\dots,i_k}\star g$ for $g\in G$.
For degree reason it is null if $k>0$. If $k=0$,  then the \hidden
product $F'_\emptyset \star g$ is given by: \eqn && H(M \times
M^g)\stackrel{{\rm tr}^G_{G\times G}}\to H(M\times
M^g)\stackrel{(i_1\times i_g)^!}\to H(M^{g}) \stackrel{\cap
e(Ex)}\to H^{M^{g}}\stackrel{p_*}\to H^(M^{g}) \eneqn and both maps
$(i_1\times i_g)^!$ and $m_*$ are the identity. Furthermore, since
$F'_\emptyset$ is transverse with every subface, the Excess bundle
is of dimension $0$. Thus $F'_\emptyset \star g=2^{n+1} g$.} The
result follows by identifying  $F_{i_1,\dots,i_k}$ with
$2^{-n-1}F'_{i_1,\dots,i_k}$ as basis element.
\end{pf}

\begin{numrmk}
It is easy to show that the \hidden coproduct is trivial. Indeed, for degree reason, only the class of $F_\emptyset$ might be non zero. Proposition~\ref{stringcoproductorbifolds} shows the \hidden coproduct is induced by the composition
$$H(S^{2n+1} )\stackrel{\sum i_g^!}\to \bigoplus H((S^{2n+1})^g) \stackrel{\cap \oplus e\big((S^{2n+1})^h\big)}\to \bigoplus H((S^{2n+1})^g).$$ Since $(S^{2n+1})^h$ is an odd dimensional sphere, its  Euler class  is 2-torsion, hence trivial by our assumption on $\kor$.
\end{numrmk}

\medskip

Since $R=(\mathbb{Z}/2\mathbb{Z})^{n+1}$ is abelian, its group algebra
  is a Frobenius algebra (see Example~\ref{prop:G} above).
\begin{prop}\label{loopsphere} Let $\kor$ be a field of characteristic
  different from $2$. There is an isomorphism of \BV-algebras as well
  as Frobenius algebras \eq\label{BVRR} H\lcom(\Lo \RR) &\cong& H\lcom (\Lo S^{2n+1})\otimes_{\kor}
  \kor[(\mathbb{Z}/2\mathbb{Z})^{n+1}].\eneq The \BV-operator on the
  right hand side is $B \otimes {\rm id}$ where $B:H\lcom({\rm L}S^{2n+1})\to H_{\scriptstyle \bullet+1}({\rm L}S^{2n+1})$ is the \BV-operator of the loop homology of $S^{2n+1}$.
\end{prop}
\begin{pf}
According to Proposition~\ref{loopdiscrete}, the free loop stack $\Lo \RR$ is presented by the groupoid  $\coprod_{g \in R} \mathcal{P}_g
 S^{2n+1}\rtimes R \toto \coprod_{g \in R} \mathcal{P}_g
 S^{2n+1}$. Hence
\eqn
 H\lcom(\LXX) &=& \big(\bigoplus_{g\in R} H\lcom(\mathcal{P}_g S^{2n+1})\big)_R.
\eneqn
Since $R$ is a subgroup of the connected Lie group $SO(2n+2)$, which acts  on $S^{2n+1}$,  for all $g\in R$ there is a continuous path $\rho: [0,1]\to
 SO(2n+2)$ connecting $g$ to the identity (that is $\rho(0)=g$,
 $\rho(1)=1$). In particular, any path $f\in P_g S^{2n+1}$ can be composed
 with the path $f(0).\rho(t)$ yielding a loop  $\Upsilon_g(f) \in {\rm
 L}S^{2n+1}$.
It is a general fact  that $\Upsilon_g\: \mathcal{P}_gS^{2n+1} \to {\rm
 L}S^{2n+1}$ is  a $G$-equivariant  homotopy equivalence (see~\cite{LUX} for details). We write
 $\Upsilon: \coprod_{g\in R} \mathcal{P}_g S^{2n+1} \to \coprod_{g\in R} {\rm L}S^{2n+1}$ for the
 map induced by the maps $\Upsilon_g$ for $g\in R$. Since the
 $G$-action on ${\rm L}M=P_e M$ is  trivial, the isomorphism~\eqref{BVRR}  follows.

\smallskip

It remains to prove that the linear isomorphism \eqref{BVRR} is an isomorphism of Frobenius algebras and \BV-algebras. To do so,
 we need the evaluation
map  $\ev_0:\Lo \RR\to \RR$ at the groupoid level. One checks   that
$\ev_0$   is represented by the maps $\ev_g\:\mathcal{P}_g
S^{2n+1}\times R \to S^{2n+1}\times R$ defined by
$\ev_g\big((f,h)\big)=(f(1),h)$. Let $f,g \in \mathcal{P}_r
S^{2n+1}\times \mathcal{P}_h S^{2n+1}$ such that $f(1)=g(1)$. The
composition of the path $f(-)$ and $g(-) \cdot h$  gives an element
$m(f,g)\in \mathcal{P}_{rh}  S^{2n+1}$. This composition induces the
stack morphism  $m: \Lo \RR\times_{\RR}\Lo \RR\to \Lo\RR$. Denote
$\tilde{m}$ the map  $\coprod_{g,h\in R} \Lo
S^{2n+1}\times_{S^{2n+1}} \Lo S^{2n+1} \stackrel{\widetilde{m}}\to
\coprod_{g\in R} \Lo S^{2n+1}$ which maps an element
$(\gamma,\gamma')\in\Lo S^{2n+1}\times_{S^{2n+1}} \Lo S^{2n+1} $ in
the component $(g,h)$
  to the element $m(\gamma,\gamma')$ in the component $gh$. Here $m$ is the usual composition of paths. The map $\Upsilon_g\: \mathcal{P}_g S^{2n+1} \to \Lo S^{2n+1}$ induces a commutative diagram of $R$-equivariant maps
{\footnotesize \eqn
 \xymatrix{\coprod_{g,h} P_g S^{2n+1}\times P_h S^{2n+1} \dto^{\coprod \Upsilon_g\times \Upsilon_h}&\coprod_{g,h} P_g S^{2n+1}\times_{S^{2n+1}} P_h S^{2n+1} \lto \rto^{\quad\quad  \quad m} \dto_{\coprod \Upsilon_g\times \Upsilon_h} & \coprod_{g} P_g S^{2n+1}\dto^{\coprod \Upsilon_g} \\
\coprod_{g,h} \Lo S^{2n+1}\times \Lo S^{2n+1} & \coprod_{g,h} \Lo S^{2n+1}\times_{S^{2n+1}} \Lo S^{2n+1}\lto \rto^{\quad \quad \quad \widetilde{m}} & \coprod_{g} \Lo S^{2n+1}.}
\eneqn}
Since $\Lo S^{2n+1}\to S^{2n+1}$, $\mathcal{P}_g S^{2n+1}\to S^{2n+1}$ are fibration,  the vertical arrows are $R$-homotopy equivalences. It follows easily that \unsure{the Gysin map $\Upsilon: \Lo\RR \to
[\coprod_{g\in R} \Lo S^{2n+1}/R]$ induces a map of algebras in homology.

The map $\Upsilon_g\:\mathcal{P}_g S^{2n+1} \to {\rm
 L}S^{2n+1} $ assemble to give a map
\eqn
(\LG\times_{\gm} \LG)_0 &\stackrel{\widetilde{\Upsilon}}\to &
\coprod_{g,h\in G} \Lo S^{2n+1}\times_M \Lo S^{2n+1}\eneqn
where  $\gm$  the
transformation groupoid $M\rtimes G\toto M$. $\tilde{\Upsilon}$
 is an homotopy equivalence since $P_g M \to \Lo S^{2n+1}$ is an homotopy
 equivalence and .Notice that $\ev_g\:P_g M\to M $ is a fibration as well as $\ev_0:\Lo S^{2n+1} \to M$.  It
 follows that the Gysin map $\gy{\Delta}$ given by the cartesian square
\eqn
&&\xymatrix{(\LG\times_{\gm} \LG)_0 \rto \dto & \coprod_{g,h\in G}\times_{} (\coprod_{g,h} P_g M\times P_h M) \dto^{\pi} \\
M\rto^{\Delta} & M\times M}
\eneqn
commutes with $\Upsilon$ i.e. $\widetilde{\Upsilon}\circ
\gy{\Delta} = \mathcal{C}_*\circ\mathcal{T}_* \circ (\Upsilon\times
\Upsilon)$ where $ \mathcal{C}_*$, $\mathcal{T}_*$ are described in
Equation~\ref{Thomcollapse}.
 It is now sufficient   to identify
$\Upsilon \circ m_*\circ \widetilde{\Upsilon}^{-1}$.   Since
$G$ acts trivially on $\coprod_{g\in G}\Lo S^{2n+1}$,}
the map
$$\cfrac{1}{|R|}\Upsilon: H\lcom\big(\big[\coprod_{g\in R} \Lo S^{2n+1}/R\big]\big) \to
H\lcom(\Lo S^{2n+1})\otimes \kor[R]$$ is a morphism of algebras. One proves
 similarly that $\cfrac{1}{|R|}\Upsilon$ is a coalgebra map.

\smallskip

 Now we need to identify the \BV-operator. Denote $\LG$ the transformation groupoid $\coprod_{g \in R} \mathcal{P}_g
 S^{2n+1}\rtimes R \toto \coprod_{g \in R} \mathcal{P}_g
 S^{2n+1}$. Since the stack $S^1$ is canonically identified with the quotient stack $[\mathbb{R}/\mathbb{Z}]$, the homology $H\lcom(S^1)$
 coincides with the homology of the groupoid $\gm':\mathbb{R}\rtimes
 \mathbb{Z}\toto \mathbb{R}$.  The $0$-dimensional simplex $(0,1)\in \mathbb{R}\rtimes
 \mathbb{Z}=\gm'_1$ defines an element in $C_0(\gm'_1)\subset C_1(\gm')$ which is the generator of $H_1(S^1)$. The  map
$\gm'\times \LG \stackrel{\theta}\to \LG$ defined, for  $(x,n)\in \mathbb{R}\times \mathbb{Z}$, $f\in \mathcal{P}_g$ and $h\in R$, by ${\theta}(x,n,f,h)(t)= f(t+x).h^n$ \unsure{(where )where, for $t\geq 0$, $f(1+t):=f(t).g$} is a groupoid morphism representing the   $S^1$-action on $\LXX$.  Since $\Upsilon(\theta((0,1),f)=f$,  $\Upsilon$ commutes with the \BV-operator.
\end{pf}
\begin{numrmk}For the sake of completeness, we recall~\cite{CJY} that,
$\sH\lcom(\Lo S^{2n+1})\cong \kor [u,v]$ with $|v|=-2n-1$ and $|u|=2n$ for $n>0$, and $\sH\lcom(\Lo S^{2n+1})\cong \kor [[u,u^{-1}]][v]$ if $n=0$. Thus
$$\sH\lcom(\Lo[S^{2n+1}/(\mathbb{Z}/2\mathbb{Z})^{n+1}])\cong \kor[(\mathbb{Z}/2\mathbb{Z})^{n+1}][u,v] \mbox{  if } n>0, \quad \mbox{ and}$$
$$\sH\lcom(\Lo [S^{1}/\mathbb{Z}/2\mathbb{Z}])\cong \kor[[u,u^{-1}]][\tau,v]/(\tau^2=1)\quad \mbox{with } |v|=1,\, |u|=0 \mbox{ if } n=0.$$
\end{numrmk}

\begin{numrmk} The stack morphism $\Phi\: \IXX \to \LXX$ of
  Section~\ref{Frobeniusmorphism} is represented at the groupoid level
  by  $\coprod_{g\in R}(S^{2n+1})^g\hookrightarrow  \coprod_{g\in R} \mathcal{P}_g S^{2n+1}$  where $x\in (S^{2n+1})^g$ is identified with  a
  constant path.  It follows easily that  the Frobenius algebra morphism is given by
$\Phi(F_\emptyset)=e$, $\Phi(F_{i_1\dots i_k})=0$ and $\Phi(g)=gv$.
\end{numrmk}

%%%%%%%%%%%%%%%%%%%%%%%%%%%%%%%%%%%%%%%%%%%%%%%%%%%%%%%%%%%%%%%%%%%%%%%%%%%%%%%%%%

%%%%%%%%%%%%%%%%%%%%%%%%%%%%%%%%%%%%%%%%%%%%%%%%%%%%%%%%%%%%%%%%%%%%%%%%%%%%%%%%%%%%%%

%%%%%%%%%%
%%%%%%%%%%%
%%%%%%%%%%

%%%%%%%

\subsection{String topology of $\Lo [*/G]$ when  $G$ is a compact Lie group}
\label{Liegroups}

Any topological group $G$ naturally defines a topological stack corresponding
to the groupoid  $G\toto \{*\}$, which is denoted by  $[*/G]$.
In this section we study  the Frobenius structures on
 the homology of its loop stack and inertia stack assuming
that   $G$ is a compact and connected Lie group.
It turns out that in this case  the two Frobenius structures obtained
are indeed isomorphic since $\IPG$ and
$\LPG$ are homotopy equivalent.
 In this section, we assume that $G$ is of dimension $d$ and
we will  work with real coefficients for  (co)homology groups.

\smallskip

First we will identify the homology groups $H\lcom(\IPG)$ and
$H\lcom(\LPG)$.

\begin{lem}\label{IPG}
The inertia stack $\IPG$ is represented by
 the transformation groupoid $G\rtimes G\toto G$,  where $G$ acts
on itself by conjugation, while the  stack $\Lambda [*/G]\times_{[*/G]}\Lambda [*/G]$ is represented by the
 groupoid $(G\times G) \rtimes G\toto G\times G$ with the diagonal
conjugacy action.
\end{lem}
%\begin{pf}
%The first assertion is immediate. The second follows from the first and the fact that the pullback of stack is represented by the pullback of groupoids.
%\end{pf}

The following result is well known \cite{Brylinski-book}.

\begin{lem}
\label{folkhomotopy}
The  map $\IPG \stackrel{\Phi}\to \LPG$ is an homotopy equivalence.
\end{lem}
\begin{pf}
Since $G$ is connected, by Proposition~\ref{prop:targetconnected},
$\LPG$ can be represented by the loop group
$\Lo G\toto \{*\}$. On the other hand,  $B\Lo G\cong \Lo BG$ is homotopy equivalent
to $EG\times_{G}G$~\cite{Ben}, \cite{Brylinski-book}. This equivalence can be seen as follows. Denote by $e$ the unit of $G$ and let $P_e G$ be the based path space of $G$, that is $P_e G$ is the set, endowed with the compact-open topology, of paths  $[0,1]\llra{f} G$   such that  $f(0)=e$.  There is an action of the loop group $LG$ on $EG\times P_e G$ given, for $(x,f)\in EG\times P_e G$ and $\gamma \in LG$, by $$(e,f)\cdot \gamma =\big(e\cdot \gamma(0),\gamma(0)^{-1}*f *\gamma\big) $$ where $*$ stands for the (pointwise) multiplication in $G$. This action is clearly free, thus $BLG \cong (EG\times P_e G) \times_{LG} \{{\rm pt}\}$.  The map $(x,f) \mapsto (x,f(1))$ induces a continuous map $\rho: (EG\times P_e G) \times_{LG} \{{\rm pt}\}\to EG\times_{G} G$.  Since $G$ is connected, for any $g \in G$, there is a path $f_g \in P_e G$ with $f_g(1)=g$. The map  $$EG\times G \ni (x,g) \mapsto (x, f_g) \in G\times P_e G$$ induces a well-defined map $\psi\: EG\times_G G \to (EG\times P_e G) \times_{LG} \{{\rm pt}\}$. It is easy that $\psi$ is independent of the choice of the $f_g$s and is a left and right inverse of $\rho$. Hence the homotopy equivalence $BLG \cong EG\times_G G$ follows. Through the isomorphism in between $\Lo [*/G]$ and $[*/LG]$ (Proposition~\ref{prop:targetconnected}), the map 
 $\Phi$ of Lemma~\ref{lem:frobeniusmorphism} is transfered to the map $\psi\: B[G/G] \to (EG\times P_e G)\times_{LG} {\rm pt}\cong B[{\rm pt}/LG]$. Hence the result. 
\end{pf}

As an immediate consequence, we have

%An obvious corollary of Proposition~\ref{folkhomotopy} above and Proposition~\ref{frobeniusmorphism} is

\begin{cor}
The map $\Phi_*\:H\lcom(\IPG) \to H\lcom(\LPG)$ is an isomorphism of Frobenius algebras.
\end{cor}
Thus it is sufficient
 to study  the Frobenius structure on
 the homology of the inertia stack $\IPG$.

%Let $M$ be a  manifold  with a  smooth $(G\times G)$-action.
% Consider $G$ as a subgroup of $G\times G$ by embedding it
%diagonally. Then $M$ is also a $G$-space and we have a  map
%  $[M/G]\to [M/G\times G]$, which is a principal $G$-bundle.
% Recall that
%the (co)homology of a transformation groupoid $[M/G]$ is the same
%than the $G$-equivariant (co)homology of $M$. In this section,
%unless otherwise stated, all groups actions are on the left.
%The loop product on $H\lcom([G/G])$ is of degree $-d$ ($d$ is the dimension of $G$). It  is dual (by universal coefficient see Remark~\ref{.....}) to a degree $-d$ coproduct  $\nabla\: H\lcom([G/G]) \to H\lcom([G/G]) \otimes H\lcom([G/G]) $ on $H([G/G])$. And similarly for the coproduct. In other words there is a dual Frobenius structure on the cohomology that we are going to compute explicitly.
%
\smallskip

According to Remark~\ref{universalcoefficient}, there is a dual Frobenius structure induced on $\big(H\com(\Lambda[*/G]),\star,\delta \big)$. We refer to   $\delta\:H\com(\Lambda[*/G])\to H\com(\Lambda[*/G])\otimes H\com(\Lambda[*/G])$ and $\star:H\com(\Lambda[*/G])\otimes H\com(\Lambda[*/G]) \to H\com(\Lambda[*/G])$ as the dual \hidden coproduct and dual \hidden product respectively. Since, it is technically easier, we will describe the Frobenius structure of $H\com(\Lambda [*/G])$.
The following result is standard  \cite{Me}.
We write $EG$ for a free $G$-space which is contractible
and $BG=EG\times_{G} *$ its classifying space so that
 $H\com([*/G])=H\com(BG)=H_G\com (*)$.

\begin{prop}
\label{cohomologyofLiegroups}
\begin{enumerate}
\item The cohomology of $G$, as a topological space,
 is
$$H\com(G)=\left(\Lambda \Gg^*\right)^{G}\cong \Lambda (y_1, y_2, \cdots, y_l)$$

\item The cohomology of $[*/G]$ is $$H\com([*/G])=\left(S^*(\Gg^*)\right)^{G}\cong S(x_1,x_2,\cdots, x_l)$$

\item The cohomology of $[G/G]$ is
 $$H\com([G/G])=\left(S^*(\Gg^*)\right)^{G}\otimes \left(\Lambda \Gg^*\right)^{G}\cong S(x_1,x_2,\cdots, x_l)\otimes \Lambda (y_1, y_2, \cdots, y_l)$$

%where $l=rank(G)$ and $deg(x_i)=2d_i$.  where $l=rank(G)$, $deg (y_i)=2d_i-1$ and $deg(x_i)=2d_i$.

\item The cohomology of $[G\times G/G]$ is \begin{eqnarray*}
H\com([G\times G/G])&=&\left(S^*(\Gg^*)\right)^{G}\otimes \left(\Lambda (\Gg^*\oplus \Gg^*)\right)^{G}\\
&\cong& S(x_1,x_2,\cdots, x_l)\otimes \Lambda (y_1, y_2, \cdots, y_l,y_1', y_2', \cdots, y_l'),
\end{eqnarray*}
%where $l=rank(G)$, $deg (y_i)=deg(y_i')=2d_i-1$ and $deg(x_i)=2d_i$.
\item The cohomology of $[G\times G/G\times G]$ is
\begin{eqnarray*}H\com([G\times G/G\times G])&=&\left(S^*(\Gg^*\oplus \Gg^*)\right)^{G}\otimes \left(\Lambda (\Gg^*\oplus \Gg^*)\right)^{G}\\
&\cong& S(x_1,x_2,\cdots, x_l,x_1',x_2',\cdots, x_l')\\
&& \otimes \Lambda (y_1, y_2, \cdots, y_l,y_1', y_2', \cdots, y_l')\end{eqnarray*}
%where $l=rank(G)$, $deg (y_i)=deg(y_i')=2d_i-1$ and $deg(x_i)=deg(x_i')=2d_i$.
\end{enumerate}
Here $l=rank(G)$, $deg (y_i)=deg(y_i')=2d_i+1$, $deg(x_i)=deg(x_i')=2d_i$
and $d_i$ are the exponents of $G$.
\end{prop}
%\begin{pf}
%It is well known that the generators $x_i,y_i$ ($i=1\dots r$) can be represented by equivariant homogeneous polynomials from $\Gg \to \mathbb{R}$ and polynomials $\eta_i\: \Gg \to \Omega^*(G)$ respectively. The Kunneth formula implies that the obvious cochain complex map $$ \big(S(\Gg^*)\otimes \Omega^*(G)\big)^{G} \otimes \big(S(\Gg^*)\otimes \Omega^*(G)\big)^{G} \to \big(S(\Gg^*\oplus \Gg^*)\otimes \Omega^*(G\times G)\big)^{G\times G}$$ is a quasi-isomorphism which gives the last statement.
%

%For the third statement, one remark that $\eta_i\otimes 1\in \big(S(\Gg^*)\otimes \Omega^*(G)\otimes \Omega^*(G)\big)^{G}\subset \Omega_G^*(G\times G)$ is a cycle which is mapped to a cycle representing the class $y_i\otimes 1 \in H^*(G)\otimes H^*(G)\cong H^*(G\times G)$. Similarly $1\otimes \eta_i $ is mapped to a generator  $1\otimes y_i$ of $H^*(G)\otimes H^*(G)\cong H^*(G\times G)$. The map $H^*_G(G\times G)\to H^*(G\times G)$ is a map of algebra. Thus Kunneth theorem applied to $H^*(G)$ yields that $H^*_G(G\times G)\to H^*(G\times G)$ is surjective and $G\times G$ is equivariantly formal. In particular its cohomology is equivalent to the $E_2$ term of its spectral sequence which gives the first statement.
%\end{pf}

To compute the Frobenius structure of $H\com(\Lambda [*/G])$, we need an
explicit construction of some Gysin maps.

Let $M$ be an oriented
manifold  with a smooth $(G\times G)$-action.
 Consider $G$ as a subgroup of $G\times G$ by
embedding it diagonally. In this way,
 $M$ becomes a $G$-space and we have a morphism of stacks
 $[M/G]\to [M/G\times G]$, which is indeed a $G$-principle
bundle. According to Section~\ref{S:Gysin}, there is a cohomology
Gysin map $\Delta_!\:H\com{[M/G]}\to H^{\scriptstyle
\bullet-d}{[M/G\times G]}$, which should be in a certain sense
fibration integration.
 %The idea to express it is
 %to make  "a change of coordinate" in order to have $G\times G$ acting on both sides and then to use equivariant space level Gysin map.

Recall that when $G$ is a compact connected Lie group,
the  cohomology of  the  quotient stack $H\com([M/G])$ with real
coefficients can be computed
 using the Cartan model $( \Omega_G(M),  d_G)$, where
$\Omega_G(M):=\left( S(\Gg^*)\otimes \Omega(M) \right)^{G}$
is the space of $G$-equivariant polynomials
 $P\:\Gg \to \Omega(M)$, and
$$d_G (P)(\xi):= d(P(\xi))-\iota_{\xi}P(\xi), \quad \forall \xi \in \Gg. $$
Here  $d$ is the de Rham differential
and $\iota_{\xi}$ is the contraction by the generating vector field of $\xi$.
Given a Lie group $K$ and a Lie subgroup
 $G\subset K$,
let $G$ act on $K$  from the right  by multiplication and
$K$ act on itself from the left by multiplication.
The submersion $K\to K/G$ is a
principal $K$-equivariant right $G$-bundle. There is an isomorphism of stacks $[M/G]\iso
[K\times_{G} M/K]$ which  induces  an isomorphism in cohomology. It is known
\cite{Me} that, on the Cartan model, this isomorphism can be described  by an induction map ${\rm Ind}^G_K:
\Omega_G(M)\to  \Omega_K(K\times_{G} M)$. Here $G$ acts
 on $K\times M $ by
$$(k,m)\cdot g=(k\cdot g, g^{-1}\cdot m). $$
 The induction map is the composition
$$\Omega_G(M)\stackrel{{\rm Pul}}\to  \Omega_{K\times G}(K\times M)
\stackrel{{\rm Car}}\to \Omega_K(K\times_{G} M),$$
 where $\Omega_G(M)
\stackrel{{\rm Pul}}\to  \Omega_{K\times G}(K\times M)$ is the
natural  pullback map,
  induced by the projections on the second factor $K\times G \to G$,
 and $ \Omega_{K\times G}(K\times M)
\stackrel{{\rm Car}}\longrightarrow \Omega_K(K\times_{G} M)$ is the
Cartan map corresponding to a $K$-invariant
connection for the $G$-bundle $K\to K/G$~\cite{Me}. We now recall the description of this map.

Let $\Theta\in \Omega^1(K)\otimes \Gg$ be a $K$-invariant connection
on the $G$-bundle $K\to K/G$.
The associated principal $G$-bundle $$G\to K\times M \to \frac{K\times M}{G}\cong K\times_{G} M $$ carries a pullback connection, denoted by the same symbol $\Theta$ .
%Here $G$ acts on $K\times M $ by $(k,m).g=(k.g, g^{-1}.m)$.
We denote $F^{\Theta}=d\Theta +\cfrac{1}{2}[\Theta,\Theta]$ its curvature,
which is an element   in $\Omega^2_K(K\times M)\otimes \Gg$.
The equivariant momentum map
$\mu^{\Theta}\in (\Gk^*\otimes \Omega^0(K))^{K}\otimes \Gg$ is defined by
 $$\xi\in \Gk \mapsto \mu^{\Theta}(\xi )=-\iota_{\xi}\Theta$$
 where $\iota_{\xi}$ is the contraction along
 $\hat{\xi}\in \cal (K)$,
 the generating vector field of $\xi$.
  Then  $F^{\Theta}+\mu^{\Theta} $ is  the equivariant curvature of $\Theta$ \cite{BGV}.
Observe that $\Omega_{K\times G}(K\times M)\cong \left( S(\Gg^*)\otimes \Omega_K(K\times M)\right)^G$ that is the space of $G$-equivariant polynomial
 functions from $\Gg$ to $\Omega_K(K\times M)$.
Hence if $x\in \Gg\otimes \Omega_K^i(K\times M)$ and $P$ is a homogeneous
 degree $q$ polynomial on $\Gg$, then by substitution of variables, we get an element $P(x)$ in $\Omega_{K}^{2q+qi}(K\times M)$.
The Cartan map $\Omega_{K\times G}(K\times M) \to \Omega_K(K\times_{G} M) $ is the composition
\eqn P\otimes \omega \in \left( S(\Gg^*)\otimes \Omega_K(K\times
  M)\right)^G &\mapsto&  P(F^{\Theta}+\mu^{\Theta})\omega\in
\Omega_K(K\times M)\\
&\mapsto& {\rm Hor}\big(P(F^{\Theta}+\mu^{\Theta})\omega \big)\in \Omega_K(K\times_{G} M), \eneqn
where $ {\rm Hor}\:\Omega(K\times M)\to \Omega(K\times_{G} M)$ is the horizontal projection with respect to $\Theta$.

If moreover $F^{\Theta}=0$ and $K\times M\to K\times_{G} M$ admits a horizontal section $\sigma\:  K\times_{G} M\to K\times M$, we have the following lemma.
\begin{lem}\label{Inductiongeneralcase}
Let $P\otimes \omega$ be an element in $\left(S(\Gg^*)\otimes \Omega(M) \right)^G\cong \Omega_G(M)$. Then, ${\rm Ind}^G_K(P\otimes \omega) \in \Omega(K\times_{G} M)$ is the $K$-equivariant polynomial  on $\Gk$ with value
 in $\Omega(K\times_{G} M)$ defined, for any $\xi \in \Gk$, by
$${\rm Ind}^G_K(P\otimes \omega)\:\xi \mapsto \sigma^*\big(P(\mu^{\Theta}(\xi)){\rm pr}_2^*(\omega)  \big). $$
\end{lem}
\begin{pf}
First of all,
 ${\rm Pul}(P\otimes \omega)\in \left( S\big((\Gk\oplus \Gg)^*\big)\otimes \Omega(K\times M)\right)^{K\times G }$ is the map $\xi\oplus y\mapsto P(y){\rm pr}_2^*(\omega)$ for any $\xi\in \Gk$ and $y\in \Gg$. Then, by hypothesis,
$${\rm Hor}\big(P(F^{\Theta}+\mu^{\Theta})\omega \big)=\sigma^*\left(
P\big(\mu^{\Theta}(\xi )\big){\rm pr}_2^*(\omega)\right) \in \left(S(\Gk^*)\otimes \Omega(K\times_{G} M)\right)^K $$ and the lemma follows.
\end{pf}

Now let $K$ be the cartesian product group $G\times G$. We view $G$ as the diagonal subgroup of $K$.
 The $K$ action on itself by left multiplication
 commutes with the right $G$-action. We have a
principal right $G$-bundle
\begin{eqnarray*}
G\longrightarrow & K(=G\times G) &\longrightarrow G \\
&(g,h) & \mapsto gh^{-1} .
\end{eqnarray*}
 The left Maurer-Cartan form $\Theta^L_{MC}\in \Omega^1(G)\otimes \Gg$
on $G$ yields a $K$-invariant
 one-form $\Theta={\rm pr}_2^*\left(\Theta^L_{MC}\right)
 \in \Omega^1(K)\otimes \Gg$ by pullback
 along the projection on the second factor.
Then  $\Theta$ is a $K (=G\times G)$-invariant connection.
Moreover it is flat, thus its equivariant curvature reduces to
 the equivariant
 momentum $\mu^{\Theta}\:\Gk=\Gg\oplus \Gg \to \Omega^0(K)\otimes \Gg$.

\begin{lem}\label{momentum}
For any $(\alpha, \beta) \in \Gk (=\Gg\oplus \Gg)$,
and $(g,h)\in K (=G\times G)$ one has
$$\mu^{\Theta}(\alpha,\beta)|_{(g,h)}=- \Ad_{h^{-1}} \beta.
 $$
\end{lem}
\begin{pf}
The generating vector field for the left $G$-action on $G$
is given, for all $\beta\in \Gg$ by
$$\hat{\beta}|_{h}=\left. \cfrac{\partial}{\partial t}\,\!\right |_{t=0}
 \exp(t \beta)h=L_{h}\left( \Ad_{h^{-1}}\beta\right). $$ It follows, for any $(g,h)\in K=G\times G$, that
\begin{eqnarray*}
\mu^{\Theta}(\alpha,\beta)|_{(g,h)}&=& -\iota_{(\hat{\alpha},\hat{\beta})}(\Theta |_{(g,h)})\\
 &=& -\iota_{\hat{\beta}}\left. \Theta^{L}_{MC}\right |_{h}\\
&=& \left . -\Theta^{L}_{MC}\right |_{1}(\Ad_{h^{-1}}\beta )=-\Ad_{h^{-1}}\beta.
\end{eqnarray*}
\end{pf}

Let $M$ be a $K$ ($=G\times G$) space. It is then a $G$-space. Thus
we have an induction map
$${\rm Ind}^G_{G\times G}\:\Omega_{G} (M) \to \Omega_{G\times G}((G\times G) \times_{G} M)\cong  \Omega_{G\times G}(G\times M).$$
The group $G\times G$ acts on $G\times M$ by
 $$(k_1,k_2)\cdot (g,m)=\big(k_1gk_2^{-1},(k_1,k_2)\cdot m\big).$$

 \begin{lem}\label{homeo}\begin{enumerate}
 \item The map
$$(G\times G) \times_{G} M\to G\times M, \ \
(k_1,k_2,m)\mapsto \big(k_1k_2^{-1}, (k_1,k_2)\cdot m\big)$$
 is a $(G\times G)$-equivariant diffeomorphism.
\item The map
$$\sigma\: G\times M \to K\times M , \ \ \sigma (g,m)= \big(g,1,(g^{-1},1)\cdot m\big)$$
is a horizontal section for the principal $G$-bundle
 $G\to K\times M \to K\times_{G} M\cong G\times M$.
 \end{enumerate}
 \end{lem}

As a consequence,
 we have an isomorphism
$$\Omega_{G\times G}((G\times G) \times_{G} M)\stackrel{\sim}{\to}
  \Omega_{G\times G}(G\times M). $$
 Thus there is an  induction map
$${\rm Ind}^G_{G\times G}\:\Omega_{G} (M) \to \Omega_{G\times G}(G\times M) .$$

To obtain the Gysin map $H\com([M/G])\to H^{\scriptstyle \bullet-d}([M/G\times G])$,
  one simply composes the  induction map ${\rm Ind}^G_{G\times G}\:H\com([M/G])\to
H^{\scriptstyle \bullet}([G\times M /G\times G])$ with
the equivariant fiber integration map  \cite{AB}
 $H^{\scriptstyle \bullet}([G\times M /G\times G])\to H^{\scriptstyle \bullet-d}([ M /G\times G])$
 over the first factor $G$.

\begin{prop}
\label{GysinforLiegroupactions}
Given a $(G\times G)$-manifold $M$, the Gysin map
 $$H\com_{G}(M)\to H^{\scriptstyle \bullet-d}_{G\times G} (M) $$
 is given,  on the Cartan model, by
the  chain map
$\Psi\:\Omega_{G} (M)\to \Omega_{G\times G}
(M)$, $\forall P\otimes \omega \in  \left(S(\Gg^*)\otimes
\Omega(M)\right)^G$,
\eq \label{eq:psi}\Psi(P\otimes \omega)&=&\left( (\xi_1,\xi_2)\mapsto
  \int_{G}P(-\xi_2)\varphi^*(\omega)
\right), \ \ \ \forall \xi_1, \xi_2\in \Gg, \eneq
 where $\varphi\:G\times M\to M$ is the map $(g,m)\stackrel{\varphi}\mapsto
 (g^{-1},1)\cdot m $, and
$\int_{G}$ stands for the
 fiber integration over the first factor $G$.
\end{prop}
\begin{pf}
The induction map ${\rm Ind}^G_{G\times G}\:\Omega_{G} (M) \to
\Omega_{G\times G}(G\times M)$ is a chain level representative of the
stacks isomorphisms
\eqn &&
[M/G] \stackrel{\sim}\longleftarrow [ G\times G \times M/G\times
  G\times G] \stackrel{\sim}\longrightarrow [ G\times G \times_{G} M/G\times
  G].
\eneqn
induced by Morita equivalences of groupoids. Thus the Gysin map
$\Delta_!\:H\com{[M/G]}\to H^{\scriptstyle \bullet-d}{[M/G\times G]}$ is the
composition of ${\rm Ind}^G_{G\times G}$ with the Gysin map
$H\com([G\times M/G\times G])\to H^{\scriptstyle \bullet -d}([M/G\times
  G])$ which, by Proposition~\ref{prop:AB} is the equivariant fiber integration.

We now need to express the induction map more explicitly.
 Recall that, for any $\alpha \in \Omega_{G} (M)$,
 ${\rm Ind}^G_{G\times G}(\alpha )\in \Omega_{G\times G}(G\times M)$.
That is,  ${\rm Ind}^G_{G\times G}(\alpha )$ is a polynomial function on
 $\Gk (=\Gg\oplus \Gg)$ valued in $\Omega(G\times M) $. Write $\varphi\: G\times M \to M$ for the composition $\varphi= {\rm pr}_2\circ \sigma $.
Thus $\varphi (g,m)= (g^{-1},1)\cdot m$.
 According to Lemma~\ref{Inductiongeneralcase}, it suffices
to   compute
 $\sigma^*\big(P(\mu^{\Theta}(\xi_1,\xi_2)){\rm pr}_2^*(\omega)  \big)$.
 By Lemma~\ref{homeo}.3 and Lemma~\ref{momentum} we find that $$\sigma^*(P(\mu^{\Theta}(\xi_1,\xi_2)))=P(-\xi_2) .$$
Now the very definition of $\varphi$ yields that for
 any $\alpha = P\otimes \omega \in \Omega_{G} (M)\cong \left(S(\Gg^*)
 \otimes \Omega(M)\right)^G$, and   $\forall \xi_1,\xi_2\in \Gg$,
$${\rm Ind}^G_{G\times G}(\alpha )(\xi_1,\xi_2)= P(-\xi_2)\varphi^*(\omega). $$
This concludes the proof.
\end{pf}

\begin{rmk}
If we identify an element of $\Omega_{G} (M)$ with a $G$-equivariant polynomial $Q\:\Gg\to \Omega(M)$, then Equation~\eqref{eq:psi} can be
written as follows.  $\forall (\xi_1,\xi_2)\in \Gk=\Gg\oplus \Gg$,
 $$\Psi(Q)(\xi_1,\xi_2)= \int_{G} \varphi^*(Q(-\xi_2)).$$
\end{rmk}

\medskip

We now go back to our special case.
% Let's recall the construcion of  the dual \hidden coproduct and  dual \hidden product in cohomomology for this case.
Denote by $m:G\times G \to G$ and $\Delta\:G\to G\times G$
the group multiplication and the
 diagonal map respectively. The diagonal map  induces a stack map $\Delta\:[G\times G/G]
 \to[G\times G/G\times G] $ and thus a Gysin map
$$\Delta_!\:H\com([G\times G/G])\to
 H^{\scriptstyle \bullet-d}([G\times G/G\times G]),$$
 which is given by Proposition~\ref{GysinforLiegroupactions}.
 Similarly the group
  multiplication $m$ induces a stack map $m:[G\times G/G]\to [G/G] $ and thus a Gysin map
 $$m_!\: H^{\scriptstyle \bullet}([G\times G/G])\to H^{\scriptstyle \bullet-d}([G/G]).$$
 Since $m$ is $G$-equivariant, this  is the usual $G$-equivariant  Gysin map
on manifolds according to Proposition \ref{prop:AB}.

Note that $H_{G}\com (G)$ is a free module over
 $H\com ([*/G])\cong S(x_1,\dots,x_l)$. In fact,
 $H\com_G(G)=H\com ([*/G])[y_1,\dots,y_l]$
(the $y_j$s are of odd degrees).
Thus elements of  $H_{G}\com (G)$ are linear combinations of
monomials $y_1^{\epsilon_1}...y_l^{\epsilon_l}$,
 where each $\epsilon_j$ is either $0$ or $1$.
Similarly $H_{G}\com (G\times G)$ is the free $H_{G}\com (G)$-module
 generated by the monomials $y_1^{\epsilon_1}\cdots
y_l^{\epsilon_l}{y_1'}^{\epsilon'_1}\dots {y'_l}^{\epsilon'_l}$.

\begin{lem}\label{mshriek}
The map $m_!$ is a  $H\com ([*/G])$ linear map defined by
$$m_!(y_1^{\epsilon_1}...y_l^{\epsilon_l}{y_1'}^{\epsilon'_1}\dots {y'_l}^{\epsilon'_l})=y_1^{\epsilon_1+\epsilon_1'-1}...y_l^{\epsilon_l+\epsilon_l'-1}$$ with the convention that $y_j^{-1}=0$.
\end{lem}
\begin{pf}
Since $m: G\times G \to G$ is $G$-equivariant, the Gysin  map
 $m_!\:H^*([G\times G/G])\to H^*([G/G])$ is a map of $H\com
 ([*/G])$-module, and by Proposition~\ref{prop:AB}, it is the equivariant fiber integration of the principal bundle $G\times G \to G$. It can be represented on the Cartan cochain complex by integration of forms, see~\cite{GuSt} for details. In particular $m_!(y_1^{\epsilon_1}...y_l^{\epsilon_l}{y_1'}^{\epsilon'_1}\dots {y'_l}^{\epsilon'_l})$ is determined by the equation
\eq \label{eq:fintm} \!\!\int_{G\times G}  \!\! \!\!\!m^*(\alpha)
\wedge (y_1^{\epsilon_1}...y_l^{\epsilon_l}{y_1'}^{\epsilon'_1}\dots
{y'_l}^{\epsilon'_l}) \!&=& \!\! \! \int_G \alpha \wedge
m_!(y_1^{\epsilon_1}...y_l^{\epsilon_l}{y_1'}^{\epsilon'_1}\dots
{y'_l}^{\epsilon'_l}) \eneq Since the volume form on $G$ and
$G\times G$ are respectively given by $y_!\dots y_l$ and
$y_1...y_l{y_1'}\dots {y'_l}$, Equation~\eqref{eq:fintm} implies
  that $m_!\: H^*(G\times G)\to  H^{*-d}(G)$ sends $y_1^{\epsilon_1}...y_l^{\epsilon_l}{y_1'}^{\epsilon'_1}\dots {y'_l}^{\epsilon'_l}$ to $ y_1^{\epsilon_1+\epsilon_1'-1}...y_l^{\epsilon_l+\epsilon_l'-1}$. This finishes the proof.
\end{pf}

The \hidden product and coproduct on $H\lcom([G/G])$,
by universal coefficient theorem,
induces a degree $-d$ coproduct
$\delta\: H\com ([G/G]) \to H\com([G/G]) \otimes H\com([G/G])$
and degree $-d$ product
$\star: H\com([G/G]) \otimes H\com([G/G])\to H\com([G/G])$
which makes $H\com([G/G])$ into a Frobenius algebra, called the
{\em dual Frobenius structure} on $H\com([G/G])$.

More explicitly, these two operations are given by
the following compositions:
\eqn \delta\: H\com([G/G]) \stackrel{m^*}\to H\com([G\times G/G])
 &\stackrel{\Delta_!} \to& H^{\scriptstyle \bullet-d}([G\times G/G \times G])\\
&\to&
 \bigoplus_{i+j=\bullet -d}H^i([G/G])
 \otimes H^j([G/G]),\eneqn
and
$$\star:H\com([G/G])\otimes H\com([G/G])\cong H\com([G\times G/G\times  G])\stackrel{\Delta^*}
\to H\com([G\times G/G])\stackrel{m_!}\to H^{\scriptstyle \bullet-d}([G/G]).$$

\begin{them}
Let $G$ be a compact connected Lie group.
The dual \hidden coproduct on $H\com([G/G])$ is trivial.
And the dual \hidden product on $H\com([G/G])$ is given as follows.
For any
$P(x_1, \dots , x_l)y_1^{\epsilon_1}...y_l^{\epsilon_l}$ and $Q(x_1, \dots ,
x_l)y_1^{\epsilon_1'}...y_l^{\epsilon_l'}$ in $H\lcom([G/G])$,
we have
 $$\left(P(x_1, \dots , x_l)y_1^{\epsilon_1}...y_l^{\epsilon_l}\right)
 \star \left(Q(x'_1, \dots , x'_l){y_1'}^{\epsilon'_1}\dots
{y'_l}^{\epsilon'_l}\right) \qquad \qquad $$
$$\qquad \qquad =(PQ)(x_1, \dots , x_l)y_1^{\epsilon_1+\epsilon_1'-1}...y_l^{\epsilon_l+\epsilon_l'-1}$$
with the convention that $y_j^{-1}=0$.
\end{them}
\begin{pf}
  On the Cartan model, by Proposition~\ref{GysinforLiegroupactions}, the
\hidden coproduct is given by the following composition of chain maps:
$$\Omega_G(G)\stackrel{p^*}\longrightarrow \Omega_{G}(G\times G)\stackrel{\Psi^*}\longrightarrow \Omega_{G\times G}(G\times G)\stackrel{\cong}\longrightarrow \Omega_G(G)\otimes \Omega_G(G) .$$
Here the last map is K{\"u}nneth formula,
and the first map
 $$p^*\:\Omega_G(G)\cong \left(S(\Gg^*)\otimes \Omega(G)\right)^G\to \Omega_{G}(G\times G)\cong \left(S(\Gg^*)\otimes \Omega(G\times G)\right)^G$$
is  $S(\Gg^*)$-linear and  given by
 $$p^*(P\otimes \omega)=P \otimes m^*(\omega), \ \ \ \forall
 P\otimes \omega \in (S(\Gg^*)\otimes \Omega(G))^G .$$
Note that the space
$\Omega_G(G)\otimes \Omega_G(G)\cong \left(S(\Gg^*)\otimes
\Omega(G)\right)^G\otimes  \left(S(\Gg^*)\otimes \Omega(G)\right)^G$ has a
$S(\Gg^*)^G$-module structure, which is
 given by multiplication on the second factor: i.e.  $\forall
Q\in S(\Gg^*)$, $P_1\otimes \omega_1 \otimes P_2\otimes \omega_2\in
\Omega_G(G)\otimes \Omega_G(G)$, one defines
$$Q\cdot (P_1\otimes \omega_1 \otimes P_2\otimes \omega_2)= \left\{(\xi_1,\xi_2)\mapsto P_1(\xi_1)\otimes \omega_1 \otimes Q(-\xi_2)P_2(\xi_2)\otimes \omega_2 \in \Omega(G)\otimes \Omega(G)\right\}.$$
By Proposition \ref{GysinforLiegroupactions},
we know that the
 Gysin map $\Psi^*\:\Omega_{G}(G\times G)\to \Omega_{G\times G}(G\times G)$
 is indeed  a $S(\Gg^*)^G$-module map.

 There are two kinds of elements $P\otimes \omega$ in
$H\com([G/G])=\big(S(\Gg^*)\big)^G \otimes \Lambda(\Gg)^{G}$.  One
consists of those where  $\omega $ is a top degree form, i.e.
 a multiple of $y_1\wedge\dots \wedge y_l$,  and the
 others are  those where
 $\omega$ corresponds to a form in $\Omega^{*<d}(G)$.
In the latter case, $\Psi (P\otimes \omega )$ vanishes
 after fiber integration for degree reasons.
In the first case, the $G$-action on $G$ is by conjugation. Since the conjugacy action is trivial in cohomology, $\int_G \varphi^*(\omega) =0$ and by
 Proposition~\ref{GysinforLiegroupactions},
 $\Psi (P\otimes \omega )$vanishes.
Hence the dual \hidden coproduct is trivial.

We now compute the dual \hidden product.
First,  by a simple computation, we know that, on the Cartan model,
the map $\Delta^*\:H_{G\times G}^*(G\times G)\to H_{G}^*(G\times G)$ is
given by \begin{align*}\Delta^*(P(x_1,\dots,x_l,y_1,\dots,y_l,& x_1',\dots,x_l\
',y_1',\dots,y_l'))\\
&=P(x_1,\dots,x_l,y_1,\dots ,y_l,x_1,\dots,x_l
,y_1',\dots,y_l'). \end{align*}
In other words,
 the map $\Delta^*$ is an algebra map that leaves the odd
degree generators $y_i, y_j'$
unchanged and send both generators $x_i$, $x_i'$ ($i=1\dots r$) to the generator $x_i$.
 By Lemma~\ref{mshriek},
one obtains that
\eqn
 (m_!\circ \Delta^* ) \big(y_1^{\epsilon_1}...y_l^{\epsilon_l},
{y_1'}^{\epsilon'_1}\dots
{y'_l}^{\epsilon'_l}\big)&=&y_1^{\epsilon_1+\epsilon_1'-1}...y_l^{\epsilon_l+\epsilon_l'-1}.
\eneqn
The dual \hidden product now follows from the explicit
$S(\Gg^*)$-module structure.
\end{pf}

\begin{numrmk}
It follows that the \hidden product on $H\lcom([G/G])$ is trivial
while the \hidden coproduct has a counit given by the fundamental
class of $G$, which is dual of the cohomology class $y_1\dots y_l$.
\end{numrmk}

\begin{numrmk}
 Let $G$ be either a compact Lie group or a discrete group. As we have seen above, T
 the stack $[*/G]$ is strongly oriented and,  according to Theorem~\ref{frobeniusloop}, $H\lcom(\Lo [*/G])$ is a Frobenius algebra. An alternative approach to string topology for $[*/G]$ has been carried out in Gruher-Salvatore~\cite{GrSa}. It would be interesting to find a precise link between the result of  this Section~\ref{Liegroups} with those of~\cite{GrSa}. Similarly, it would be interesting to find the connection between the construction relate  of our \hidden product (Theorem~\ref{th:stringproduct}) with that of Abbaspour-Cohen-Gruher \cite{ACG} for Poincar\'e duality groups. Another approach to the $BV$-algebra structure for $\Lo[*/G]$ was recently studied by Chataur-Menichi~\cite{CM}. Their results seem to agree with ours.
\end{numrmk}

%%%%%%%%%%%%%%%%%%%%%%%%%%%%%%%%%%%%%%%%%%%%%%%%%%%%%%%%%%%%%%%%%%%%%%%%%%%%%%%%

\notsure{

\appendix 

\section{Categories fibered in groupoids}
\label{S:FiberedCat}
The formalism of categories fibered in groupoids provides a
convenient framework for working with lax groupoid-valued functors.
In this section, we recall some basic facts about categories fibered
in groupoids.

Let $\sfT$ be a fixed category. An example to keep in mind is
$\sfT=\Top$, the category of topological spaces. A {\bf category
fibered in groupoids} over $\sfT$ is a category $\X$ together with a
functor $\pi \: \X \to \sfT$ satisfying the following properties:

\begin{itemize}
  \item[($\mathbf{i}$)] For every arrow $f \: V \to U$ in $\sfT$,
   and for every object $X$ in $\X$ such that $\pi(X)=U$, there is
   an arrow $F \: Y \to X$ in $\X$ such that $\pi(F)=f$.

  \item[($\mathbf{ii}$)] Given a commutative triangle in $\sfT$,
    and a partial lift for it to $\X$ as in the diagram
     $$\xymatrix@C=12pt@R=0pt@M=8pt{ Y \ar[rd]^{F} & & & V \ar[rd]^{f} \ar[dd]_h & \\
         & X & \llra{\pi} & & U\\
        Z \ar[ru]_{G} & & & W \ar[ru]_{g} & } $$
  there is a unique morphism $H \: Y \to Z$ such that the triangle commutes and
  $\pi(H)=h$.

\end{itemize}

We will often drop the base functor $\pi$ from the notation and
denote  a fibered category $\pi \: \X \to \sfT$ by $\X$.

For a fixed object $T \in \sfT$, we let $\X(T)$ denote the category
of objects $X \in \X$ such that $\pi(X)=T$. Morphisms in $\X(T)$ are
morphisms $f \: X \to Y$ in $\X$ such that $\pi(f)=\id_T$.

It is easy to see that $\X(T)$ is a groupoid. This groupoid is
sometimes called the {\em fiber of $\X$ over $T$}. It is also called
the groupoids of {\em $T$-points} of $\X(T)$.

\begin{ex}{\label{E:fibered}}\par\noindent
  \begin{itemize}
  \item[$\mathbf{1.}$]  Let $T=\Top$, and let $G$ be a topological
  group. Let $\mathcal{B}G$ be the category of all principal $G$-bundles $P \to
  T$.
  A morphisms in $\mathcal{B}G$ is a $G$-equivariant cartesian
  diagram
      $$\xymatrix@C=12pt@R=10pt@M=8pt{P' \dto \rto & P
      \dto\\ T' \rto  & T}$$
  The base functor $\mathcal{B}G \to \Top$ is the forgetful functor
  that sends $P \to T$ to $T$.
  Observe that $\mathcal{B}G(T)$ is the groupoid of principal
  $G$-bundles over $T$.

  \item[$\mathbf{2.}$]  Let $\sfT=\Top$, and let $X$ be a topological
  space. Let $\X$ be the category of continuous maps $T \to X$.
  A morphism in $\X$ is a commutative triangle
         $$\xymatrix@C=6pt@R=10pt@M=8pt{T' \ar[rd] \ar[rr] & & T
         \ar[ld] \\ & X & }$$
  The forgetful functor that sends $T \to X$ to $T$ makes $\X$ a
  category fibered in groupoids over $\Top$.

  The groupoids $\X(T)$ is in fact equivalent to a set, namely,
  the set of continuous maps $T \to X$ (i.e., the set of  $T$-points
  of $X$).
\end{itemize}
\end{ex}

\begin{rem}
  There are
  two ways of thinking of a fibered category $\X \to \sfT$. One is
  to think of it as a device for cataloguing  the objects parameterized
  by a {\em moduli problem} over $\sfT$. In this case, an object $X \in
  \X(T)$ is viewed as a ``family parameterized by $T$.''

  The second point of view is to think of $\X$ as some kind of
  a {\em space}. In this case, an object in $\X(T)$ is simply thought of as
  a $T$-valued point of $\X$, that is, a map from $T$ to $\X$.

  The Yoneda type Lemma \ref{L:Yoneda} clarifies this dual point of view.
\end{rem}

\begin{rem}{\label{R:fibered}}\par\noindent
 \begin{itemize}
    \item[$\mathbf{1.}$]
     Conditions ($\mathbf{i}$) and ($\mathbf{ii}$) imply that, for
     every morphism $f \: T' \to T$ in $\sfT$, every object $X \in
     \X(T)$ has a ``pull-back'' $f^*(X)$ in $\X(T')$. The pull-back
     is  unique up to a unique isomorphism. We sometimes denote
     $f^*(X)$ by $X|_{T'}$.

    \item[$\mathbf{2.}$] The pull-back functors $f^*$ (whose definition
    involves making some choices) give rise to a lax groupoid-valued
    functor $T \mapsto \X(T)$. Conversely, given a lax groupoid-valued functor on
    $\sfT$, it is possible to construct a category fibered in groupoids over
    $\sfT$ via the so-called Grothendieck construction.
 \end{itemize}
\end{rem}

\subsection{The 2-category of fibered categories}

Categories fibered in groupoids over $\sfT$ form a 2-category. Let
us explain how this works.

A {\em morphism} $f \: \X \to \Y$ of fibered categories is a functor
$f \: \X \to \Y$ between the underlying categories such that
$\pi_{\Y}\circ f =\pi_{\X}$. Given two such morphisms  $f,g \: \X
\to \Y$, a {\em 2-morphism} $\varphi \: f \Rightarrow g$ between
them is a natural transformation of functors $\varphi$ from $f$ to
$g$ such that the composition $\pi_{\Y}\circ\varphi$ is the identity
transformation from $\pi_{\X}$ to itself.

With morphisms and 2-morphisms as above, categories fibered in
groupoids over $\sfT$ form a 2-category $\mathfrak{Fib}_{\sfT}$. The
2-morphisms in $\mathfrak{Fib}_{\sfT}$ are automatically invertible.

The construction in Example \ref{E:fibered}.2 can be performed in
any category $\sfT$ and it gives rise to a functor $\sfT \to
\mathfrak{Fib}_{\sfT}$. From now on, we will use the same notation
for an object $T$ in $\sfT$ and for its corresponding category
fibered in groupoids.

We have the following Yoneda-type lemma.

\begin{lem}[Yoneda lemma]{\label{L:Yoneda}}
  Let $\X$ be a category fibered in groupoids over $\sfT$, and let $T$
  be an object in $\sfT$. Then, the natural functor
    $$\Hom_{\mathfrak{Fib}_{\sfT}}(T,\X) \to \X(T)$$
  is an equivalence of groupoids.
\end{lem}

This lemma implies that the functor $\sfT \to \mathfrak{Fib}_{\sfT}$
is fully faithful. That is, we can think of the category $\sfT$ as a
full subcategory of $\mathfrak{Fib}_{\sfT}$. For this reason, in the
sequel we quite often do not distinguish between an object $T$ and
the fibered category associated to it.

\subsection{Descent condition}

To simplify the exposition, and to avoid the discussion of
Grothendieck topologies, we will assume from now on that
$\sfT=\Top$.

We say that a category $\X$ fibered in groupoids over $\sfT$ is a
{\bf stack}, if the following two conditions are satisfied:

\begin{itemize}
  \item[($\mathbf{i}$)] ({\em Gluing morphisms.}) Given two objects
  $X$  and $Y$ in $\X$ over a fixed topological spaces $T$,
  morphisms between them form a sheaf. That is, the presheaf of sets on
  $T$ defined by
   $$U \mapsto \Hom_{\X(U)}(X|_U,Y|_U)$$
  is a sheaf.

  \item[($\mathbf{ii}$)] ({\em Gluing objects.})  Let $T$ be a
  topological space, and let $\{U_i\}$ be an open covering of $T$.
  Assume we are given objects $X_i \in \X(U_i)$, together with
  isomorphisms $\varphi_{ij} \: X_j|_{U_i\cap U_j} \to X_i|_{U_i\cap
  U_j}$ in $\X(U_i\cap U_j)$ which satisfy the cocycle condition
    $$\varphi_{ij}\circ\varphi_{jk}=\varphi_{ik}$$
  on $U_i\cap U_j\cap U_k$ for every triple of indices $i$, $j$ and $k$. Then, there is an
  object $X$ over $T$, together with isomorphisms $\varphi_i \:
  X|_{U_i} \to X_i$ such that  $\varphi_{ij}\circ\varphi_i=\varphi_j$.
\end{itemize}

The data given in ($\mathbf{ii}$) is usually called a {\em gluing
data} or a {\em descent data}. It follows from ($\mathbf{i}$) that
the object $X$ in ($\mathbf{ii}$) is unique up to a unique
isomorphism.

Stacks over $\sfT$ form a full sub 2-category of $\mathfrak{Fib}$.

\begin{ex}{\label{E:descent}}\par\noindent
  \begin{itemize}
   \item[$\mathbf{1.}$]  The fibered category $\mathcal{B}G$ of Example
      \ref{E:fibered}.1 is a stack. This is because one can glue principal
       $G$-bundles over a fixed space $T$ using a gluing
       data (and the same thing is true for morphisms of principal
       $G$-bundles as well).

   \item[$\mathbf{2.}$] The fibered category $\X$ of Example
      \ref{E:fibered}.2 is a stack. This is because, given a
      collection of continuous maps $f_i \: U_i \to X$ which are
      equal
      over the intersections $U_i \cap U_j$, we can uniquely glue
      them to a continuous map $f \: T \to X$.
\end{itemize}
\end{ex}

Note that the cocycle condition over triple intersections does not
appear in Example \ref{E:descent}.2. The reason for this is that the
fiber  groupoids $\X(U)$ are equivalent to sets. That is, if there
is a morphisms between two objects in $\X(U)$ it has to be unique.

In view of Example \ref{E:descent}.2 (and Lemma \ref{L:Yoneda}), the
descent condition for a stack $\X$ can be interpreted as follows.
Let  $T$ be a topological space and $\{U_i\}$ an open covering  of
$T$. Assume we are given morphisms $f_i \: U_i \to \X$, together
with 2-isomorphisms $\varphi_{ij} \: f_j|_{U_i\cap U_j} \Rightarrow
f_i|_{U_i\cap U_j}$, satisfying the cocycle condition
$\varphi_{ij}\circ\varphi_{jk}=\varphi_{ik}$. (This should be
thought of as saying that $\varphi_{ij}$ are ``identifying $f_i$ and
$f_j$ along $U_i\cap U_j$.'') Then, we can glue $f_i$ to a global
map $f \: T \to \X$ whose restriction $f|_{U_i}$ to $U_i$ is
identified to $f_i$ via a 2-isomorphism $\varphi_i \: f|_{U_i}
\Rightarrow f_i$.

\subsection{Quotient stacks}

To any topological groupoid $\mathbb{X}=[X_1\sst{}X_0]$ one can
associate a stack $[X_0/X_1]$ called the {\em quotient stack} of the
groupoid. A quick definition for this quotient stack is as follows.
By definition, $[X_0/X_1]$ is the stack associated to the (fibered
category associated to the) presheaf of groupoids
  $$ T \ \mapsto \ [X_1(T) \toto X_0(T)].$$
Since we have not discussed the stack associated to a category
fibered in groupoids, we give an alternative description of
 $[X_0/X_1]$ in terms of principal bundles.

We only describe the case when $\mathbb{X}$ is the action groupoid
$[X\times G \sst{} X]$ of the action of a topological group $G$ on a
topological space $X$ and refer the reader for the general case to
(\cite{Noohi}, $\S$\;12). In the case of a group action, the
quotient stack is denoted by $[X/G]$.

For a topological space $T$, the groupoid $[X/G](T)$ of $T$-points
of $[X/G]$ is the groupoid of pairs $(P,\varphi)$, where $P$ is a
principal $G$-bundle over $T$, and $\varphi \: P \to X$ is a
$G$-equivariant map.  The morphisms in $[X/G](T)$ are
$G$-equivariant morphisms $f \: P' \to P$ such that
$\varphi'=\varphi\circ f$.

It is easy to verify that $[X/G]$ is a stack. When $X$ is a point,
the quotient stack $[*/G]$ coincides with $\mathcal{B}G$ of Example
\ref{E:fibered}.1 and is called the {\em classifying stack} of $G$.
Remark that, by Lemma \ref{L:Yoneda}, the groupoid
$\Hom(T,\mathcal{B}G)$ of morphisms from $T$ to $\mathcal{B}G$ is
equivalent to the groupoid of principal $G$-bundles over $T$.

\section{Generalized Fulton-MacPherson bivariant theories}
\label{S:Generalizedbivariant}
In this section we recall the axioms of a
Fulton-MacPherson bivariant theory. Our theory
is slightly more general than the original approach of Fulton-MacPherson
in the following ways:

\begin{itemize}
  \item[-] Since we need to work with stacks, the underlying category of our
  theory is indeed a 2-category. All fiber products and commutative diagrams
  should be interpreted in the 2-categorical sense.
  Our bivariant groups, however, will be equal for 2-isomorphic moprhisms.
  \item[-] Fulton-MacPherson have a notion of a `confined morphism' 
  (along which you can push forward bivariant classes) while
  we believe it is more natural to have `confined triangles'.
  \item[-] Product of bivariant classes are only partially defined.
\end{itemize}

In the context of this paper, these differences, however, are not crucial
and give us the right amount of generality to define the desired Gysin maps.

% --------------------------------------------
\subsection{The underlying (2-)category}

 The underlying category of a generalized bivariant theory is
 a category $\sfC$ (rather, 2-category) with fiber products and a final object. The
 category $\sfC$ is equipped with
 the following structure:
 \begin{itemize}
 \item A class of commutative triangles called {\bf confined triangles}
       $$\xymatrix@=16pt@M=8pt{  X \ar[rd]_u\ar[rr]^f  & & Y \ar[ld]^v \\
                         & S &     }$$
          We usually write this triangle as $X\llra{f}Y\llra{v}S$.
          We sometime refer to the above triangle as {\em a
          morphism $f \: X \to Y$ confined relative to $S$}.

 \item  a class of squares called {\bf independent squares}
        $$\xymatrix@=16pt@M=8pt{  X' \ar[d]_{f'} \ar[r]^{g'} & X \ar[d]^f \\
                         Y' \ar[r]_g &    Y }$$
           Note: we will distinguish the above square from its transpose,
           so the transpose of an independent square may not be
           independent.
\item a class of morphisms  called {\bf adequate}.
 \end{itemize}

We require the following axioms to be satisfied:

\begin{itemize}
   \item[$\mathbf{A1.}$] A triangle $X\llra{f}X\llra{v}Z$ in
     which $f=\id_X$ is the identity map is confined.

   \item[$\mathbf{A2.}$] If the inside triangles in
     $$\xymatrix@=16pt@M=8pt{  X \ar[rd]_u\ar[r]^f  &
                                  Y \ar[d]^v \ar[r]^g &  Z\ar[ld]^w \\
                         & S &     }$$
      are confined, then so is the outside triangle.

   \item[$\mathbf{B1.}$] Any commutative square in  which the top
     and the bottom morphisms arethe  identity maps is independent.

   \item[$\mathbf{B2.}$] Any square obtained from juxtaposition
   (vertical, or horizontal) of independent squares is independent.

   \item[$\mathbf{C.}$]  If in the commutative diagram
       $$\xymatrix@=16pt@M=8pt{  X' \ar[d]_{g'} \ar[r]^{f'} &
                 Y' \ar[d]^g   \ar[r]^{v'} & S' \ar[d]    \\
                         X \ar[r]_f &   Y  \ar[r]_v & S}$$
       the left square (or its transpose) is independent and
       $f$ is confined relative to $S$, then $f'$ is confined
       relative to $S'$.
\item[$\mathbf{D.}$] All isomorphisms are adequate.

\end{itemize}

\begin{lem}{\label{L:confined1}}
 Given
  $$\xymatrix@C=14pt{X \ar[r]^{f}
                & Y  \ar[r]  & Z \ar[r] & W}$$
  if $f$ is confined relative to $W$ then $f$ is confined relative to
  $Z$.
\end{lem}

\begin{proof}
  Use Axioms $\mathbf{B1}$ and $\mathbf{C}$.
\end{proof}

% --------------------------------------------
\subsection{Axioms for a bivariant theory}
\label{Axioms}

A bivariant theory $T$ on such a category $\sfC$ assigns to every
morphism $f \: X \to Y$ in $\sfC$ a graded abelian group
$T(X\llra{f}Y)$, or $T(f)$ for short. We denote the $i^{th}$ graded
component, $i \in \bbZ$, of $T$ by $T^i$. We sometimes denote an
element $\alpha \in T(X\llra{f}Y)$ by
  $$\lllra{X}{f}{\al}{Y}.$$

 The functor $T$ support
three types of operations:

\begin{itemize}
\item {\em Product.} For \comment{I've added the adequate condition} every $f \: X \to Y$ and adequate
   $g \: Y \to Z$, there is a product
     $$T^i(X\llra{f}Y)\otimes   T^j(Y\llra{g}Z) \llra{\cdot} T^{i+j}(X\llra{g\circ f}Z).$$

\item {\em Pushforward.} Given a confined triangle
      $$\xymatrix@=16pt@M=8pt{  X \ar[rd]_u\ar[rr]^f  & & Y \ar[ld]^v \\
                         & S &     }$$
      there is a pushforward homomorphism
        $$f_* \: T^i(X\llra{u}S) \lra  T^i(Y\llra{v}S).$$
\item{\em Pullback.} For every independent square
      $$\xymatrix@=16pt@M=8pt{  X' \ar[d]_{f'} \ar[r]^{g'} & X \ar[d]^f \\
                         Y' \ar[r]_g &    Y }$$
      there is a pullback homomorphism
        $$g^* \: T^i(X\llra{f}Y) \lra  T^i(X'\llra{f'}Y').$$
  (Observe the abuse of notation.)
\end{itemize}

These operations should satisfy the following compatibility axioms:

\begin{itemize}
  \item[$\mathbf{A1.}$] {\em Product is associative.}
    Given a diagram
       $$\xymatrix@C=14pt{X \ar[r]^{f}_*+[o][F-]{\al} & Y
          \ar[r]^{g}_*+<4pt>[o][F-]{\be} & Z \ar[r]^{h}_*+<4pt>[o][F-]{\ga} & W}$$
 where $g$, $h$ and $h\circ g$ are adequate,   we have
      $$(\al\cdot\be)\cdot\ga=\al\cdot(\be\cdot\ga)$$
   in $T(h\circ g\circ f)$.

  \item[$\mathbf{A2.}$] {\em Pushforward is functorial.}
   If the triangles in
     $$\xymatrix@=16pt@M=8pt{  X \ar[rd]_u\ar[r]^f  &
                                  Y \ar[d]^v \ar[r]^g &  Z\ar[ld]^w \\
                         & S &     }$$
      are confined, then
        $$(g\circ f)_* = g_*\circ f_* \: T^i(X\llra{u}S) \lra  T^i(Z\llra{w}S).$$

  \item[$\mathbf{A3.}$] {\em Pullback is functorial.}
     If the squares in
    $$\xymatrix@=16pt@M=8pt{  X'' \ar[d]_{f''} \ar[r]^{h'}  &
    X'\ar[d]_{f'} \ar[r]^{g'} & X \ar[d]_{f}\\
                  Y'' \ar[r]_h      & Y' \ar[r]_g & Y    }$$
     are independent, then
      $$(g\circ h)^* = h^*\circ g^* \: T^i(X\llra{f}Y) \lra  T^i(X''\llra{f''}Y'').$$

  \item[$\mathbf{A12.}$] {\em Product and pushforward commute.}
         Given
           $$\xymatrix@C=14pt{X \ar[r]^{f} \ar@/_1pc/ [rr]_*+[o][F-]{\al}
                & Y  \ar[r]^{g} & Z \ar[r]^{h}_*+<4pt>[o][F-]{\be} & W}$$
         with $f$ confined relative to $W$ and $h$ adequate, we have
               $$f_*(\al\cdot\be)=f_*(\al)\cdot\be$$
         in $T(h\circ g)$.

  \item[$\mathbf{A13.}$] {\em Product and pullback commute.}
     Given
        $$\xymatrix@=16pt@M=8pt{  X' \ar[d]_{f'} \ar[r]^{h''}  &
                           X\ar[d]_{f}^*+<4pt>[o][F-]{\al} \\
         Y' \ar[r]^{h'} \ar[d]_{g'} & Y \ar[d]_g^*+<3pt>[o][F]{\be}  \\
                        Z' \ar[r]^h & Z   }$$
         with independent squares, $g$ and $g'$ adequate, we have
                $$h^*(\al\cdot\be)=h'^*(\al)\cdot h^*(\be)$$
         in $T(g'\circ f')$.

  \item[$\mathbf{A23.}$] {\em Pushforward and pullback commute.}
       Given
                $$\xymatrix@=18pt@M=8pt{  X' \ar[d]_{f'} \ar[r]^{h''}  &
                           X \ar[d]_f \ar@/^1pc/[dd]^*+<4pt>[o][F-]{\al} \\
         Y' \ar[r]^{h'} \ar[d]_{g'} & Y \ar[d]_g \\
                        Z' \ar[r]^h & Z   }$$
         with independent squares and $f$ confined relative to $Z$, we have
               $$f'_*(h^*\al)=h^*f_*(\al)$$
         in $T(g')$.

  \item[$\mathbf{A123.}$] {\em Projection formula.} Given
         $$\xymatrix@=22pt@M=8pt{
             X'\ar[d]_{f'} \ar[r]^{g'} & X \ar[d]_{f}^*+<4pt>[o][F-]{\al} & \\
                Y' \ar[r]^g \ar@/_1pc/[rr]_*+<4pt>[o][F-]{\be}& Y \ar[r]^h &   Z}$$
         with independent square, $g$ adequate and confined relative to
             $Z$ and $h\circ g$ adequate,
         we have
            $$\al\cdot g_*(\be)=g'_*(g^*\al\cdot\be)$$
         in $T(h\circ f)$.
\end{itemize}

We say a bivariant theory $T$ has  {\bf unital} if for every $X \in
\sfC$ there is an element $1_X \in T^0(X\llra{id}X)$ with the
following properties:
\begin{itemize}
   \item  For every $f \: W \to X$ and every $\al \in
          T(W\llra{f}X)$, we have $\al\cdot 1_X=\al$.
   \item  For every $g \: X \to Y$ and every $\be \in
          T(X\llra{g}Y)$, we have $1_X\cdot \be=\be$.
   \item   For every $g\: X' \to X$, we have $g^*(1_X)=1_{X'}$.
\end{itemize}

A bivariant theory $T$ is called {\bf skew-commutative}
(respectively, {\bf commutative}), if for any square
  $$\xymatrix@=22pt@M=8pt{
             X'\ar[d]_{f'} \ar[r]^{g'} & X \ar[d]_{f}^*+<4pt>[o][F-]{\al}  \\
                Y' \ar[r]^g_*+<4pt>[o][F-]{\be}& Y }$$
   that is independent or its transpose is independent, $g$ and $f$
             are adequate, we
   have
      $$g^*(\al)\cdot\be=(-1)^{\deg(\al)\deg(\be)}f^*(\be)\cdot\al$$
   (respectively, $g^*(\al)\cdot\be=f^*(\be)\cdot\al$).

\medskip

Note that we don't assume the class of adequate morphisms to be
closed; that is, if $f,g$ are adequate, $g\circ f$ might not be
adequate. However, in practice, it is convenient to specify  a
(large)  closed subclass of adequate maps, called the {\bf strongly
adequate} morphisms. In particular, the product of bivariant classes
are always defined and associative on the subclass  of strongly
adequate morphisms.

\bigskip

Using the definitions of Section~\ref{S:Bivariant} and results of
Sections~\ref{S:Thom}, \ref{S:Proper}, \ref{S:Techinical}, it is
straightforward to prove
\begin{thm}
The bivariant theory of Section~\ref{S:Bivariant} is a generalized
Fulton-MacPherson bivariant theory.
\end{thm}
Note that, in view of Lemma~\ref{L:superproper} and Example~\ref{E:super}.1, 
we can choose the class of strongly adequate morphisms to be
the class of \superproper maps.

}

%%%%%%%%%%%%%%%%%%%%%%%%%%%%%%%%%%%%%%%%%%%%%%%%%%%%%%%%%%%%%%%%
%%%%%%%%%%%%%%%%%%%%%%%%%%%%%%%%%%%%%%%%%%%%%%%%%%%%%%%%%%%%%%
%%%%%%%%%%%%%%%%%%%%%%%%%%%%%%%%%%%%%%%%%%%%%%%%%%%%%%%%%%%%%%

\end{document}